\documentclass[11pt,twoside]{article}
\usepackage[hmargin=0.8in,vmargin=0.9in]{geometry}
\usepackage{fancyhdr}
\usepackage{graphicx}
\usepackage{subfigure}
\usepackage{amssymb}
\usepackage{amsmath}
\usepackage{amsfonts}
\usepackage{theorem}
\usepackage{mathrsfs}
\usepackage{mathtools}
\usepackage{bm}
\usepackage{color}
\usepackage{setspace}
\usepackage{exscale}
\usepackage{relsize}
\usepackage{epstopdf}
\usepackage{float}
\usepackage{cite}
\usepackage{hyperref}

\usepackage{color}

\usepackage{setspace}

\usepackage{booktabs,multirow} 
\usepackage{array} 
\usepackage{paralist} 
\usepackage{verbatim} 
\usepackage{subfigure} 

\theoremstyle{plain}			
\newtheorem{thm}{Theorem}[section]

{\theorembodyfont{\rmfamily}}

\newenvironment{DA}{{\flushleft \bf Declarations:}}{}

\allowdisplaybreaks[1]

\numberwithin{equation}{section}
\numberwithin{figure}{section}
\numberwithin{table}{section}

\newcommand\eref[1]{(\ref{#1})}

\newcommand*\xbar[1]{%
  \hbox{%
    \vbox{%
      \hrule height 0.5pt 
      \kern0.4ex
      \hbox{%
        \kern-0.05em
        \ensuremath{#1}%
        \kern-0.00em
      }%
    }%
  }%
}

\newcommand{\mF}{\bm{F}}

\newcommand{\mG}{\bm{G}}

\newcommand{\mH}{\bm{H}}
\newcommand{\mU}{\bm{U}}

\newcommand{\mo}{\bm{0}}

\newcommand{\dt}{\Delta t}
\newcommand{\dx}{\Delta x}
\newcommand{\dy}{\Delta y}

\newcommand{\hf}{{\frac{1}{2}}}

\newcommand{\jph}{{j+\frac{1}{2}}}
\newcommand{\jmh}{{j-\frac{1}{2}}}
\newcommand{\kph}{{k+\frac{1}{2}}}
\newcommand{\kmh}{{k-\frac{1}{2}}}

\newcommand{\ajphp}{{a_{j+\frac{1}{2}}^+}}
\newcommand{\ajphm}{{a_{j+\frac{1}{2}}^-}}

\newtheorem{rmk}[thm]{Remark}
\title{A Comparative Study of Low-Dissipation Numerical Schemes for Hyperbolic Conservation Laws}
\author{Shaoshuai Chu\thanks{Department of Mathematics, RWTH Aachen University, 52056, Aachen, Germany;
{\tt chu@igpm.rwth-aachen.de}}~ and ~Michael Herty\thanks{Department of Mathematics, RWTH Aachen University, 52056, Aachen, Germany; Department of Mathematics and Applied Mathematics, University of Pretoria, Private Bag X20, Hatfield 0028, South Africa; {\tt herty@igpm.rwth-aachen.de}}}
\date{}
\begin{document}
\date{}
\maketitle
\begin{abstract}
This work provides a comparative assessment of several low-dissipation numerical schemes for hyperbolic conservation laws, highlighting their performance relative to the classical Harten–Lax–van Leer (HLL) schemes. The schemes under consideration include the classical
Harten–Lax–van Leer–Contact (HLLC), the recently proposed TV flux splitting, the low-dissipation Central-Upwind (LDCU), and the local characteristic decomposition-based Central-Upwind (LCDCU) schemes. These methods are extended to higher orders of accuracy, up to the fifth order, within both finite-volume and finite-difference frameworks. A series of numerical experiments for the one- and two-dimensional Euler equations of gas dynamics are performed to evaluate the accuracy, robustness, and computational efficiency of the studied schemes. The comparison highlights the trade-offs between resolution of contact and shear waves, robustness in the presence of shocks, and computational cost. The investigated low-dissipation schemes show comparable levels of numerical dissipation, with only subtle differences appearing in selected benchmark problems. The results provide practical guidance for selecting efficient low-dissipation solvers for the simulation of complex compressible flows.
\end{abstract}



\section{Introduction}

This paper focuses on the numerical solution of hyperbolic conservation laws, which in one and two spatial dimensions take the form
\begin{equation}\label{1.1}
  \mU_t + \mF(\mU)_x = \mo,
\end{equation}
and
\begin{equation}\label{1.2}
  \mU_t + \mF(\mU)_x + \mG(\mU)_y = \mo,
\end{equation}
where $t$ denotes time, $x$ and $y$ are spatial coordinates, $\mU \in \mathbb{R}^d$ is the vector of conserved variables, and $\mF, \mG: \mathbb{R}^d \to \mathbb{R}^d$ are the flux functions in the $x$- and $y$-directions, respectively.

It is well-known that the systems \eref{1.1} and \eref{1.2} may involve complex wave patterns, including shocks, rarefactions, and contact discontinuities, even when the initial data are very smooth, which presents significant challenges for numerical methods. Traditional approaches, particularly first-order methods, often fail to capture these discontinuities accurately and require extremely fine meshes to produce reliable results.

Since the pioneering works \cite{Fri,Lax,Godunov59}, numerous methods have been developed to solve hyperbolic systems \eref{1.1} and \eref{1.2}; for a detailed overview, see, e.g., the monographs and review papers \cite{Hesthaven18,Leveque02,KLR20,Shu20,Toro2009,BAF} and the references therein. The Harten--Lax--van Leer (HLL) solver from  \cite{HLL1983} is one of the most widely used approximate Riemann solvers for hyperbolic conservation laws. It replaces the exact Riemann solution, which may involve multiple nonlinear waves, by a simple two-wave model. Specifically, the solution is approximated by a single intermediate state separated from the left and right data by two waves propagating with estimated minimum and maximum signal velocities. The HLL flux is robust and, under suitable signal-speed estimates and a CFL restriction, positivity preserving, but the main drawback is that intermediate waves such as contact and shear discontinuities are excessively smeared, since they are not explicitly represented in the two-wave model. To overcome this limitation, the Harten–Lax–van Leer–Contact (HLLC) solver was introduced in \cite{TSS1994}, which extends the HLL framework by restoring the missing contact wave. The HLLC scheme thus consists of three waves, separating the solution into two intermediate states. The presence of the contact wave allows the solver to resolve contact discontinuities and shear layers more sharply, while maintaining the robustness and efficiency of the HLL scheme. Owing to this balance between accuracy and robustness, the HLLC flux has become one of the most popular approximate Riemann solvers for the compressible Euler equations and related hyperbolic systems.

Another important class of schemes is based on flux vector splitting. In the flux splitting methods, one can decompose fluxes into components corresponding to different wave families, such as shock, contact, and expansion waves, which makes them fundamental in computational fluid dynamics. This decomposition enhances the resolution of shock waves and discontinuities, which is crucial for accurate simulations of compressible flows; see, e.g., \cite{ATvL1986, ATR1989, SW1981, VL1982, VL1982a, Lious1993, Lious1996, Lious1998, Lious2006,DT2011,  GVM1999, ZB1993}. However, the classical flux splitting schemes (e.g., \cite{ATvL1986, ATR1989, SW1981, VL1982, VL1982a}) often struggle to resolve intermediate characteristic fields, leading to excessive dissipation or numerical artifacts, even though they are effective. To overcome these limitations, more advanced methods, such as the advection upstream splitting method (AUSM), have been developed (see, e.g., \cite{Lious1993}), offering better handling of contact waves and improved resolution of wave patterns. AUSM has since garnered considerable attention and undergone refinement, with further developments \cite{Lious1996, Lious1998, Lious2006}. Additionally, a flux splitting approach similar to AUSM was proposed in \cite{ZB1993}, and subsequent advancements have been recorded in works such as \cite{DT2011, LLN2012, GVM1999, KSFW2011}, further enriching the field of computational fluid dynamics. Recently, a new flux splitting method, known as the TV splitting scheme for the one-dimensional (1-D) Euler equations of gas dynamics, was introduced  in \cite{TV2012} and later extended to higher dimensions and orders in \cite{TCL2015,CHT25}. The proposed schemes are characterized by their simplicity, robustness, and accuracy, offering significant improvements over existing flux splitting methods. In particular, they effectively capture contact and shear waves while precisely preserving isolated stationary contacts. Beyond their applications to high-dimensional Euler equations, it has also been extended to other systems, including magnetohydrodynamics and shallow water equations; see, e.g., \cite{BMT2016,DBTF2019,TCVS2022,TT2017}.

While flux splitting methods such as TV splitting focus on decomposing fluxes along characteristic fields, an alternative direction is offered by central schemes. These methods avoid explicit Riemann solvers altogether and instead compute fluxes through staggered or modified control volumes \cite{Nessyahu90,Levy99,Liu98,Lie03a,Jiang98,Arminjon95}. Central schemes are easy to implement and broadly applicable, but they tend to suffer from relatively large numerical dissipation. To overcome this drawback, central-upwind (CU) schemes were developed in \cite{Kurganov01,Kurganov00}, in which control volumes are adapted to local wave speeds. CU schemes retain the simplicity of central methods while achieving significantly improved accuracy. Further refinements include anti-diffusion corrections \cite{Kurganov07}, dissipation switches \cite{Kurganov21a}, and more accurate wave-speed estimates \cite{Garg21}. Recently, two refinements of CU schemes have been proposed to further reduce dissipation. Low-dissipation CU (LDCU) schemes \cite{KX_22,CKX_24} introduce a novel projection step based on subcell resolution, which sharply approximates contact waves while ensuring non-oscillatory behavior of the projected solution; see also, e.g., \cite{CKX22,CKX24,KLX2025}. This significantly decreases dissipation in contact- and shear-dominated regions. In parallel, LCD-based CU (LCDCU) schemes \cite{CCHKL_22} apply local characteristic decomposition to the diffusion term, which also leads to substantially improved resolution by further suppressing excessive numerical dissipation; see also, e.g., \cite{CHK25,CK2023}.

As mentioned above, the HLLC, TV splitting, LDCU, and LCDCU schemes are low-dissipation numerical schemes for solving the hyperbolic conservation laws \eref{1.1} and \eref{1.2}, and it is instructive to check the dissipations contained in these schemes by comparing the resolution of the computed numerical results in a series of benchmarks. 
To this end, we carry out a systematic comparative study of the first-, second-, third-, and fifth-order HLL, HLLC, TV splitting, LCDCU, and LDCU schemes, focusing on their robustness, accuracy, and efficiency in resolving complex wave interactions in hyperbolic conservation laws. The second-order extension employs a piecewise linear interpolant for reconstructing one-sided point values, while the third- and fifth-order schemes are formulated within the finite-difference (FD) A-WENO framework. This framework effectively generalizes low-order finite-volume (FV) schemes to higher-order FD ones, especially in multidimensional settings, through its simple dimension-by-dimension reconstruction; see, e.g., \cite{liu17,Jiang13,Jiang96}.

The rest of this paper is organized as follows. In \S\ref{sec2}, we begin with a brief overview of the first-order HLL scheme for the 1-D Euler equations of gas dynamics and then extend it to second-, third-, and fifth-order accuracy. The corresponding HLLC, TV, LDCU, and LCDCU schemes are subsequently reviewed. In \S\ref{sec3}, we describe the two-dimensional (2-D) extensions of these schemes in a dimension-by-dimension manner. Finally, in \S\ref{sec4}, we present a number of 1-D and 2-D numerical results to compare their performance.

\section{One-Dimensional Schemes}\label{sec2}
In this section, we consider the 1-D Euler equations of gas dynamics and briefly review the first-, second-, third-, and fifth-order 1-D HLL, HLLC, TV, LDCU, and LCDCU schemes, respectively. The 1-D 
Euler equations of gas dynamics read as \eref{1.1} with 
\begin{equation}\label{2.1}
 \bm U:=(\rho,\rho u,E)^\top, \quad {\rm and} \quad \bm F(\bm U)=(\rho u,\rho u^2+p,u(E+p))^\top,
\end{equation}
where $\rho$, $u$, $p$, and $E$ are the density, velocity, pressure, and total energy, respectively. The system is completed through the following equations of state (EOS): 
\begin{equation}\label{2.2}
p=(\gamma-1)\Big[E-\hf\rho u^2\Big],
\end{equation}
where the parameter $\gamma$ represents the specific heat ratio.

\subsection{One-Dimensional HLL Schemes}\label{sec2.1}
We first briefly review the 1-D HLL scheme from \cite{HLL1983} and show its high-order extensions.

\subsubsection{One-Dimensional First-Order HLL Scheme}
Assume the computational domain is covered with uniform cells $C_j:=[x_\jmh,x_\jph]$ with $x_\jph-x_\jmh\equiv\dx$ centered at $x_j=(x_\jmh+x_\jph)/2$, $\,j=1,\ldots,N$, and the cell average values 
\begin{equation*}
  \xbar\mU_j(t):\approx\frac{1}{\dx}\int\limits_{C_j}\mU(x,t)\,{\rm d}x
\end{equation*}
are available at a certain time level $t$. The computed cell averages $\xbar \mU_j$ of the 1-D system \eref{1.1} are evolved in time by solving the following semi-discrete system of ordinary differential equations (ODEs):
\begin{equation}
\frac{{\rm d}\xbar \mU_j}{{\rm d}t}=-\frac{\bm{{\cal F}}^{\rm FV}_\jph-\bm{{\cal F}}^{\rm FV}_\jmh}{\dx},
\label{2.3}
\end{equation}
where $\bm{{\cal F}}^{\rm FV}_\jph=\bm{{\cal F}}^{\rm FV}_\jph\big(\bm U_\jph^-,\bm U_\jph^+\big)$ is the numerical flux, defined by
\begin{equation}\label{2.4}
\resizebox{0.93\linewidth}{!}{$
\bm{{\cal F}}^{\rm FV}(\mU^-_\jph,\mU^+_\jph) =
\begin{cases}
\mF(\mU^-_\jph), & a^-_\jph \geq 0, \\[8pt]
\dfrac{a^+_\jph \mF(\mU^-_\jph) - a^-_\jph \mF(\mU^+_\jph) +a^+_\jph a^-_\jph (\mU^+_\jph -\mU^-_\jph)}{a^+_\jph - a^-_\jph}, & a^-_\jph \leq 0 \leq a^+_\jph, \\[12pt]
\mF(\mU^+_\jph), & a^+_\jph \leq 0 ,
\end{cases}
$}
\end{equation}
where $\mU^\pm_\jph$ are the left/right-sided point values of $\mU$ at the cell interfaces $x_\jph$. For the first order scheme, we take $\mU^+_\jph=\mU_{j+1}$ and $\mU^-_\jph=\mU_{j}$. The one-sided local speeds of propagation $a^\pm_\jph$ can be estimated by
\begin{equation}\label{2.4a}
a^+_\jph=\max\Big\{u^+_\jph+c^+_\jph,u^-_\jph+c^-_\jph\Big\},\quad a^-_\jph=\min\Big\{u^+_\jph-c^+_\jph,u^-_\jph-c^-_\jph\Big\},
\end{equation}
where
\begin{equation*}
u^\pm_\jph=\frac{(\rho u)_\jph^\pm}{\rho^\pm_\jph},\quad p^\pm_\jph=(\gamma-1)\Big[E^\pm_\jph-\hf\rho^\pm_\jph\big(u^\pm_\jph\big)^2\Big],
\quad c^\pm_\jph=\sqrt{\frac{\gamma p^\pm_\jph}{\rho^\pm_\jph}}.
\end{equation*}

\subsubsection{One-Dimensional Second-Order HLL Scheme}
We now extend the first-order HLL scheme introduced in \S \ref{sec2.1} to the second order of accuracy. The resulting scheme \eref{2.3}--\eref{2.4a} achieves second-order accuracy provided that the one-sided point values $\mU^\pm_\jph$, used to compute the numerical flux \(\bm{{\cal F}}^{\rm FV}_\jph\), are second-order accurate. To this end, we approximate $\mU^\pm_\jph$ using the generalized minmod reconstruction \cite{lie03,Nessyahu90,Sweby84} with $\theta=1.3$

\subsubsection{One-Dimensional Third-Order HLL Scheme}
In this section, we extend the HLL scheme to the third-order accuracy in the framework of the FD A-WENO scheme introduced in \cite{Jiang13} (see also \cite{liu17}), which has been proven to be a powerful tool for generalizing low-order FV schemes to higher-order FD ones.

Following \cite{Jiang13}, the point values $\mU_j$ are evolved in time by solving the following system of ODEs:
\begin{equation}
\frac{{\rm d}\mU_j}{{\rm d}t}=-\frac{{\mH_\jph}-{\mH_\jmh}}{\dx},
\label{2.9}
\end{equation}
where ${\mH_\jph}$ is the (third-order accurate) numerical flux defined by
\begin{equation}
{\mH_\jph}=\bm{{\cal F}}^{\rm FV}_\jph-\frac{1}{24}(\dx)^2(\mF_{xx})_\jph.
\end{equation}
Here, $\bm{{\cal F}}^{\rm FV}_\jph$ is the FV numerical flux as in \eref{2.4} and $(\mF_{xx})_\jph$ is the higher-order correction term used to increase the order of the numerical flux. The correction term $(\mF_{xx})_\jph$ can be approximated with the help of the FV numerical fluxes $\bm{{\cal F}}^{\rm FV}_\jph$:
\begin{equation}
(\mF_{xx})_\jph=\frac{1}{(\dx)^2}\big[\bm{{\cal F}}^{\rm FV}_\jmh-2\bm{{\cal F}}^{\rm FV}_\jph+\bm{{\cal F}}^{\rm FV}_{j+\frac{3}{2}}\big],
\end{equation}
which has been proved to be more efficient than the old version
\begin{equation}
(\mF_{xx})_\jph=\frac{1}{2(\dx)^2}\big[{\mF}_{j-1}-{\mF}_{j}-{\mF}_{j+1}+{\mF}_{j+2}\big],
\end{equation}
where $\mF_j=\mF(\mU_j)$, while affecting neither the accuracy nor the quality of resolution. The resulting scheme is third-order once the one-sided point values $\mU^\pm_\jph$ employed to compute the numerical flux $\bm{{\cal F}}^{\rm FV}_\jph$ are third-order accurate. This can be done by implementing a certain nonlinear limiting procedure like the third-order WENO-type interpolation (see, e.g., \cite{CH_third,GRR2017,GRR2018,LXDYG2022}) applied to the local characteristic variables; see \cite[Appendix A]{CHT25} for a detailed explanation.

\subsubsection{One-Dimensional Fifth-Order HLL Scheme}
According to \cite{Jiang13}, to achieve fifth-order accuracy, the point values $\mU_j$ are evolved in time by solving the system \eref{2.9} with the (fifth-order accurate) numerical flux 
\begin{equation}
{\mH_\jph}=\bm{{\cal F}}^{\rm FV}_\jph-\frac{1}{24}(\dx)^2(\mF_{xx})_\jph+\frac{7}{5760}(\dx)^4(\mF_{xxxx})_\jph,
\end{equation}
where $\bm{{\cal F}}^{\rm FV}_\jph$ is the FV numerical flux as in \eref{2.4}, $(\mF_{xx})_\jph$ and $(\mF_{xxxx})_\jph$ are the higher-order correction terms computed by the fourth- and second-order accurate FDs, respectively; see, e.g., \cite{CKX23,CCK23_Adaptive}. Here, we have used the following
higher-order correction terms from \cite{CKX23}:
\begin{equation}\label{2.11aa}
\begin{aligned}
&(\mF_{xx})_\jph=\frac{1}{12(\dx)^2}\Big[-\bm{{\cal F}}^{\rm FV}_{j-\frac{3}{2}}+16\bm{{\cal F}}^{\rm FV}_\jmh-30\bm{{\cal F}}^{\rm FV}_\jph
+16\bm{{\cal F}}^{\rm FV}_{j+\frac{3}{2}}-\bm{{\cal F}}^{\rm FV}_{j+\frac{5}{2}}\Big],\\
&(\mF_{xxxx})_\jph=\frac{1}{(\dx)^4}\Big[\bm{{\cal F}}^{\rm FV}_{j-\frac{3}{2}}-4\bm{{\cal F}}^{\rm FV}_\jmh+6\bm{{\cal F}}^{\rm FV}_\jph-
4\bm{{\cal F}}^{\rm FV}_{j+\frac{3}{2}}+\bm{{\cal F}}^{\rm FV}_{j+\frac{5}{2}}\Big].
\end{aligned}
\end{equation}
In order to ensure the resulting scheme is fifth order, the one-sided point values $\mU^\pm_\jph$ employed to compute the numerical flux $\bm{{\cal F}}^{\rm FV}_\jph$ need to be at least fifth-order accurate. This can be done by using a certain nonlinear limiting procedure like the fifth-order WENO-Z interpolation from \cite{Jiang13,liu17,wang18} applied to the local characteristic variables; see \cite[Appendix B]{CHT25} for details.

\subsection{One-Dimensional HLLC Schemes}
According to \cite{TSS1994}, the 1-D first-, second-, third-, and fifth-order HLLC schemes can be obtained similarly to the corresponding HLL schemes by replacing the numerical fluxes $\bm{{\cal F}}^{\rm FV}_\jph$ in \eref{2.4} by 
\begin{equation}
\bm{{\cal F}}^{\rm FV}_\jph(\mU^-_\jph,\mU^+_\jph) =
\begin{cases}
\mF(\mU^-_\jph), & a^-_\jph \geq 0, \\[8pt]
\mF(\mU^-_\jph) + a^-_\jph \big(\mU_\jph^{*,-} - \mU^-_\jph \big), & a^-_\jph \leq 0 \leq a^*_\jph, \\[8pt]
\mF(\mU^+_\jph) + a^+_\jph \big(\mU_\jph^{*,+} - \mU^+_\jph \big), & a^*_\jph \leq 0 \leq a^+_\jph, \\[8pt]
\mF(\mU^+_\jph), & a^+_\jph \leq 0 ,
\end{cases}
\end{equation}
where 
\begin{equation}
\mU_\jph^{*,\pm}
=
\rho^\pm_\jph
\frac{a^\pm_\jph-u^\pm_\jph}
     {a^\pm_\jph-a^*_\jph}
\begin{pmatrix}
1 \\[4pt]
a^*_\jph \\[4pt]
\displaystyle
\frac{E^\pm_\jph}{\rho^\pm_\jph}
+
\big(a^*_\jph-u^\pm_\jph\big)
\left(
a^*_\jph+
\frac{p^\pm_\jph}
     {\rho^\pm_\jph\big(a^\pm_\jph-u^\pm_\jph\big)}
\right)
\end{pmatrix},
\end{equation}
and 
\begin{equation}
a^*_\jph= \frac{p^+_\jph -p^-_\jph + \rho^-_\jph u^-_\jph (a^-_\jph-u^-_\jph) - \rho^+_\jph u^+_\jph (a^+_\jph-u^+_\jph)}{\rho^-_\jph (a^-_\jph-u^-_\jph)-\rho^+_\jph (a^+_\jph-u^+_\jph)}.
\end{equation}

\subsection{One-Dimensional TV Splitting Schemes}
According to \cite{TV2012,CHT25}, the 1-D first-, second-, third-, and fifth-order TV schemes can be obtained by replacing the numerical fluxes $\bm{{\cal F}}^{\rm FV}_\jph$ in \eref{2.4} by 
\begin{equation}
\bm{{\cal F}}^{\rm FV}_\jph\big(\bm U_\jph^-,\bm U_\jph^+\big)= \bm{{\cal F}}^A_\jph\big(\bm U_\jph^-,\bm U_\jph^+\big) + \bm{{\cal F}}^P_\jph \big(\bm U_\jph^-,\bm U_\jph^+\big).
\end{equation}
Here,  $\bm{{\cal F}}^A_\jph\big(\bm U_\jph^-,\bm U_\jph^+\big)$ is the advection flux given by 
\begin{equation}
  \bm{{\cal F}}^A_\jph\big(\bm U_\jph^-,\bm U_\jph^+\big)=\begin{cases}
                                                          u^*_\jph\begin{pmatrix} \rho^-_\jph \\[0.8ex] 
                                                          (\rho u)^-_\jph\\[0.8ex]\dfrac{1}{2} \rho^-_\jph (u^-_\jph)^2\end{pmatrix}, & \mbox{if } u^*_\jph \ge 0, \\[8.ex]
                                                          u^*_\jph \begin{pmatrix} \rho^+_\jph \\[0.8ex] (\rho u)^+_\jph\\[0.8ex]\dfrac{1}{2} \rho^+_\jph (u^+_\jph)^2\end{pmatrix}, & \mbox{otherwise},
                                                  \end{cases}
\end{equation}
and  $\bm{{\cal F}}^P_\jph \big(\bm U_\jph^-,\bm U_\jph^+\big)$ is the pressure flux given by 
\begin{equation}
  \bm{{\cal F}}^P_\jph \big(\bm U_\jph^-,\bm U_\jph^+\big)= \begin{pmatrix} 0,\, p^*_\jph,\, \dfrac{\gamma  u^*_\jph  p^*_\jph}{\gamma-1} \end{pmatrix}^\top,
\end{equation}
where 
\begin{equation}
 \begin{aligned}
  u^*_\jph &= \frac{C^+_\jph u^+_\jph - C^-_\jph u^-_\jph} { C^+_\jph-C^-_\jph}-\frac{2}{ C^+_\jph-C^-_\jph} (p^+_\jph - p^-_\jph), \\
  p^*_\jph &= \frac{C^+_\jph p^-_\jph - C^-_\jph p^+_\jph} { C^+_\jph-C^-_\jph}+\frac{C^+_\jph C^-_\jph}{ 2(C^+_\jph-C^-_\jph)} (u^+_\jph - u^-_\jph),\\
  C^\pm_\jph &=\rho^\pm_\jph \bigg(u^\pm_\jph \pm \sqrt{\big(u^\pm_\jph \big)^2+4 \big(c^\pm_\jph\big)^2}\bigg).
 \end{aligned}
\end{equation}

\subsection{One-Dimensional LDCU Schemes}
According to \cite{KX_22,CKX_24,CKX22}, the 1-D first-, second-, third-, and fifth-order LDCU schemes can be obtained by replacing the numerical fluxes $\bm{{\cal F}}^{\rm FV}_\jph$ in \eref{2.4} by 
\begin{equation}
\bm{{\cal F}}^{\rm FV}_\jph\big(\bm U_\jph^-,\bm U_\jph^+\big)=\frac{\ajphp\mF\big(\mU^-_\jph\big)-\ajphm \mF\big(\mU^+_\jph\big)}{\ajphp-\ajphm}+
\frac{\ajphp\ajphm}{\ajphp-\ajphm}\Big(\mU^+_\jph-\mU^-_\jph\Big)+\bm q_\jph,
\end{equation}
and the anti-diffusion term $\bm q_\jph$ is given by \cite{CKX_24}.

Finally, the one-sided local speeds of propagation $a^\pm_\jph$ are slightly different from \eref{2.4a}, and are defined by 
\begin{equation}
a^+_\jph=\max\Big\{u^+_\jph+c^+_\jph,u^-_\jph+c^-_\jph,0\Big\},\quad a^-_\jph=\min\Big\{u^+_\jph-c^+_\jph,u^-_\jph-c^-_\jph,0\Big\}.
\end{equation}

\subsection{One-Dimensional LCDCU Schemes}
According to \cite{CCHKL_22,CCK23_Adaptive}, the 1-D first-, second-, third-, and fifth-order LCDCU schemes can be obtained by  replacing the numerical fluxes $\bm{{\cal F}}^{\rm FV}_\jph$ in \eref{2.4} by 
\begin{equation*}
\resizebox{\linewidth}{!}{$
\bm{{\cal F}}^{\rm FV}_\jph\big(\bm U_\jph^-,\bm U_\jph^+\big)=R_\jph P^{\rm LCD}_\jph R^{-1}_\jph\mF^-_\jph+R_\jph M^{\rm LCD}_\jph R^{-1}_\jph\mF^+_\jph+
R_\jph Q^{\rm LCD}_\jph R^{-1}_\jph\big(\mU^+_\jph-\mU^-_\jph\big).
$}
\end{equation*}
The diagonal matrices $P^{\rm LCD}_\jph$, $M^{\rm LCD}_\jph$, and $Q^{\rm LCD}_\jph$ are given by
\begin{equation*}
\begin{aligned}
&P^{\rm LCD}_\jph={\rm diag}\big((P^{\rm LCD}_1)_\jph,\ldots,(P^{\rm LCD}_d)_\jph\big),\quad
M^{\rm LCD}_\jph={\rm diag}\big((M^{\rm LCD}_1)_\jph,\ldots,(M^{\rm LCD}_d)_\jph\big),\\
&Q^{\rm LCD}_\jph={\rm diag}\big((Q^{\rm LCD}_1)_\jph,\ldots,(Q^{\rm LCD}_d)_\jph\big)
\end{aligned}
\end{equation*}
with
\begin{align*}
&\hspace*{-0.2cm}\big((P^{\rm LCD}_i)_\jph,(M^{\rm LCD}_i)_\jph,(Q^{\rm LCD}_i)_\jph\big)\label{2.7}\\
&=\left\{\begin{aligned}
&\frac{1}{(\lambda^+_i)_\jph-(\lambda^-_i)_\jph}\Big((\lambda^+_i)_\jph,-(\lambda^-_i)_\jph,(\lambda^+_i)_\jph(\lambda^-_i)_\jph\Big)&&
\mbox{if}~(\lambda^+_i)_\jph-(\lambda^-_i)_\jph> \varepsilon_0,\\[0.5ex]
&\Big(\hf,\hf,0\Big)&& \mbox{otherwise},
\end{aligned}\right.\nonumber
\end{align*}
where
\begin{equation*}
(\lambda^+_i)_\jph=\max\left\{\lambda_i\big(A(\mU^-_\jph)\big),\,\lambda_i\big(A(\mU^+_\jph)\big),\,0\right\}, \,\,\,
(\lambda^-_i)_\jph=\min\left\{\lambda_i\big(A(\mU^-_\jph)\big),\,\lambda_i\big(A(\mU^+_\jph)\big),\,0\right\},
\end{equation*}
for $i=1,\ldots,d$. Here, $A(\mU)=\frac{\partial\mF}{\partial\mU}(\mU)$ is the Jacobian, $\lambda_1(A(\mU)) \le \ldots \lambda_d(A(\mU))$ are its eigenvalues, and the matrices $R_\jph$ and $R^{-1}_\jph$ are such that
$R^{-1}_\jph\widehat A_\jph R_\jph$ is diagonal, where $\widehat A_\jph=A(\widehat\mU_\jph)$ and $\widehat\mU_\jph$ is either a simple
average $(\xbar\mU_j+\xbar\mU_{j+1})/2$ or another type of average of the $\xbar\mU_j$ and $\xbar\mU_{j+1}$ states. Finally, $\varepsilon_0$  is a very small desingularization constant, which is taken $\varepsilon_0=10^{-16}$ in all of the
numerical examples reported in \S\ref{sec4}.

\section{Two-Dimensional Schemes}\label{sec3}
We now present a 2-D overview of the HLL, HLLC, TV, LDCU, and LCDCU schemes for the 2-D Euler equations of gas dynamics, which read as 
\begin{equation}
\mU_t+\mF(\mU)_x+\mG(\mU)_y=\bm0,
\label{3.1}
\end{equation}
with $\bm U:=(\rho,\rho u,\rho v,E)^\top$, $\bm F(\bm U)=(\rho u,\rho u^2+p,\rho uv,u(E+p))^\top$, and $\bm G(\bm U)=(\rho v,\rho uv,\rho v^2+p,v(E+p))^\top$. Here, $v$ is the $y$-velocity and the rest of the notation is the same as in the 1-D case \eref{2.1}--\eref{2.2}. The system is completed through the following
EOS:
\begin{equation}\label{3.2}
p=(\gamma-1)\Big[E-\frac{\rho}{2}(u^2+v^2)\Big].
\end{equation}

\subsection{Two-Dimensional HLL Schemes}
As in the 1-D case, we first briefly review the first-order 2-D HLL scheme from \cite{HLL1983} and show its high-order extensions.

\subsubsection{Two-Dimensional First-Order HLL Scheme}
Supposing that the computational domain is covered with uniform cells $C_{j,\,k}:=[x_\jmh,x_\jph]\times[y_\kmh,y_\kph]$ centered at
$(x_j,y_k)=\big((x_\jmh+x_\jph)/2$, $(y_\kph+y_\kmh)/2\big)$ with $x_\jph-x_\jmh\equiv\dx$ and $y_\kph-y_\kmh\equiv\dy$ for all $j,k$, we 
assume that the cell averages 
\begin{equation*}
  \xbar\mU_{j,k}(t):\approx\frac{1}{\dx \dy}\int\limits_{C_{j,k}}\mU(x,y,t)\,{\rm d}x{\rm d}y
\end{equation*}
are available at a certain time level $t$. The cell averages $\xbar \mU_{j,k}$ are then evolved in time by numerically solving the following system of ODEs:
\begin{equation}\label{3.3}
\begin{aligned}
\frac{{\rm d}\xbar \mU_{j,k}}{{\rm d}t}=-\frac{\bm{{\cal F}}^{\rm FV}_{\jph,k}-\bm{{\cal F}}^{\rm FV}_{\jmh,k}}{\dx}-\frac{\bm{{\cal G}}^{\rm FV}_{j,\kph}-\bm{{\cal G}}^{\rm FV}_{j,\kmh}}{\dy}.
\end{aligned}
\end{equation}
Here, $\bm{{\cal F}}^{\rm FV}_{\jph,k}\big(\bm U_{\jph,k}^-,\bm U_{\jph,k}^+\big)$ and $\bm{{\cal G}}^{\rm FV}_{j,\kph}\big(\bm U_{j,\kph}^-,\bm U_{j,\kph}^+\big)$ are the numerical fluxes, defined by
\begin{align}\label{3.4}
\resizebox{0.93\linewidth}{!}{$
\bm{{\cal F}}^{\rm FV}(\mU^-_{\jph,k},\mU^+_{\jph,k}) =
\begin{cases}
\mF(\mU^-_{\jph,k}), & a^-_{\jph,k} \geq 0, \\[8pt]
\dfrac{a^+_{\jph,k} \mF(\mU^-_{\jph,k}) - a^-_{\jph,k} \mF(\mU^+_{\jph,k}) +a^+_{\jph,k} a^-_{\jph,k} (\mU^+_{\jph,k} -\mU^-_{\jph,k})}{a^+_{\jph,k} - a^-_{\jph,k}}, & a^-_{\jph,k} \leq 0 \leq a^+_{\jph,k}, \\[12pt]
\mF(\mU^+_{\jph,k}), & a^+_{\jph,k} \leq 0,
\end{cases}
$}\\[2.ex]
\resizebox{0.93\linewidth}{!}{$ \label{3.5}
\bm{{\cal G}}^{\rm FV}(\mU^-_{j,\kph},\mU^+_{j,\kph}) =
\begin{cases}
\mG(\mU^-_{j,\kph}), & b^-_{j,\kph} \geq 0, \\[8pt]
\dfrac{b^+_{j,\kph} \mG(\mU^-_{j,\kph}) - b^-_{j,\kph} \mG(\mU^+_{j,\kph}) +b^+_{j,\kph} b^-_{j,\kph} (\mU^+_{j,\kph} -\mU^-_{j,\kph})}{b^+_{j,\kph} - b^-_{j,\kph}}, & b^-_{j,\kph} \leq 0 \leq b^+_{j,\kph}, \\[12pt]
\mG(\mU^+_{j,\kph}), & b^+_{j,\kph} \leq 0,
\end{cases}
$}
\end{align}
where $\mU^\pm_{\jph,k}$ and $\mU^\pm_{j,\kph}$ are the left/right-sided point values of $\mU$ at the cell interfaces $(x_\jph,y_k)$ and $(x_j,y_\kph)$, respectively. In the first-order scheme,  we take $\mU^+_{\jph,k}=\mU_{j+1,k}$, $\mU^+_{j,\kph}=\mU_{j,k+1}$, and $\mU^-_{\jph,k}=\mU^-_{j,\kph}=\mU_{j,k}$. Here, the one-sided local speeds of propagation  $a^\pm_{\jph,k}$ and $b^\pm_{j,\kph}$ can be estimated in a ``dimension-by-dimension" manner to the $x$- and $y$-directions respectively as in the 1-D case. Here, we omit the details for the sake of brevity.

\subsubsection{Two-Dimensional Second-Order HLL Scheme}
As in the 1-D case, the resulting scheme \eref{3.3}--\eref{3.5} is second-order accurate once the one-sided point values $\mU^\pm_{\jph,k}$ and $\mU^\pm_{j,\kph}$ employed to compute the numerical fluxes \eref{3.4}--\eref{3.5} are second order. To this end, we approximate $\mU^\pm_{\jph,k}$ and $\mU^\pm_{j,\kph}$ by the generalized minmod reconstruction \cite{lie03,Nessyahu90,Sweby84} with $\theta=1.3$ in the $x$- and $y$- directions, respectively.

\subsubsection{Two-Dimensional Third-Order HLL Scheme}
Following \cite{Jiang13}, the point values $\mU_{j,k}$ are evolved in time by solving the following system of ODEs:
\begin{equation}\label{3.6}
\frac{{\rm d}\mU_{j,k}}{{\rm d}t}=-\frac{{\mH_{\jph,k}}-{\mH_{\jmh,k}}}{\dx}-\frac{{\mH_{j,\kph}}-{\mH_{j,\kmh}}}{\dy},
\end{equation}
where the numerical fluxes $\mH_{\jph,k}$ and $\mH_{j,\kph}$ are defined by
\begin{equation*}
{\mH_{\jph,k}}=\bm{{\cal F}}^{\rm FV}_{\jph,k}-\frac{1}{24}(\dx)^2(\mF_{xx})_{\jph,k}, \quad  {\mH_{j,\kph}}=\bm{{\cal G}}^{\rm FV}_{j,\kph}-\frac{1}{24}(\dy)^2(\mG_{yy})_{j,\kph}.
\end{equation*}
Here, $\bm{{\cal F}}^{\rm FV}_{\jph,k}$ and ${\mG^{\rm FV}_{j,\kph}}$ are the FV numerical fluxes as in \eref{3.4}--\eref{3.5}, $(\mF_{xx})_{\jph,k}$ and $(\mG_{yy})_{j,\kph}$ are the higher-order correction terms computed by the numerical fluxes
\begin{equation}\label{3.7}
\begin{aligned}
(\mF_{xx})_{\jph,k}&=\frac{1}{(\dx)^2}\big[\bm{{\cal F}}^{\rm FV}_{\jmh,k}-2\bm{{\cal F}}^{\rm FV}_{\jph,k}+\bm{{\cal F}}^{\rm FV}_{j+\frac{3}{2},k}\big],\\
(\mG_{yy})_{j,\kph}&=\frac{1}{(\dy)^2}\big[\bm{{\cal G}}^{\rm FV}_{j,\kmh}-2\bm{{\cal G}}^{\rm FV}_{j,\kph}+\bm{{\cal G}}^{\rm FV}_{j,k+\frac{3}{2}}\big].
\end{aligned}
\end{equation}
To ensure the resulting scheme is third-order accurate, the one-sided point values $\mU^\pm_{\jph,k}$ and $\mU^\pm_{j,\kph}$ are also computed using third-order WENO-type interpolation applied to the local characteristic variables. Note that this can be done in a ``dimension-by-dimension" manner as in the 1-D case;  we therefore omit the details for the sake of brevity.  

\subsubsection{Two-Dimensional Fifth-Order HLL Scheme}
According to \cite{Jiang13}, the point values $\mU_{j,k}$ are evolved in time by solving the system of ODEs \eref{3.6} with the following numerical fluxes ${\mH_{\jph,k}}$ and ${\mH_{j,\kph}}$:
\begin{equation*}
\begin{aligned}
& {\mH_{\jph,k}}=\bm{{\cal F}}^{\rm FV}_{\jph,k}-\frac{1}{24}(\dx)^2(\mF_{xx})_{\jph,k}+\frac{7}{5760}(\dx)^4(\mF_{xxxx})_{\jph,k},\\
& {\mH_{j,\kph}}=\bm{{\cal G}}^{\rm FV}_{j,\kph}-\frac{1}{24}(\dy)^2(\mG_{yy})_{j,\kph}+\frac{7}{5760}(\dy)^4(\mG_{yyyy})_{j,\kph}.
\end{aligned}
\end{equation*}
Here, $\bm{{\cal F}}^{\rm FV}_{\jph,k}$ and ${\mG^{\rm FV}_{j,\kph}}$ are the FV numerical fluxes as in \eref{3.4}--\eref{3.5}, $(\mF_{xx})_{\jph,k}$,
$(\mF_{xxxx})_{\jph,k}$, $(\mG_{yy})_{j,\kph}$, $(\mG_{yyyy})_{j,\kph}$ are approximations of the second- and fourth-order spatial derivatives of $\mF$ at $(x,y)=(x_\jph,y_k)$ and $\mG$ at $(x,y)=(x_j,y_\kph)$, respectively. The higher-order correction terms are obtained from the one-dimensional formulas in \eref{2.11aa} by applying the same finite-difference approximations dimension by dimension.

To achieve fifth-order accuracy, the one-sided point values $\mU^\pm_{\jph,k}$ and  $\mU^\pm_{j,\kph}$ employed to compute the numerical flux $\bm{{\cal F}}^{\rm FV}_{\jph,k}$ and $\bm{{\cal G}}^{\rm FV}_{j,\kph}$ need to be at least fifth-order accurate. This can also be done in a ``dimension-by-dimension" manner as in the 1-D case; we therefore omit the details for the sake of brevity.

\subsection{Two-Dimensional HLLC Schemes}
 According to \cite{TSS1994}, the two-dimensional first-, second-, third-, and fifth-order HLLC schemes are obtained by replacing the numerical fluxes in \eref{3.4}--\eref{3.5} with the corresponding HLLC fluxes. The $x$-direction numerical flux $\bm{{\cal F}}^{\rm FV}_{\jph,k}$ is given by

\begin{align*}
\bm{{\cal F}}^{\rm FV}_{\jph,k}(\mU^-_{\jph,k},\mU^+_{\jph,k}) =
\begin{cases}
\mF(\mU^-_{\jph,k}), & a^-_{\jph,k} \geq 0, \\[8pt]
\mF(\mU^-_{\jph,k}) + a^-_{\jph,k} \big(\mU_{\jph,k}^{*,-} - \mU^-_{\jph,k} \big), & a^-_{\jph,k} \leq 0 \leq a^*_{\jph,k}, \\[8pt]
\mF(\mU^+_{\jph,k}) + a^+_{\jph,k} \big(\mU_{\jph,k}^{*,+} - \mU^+_{\jph,k} \big), & a^*_{\jph,k} \leq 0 \leq a^+_{\jph,k}, \\[8pt]
\mF(\mU^+_{\jph,k}), & a^+_{\jph,k} \leq 0 ,
\end{cases}
\end{align*}
where 
\begin{equation}\resizebox{0.93\linewidth}{!}{$
\mU_{\jph,k}^{*,\pm}
=
\rho^\pm_{\jph,k}
\frac{a^\pm_{\jph,k}-u^\pm_{\jph,k}}
     {a^\pm_{\jph,k}-a^*_{\jph,k}}
\begin{pmatrix}
1 \\[4pt]
a^*_{\jph,k} \\[4pt]
v^\pm_{\jph,k} \\[4pt]
\displaystyle
\frac{E^\pm_{\jph,k}}{\rho^\pm_{\jph,k}}
+
\big(a^*_{\jph,k}-u^\pm_{\jph,k}\big)
\left(
a^*_{\jph,k}
+
\frac{p^\pm_{\jph,k}}
{\rho^\pm_{\jph,k}
 \big(a^\pm_{\jph,k}-u^\pm_{\jph,k}\big)}
\right)
\end{pmatrix},\\[2.ex]
$}
\end{equation}
and 
$$
a^*_{\jph,k}= \frac{p^+_{\jph,k} -p^-_{\jph,k} + \rho^-_{\jph,k} u^-_{\jph,k} (a^-_{\jph,k}-u^-_{\jph,k}) - \rho^+_{\jph,k} u^+_{\jph,k} (a^+_{\jph,k}-u^+_{\jph,k})}{\rho^-_{\jph,k} (a^-_{\jph,k}-u^-_{\jph,k})-\rho^+_{\jph,k} (a^+_{\jph,k}-u^+_{\jph,k})}.
$$
The corresponding $y$-direction numerical flux $\bm{{\cal G}}^{\rm FV}_{j,\kph}$ is obtained analogously by interchanging the roles of the $x$- and $y$-directions. We therefore omit the details here.

\subsection{Two-Dimensional TV Splitting Schemes}
According to \cite{CHT25,TCL2015}, the 2-D first-, second-, third-, and fifth-order TV schemes can be obtained by  replacing the numerical fluxes $\bm{{\cal F}}^{\rm FV}_{\jph,k}$ and $\bm{{\cal G}}^{\rm FV}_{j,\kph}$ in \eref{3.4}--\eref{3.5} by 
\begin{equation*}
\begin{aligned}
&\bm{{\cal F}}^{\rm FV}_{\jph,k}\big(\bm U_{\jph,k}^-,\bm U_{\jph,k}^+\big)= \bm{{\cal F}}^A_{\jph,k}\big(\bm U_{\jph,k}^-,\bm U_{\jph,k}^+\big) + \bm{{\cal F}}^P_{\jph,k} \big(\bm U_{\jph,k}^-,\bm U_{\jph,k}^+\big),\\
&\bm{{\cal G}}^{\rm FV}_{j,\kph}\big(\bm U_{j,\kph}^-,\bm U_{j,\kph}^+\big)= \bm{{\cal G}}^A_{j,\kph}\big(\bm U_{j,\kph}^-,\bm U_{j,\kph}^+\big) + \bm{{\cal G}}^P_{j,\kph} \big(\bm U_{j,\kph}^-,\bm U_{j,\kph}^+\big).
\end{aligned}
\end{equation*}

Here, $\bm{{\cal F}}^A_{\jph,k}\big(\bm U_{\jph,k}^-,\bm U_{\jph,k}^+\big)$ is the $x$-direction advection flux given by 
\begin{equation*}
  \bm{{\cal F}}^A_{\jph,k}\big(\bm U_{\jph,k}^-,\bm U_{\jph,k}^+\big)=\begin{cases}
                                                          u^*_{\jph,k}\begin{pmatrix} \rho^-_{\jph,k} \\[0.8ex]
                                                           (\rho u)^-_{\jph,k}\\[0.8ex]
                                                           (\rho v)^-_{\jph,k}\\[0.8ex]
                                                          \dfrac{1}{2} \rho^-_\jph \Big[ (u^-_{\jph,k})^2+(v^-_{\jph,k})^2 \Big]\end{pmatrix}, & \mbox{if } u^*_{\jph,k} \ge 0, \\[12.ex]
                                                          u^*_{\jph,k} \begin{pmatrix} \rho^+_{\jph,k} \\[0.8ex]
                                                           (\rho u)^+_{\jph,k}\\[0.8ex]
                                                           (\rho v)^+_{\jph,k}\\[0.8ex]
                                                          \dfrac{1}{2} \rho^+_\jph \Big[ (u^+_{\jph,k})^2+(v^+_{\jph,k})^2 \Big]\end{pmatrix}, & \mbox{otherwise},
                                                  \end{cases}
\end{equation*}
and $\bm{{\cal F}}^P_{\jph,k} \big(\bm U_{\jph,k}^-,\bm U_{\jph,k}^+\big)$ is the $x$-direction pressure flux given by 
\begin{equation*}
  \bm{{\cal F}}^P_{\jph,k} \big(\bm U_{\jph,k}^-,\bm U_{\jph,k}^+\big)= \begin{pmatrix} 0,\,  p^*_{\jph,k},\, 0,\, \dfrac{\gamma  u^*_{\jph,k}\,  p^*_{\jph,k}}{\gamma-1} \end{pmatrix}^\top,
\end{equation*}
where 
\begin{equation*}
 \begin{aligned}
  u^*_{\jph,k} &= \frac{C^+_{\jph,k} u^+_{\jph,k} - C^-_{\jph,k} u^-_{\jph,k}} {C^+_{\jph,k}-C^-_{\jph,k}}-\frac{2}{C^+_{\jph,k}-C^-_{\jph,k}} \big(p^+_{\jph,k} - p^-_{\jph,k}\big), \\
  p^*_{\jph,k} &= \frac{C^+_{\jph,k} p^-_{\jph,k} - C^-_{\jph,k} p^+_{\jph,k}} {C^+_{\jph,k}-C^-_{\jph,k}}+\frac{C^+_{\jph,k} C^-_{\jph,k}}{ 2(C^+_{\jph,k}-C^-_{\jph,k})} \big(u^+_{\jph,k} - u^-_{\jph,k}\big),\\
   C^\pm_{\jph,k}&=\rho^\pm_{\jph,k} \bigg(u^\pm_{\jph,k} \pm \sqrt{\big(u^\pm_{\jph,k} \big)^2+4 \big(c^\pm_{\jph,k}\big)^2}\bigg).
 \end{aligned}
\end{equation*}
Similarly, the detailed formulas for the $y$-direction advection and pressure fluxes $\bm{{\cal G}}^{A}_{j,\kph}\big(\bm U_{j,\kph}^-,\bm U_{j,\kph}^+\big)$ and $\bm{{\cal G}}^{P}_{j,\kph}\big(\bm U_{j,\kph}^-,\bm U_{j,\kph}^+\big)$ are obtained by interchanging the roles of the $x$- and $y$-directions. We therefore omit them here.

\subsection{Two-Dimensional LDCU Schemes}
According to \cite{KX_22,CKX_24,CKX22}, the 2-D first-, second-, third-, and fifth-order LDCU schemes can be obtained by  replacing the numerical fluxes $\bm{{\cal F}}^{\rm FV}_{\jph,k}$ and $\bm{{\cal G}}^{\rm FV}_{j,\kph}$ in \eref{3.4}--\eref{3.5} by 
\begin{equation*}
\begin{aligned}
\bm{{\cal F}}^{\rm FV}_{\jph,k}\big(\bm U_{\jph,k}^-,\bm U_{\jph,k}^+\big)&=\frac{a^+_{\jph,k}\mF\big(\bm U^-_{\jph,k}\big)-a^-_{\jph,k}\mF\big(\bm U^+_{\jph,k}\big)}
{a^+_{\jph,k}-a^-_{\jph,k}}\\
&+\frac{a^+_{\jph,k}a^-_{\jph,k}}{a^+_{\jph,k}-a^-_{\jph,k}}\Big(\mU^+_{\jph,k}-\mU^-_{\jph,k}\Big)+
\bm q^x_{\jph,k},\\
\bm{{\cal G}}^{\rm FV}_{j,\kph}\big(\bm U_{j,\kph}^-,\bm U_{j,\kph}^+\big)&=\frac{b^+_{j,\kph}\mG\big(\bm U^-_{j,\kph}\big)-b^-_{j,\kph}\mG\big(\bm U^+_{j,\kph}\big)}
{b^+_{j,\kph}-b^-_{j,\kph}}\\
&+\frac{b^+_{j,\kph}b^-_{j,\kph}}{b^+_{j,\kph}-b^-_{j,\kph}}\Big(\mU^+_{j,\kph}-\mU^-_{j,\kph}\Big)+
\bm q^y_{j,\kph},
\end{aligned}
\end{equation*}
where $\bm q^x_{\jph,k}$ and $\bm q^y_{j,\kph}$ are the $x$- and $y$-direction anti-diffusion terms, which are given by \cite{CKX_24}. The one-sided local speeds of propagation $a^\pm_{\jph,k}$ and $b^\pm_{j,\kph}$ are defined by 
\begin{equation}
\resizebox{0.9\linewidth}{!}{$
\begin{aligned}
&a^+_{\jph,k}=\max\Big\{u^+_{\jph,k}+c^+_{\jph,k},u^-_{\jph,k}+c^-_{\jph,k},0\Big\},~\,
a^-_{\jph,k}=\min\Big\{u^+_{\jph,k}-c^+_{\jph,k},u^-_{\jph,k}-c^-_{\jph,k},0\Big\},\\
&b^+_{j,\kph}=\max\Big\{v^+_{j,\kph}+c^+_{j,\kph},v^-_{j,\kph}+c^-_{j,\kph},0\Big\},\quad
b^-_{j,\kph}=\min\Big\{v^+_{j,\kph}-c^+_{j,\kph},v^-_{j,\kph}-c^-_{j,\kph},0\Big\}.\\
\end{aligned}
$}
\end{equation}

\subsection{Two-Dimensional LCDCU Schemes}\label{sec3.5}
According to \cite{CCHKL_22,CCK23_Adaptive}, the 2-D first-, second-, third-, and fifth-order LCDCU schemes can be obtained by  replacing the numerical fluxes $\bm{{\cal F}}^{\rm FV}_{\jph,k}$ and $\bm{{\cal G}}^{\rm FV}_{j,\kph}$ in \eref{3.4}--\eref{3.5} by 
\begin{align*}
\resizebox{\linewidth}{!}{$
\bm{{\cal F}}^{\rm FV}_{\jph,k}\big(\bm U_{\jph,k}^-,\bm U_{\jph,k}^+\big)=R_\jph P^{\rm LCD}_{\jph,k} R^{-1}_{\jph,k}\mF^-_{\jph,k}+R_{\jph,k} M^{\rm LCD}_{\jph,k} R^{-1}_{\jph,k}\mF^+_{\jph,k}+
R_{\jph,k} Q^{\rm LCD}_{\jph,k} R^{-1}_{\jph,k}\big(\mU^+_{\jph,k}-\mU^-_{\jph,k}\big),
$}\\[1.ex]
\resizebox{\linewidth}{!}{$
\bm{{\cal G}}^{\rm FV}_{j,\kph}\big(\bm U_{j,\kph}^-,\bm U_{j,\kph}^+\big)=R_{j,\kph} P^{\rm LCD}_{j,\kph} R^{-1}_{j,\kph}\mG^-_{j,\kph}+R_{j,\kph} M^{\rm LCD}_{j,\kph} R^{-1}_{j,\kph}\mG^+_{j,\kph}+
R_{j,\kph} Q^{\rm LCD}_{j,\kph} R^{-1}_{j,\kph}\big(\mU^+_{j,\kph}-\mU^-_{j,\kph}\big).
$}
\end{align*}
The matrices $R_{\jph,k}$, $R^{-1}_{\jph,k}$ and $R_{j,\kph}$, $R^{-1}_{j,\kph}$ are the matrices such that $R^{-1}_{\jph,k}\widehat A_{\jph,k}R_{\jph,k}$ and $R^{-1}_{j,\kph}\widehat B_{j,\kph}R_{j,\kph}$ are diagonal. Here, $\widehat A_{\jph,k}=A(\widehat\mU_{\jph,k})$, $\widehat B_{j,\kph}=B(\widehat\mU_{j,\kph})$ with $A(\mU)=\frac{\partial\mF(\mU)}{\partial\mU}$, $B(\mU)=\frac{\partial\mG(\mU)}{\partial\mU}$, and $\widehat\mU_{\jph,k}$, $\widehat\mU_{j,\kph}$ are either simple averages $(\,\xbar\mU_{j,k}+\xbar\mU_{j+1,k})/2$, $(\,\xbar\mU_{j,k}+\xbar\mU_{j,k+1})/2$ or another type of averages of $\xbar\mU_{j,k}$, $\xbar\mU_{j+1,k}$ and $\xbar\mU_{j,k}$, $\xbar\mU_{j,k+1}$ states, respectively.  The diagonal matrices $P^{\rm LCD}_{\jph,k}$, $M^{\rm LCD}_{\jph,k}$, $Q^{\rm LCD}_{\jph,k}$,
$P^{\rm LCD}_{j,\kph}$, $M^{\rm LCD}_{j,\kph}$, and $Q^{\rm LCD}_{j,\kph}$ are defined by
\begin{equation*}
\begin{aligned}
&P^{\rm LCD}_{\jph,k}={\rm diag}\left(\big(P^{\rm LCD}_1\big)_{\jph,k},\ldots,\big(P^{\rm LCD}_d\big)_{\jph,k}\right),&&
P^{\rm LCD}_{j,\kph}={\rm diag}\left(\big(P^{\rm LCD}_1\big)_{j,\kph},\ldots,\big(P^{\rm LCD}_d\big)_{j,\kph}\right),\\
&M^{\rm LCD}_{\jph,k}={\rm diag}\left(\big(M^{\rm LCD}_1\big)_{\jph,k},\ldots,\big(M^{\rm LCD}_d\big)_{\jph,k}\right),&&
M^{\rm LCD}_{j,\kph}={\rm diag}\left(\big(M^{\rm LCD}_1\big)_{j,\kph},\ldots,\big(M^{\rm LCD}_d\big)_{j,\kph}\right),\\
&Q^{\rm LCD}_{\jph,k}={\rm diag}\left(\big(Q^{\rm LCD}_1\big)_{\jph,k},\ldots,\big(Q^{\rm LCD}_d\big)_{\jph,k}\right),&&
Q^{\rm LCD}_{j,\kph}={\rm diag}\left(\big(Q^{\rm LCD}_1\big)_{j,\kph},\ldots,\big(Q^{\rm LCD}_d\big)_{j,\kph}\right),
\end{aligned}
\end{equation*}
where
$$
\begin{aligned}
&\hspace*{-1.0cm}\left((P^{\rm LCD}_i)_{\jph,k},(M^{\rm LCD}_i)_{\jph,k},(Q^{\rm LCD}_i)_{\jph,k}\right)\\
&=\left\{\begin{aligned}
&\frac{1}{\Delta(\lambda_i)_{\jph,k}}
\left((\lambda^+_i)_{\jph,k},-(\lambda^-_i)_{\jph,k},(\lambda^+_i)_{\jph,k}(\lambda^-_i)_{\jph,k}\right)
&&\mbox{if}~\Delta(\lambda_i)_{\jph,k}>\varepsilon_0,\\
&\Big(\frac{1}{2}, \frac{1}{2}, 0\Big)&&\mbox{otherwise},
\end{aligned}\right.\\
&\hspace*{-1.0cm}\left((P^{\rm LCD}_i)_{j,\kph},(M^{\rm LCD}_i)_{j,\kph},(Q^{\rm LCD}_i)_{j,\kph}\right)\\
&=\left\{\begin{aligned}
&\frac{1}{\Delta(\mu_i)_{j,\kph}}\left((\mu^+_i)_{j,\kph},-(\mu^-_i)_{j,\kph},(\mu^+_i)_{j,\kph}(\mu^-_i)_{j,\kph}\right)
&&\mbox{if}~\Delta(\mu_i)_{j,\kph}>\varepsilon_0,\\
&\Big(\frac{1}{2}, \frac{1}{2}, 0\Big)&&\mbox{otherwise}.
\end{aligned}\right.
\end{aligned}
$$
Here, $\Delta(\lambda_i)_{\jph,k}:=(\lambda^+_i)_{\jph,k}-(\lambda^-_i)_{\jph,k}$,
$\Delta(\mu_i)_{j,\kph}:=(\mu^+_i)_{j,\kph}-(\mu^-_i)_{j,\kph}$, and
\begin{equation*}
\begin{aligned}
(\lambda^+_i)_{\jph,k}=\max\big\{\lambda_i\big(A(\mU^-_{\jph,k})\big),\,\lambda_i\big(A(\mU^+_{\jph,k})\big),\,0\big\},\\
(\lambda^-_i)_{\jph,k}=\min\big\{\lambda_i\big(A(\mU^-_{\jph,k})\big),\,\lambda_i\big(A(\mU^+_{\jph,k})\big),\,0\big\},\\
(\mu^+_i)_{j,\kph}=\max\big\{\lambda_i\big(B(\mU^-_{j,\kph})\big),\,\lambda_i\big(B(\mU^+_{j,\kph})\big),\,0\big\},\\
(\mu^-_i)_{j,\kph}=\min\big\{\lambda_i\big(B(\mU^-_{j,\kph})\big),\,\lambda_i\big(B(\mU^+_{j,\kph})\big),\,0\big\},
\end{aligned}
\end{equation*}
where $\lambda_i$ and $\mu_i$ are the eigenvalues of the Jacobians $A(\mU)$ and $B(\mU)$: $\lambda_1(A)\le\ldots\le\lambda_d(A)$ and
$\mu_1(B)\le\ldots\le\mu_d(B)$, respectively.
\section{Numerical Examples} \label{sec4}
In this section, we test the studied first-, second-, third-, and fifth-order schemes on several numerical examples and compare their performance. For the sake of brevity, these schemes will be referred to as the 1-Order, 2-Order, 3-Order, and 5-Order schemes, respectively.

We numerically integrate the ODE systems \eref{2.3}, \eref{2.9}, \eref{3.3}, and \eref{3.6} by the three-stage third-order strong stability preserving Runge-Kutta (SSP RK3) method (see, e.g., \cite{Gottlieb11,Gottlieb12}) and use the CFL number 0.45.

\subsection{One-Dimensional Examples}
We begin with the 1-D Euler equations of gas dynamics \eref{1.1}, \eref{2.1}--\eref{2.2}. In  Examples 1--8, we take the specific heat ratio $\gamma=1.4$.

\subsubsection*{Example 1---1-D Accuracy Test}
In the first example taken from \cite{CHT25}, we consider the system \eref{1.1},\eref{2.1}--\eref{2.2} subject to the following periodic initial
conditions,
\begin{equation*}
\rho(x,0)=1+\frac{1}{5}\sin(2\pi x),\quad u(x,0)\equiv1,\quad p(x,0)\equiv1.
\end{equation*}
The exact solution of this initial value problem is given by
$$
\rho(x,t)=1+\frac{1}{5}\sin\left[2\pi(x-t)\right],\quad u(x,t)\equiv1,\quad p(x,t)\equiv1.
$$

We first compute the numerical solution on the computational domain $[-1,1]$ until the final time $t=0.1$ by the 1-Order, 2-Order, 3-Order, and 5-Order schemes on a sequence of uniform meshes: 100, 200, and 400, measure the $L^1$-errors, and then compute the corresponding experimental convergence rates for the density. The obtained results are presented in Table \ref{tab1},  where one can clearly see that the expected order of accuracy is achieved for the studied schemes. At the same time, one can see that the four low-dissipation schemes are more accurate than the HLL counterparts.

\begin{table}[ht!]
\centering
\begin{tabular}{|c|cc|cc|cc|cc|cc|cc|}
\hline
\multirow{2}{2em}{Mesh}&\multicolumn{2}{c|}{HLL, 1-Order}&\multicolumn{2}{c|}{HLL, 2-Order}&\multicolumn{2}{c|}{HLL, 3-Order}&\multicolumn{2}{c|}{HLL, 5-Order}\\
\cline{2-9}&Error&Rate&Error&Rate&Error&Rate&Error&Rate\\
\hline
$100$&9.91e-03 &---   &1.01e-03 &---  &2.37e-05&--- &7.01e-08&---\\
$200$&4.98e-03 &0.992 &2.46e-04 &2.04 &2.96e-06&3.00&2.20e-09&5.00\\
$400$&2.50e-03 &0.996 &5.98e-05 &2.04 &3.70e-07&3.00&6.86e-11&5.00\\
\hline
\multirow{2}{2em}{Mesh}&\multicolumn{2}{c|}{HLLC, 1-Order}&\multicolumn{2}{c|}{HLLC, 2-Order}&\multicolumn{2}{c|}{HLLC, 3-Order}&\multicolumn{2}{c|}{HLLC, 5-Order}\\
\cline{2-9}&Error&Rate&Error&Rate&Error&Rate&Error&Rate\\
\hline
$100$&7.99e-03 &---   &9.40e-04 &---  &1.99e-05&--- &5.90e-08&---\\
$200$&4.03e-03 &0.990 &2.24e-04 &2.07 &2.49e-06&3.00&1.85e-09&5.00\\
$400$&2.02e-03 &0.995 &5.51e-05 &2.03 &3.12e-07&3.00&5.78e-11&5.00\\
\hline
\multirow{2}{2em}{Mesh}&\multicolumn{2}{c|}{TV, 1-Order}&\multicolumn{2}{c|}{TV, 2-Order}&\multicolumn{2}{c|}{TV, 3-Order}&\multicolumn{2}{c|}{TV, 5-Order}\\
\cline{2-9}&Error&Rate&Error&Rate&Error&Rate&Error&Rate\\
\hline
$100$&7.99e-03 &---   &9.40e-04 &---  &1.99e-05&--- &5.90e-08&---\\
$200$&4.03e-03 &0.990 &2.24e-04 &2.07 &2.49e-06&3.00&1.85e-09&5.00\\
$400$&2.02e-03 &0.995 &5.51e-05 &2.03 &3.12e-07&3.00&5.78e-11&5.00\\
\hline
\multirow{2}{2em}{Mesh}&\multicolumn{2}{c|}{LDCU, 1-Order}&\multicolumn{2}{c|}{LDCU, 2-Order}&\multicolumn{2}{c|}{LDCU, 3-Order}&\multicolumn{2}{c|}{LDCU, 5-Order}\\
\cline{2-9}&Error&Rate&Error&Rate&Error&Rate&Error&Rate\\
\hline
$100$&7.99e-03 &---   &9.40e-04 &---  &1.99e-05&--- &5.90e-08&---\\
$200$&4.03e-03 &0.990 &2.24e-04 &2.07 &2.49e-06&3.00&1.85e-09&5.00\\
$400$&2.02e-03 &0.995 &5.51e-05 &2.03 &3.12e-07&3.00&5.78e-11&5.00\\
\hline
\multirow{2}{2em}{Mesh}&\multicolumn{2}{c|}{LCDCU, 1-Order}&\multicolumn{2}{c|}{LCDCU, 2-Order}&\multicolumn{2}{c|}{LCDCU, 3-Order}&\multicolumn{2}{c|}{LCDCU, 5-Order}\\
\cline{2-9}&Error&Rate&Error&Rate&Error&Rate&Error&Rate\\
\hline
$100$&7.99e-03 &---   &9.40e-04 &---  &1.99e-05&--- &5.90e-08&---\\
$200$&4.03e-03 &0.990 &2.24e-04 &2.07 &2.49e-06&3.00&1.85e-09&5.00\\
$400$&2.02e-03 &0.995 &5.51e-05 &2.03 &3.12e-07&3.00&5.78e-11&5.00\\
\hline 
\end{tabular}
\caption{\sf Example 1: The $L^1$-errors and experimental convergence rates for the density $\rho$ computed by the 1-Order, 2-Order, 3-Order, and 5-Order schemes.\label{tab1}}
\end{table}

\begin{rmk}
We stress that in order to achieve the fifth order of accuracy for the 5-Order scheme, we use smaller time steps with $\dt \sim (\dx)^{\frac{5}{3}}$ to balance the spatial and temporal errors.
\end{rmk}

\subsubsection*{Example 2---Moving Contact Wave} In the second example, we consider the moving contact discontinuity problem from \cite{Kurganov07} with the following initial conditions:
\begin{equation*}
(\rho, u,p)(x,0)=\begin{cases}
(1.4,0.1,1),&x<0.5,\\
(1,0.1,1),  &x>0.5,
\end{cases}
\end{equation*}
which is considered on the interval $[0,1]$ with the free boundary conditions at both ends.

\begin{figure}[ht!]
\centerline{\includegraphics[trim=1.0cm 0.3cm 0.9cm 0.8cm, clip, width=4.cm]{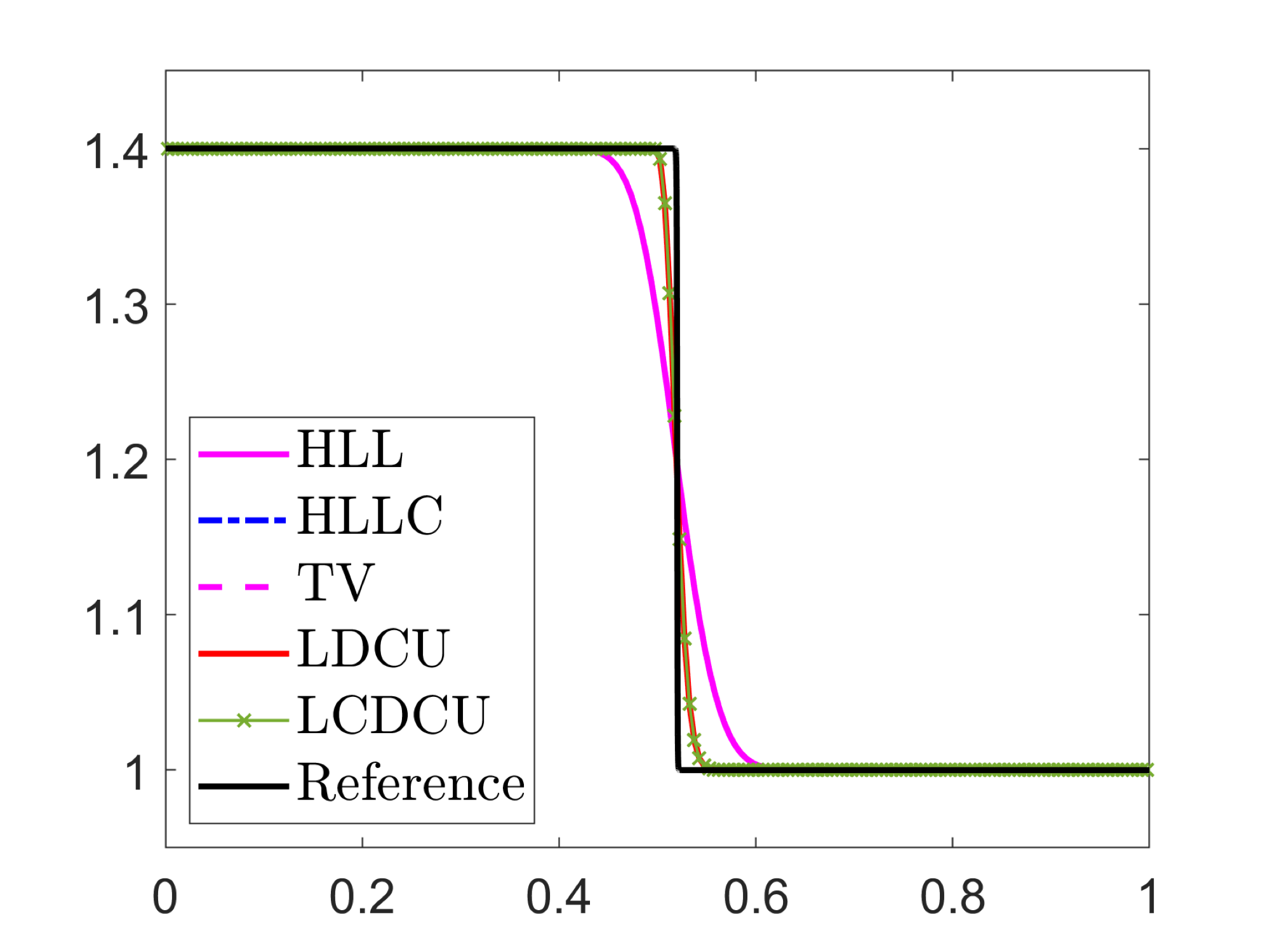}
            \includegraphics[trim=1.0cm 0.3cm 0.9cm 0.8cm, clip, width=4.cm]{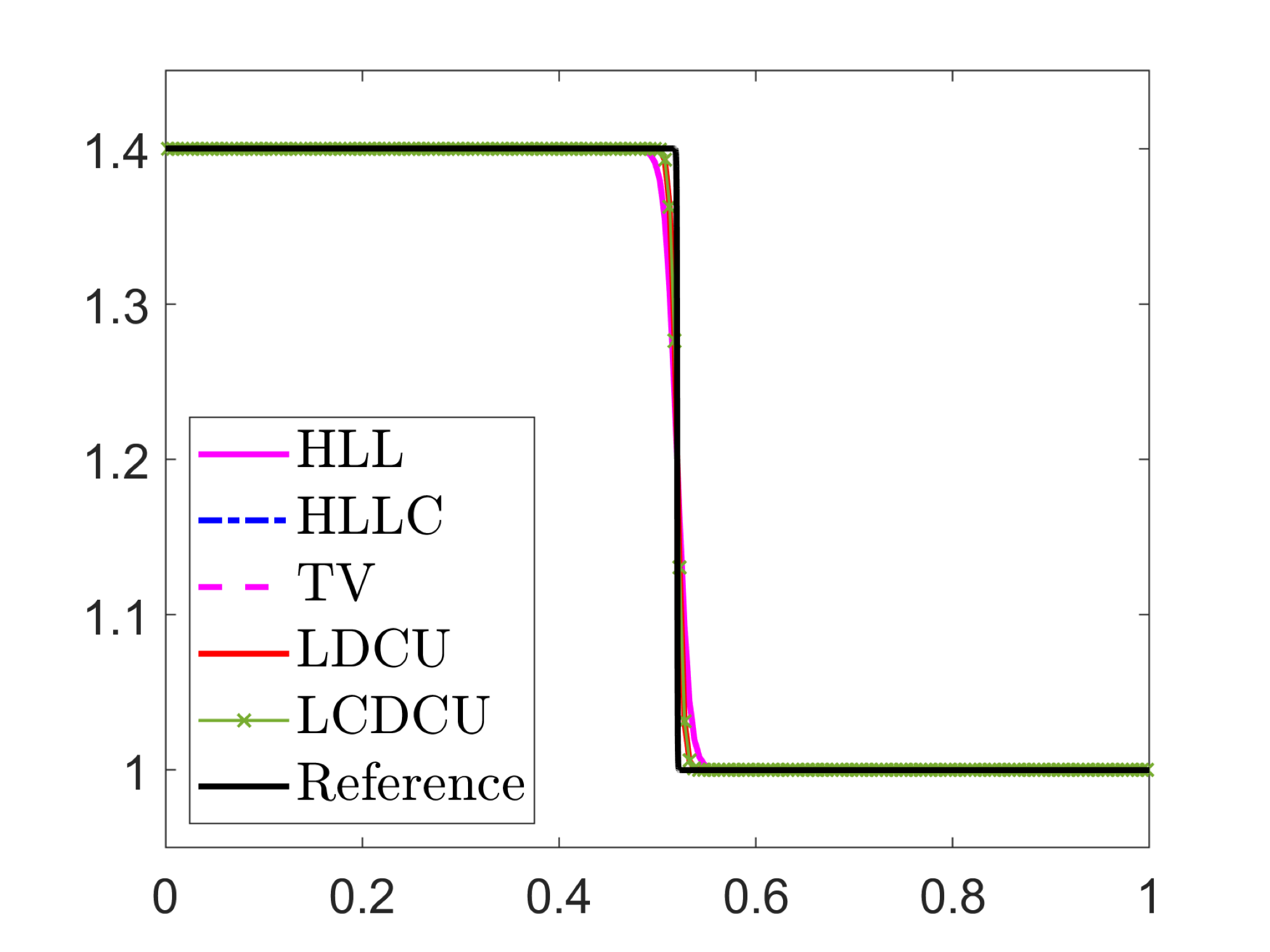}
            \includegraphics[trim=1.0cm 0.3cm 0.9cm 0.8cm, clip, width=4.cm]{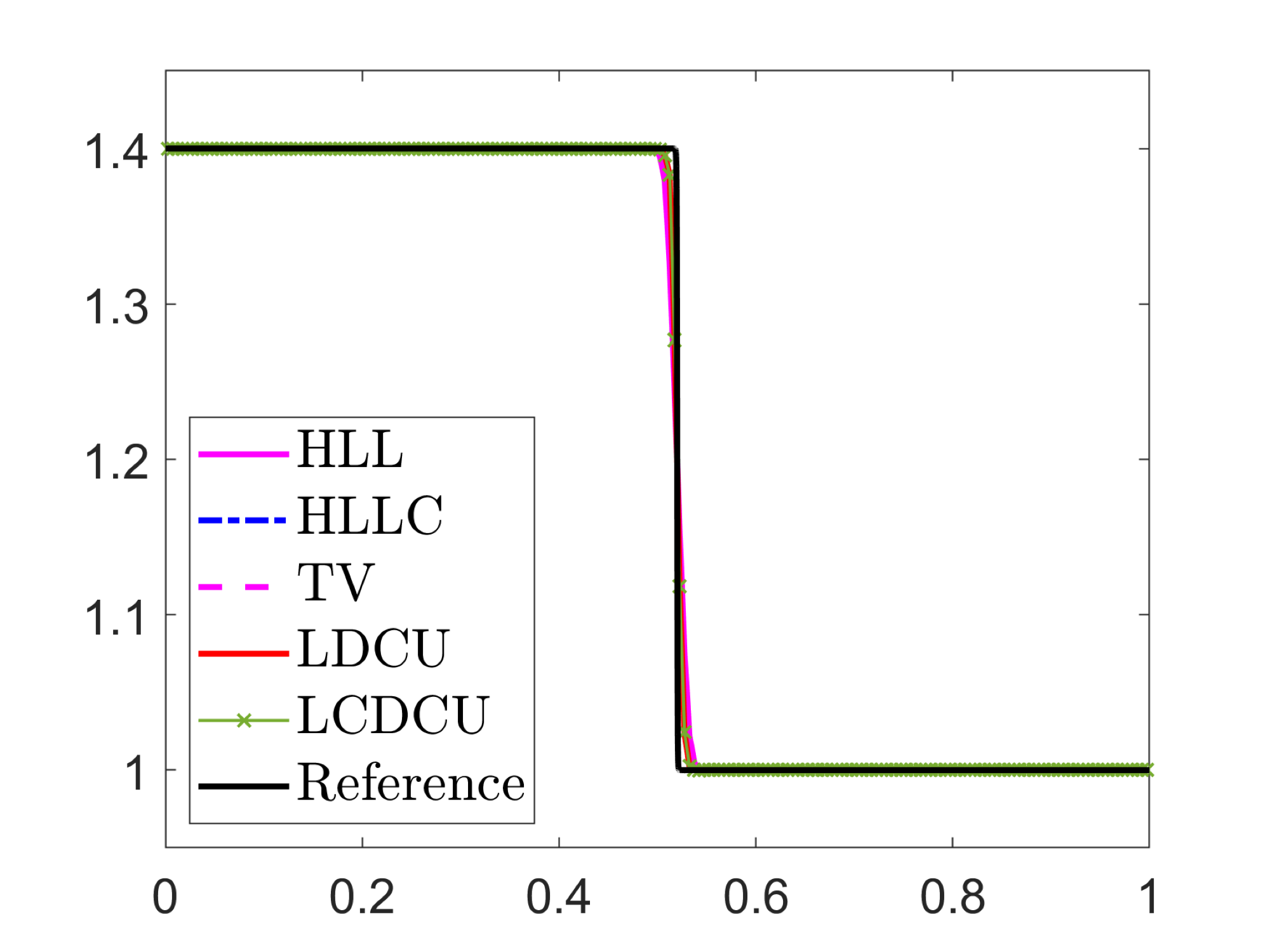}
            \includegraphics[trim=1.0cm 0.3cm 0.9cm 0.8cm, clip, width=4.cm]{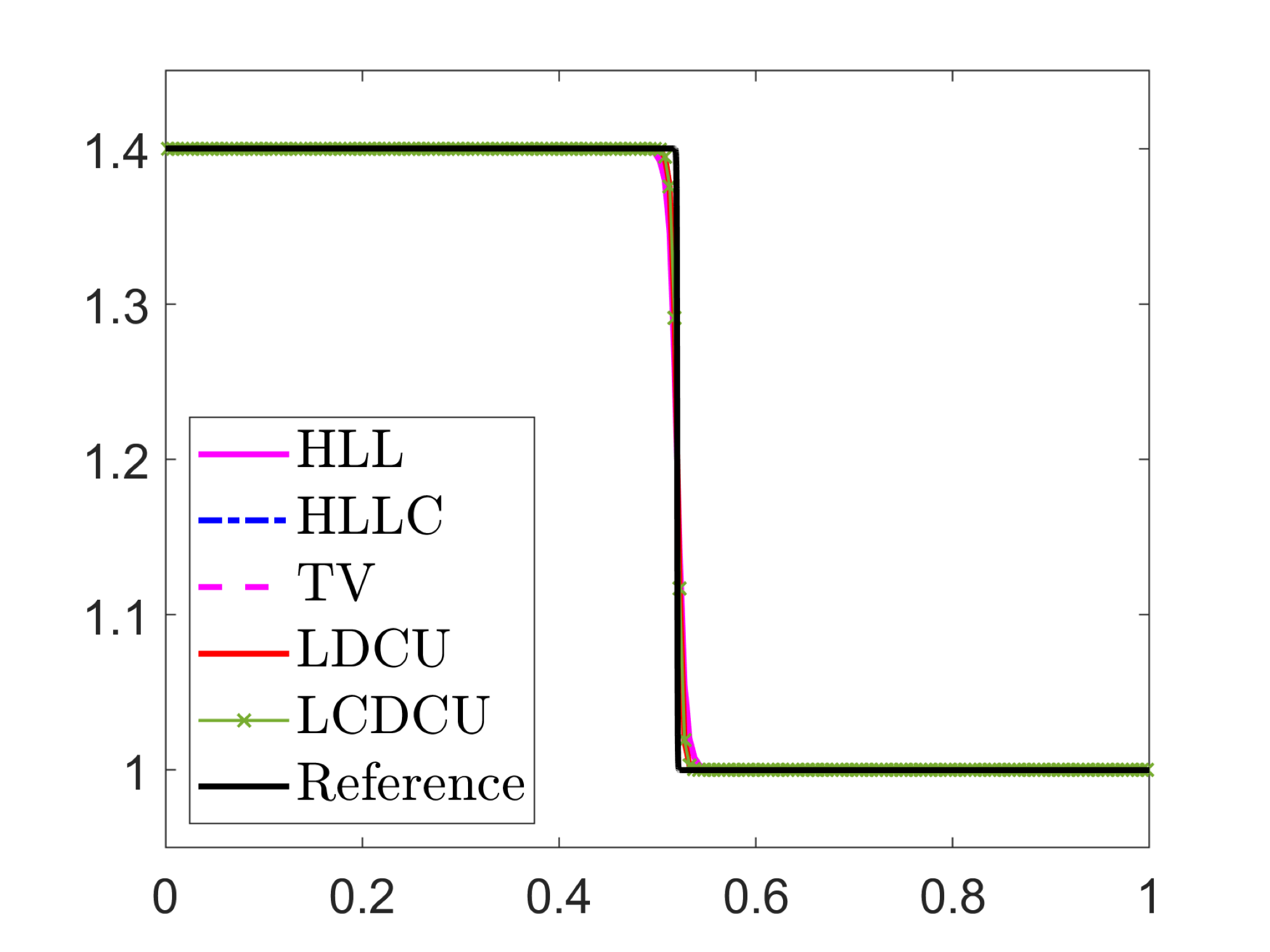}}
\vskip 12pt
\centerline{\includegraphics[trim=1.0cm 0.3cm 0.9cm 0.8cm, clip, width=4.cm]{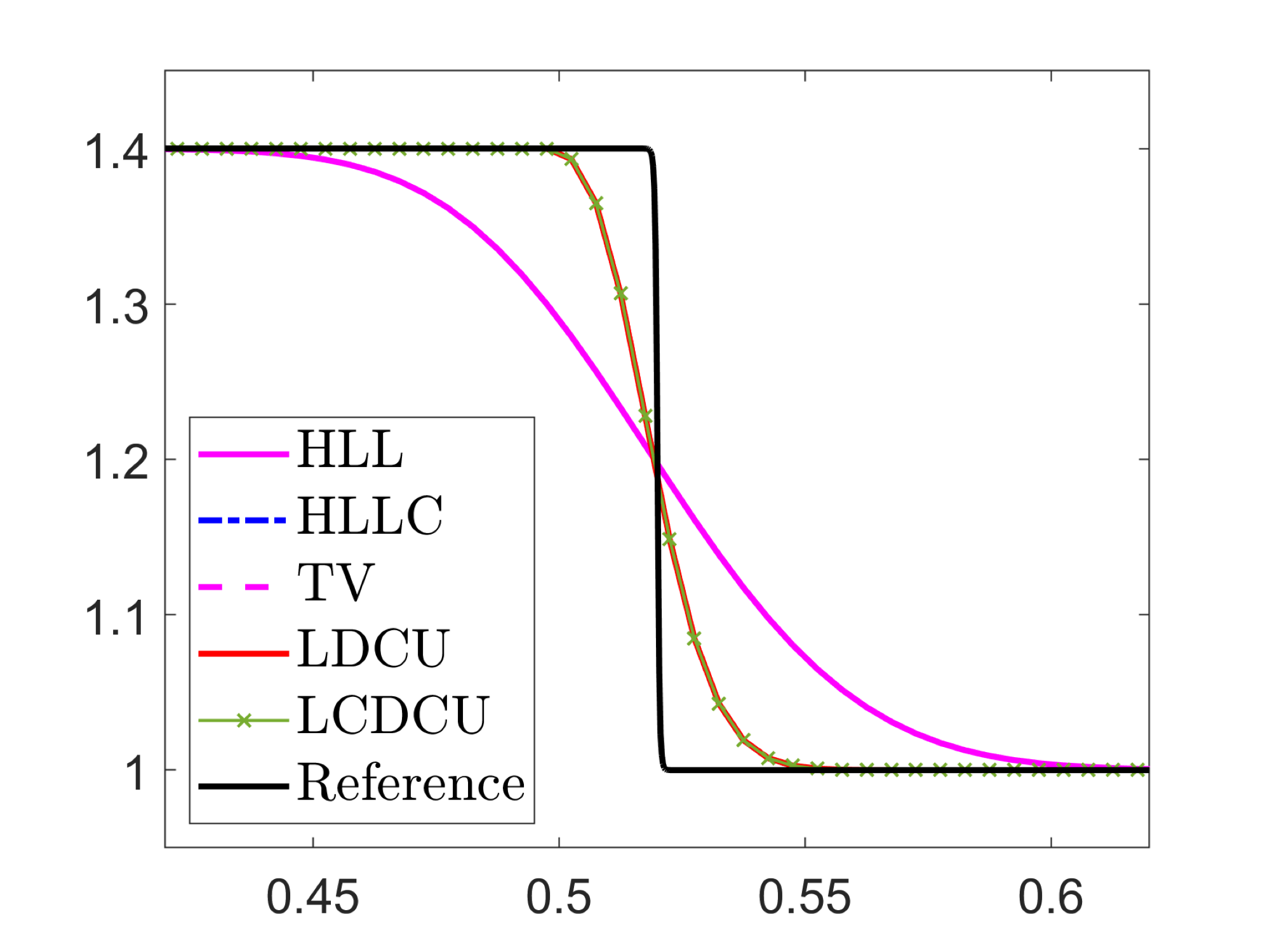}
            \includegraphics[trim=1.0cm 0.3cm 0.9cm 0.8cm, clip, width=4.cm]{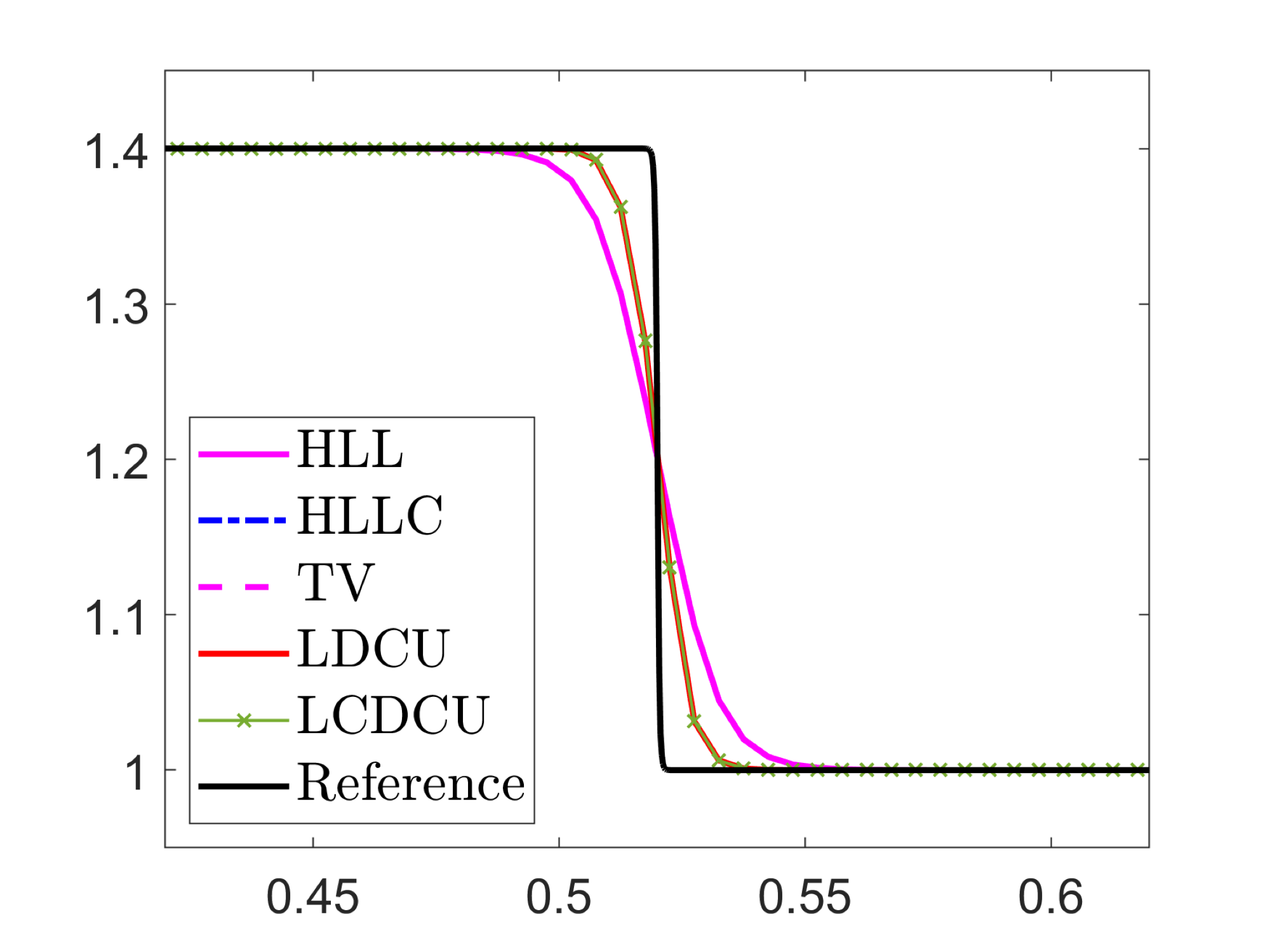}
            \includegraphics[trim=1.0cm 0.3cm 0.9cm 0.8cm, clip, width=4.cm]{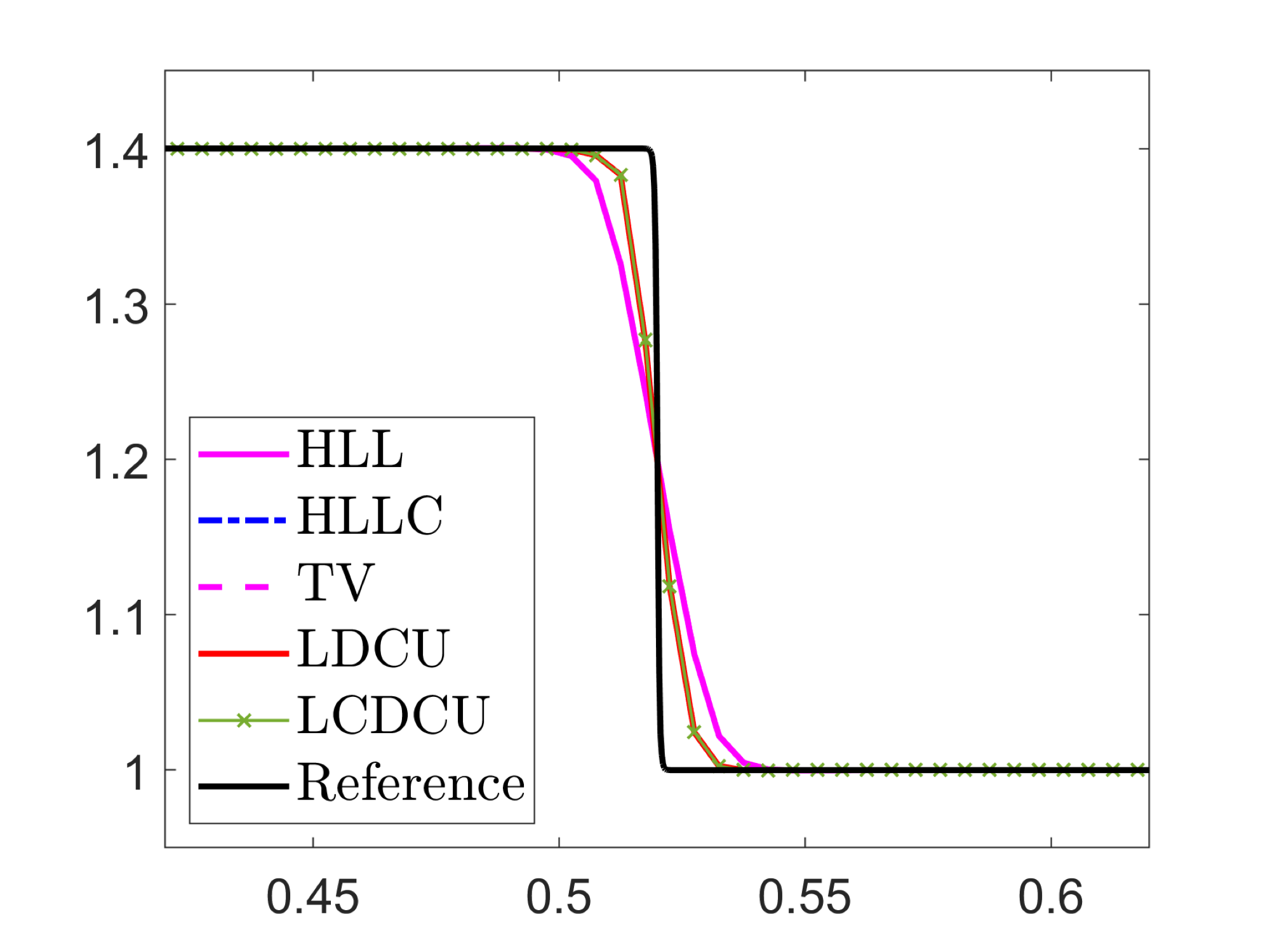}
            \includegraphics[trim=1.0cm 0.3cm 0.9cm 0.8cm, clip, width=4.cm]{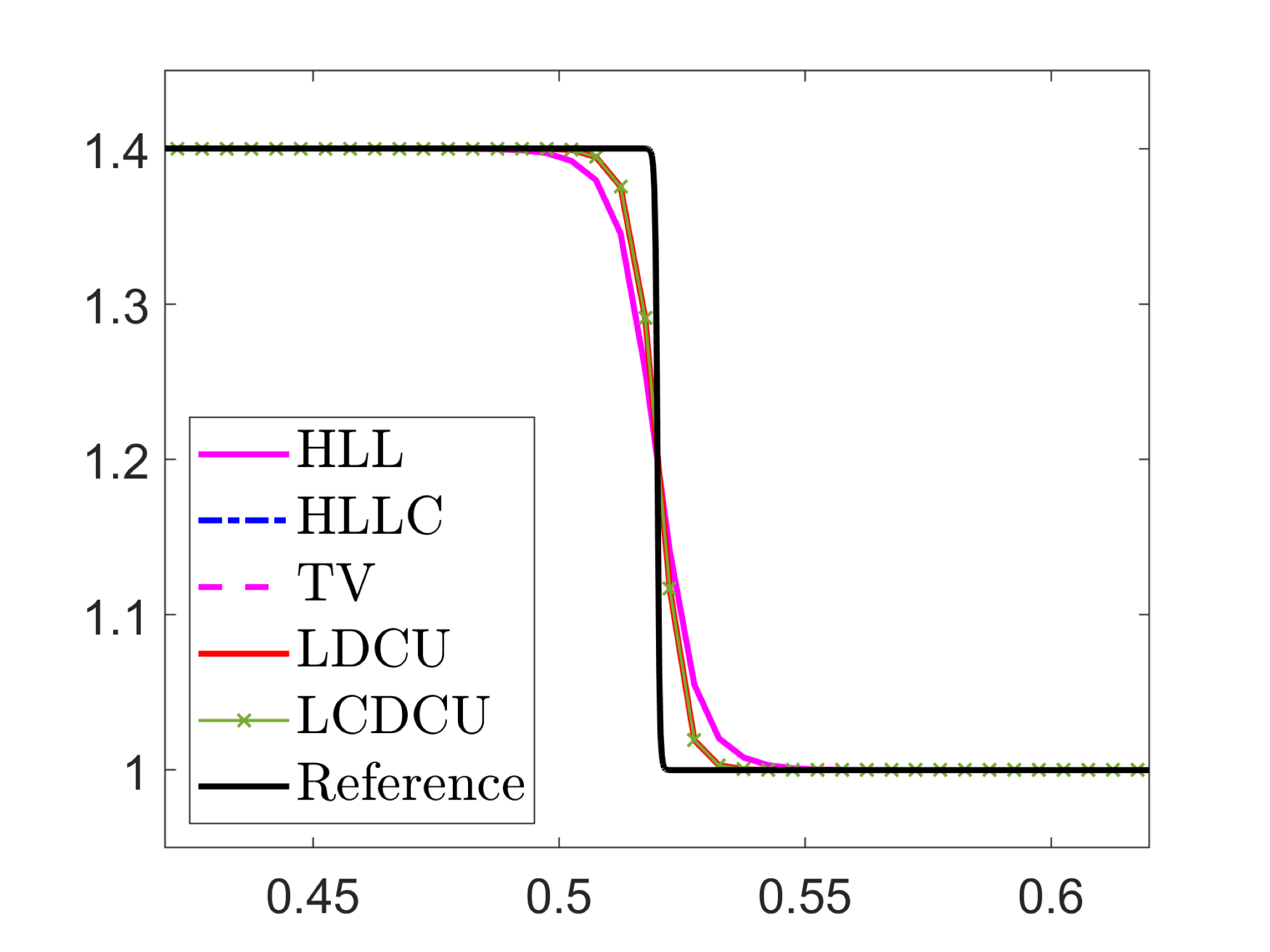}}
\caption{\sf Example 2: Density $\rho$ computed by the 1-Order, 2-Order, 3-Order and 5-Order schemes (top row) and zoom at $[0.42, 0.62]$ (bottom row).\label{fig2}}
\end{figure}

We compute the numerical solutions until the final time $t=0.2$ by the 1-Order, 2-Order, 3-Order, and 5-Order schemes on a uniform mesh with $\dx=1/200$ and then plot them in Figure \ref{fig2} together with the exact solution. As one can see, the numerical results computed by four low-dissipation schemes coincide and are better than those computed by the corresponding HLL schemes, even when they are extended to high orders. 

\subsubsection*{Example 3---Stationary Contact Wave, Traveling Shock and Rarefaction Wave}
In the third example, we consider the strong shocks interaction problem proposed in \cite{Woodward88}. The initial conditions,
\begin{equation*}
(\rho,u,p)(x,0)=\begin{cases}
(1,-19.59745,1000)&\mbox{if}~x<0.8,\\
(1,-19.59745,0.01)&\mbox{otherwise},\\
\end{cases}
\end{equation*}
are prescribed in the computational domain $[-5,5]$, in which free boundary conditions are implemented. 

\begin{figure}[ht!]
\centerline{\includegraphics[trim=0.8cm 0.3cm 1.4cm 0.6cm, clip, width=4.cm]{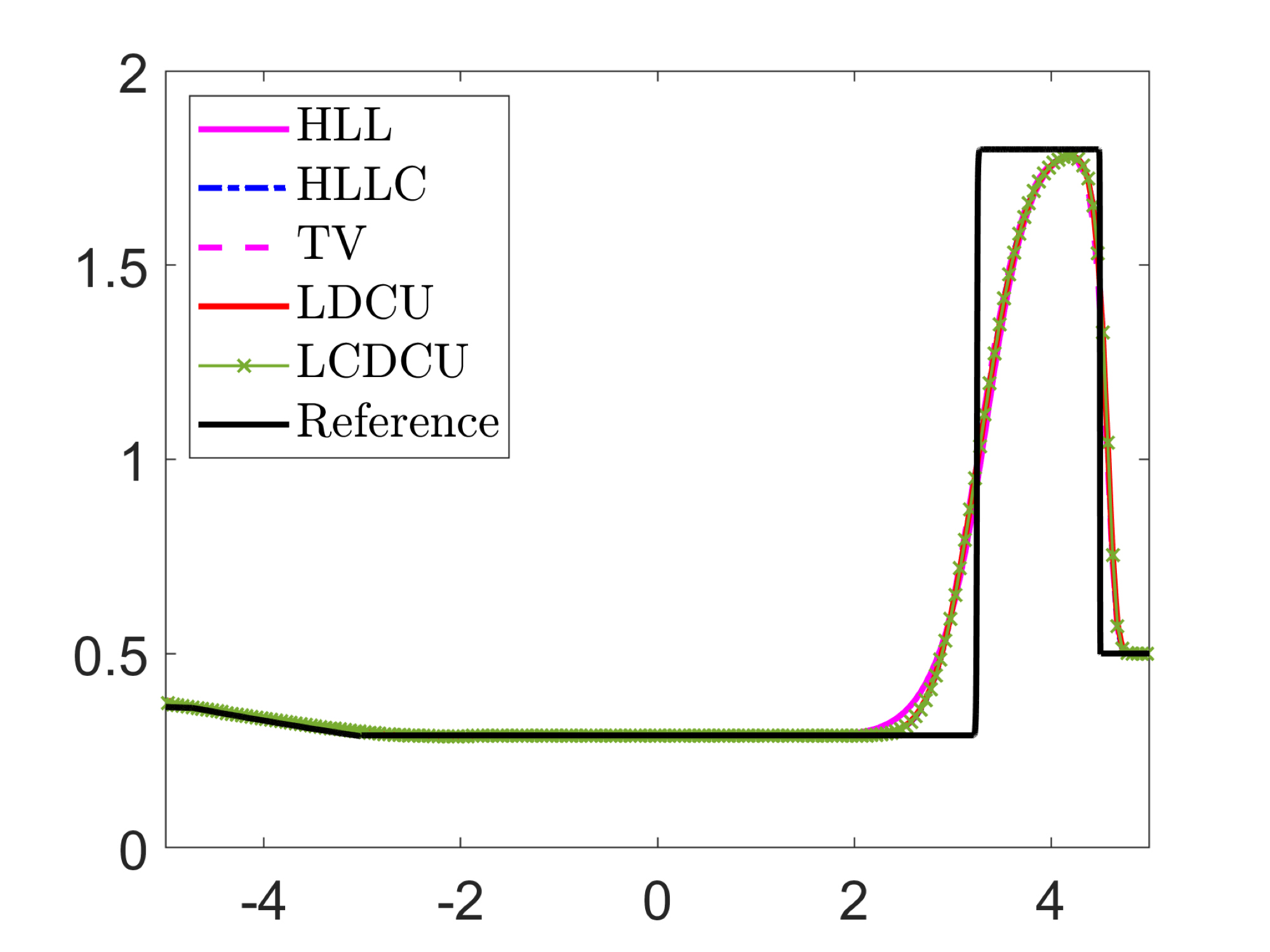}\hspace{0.2cm}
            \includegraphics[trim=0.8cm 0.3cm 1.4cm 0.6cm, clip, width=4.cm]{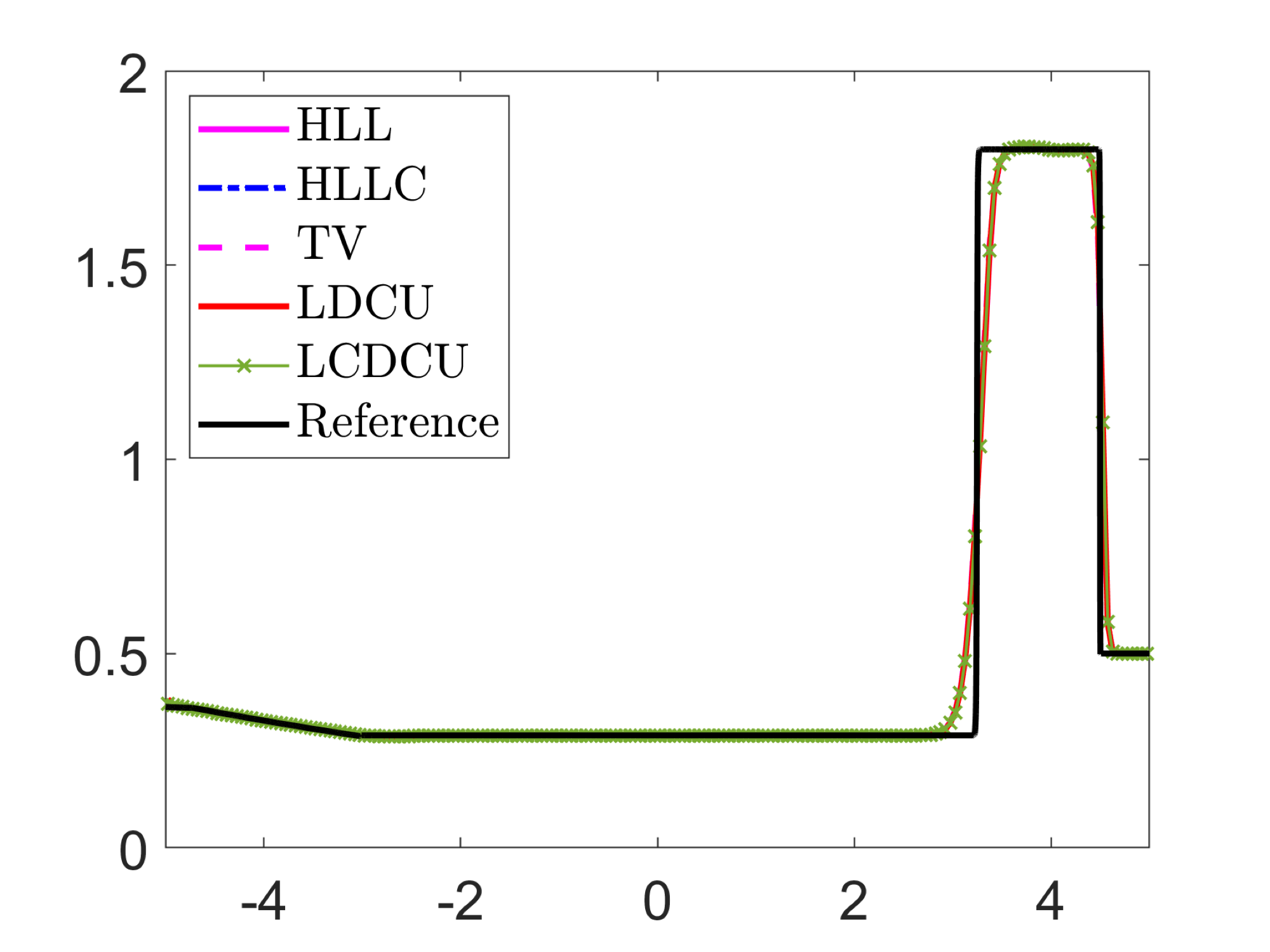}\hspace{0.2cm}
            \includegraphics[trim=0.8cm 0.3cm 1.4cm 0.6cm, clip, width=4.cm]{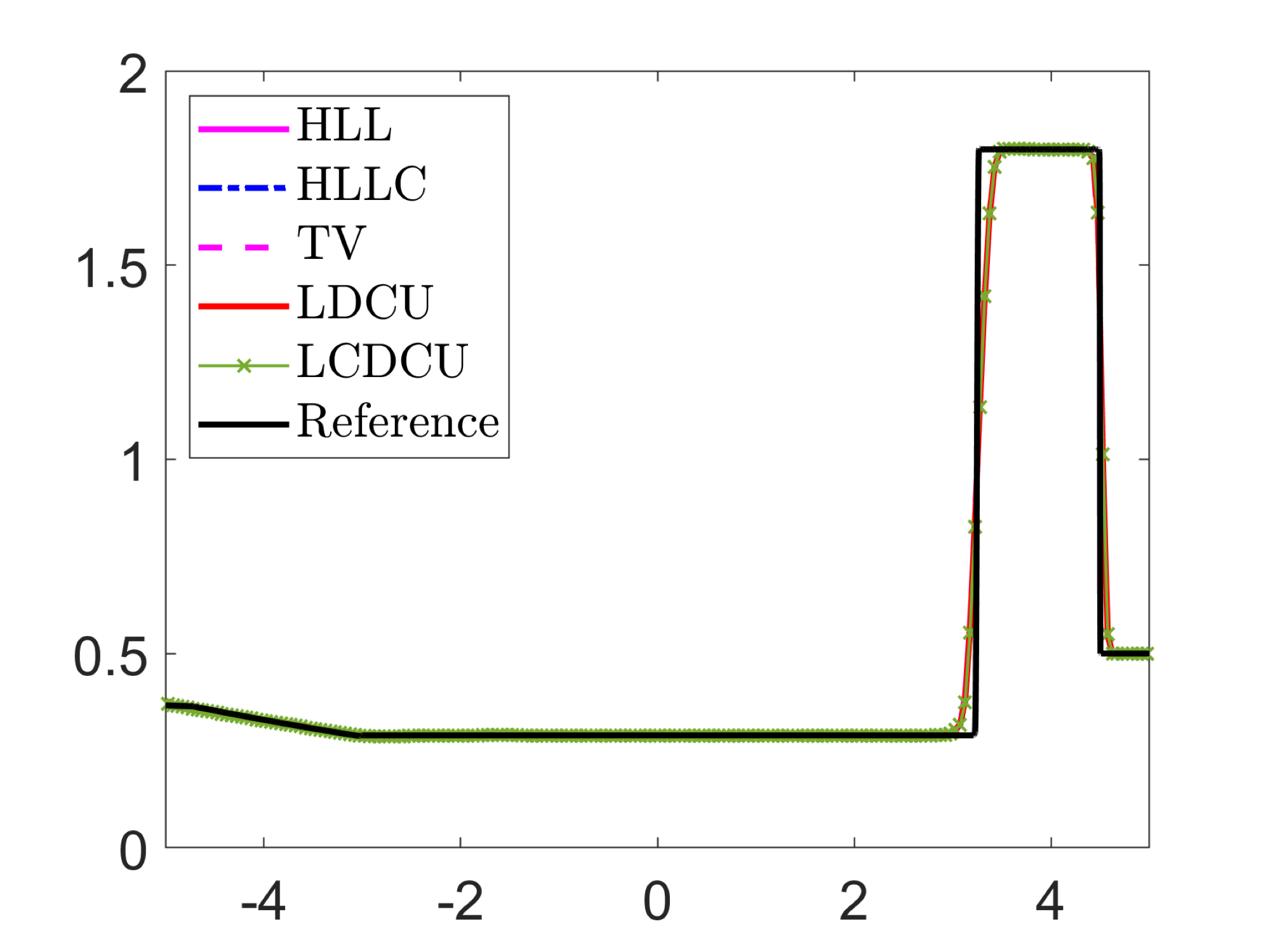}\hspace{0.2cm}
            \includegraphics[trim=0.8cm 0.3cm 1.4cm 0.6cm, clip, width=4.cm]{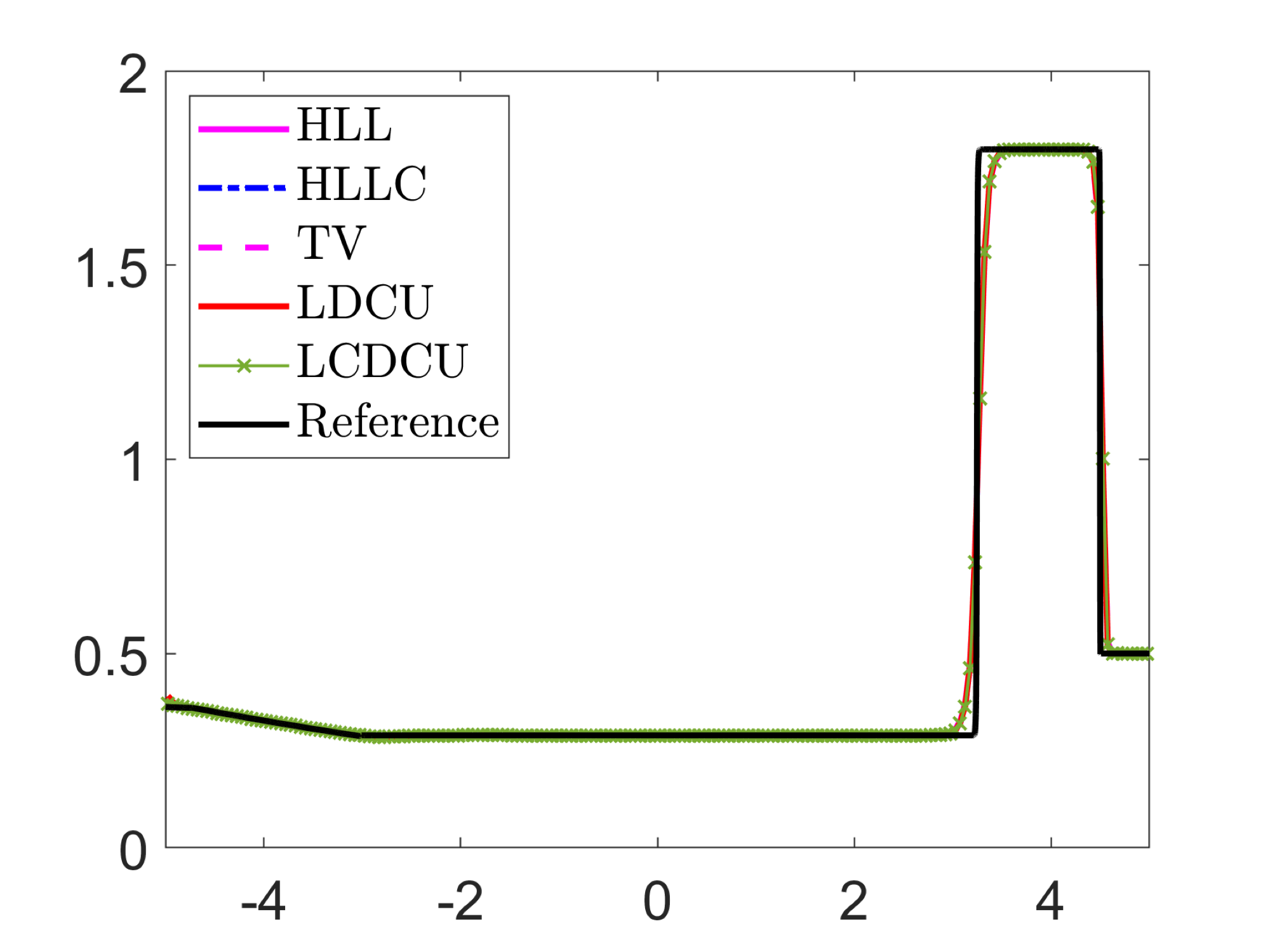}}
\vskip 12pt
\centerline{\includegraphics[trim=0.8cm 0.3cm 1.4cm 0.6cm, clip, width=4.cm]{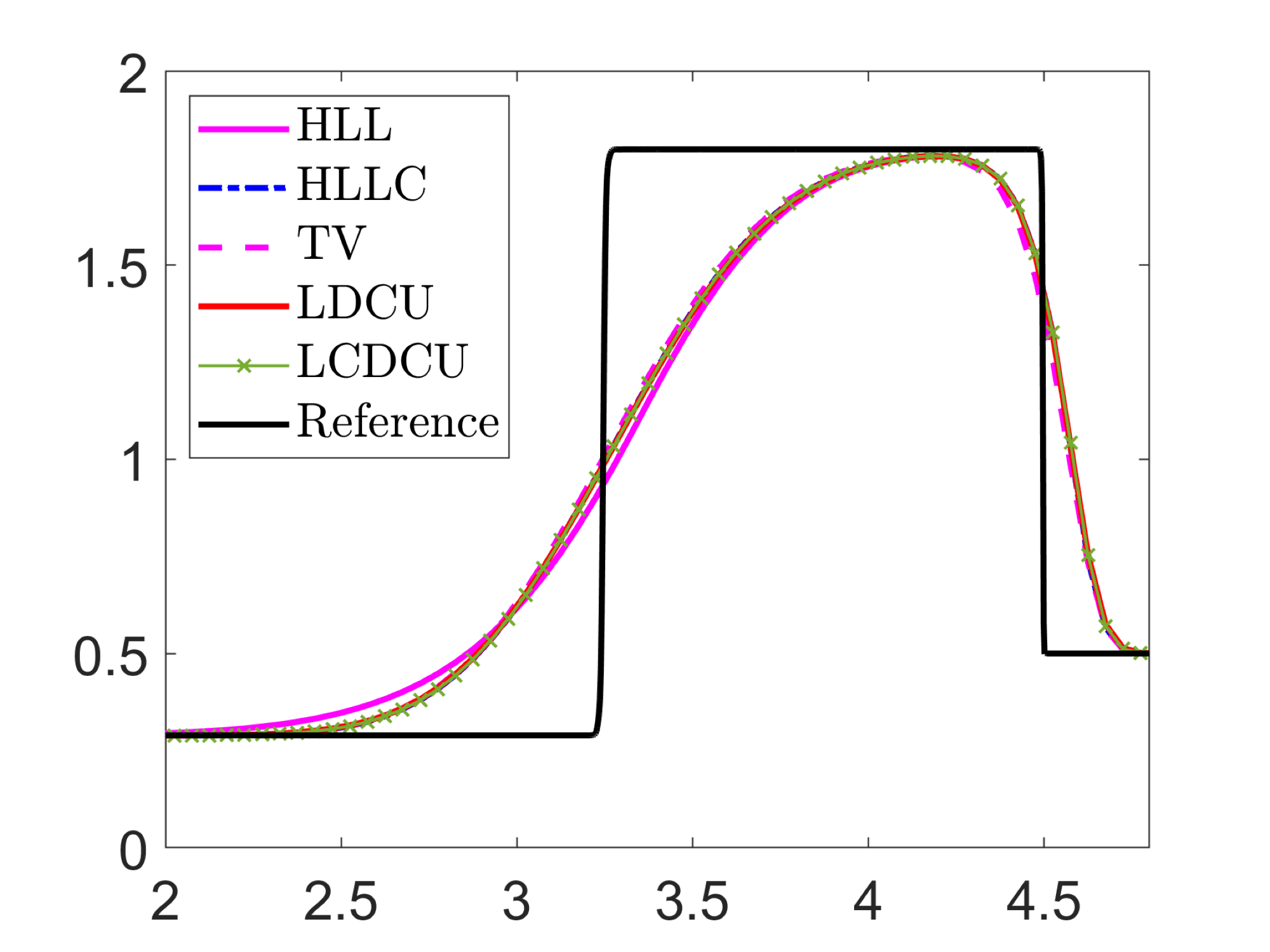}\hspace{0.2cm}
            \includegraphics[trim=0.8cm 0.3cm 1.4cm 0.6cm, clip, width=4.cm]{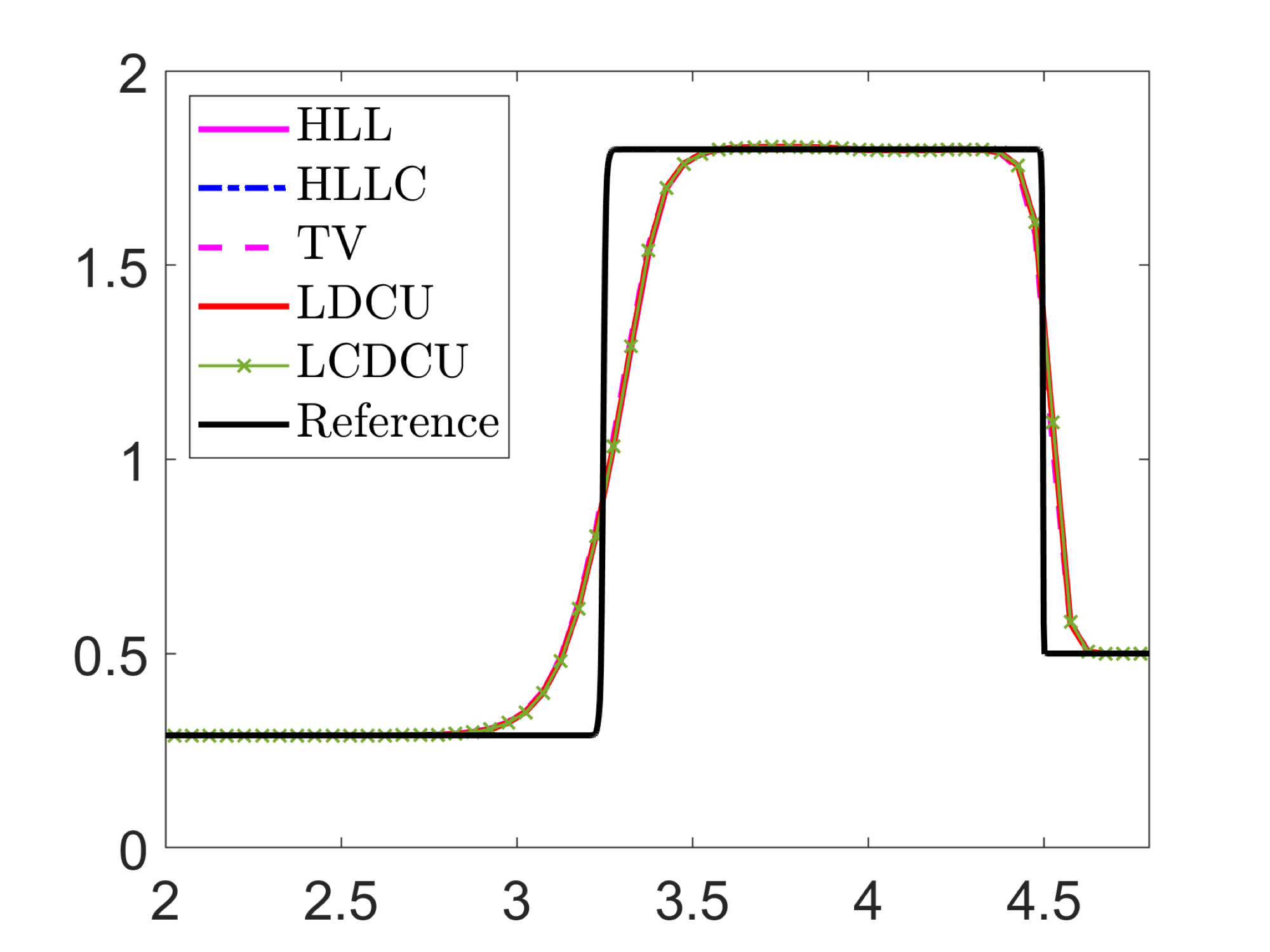}\hspace{0.2cm}
            \includegraphics[trim=0.8cm 0.3cm 1.4cm 0.6cm, clip, width=4.cm]{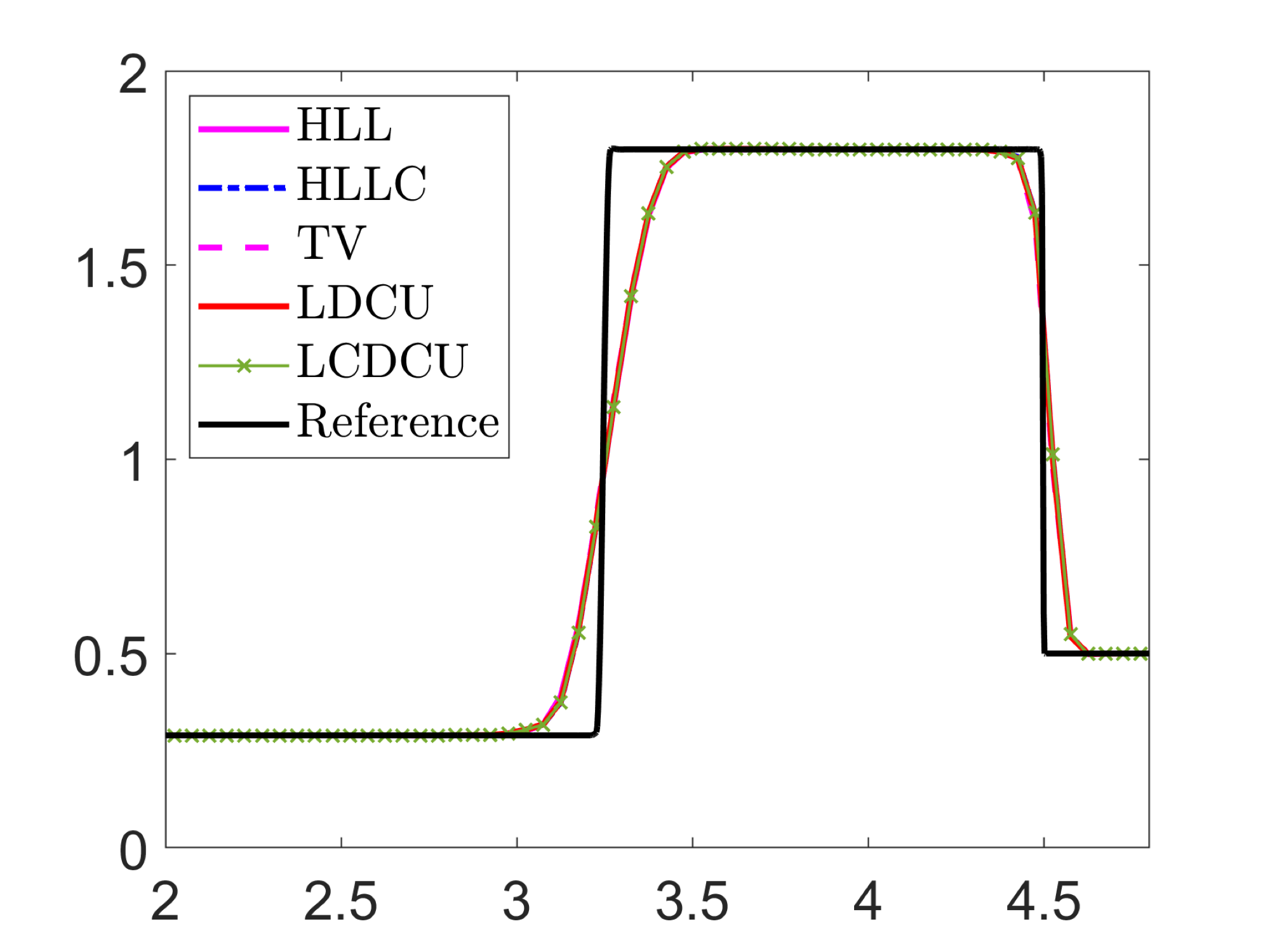}\hspace{0.2cm}
            \includegraphics[trim=0.8cm 0.3cm 1.4cm 0.6cm, clip, width=4.cm]{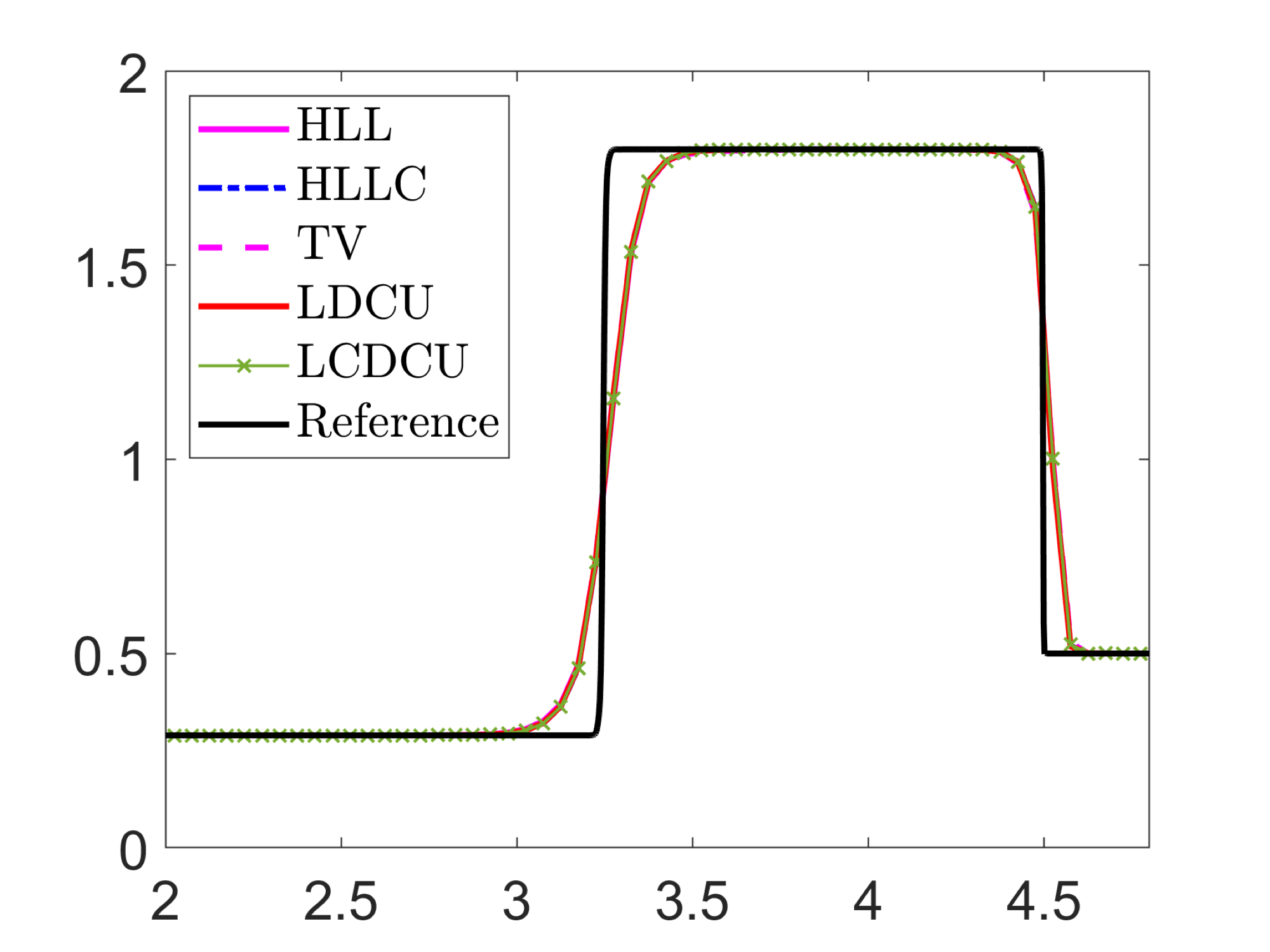}}
\caption{\sf Example 3: Density $\rho$ computed by the 1-Order, 2-Order, 3-Order, and 5-Order schemes (top row) and zoom at $[2, 4.8]$ (bottom row).
\label{fig3a}}
\end{figure}

We compute the numerical solutions until the final time $t=0.03$ using the studied 1-Order, 2-Order, 3-Order, and 5-Order schemes on a uniform mesh with $\dx=1/20$. The obtained numerical results are plotted in Figure \ref{fig3a}, as well as the reference solution, which is obtained by the HLL scheme on a much finer mesh with $\dx=1/400$, showing that the numerical results obtained with the four low-dissipation schemes are consistent with each other and outperform those produced by the corresponding HLL schemes, especially near the contact wave; see Figures \ref{fig3a} (bottom right), where we zoom at the neighborhood of the contact wave. When the studied schemes are extended to higher orders, the differences are very limited. 

\subsubsection*{Example 4---Shock Bubble Wave} 
In the fourth example, we consider the ``shock-bubble'' interaction problem taken from \cite{KX_22}. The initial data, given by
\begin{equation*}
(\rho, u,p)(x,0)=\begin{cases}
(13.1538,0,1)&\mbox{if}~|x|<0.25,\\
(1.3333,-0.3535,1.5)&\mbox{if}~x>0.75,\\
(1,0,1)&\mbox{otherwise},
\end{cases}
\end{equation*}
correspond to a left-moving shock, initially located at $x=0.75$, and a bubble of radius 0.25, initially located at the origin. We implement solid wall boundary conditions on the left and free boundary conditions on the right of the computational domain $[-1,1]$.

\begin{figure}[ht!]
\centerline{\includegraphics[trim=1.0cm 0.3cm 1.3cm 0.8cm, clip, width=4.cm]{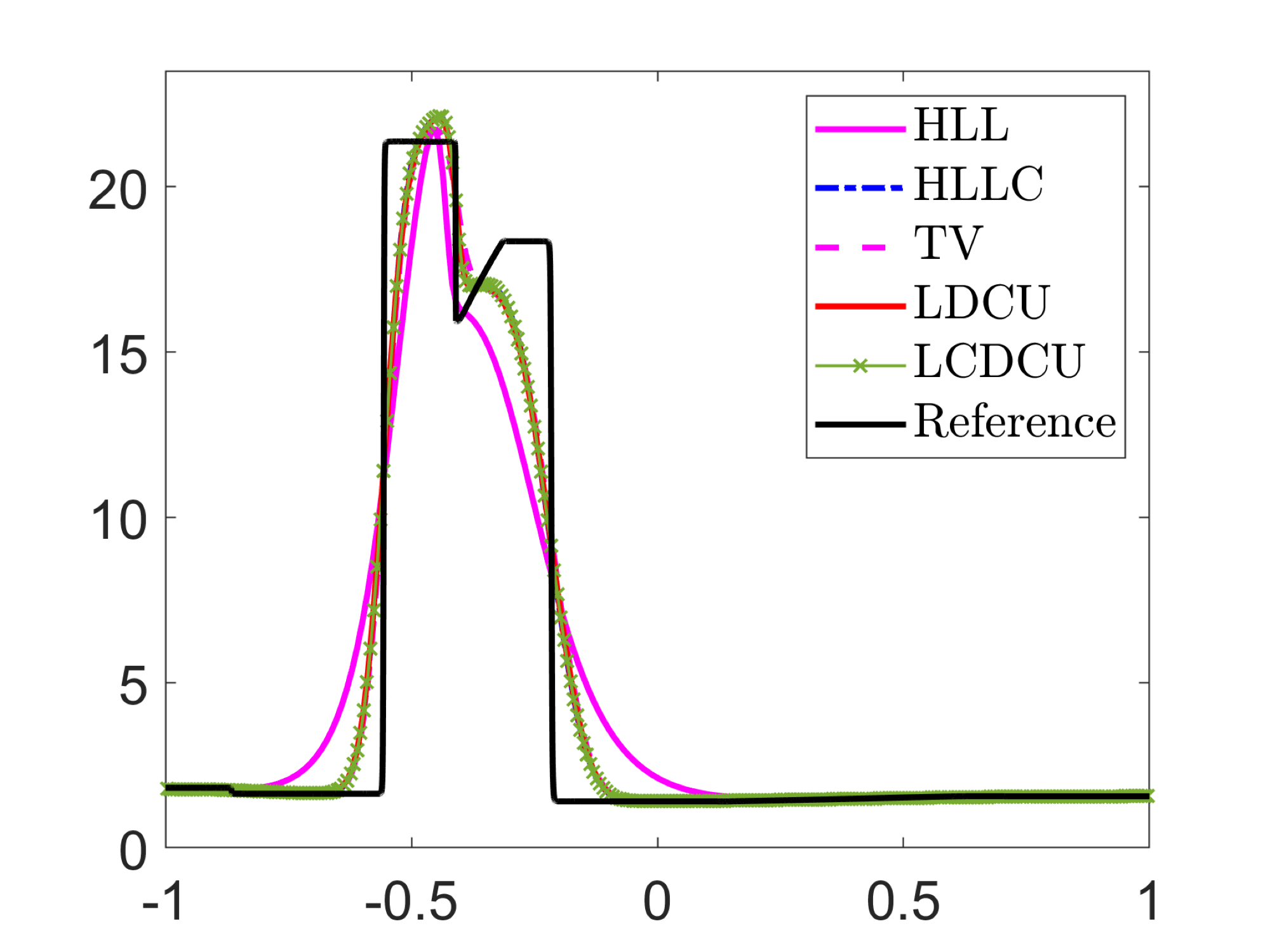}\hspace{0.2cm}
            \includegraphics[trim=1.0cm 0.3cm 0.9cm 0.8cm, clip, width=4.cm]{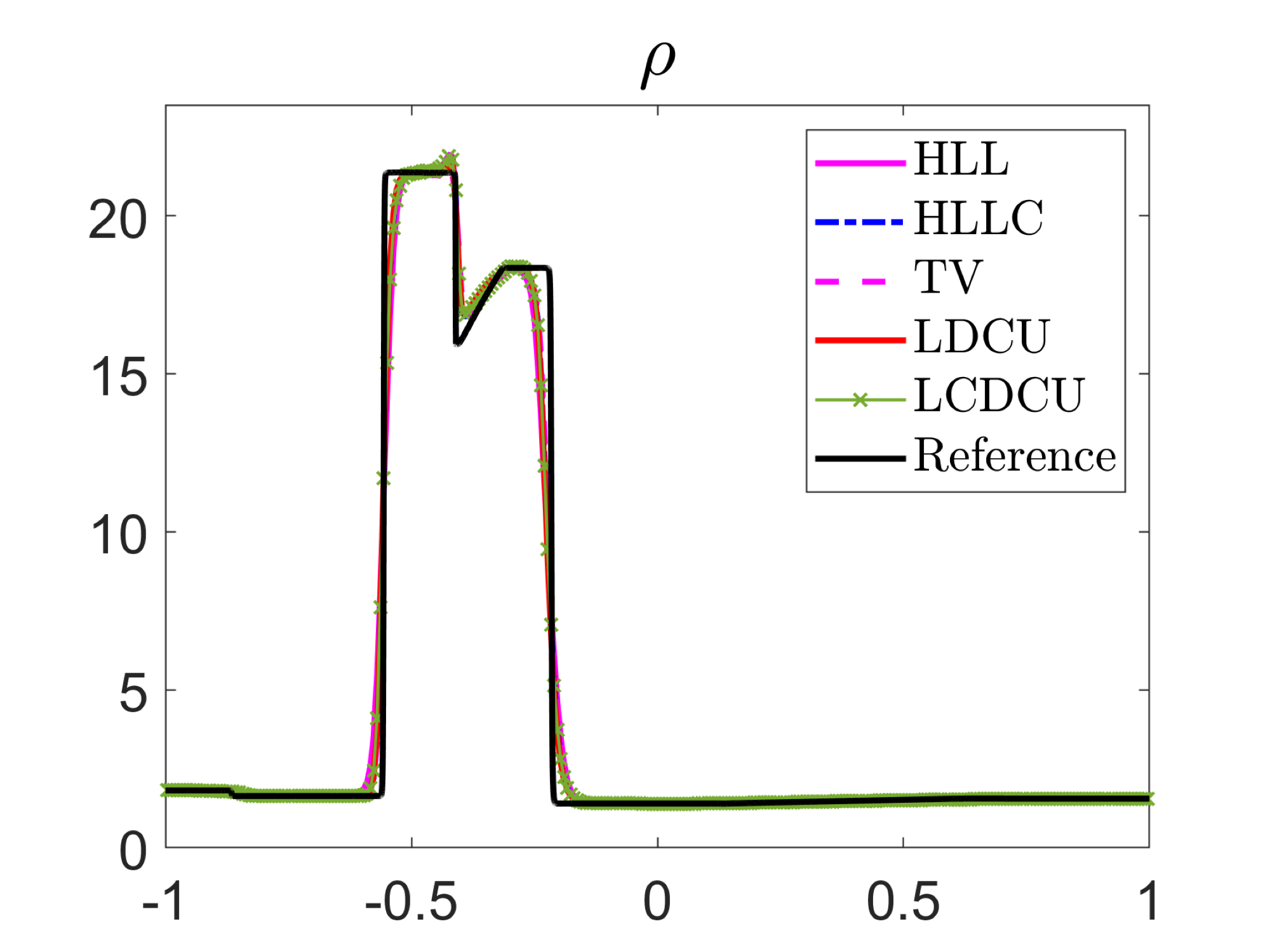}\hspace{0.2cm}
            \includegraphics[trim=1.0cm 0.3cm 0.9cm 0.8cm, clip, width=4.cm]{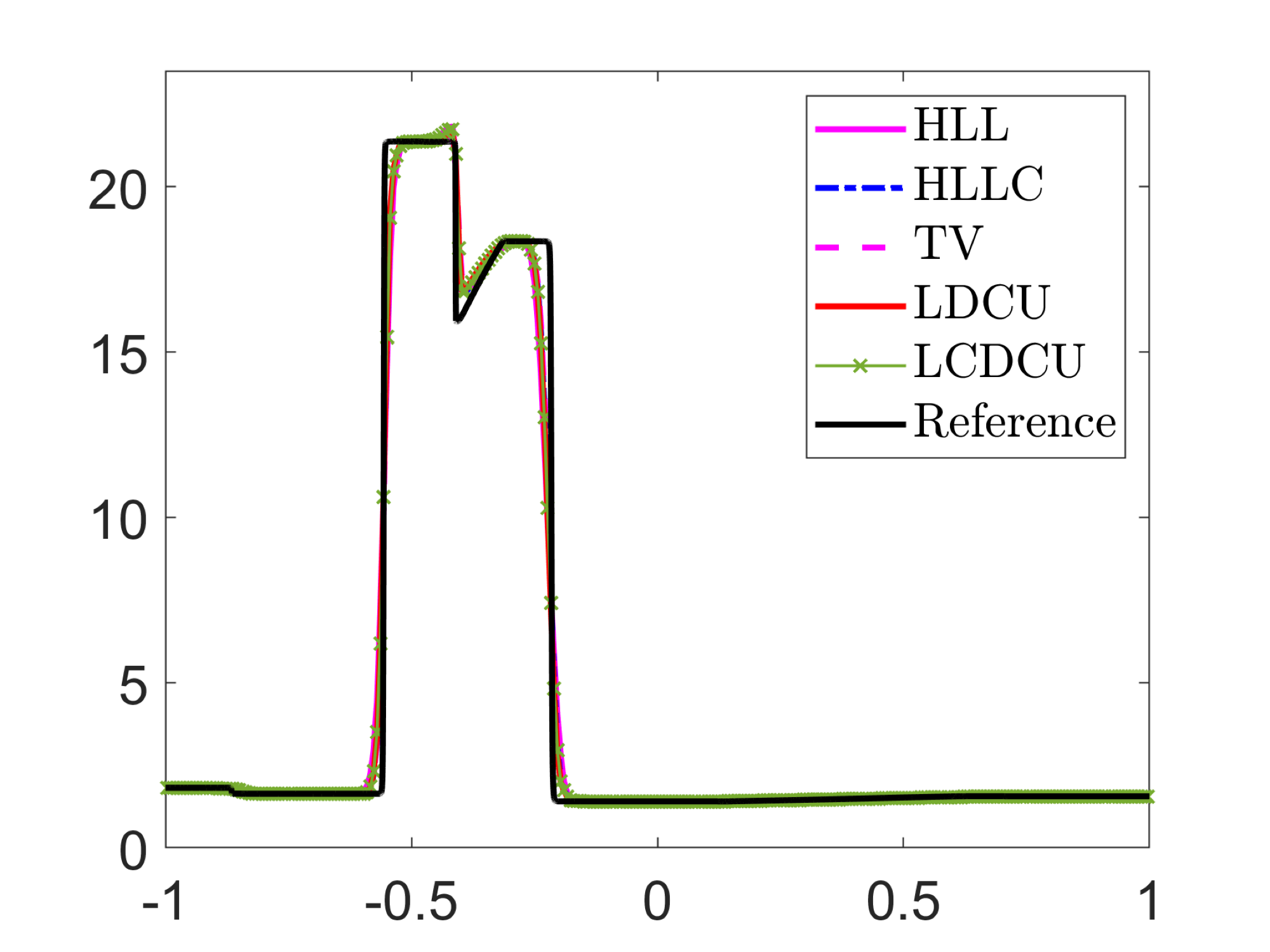}\hspace{0.2cm}
            \includegraphics[trim=1.0cm 0.3cm 0.9cm 0.8cm, clip, width=4.cm]{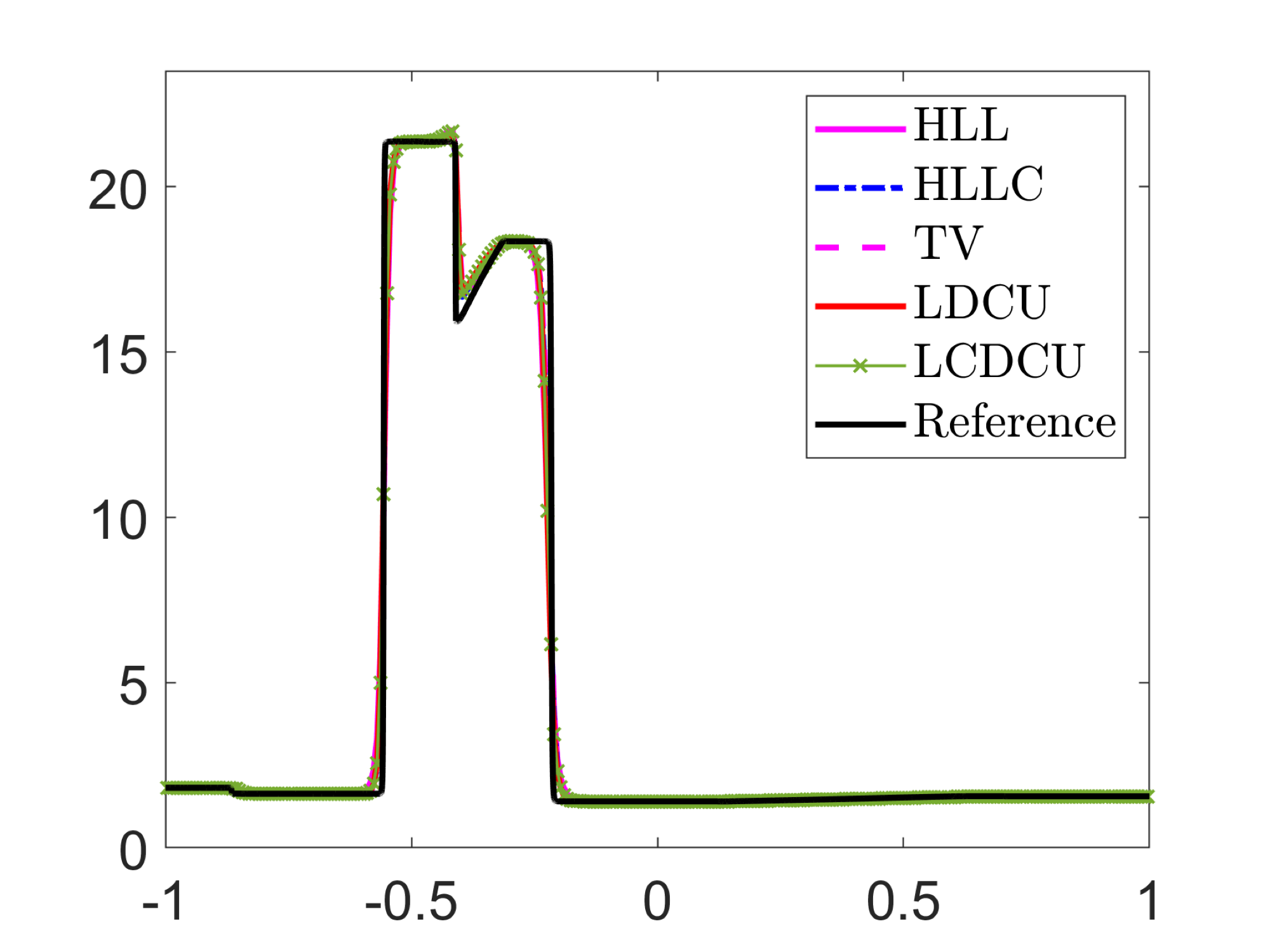}}
\vskip 12pt
\centerline{\includegraphics[trim=1.0cm 0.3cm 0.9cm 0.8cm, clip, width=4.cm]{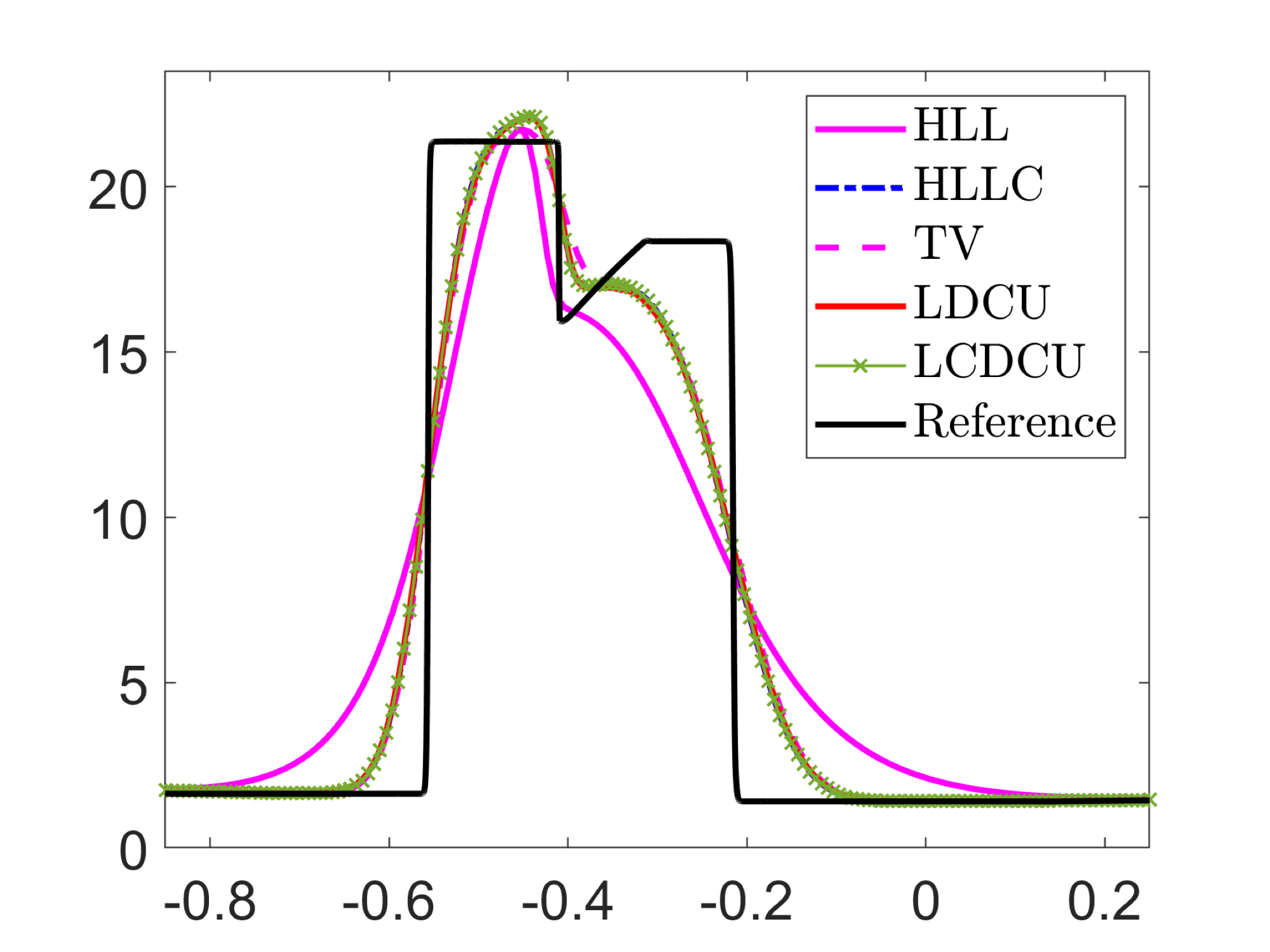}\hspace{0.2cm}
            \includegraphics[trim=1.0cm 0.3cm 0.9cm 0.8cm, clip, width=4.cm]{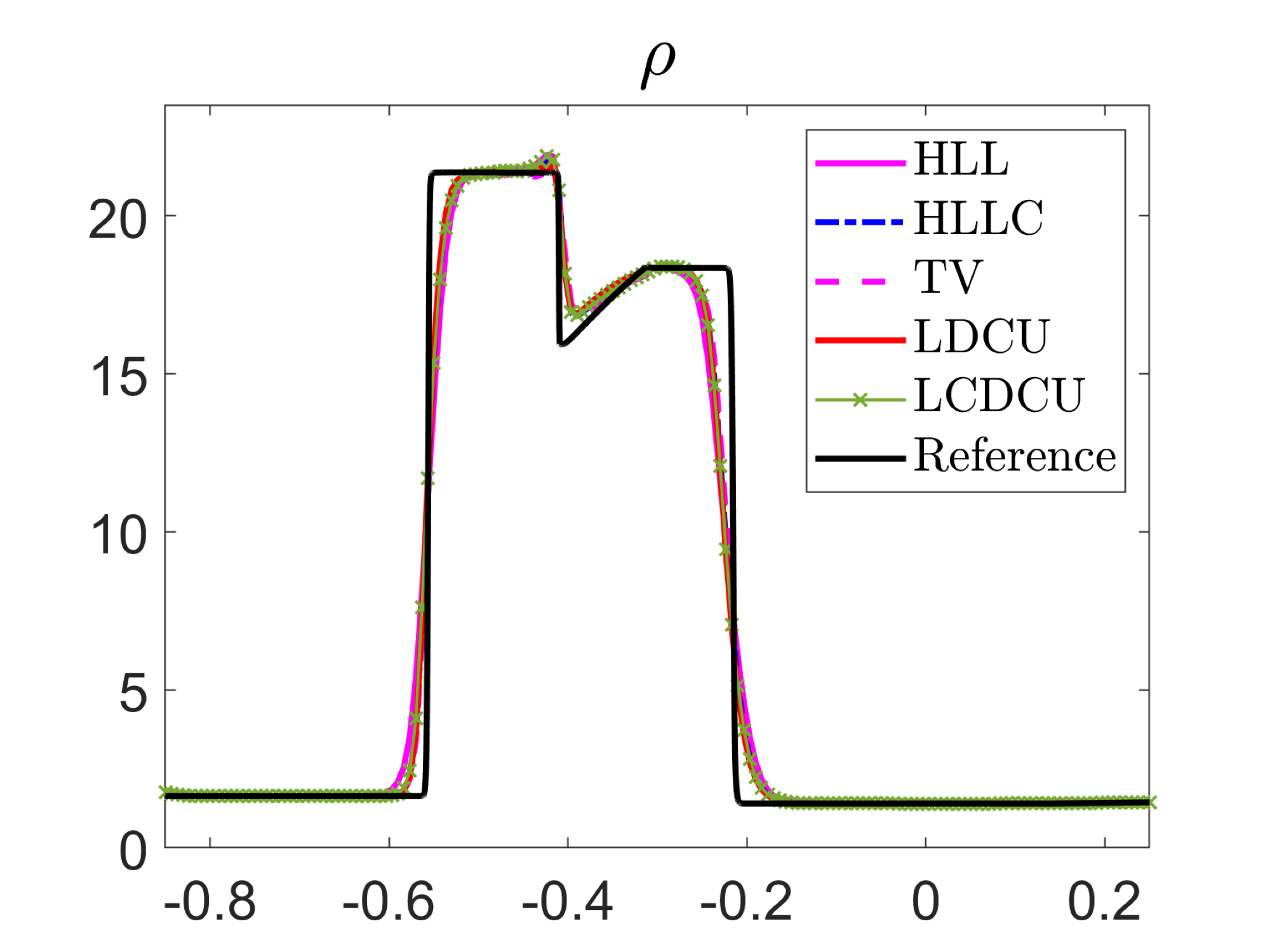}\hspace{0.2cm}
            \includegraphics[trim=1.0cm 0.3cm 0.9cm 0.8cm, clip, width=4.cm]{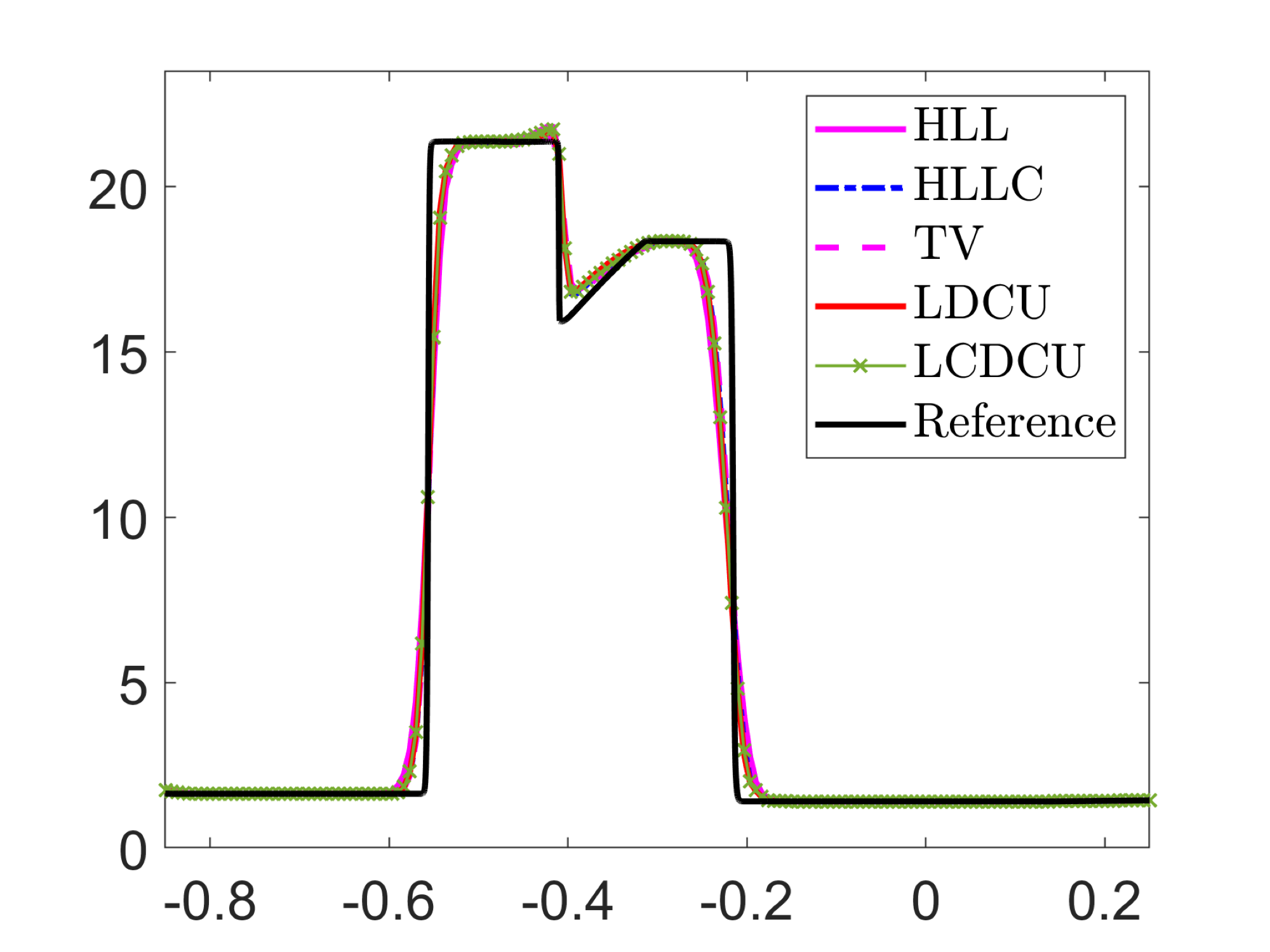}\hspace{0.2cm}
            \includegraphics[trim=1.0cm 0.3cm 0.9cm 0.8cm, clip, width=4.cm]{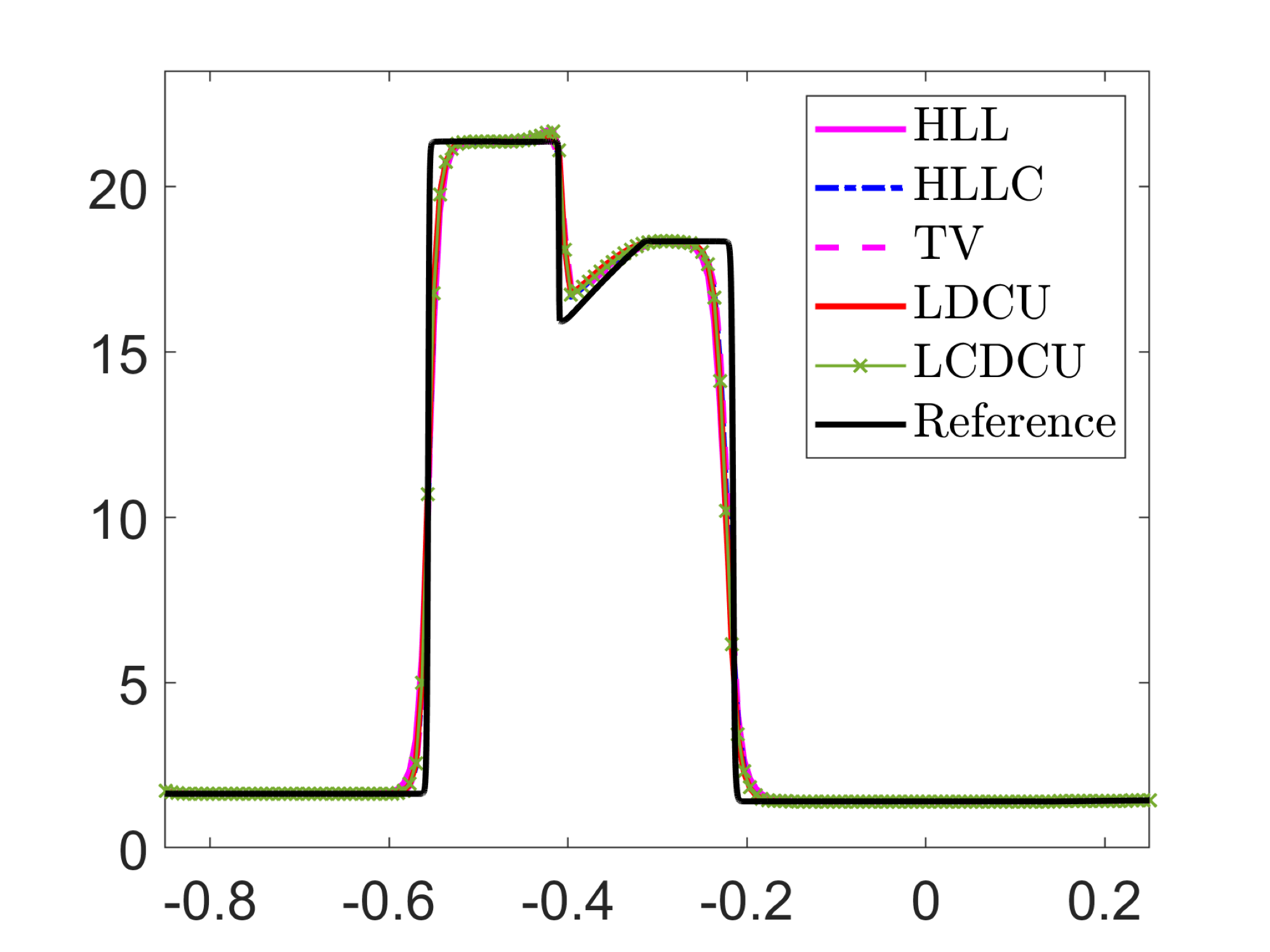}}
\caption{\sf Example 4: Density $\rho$ computed by the 1-Order, 2-Order, 3-Order, and 5-Order schemes (top row) and zoom at $[-0.85,0.25]$ (bottom row).
\label{fig4a}}
\end{figure}

We compute the numerical solution until the final time $t=3$ using the studied 1-Order, 2-Order, 3-Order, and 5-Order schemes on the uniform mesh with $\dx=1/100$. The obtained numerical results are presented in Figure \ref{fig4a} along with the reference solution computed by the HLL schemes on a much finer mesh with $\dx=1/2000$. One can notice that there are no clear differences in the numerical results computed by the four low-dissipation schemes, but the results computed by the low-dissipation schemes are better than the HLL scheme, especially in the 1-Order and 2-Order results. 

\subsubsection*{Example 5---Shock-Density Wave Interaction Problem}
In this example taken from \cite{SO89}, we consider the shock-density wave interaction problem with the following initial data,
\begin{equation*}
(\rho,u,p)\Big|_{(x,0)}=\begin{cases}\bigg(\dfrac{27}{7},\dfrac{4\sqrt{35}}{9},\dfrac{31}{3}\bigg),&x<-4,\\[0.8ex]
(1+0.2\sin(5x),0,1),&x>-4,
\end{cases}
\end{equation*}
prescribed in the computational domain $[-5,5]$ subject to the free boundary conditions.

\begin{figure}[ht!]
\centerline{\includegraphics[trim=0.8cm 0.3cm 0.9cm 0.8cm, clip, width=4.5cm]{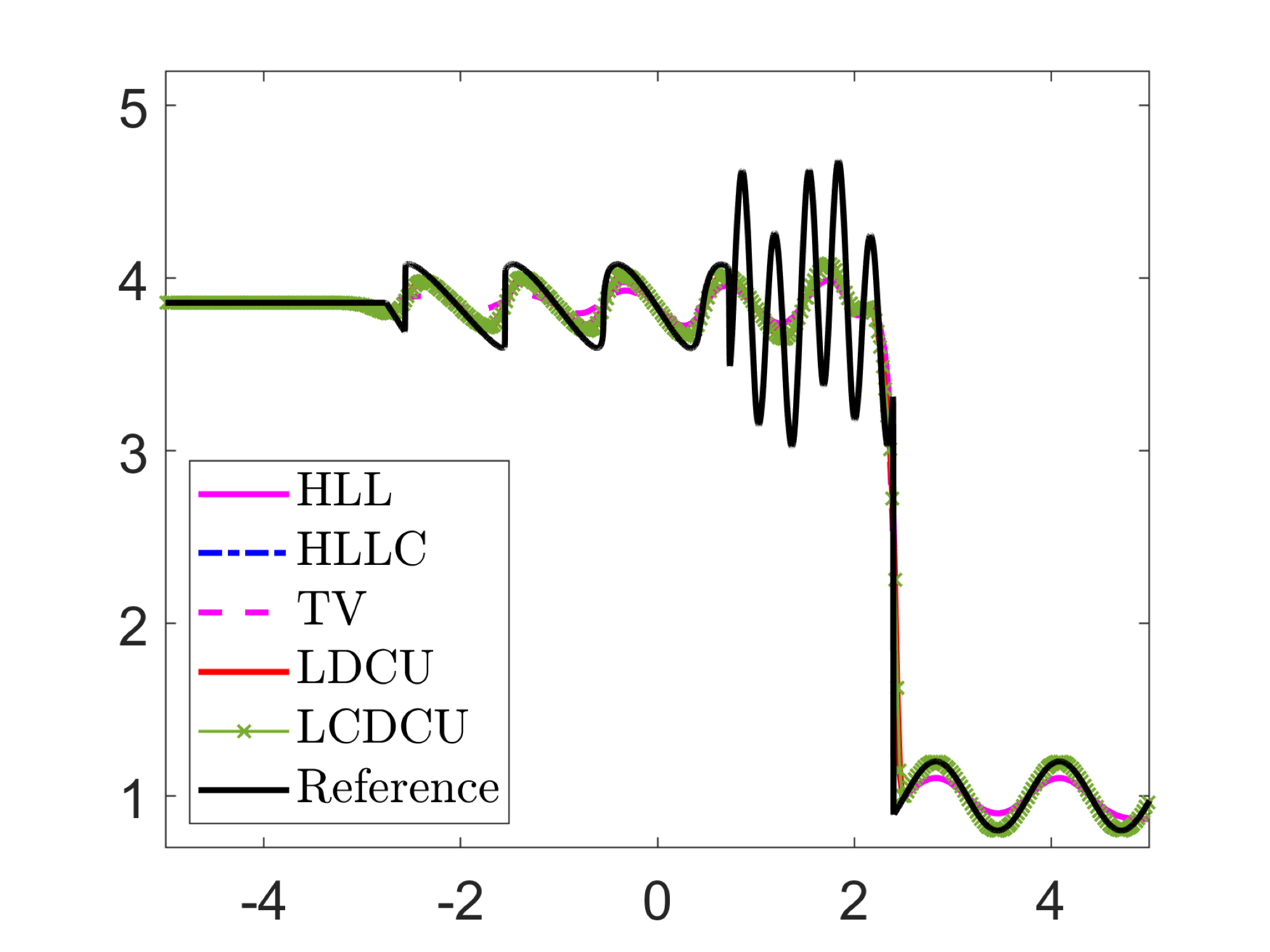}\hspace{0.5cm}
            \includegraphics[trim=0.8cm 0.3cm 0.9cm 0.8cm, clip, width=4.5cm]{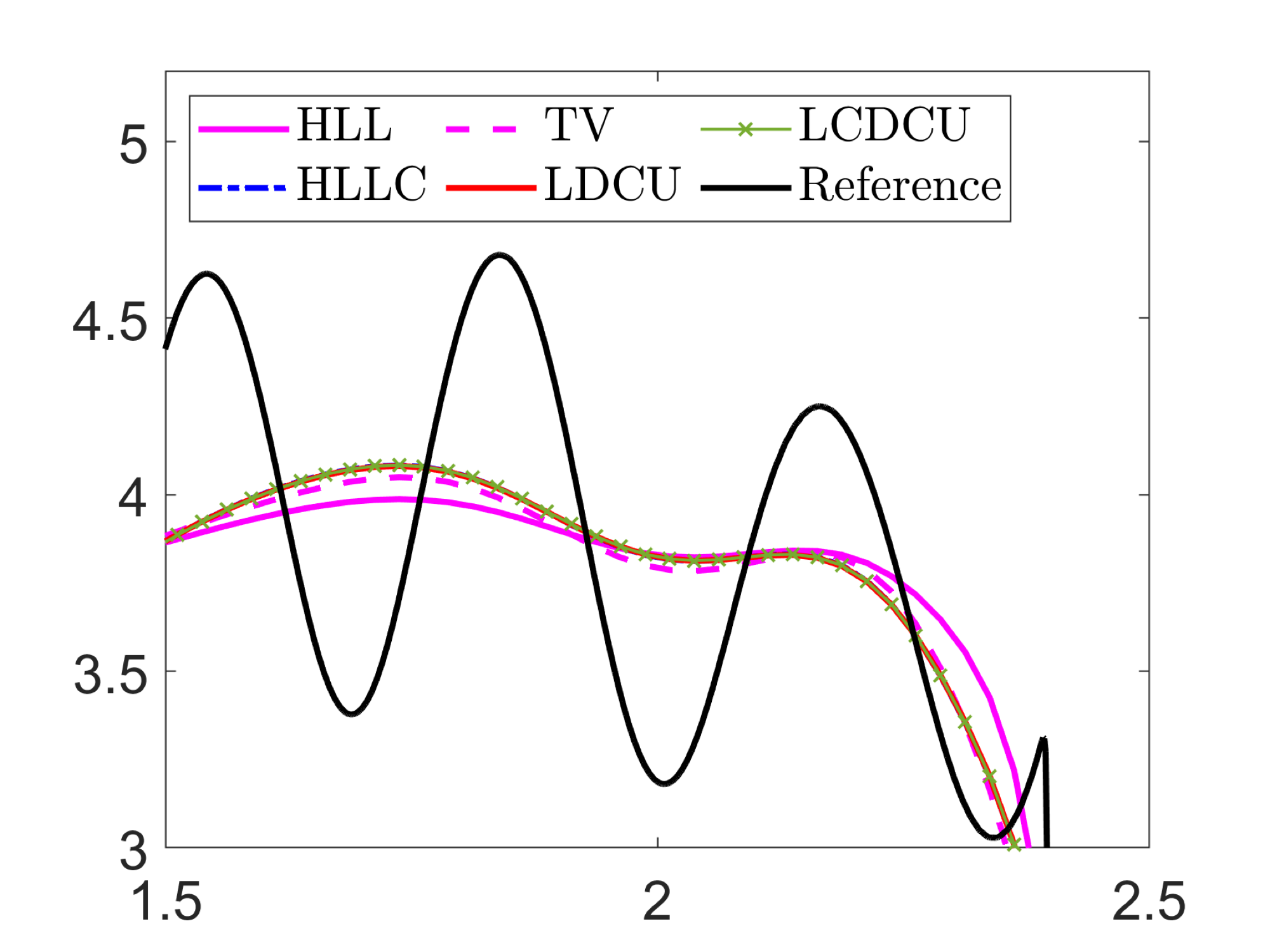}\hspace{0.5cm}
            \includegraphics[trim=0.8cm 0.3cm 0.9cm 0.8cm, clip, width=4.5cm]{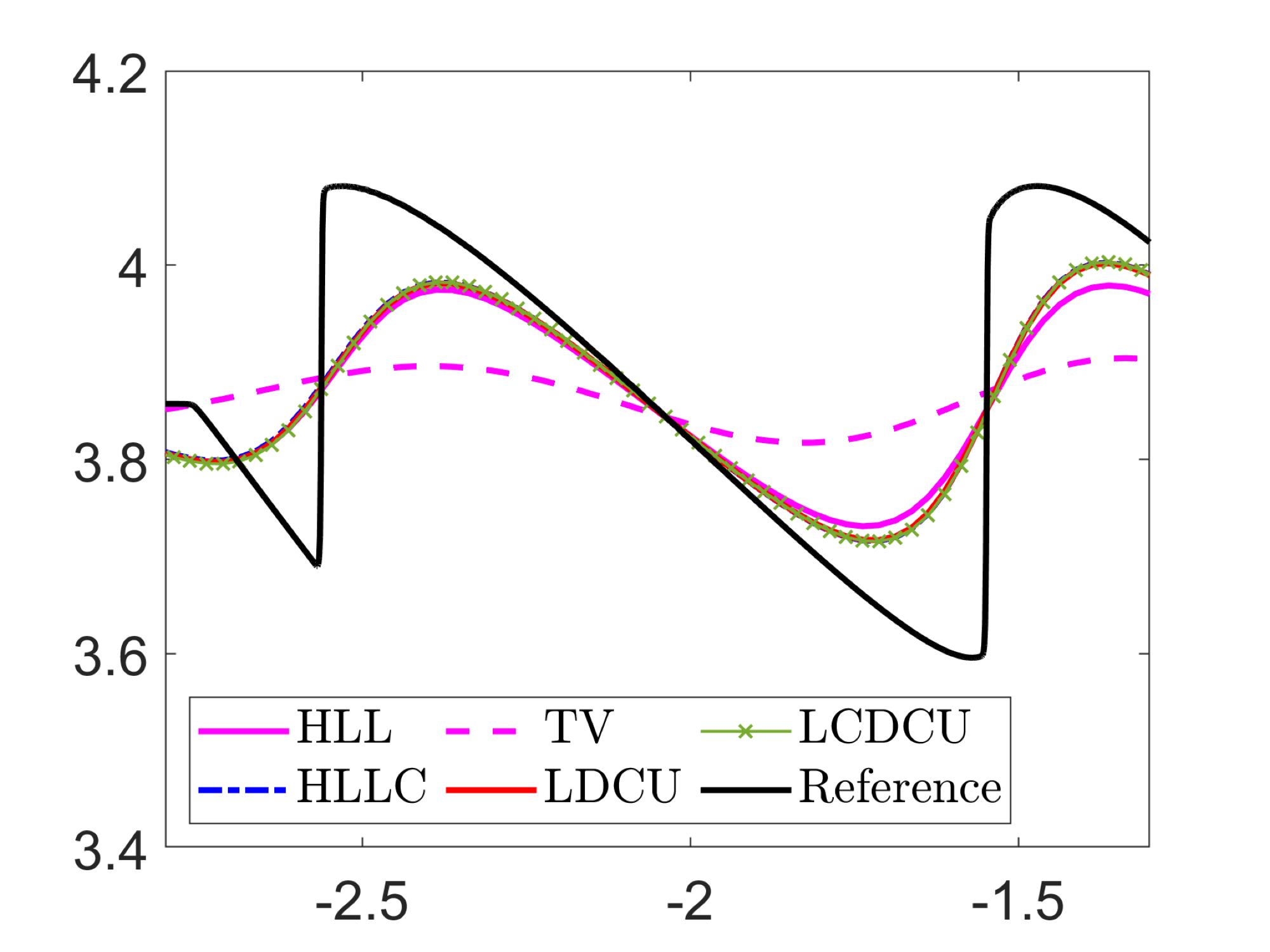}}
\vskip 12pt
\centerline{\includegraphics[trim=0.8cm 0.3cm 0.9cm 0.8cm, clip, width=4.5cm]{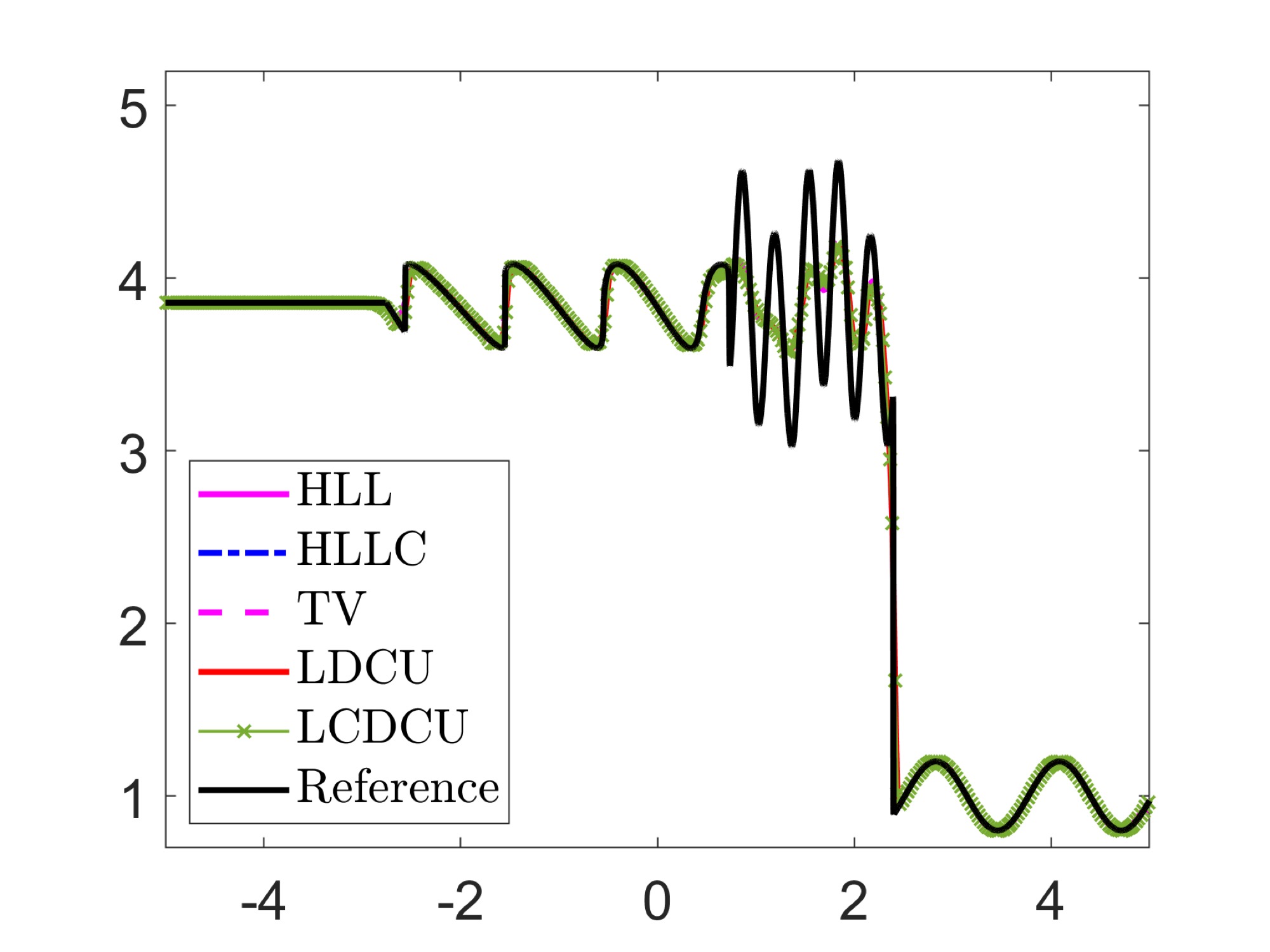}\hspace{0.5cm}
            \includegraphics[trim=0.8cm 0.3cm 0.9cm 0.8cm, clip, width=4.5cm]{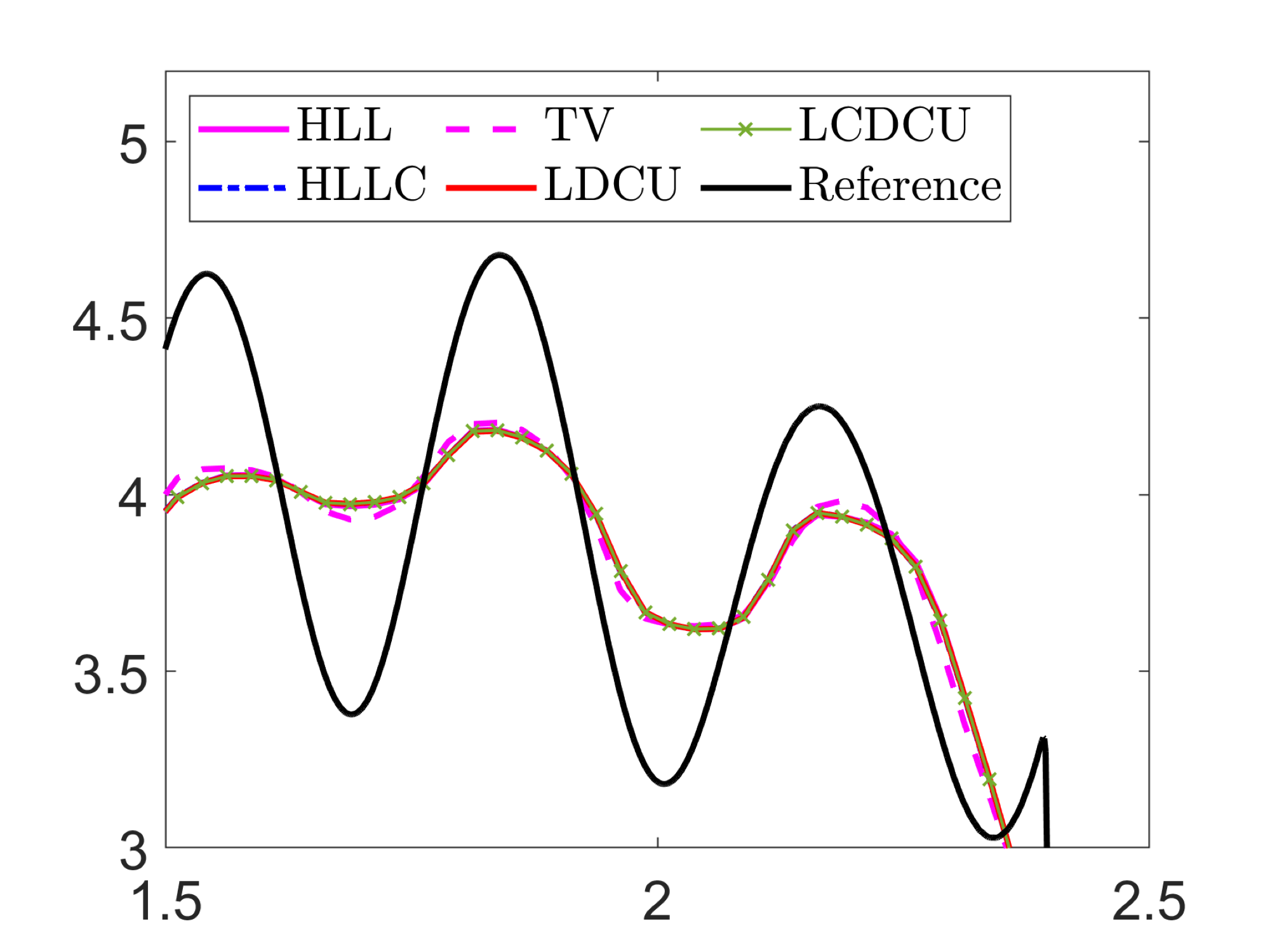}\hspace{0.5cm}
            \includegraphics[trim=0.8cm 0.3cm 0.9cm 0.8cm, clip, width=4.5cm]{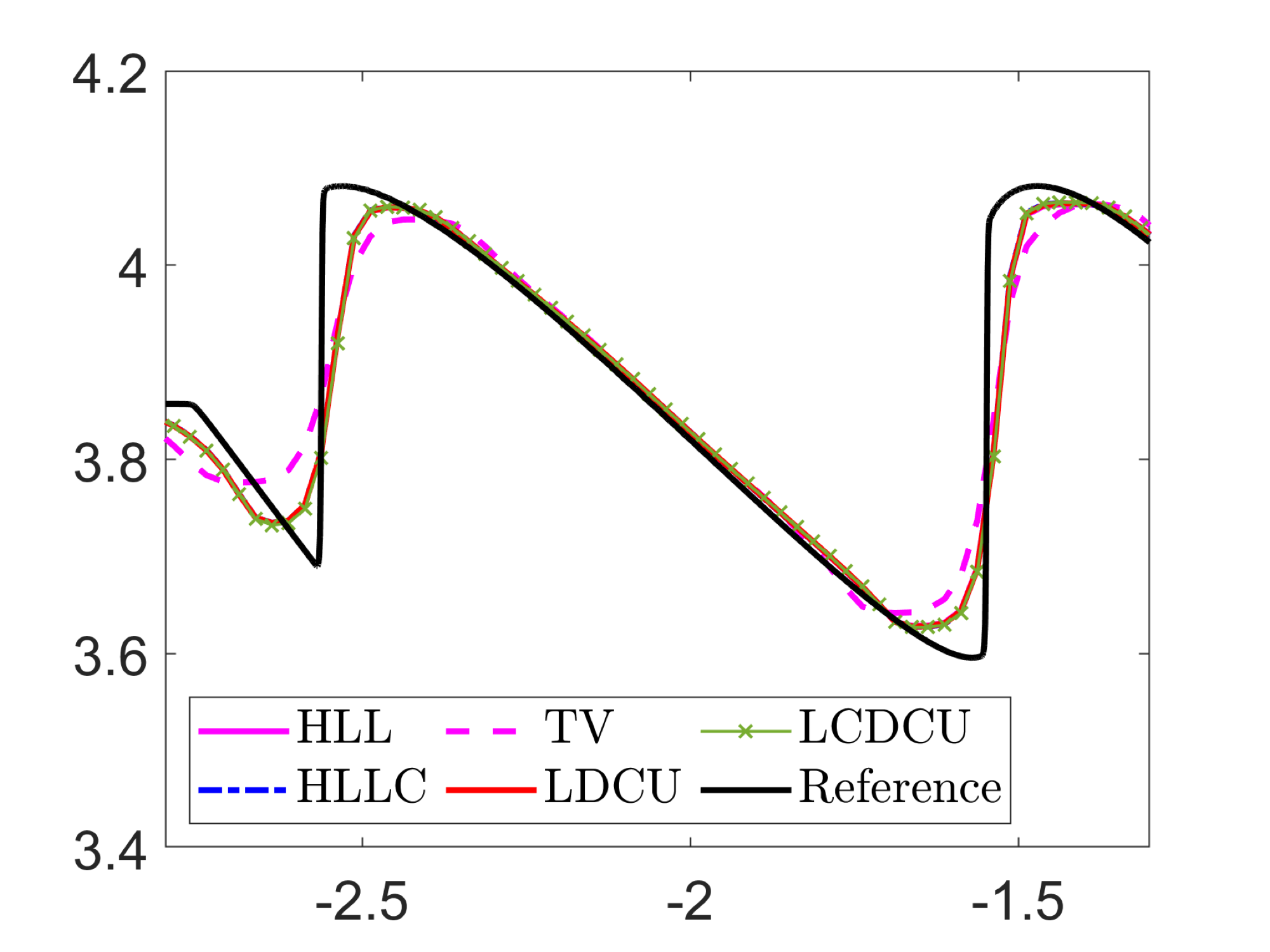}}
\caption{\sf Example 5: Density $\rho$ computed by the 1-Order and 2-Order schemes (left column) and zoom at $[1.5,2.5]$ (middle column), and  $[-2.8,-1.3]$ (right column).
\label{fig5a}}
\end{figure}

\begin{figure}[ht!]
\centerline{\includegraphics[trim=0.8cm 0.3cm 0.9cm 0.8cm, clip, width=4.5cm]{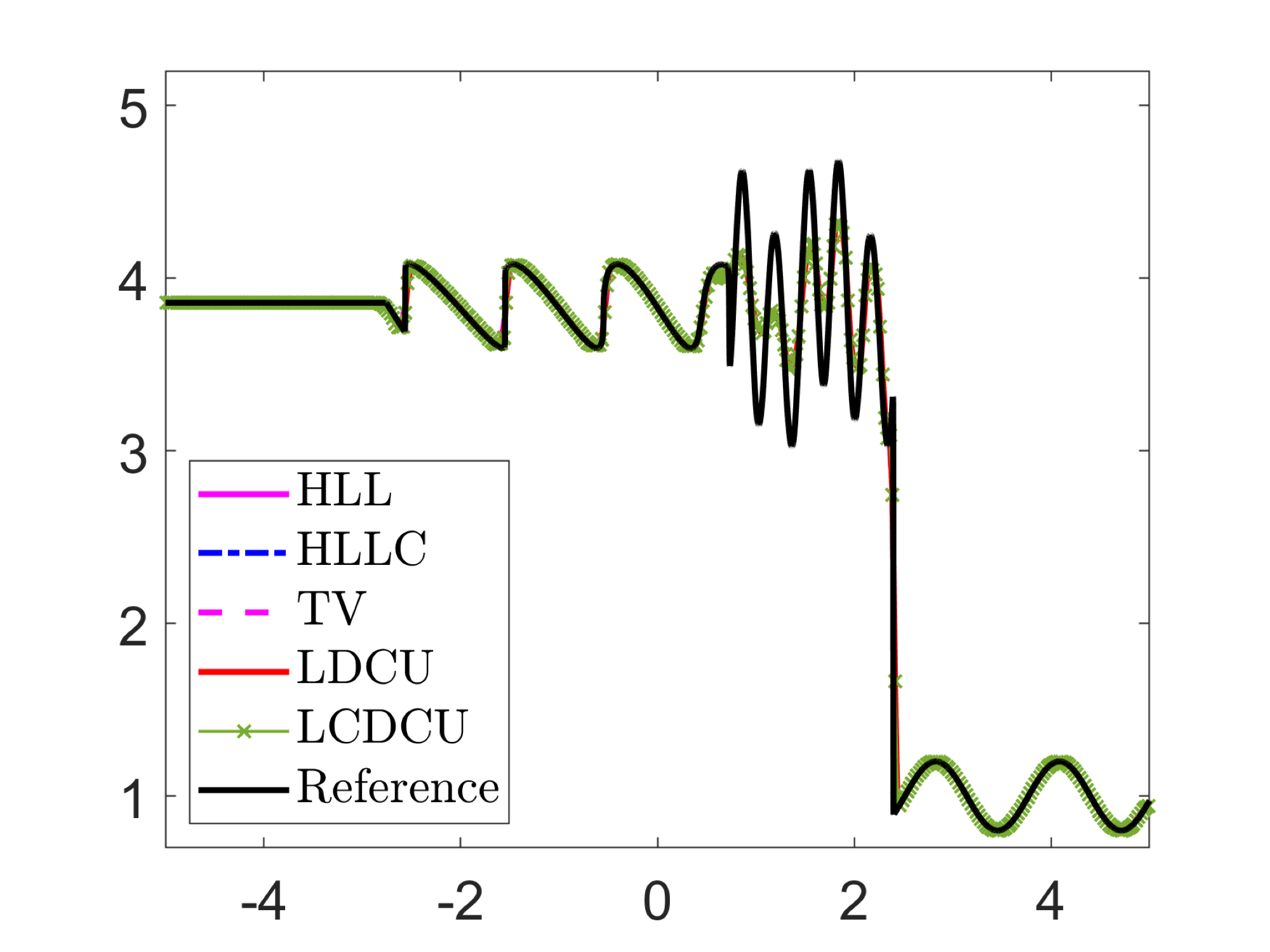}\hspace{0.5cm}
            \includegraphics[trim=0.8cm 0.3cm 0.9cm 0.8cm, clip, width=4.5cm]{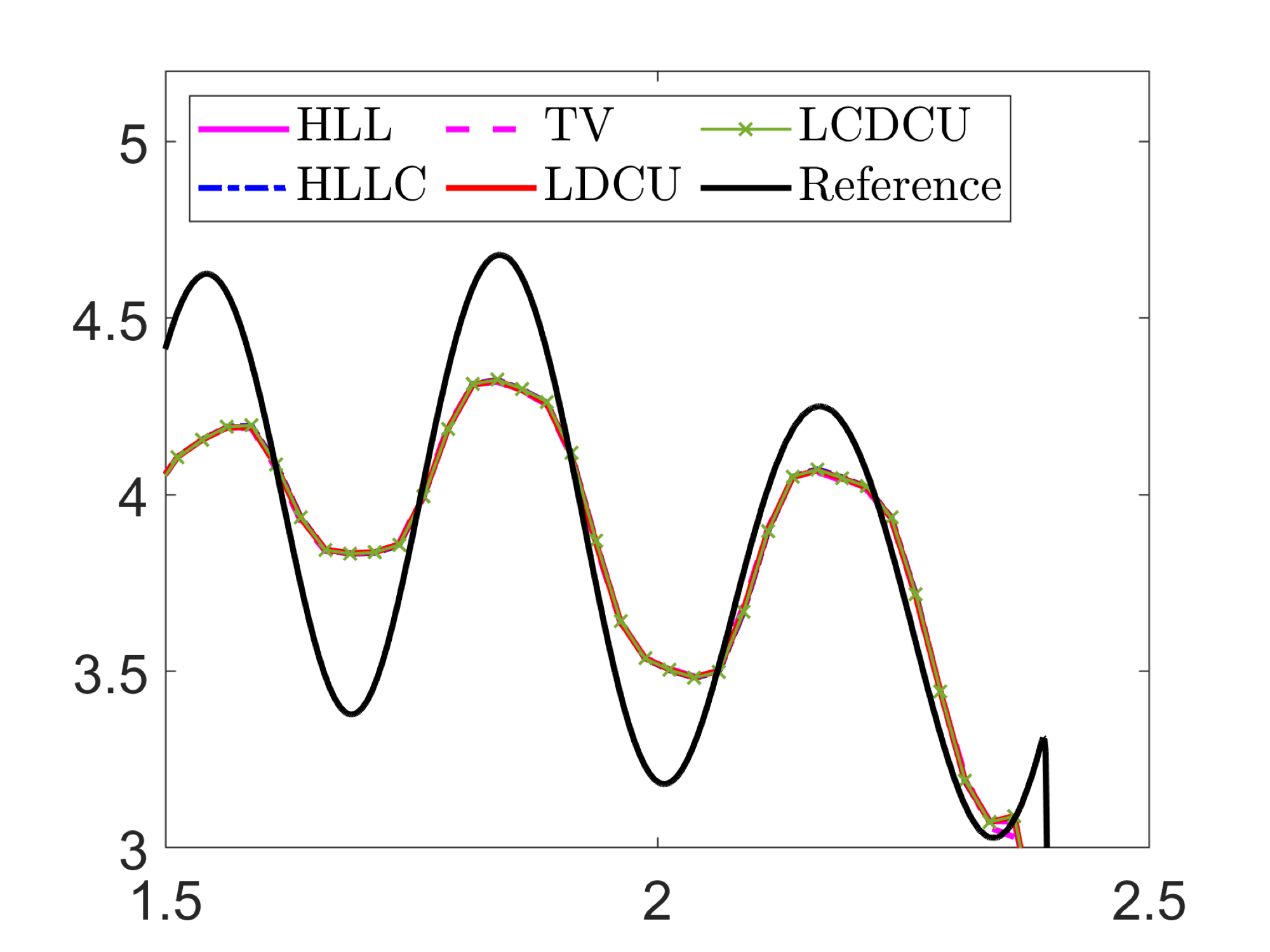}\hspace{0.5cm}
            \includegraphics[trim=0.8cm 0.3cm 0.9cm 0.8cm, clip, width=4.5cm]{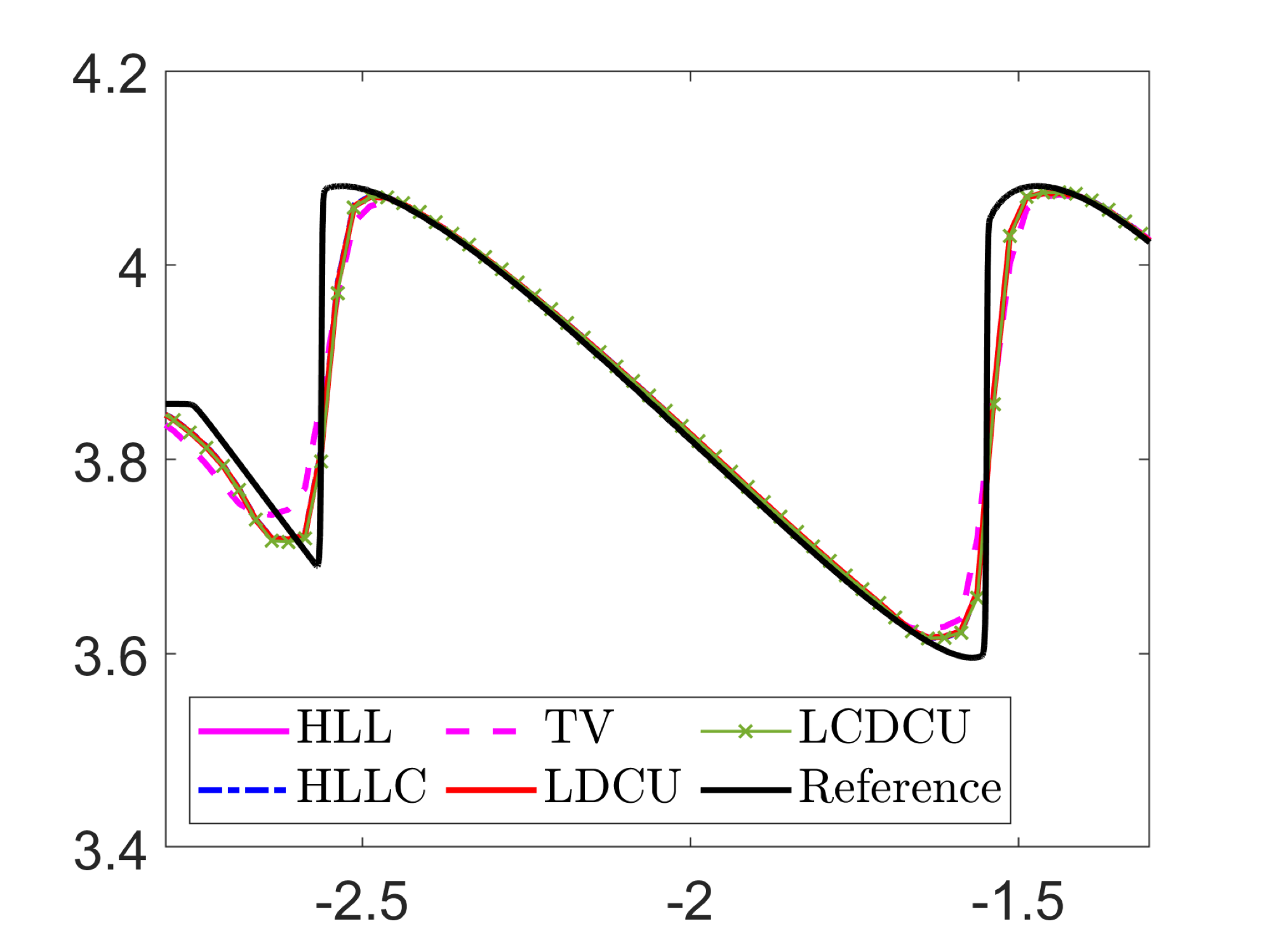}}
            \vskip 12pt
\centerline{\includegraphics[trim=0.8cm 0.3cm 0.9cm 0.8cm, clip, width=4.5cm]{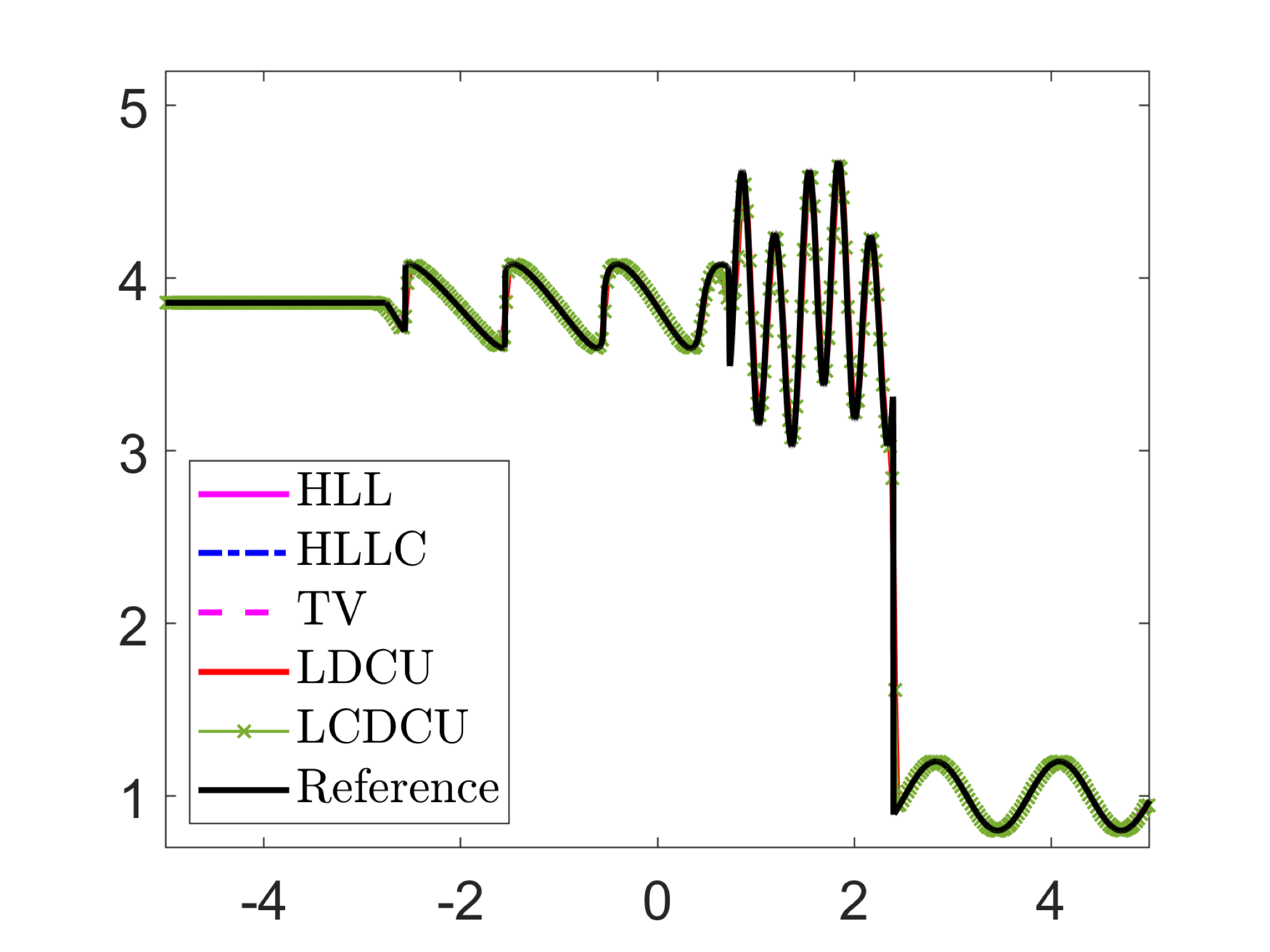}\hspace{0.5cm}
            \includegraphics[trim=0.8cm 0.3cm 0.9cm 0.8cm, clip, width=4.5cm]{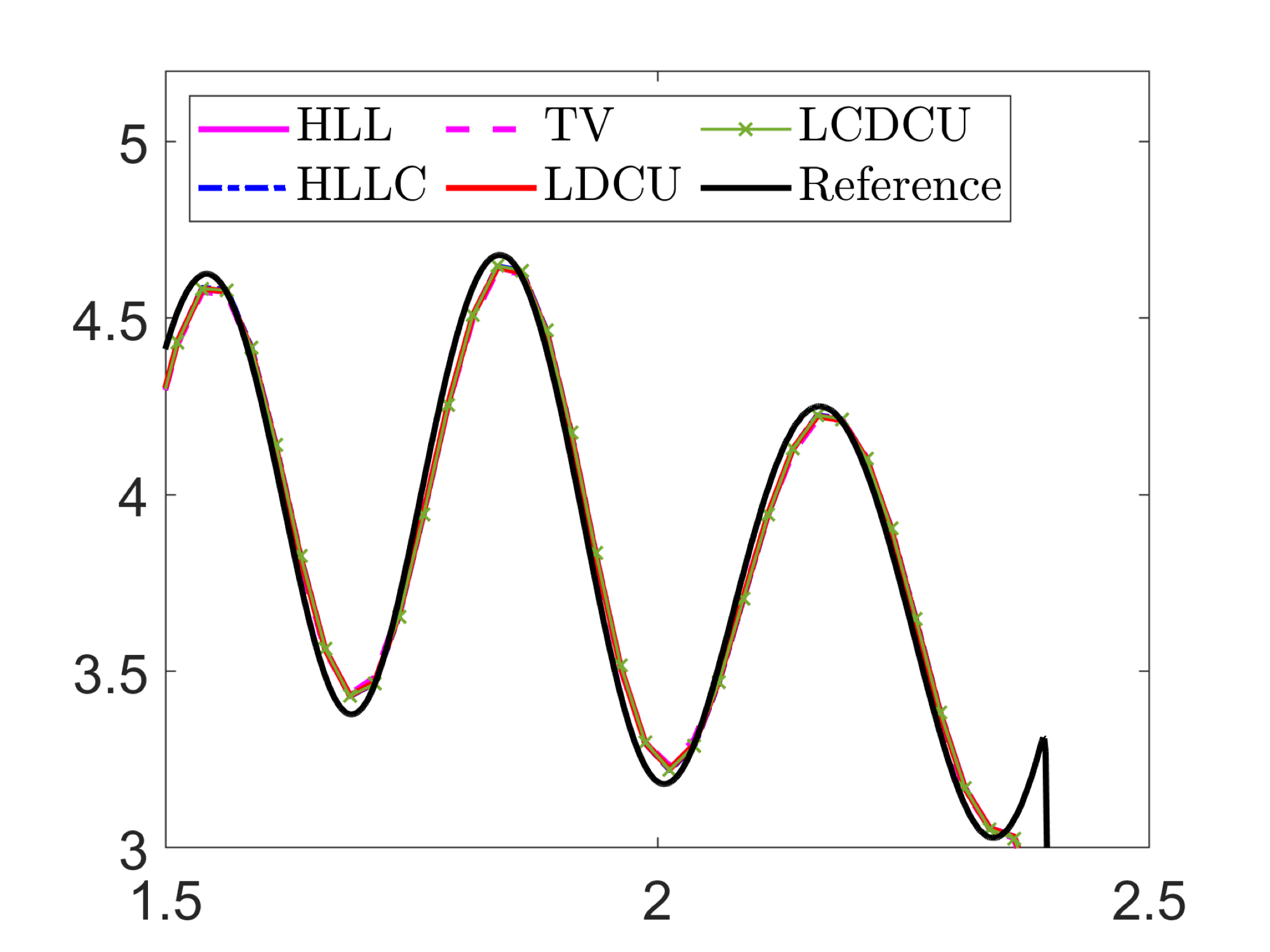}\hspace{0.5cm}
            \includegraphics[trim=0.8cm 0.3cm 0.9cm 0.8cm, clip, width=4.5cm]{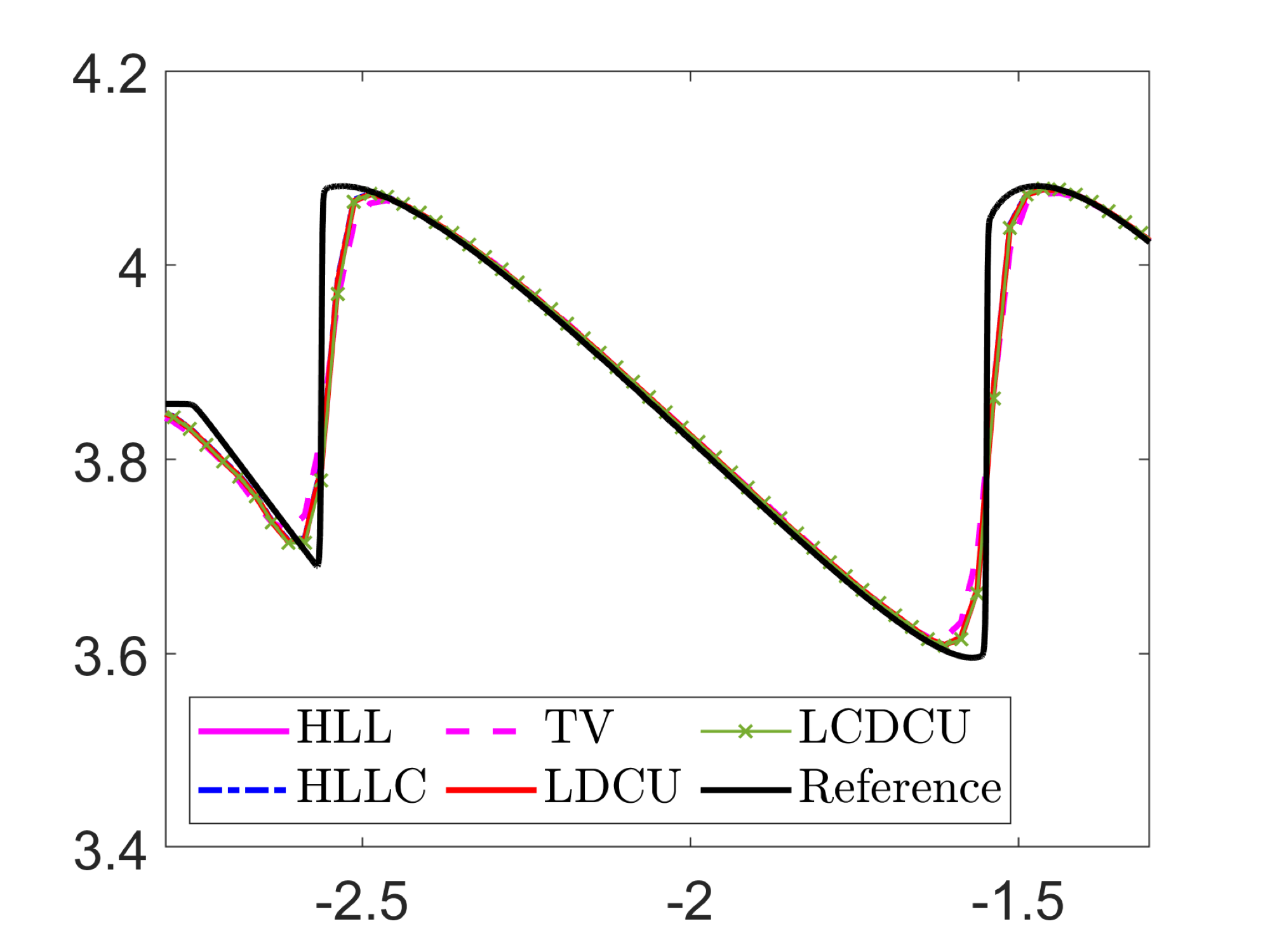}}
\caption{\sf Example 5: Density $\rho$ computed by the 3- and 5-Order schemes (left column) and zoom at $[1.5,2.5]$ (middle column), and  $[-2.8,-1.3]$ (right column).
\label{fig5}}
\end{figure}

We compute the numerical solutions until the final time $t=5$ by the 1-Order, 2-Order, 3-Order, and 5-Order schemes on a uniform mesh of 400 cells, and present the obtained numerical results in Figures \ref{fig5a}--\ref{fig5} together with the reference solution computed by the HLL scheme on a much finer mesh of 8000 cells. It can be seen from this example that the TV schemes exhibit higher dissipation than the other schemes, including the HLL scheme, near the shock waves even for the high-order schemes. However, the TV schemes achieve better resolution than HLL schemes in the smooth parts of the results. 

\subsubsection*{Example 6---Shock-Entropy Wave Interaction Problem}
In this 1-D example, we consider the shock-entropy problem taken from \cite{Shu88}. The initial conditions,
\begin{equation*}
(\rho,u,p)(x,0)=\begin{cases}
(1.51695,0.523346,1.805),&x<-4.5,\\
(1+0.1\sin(20x),0,1),&x>-4.5,
\end{cases}
\end{equation*}
correspond to a forward-facing shock wave of Mach 1.1 interacting with high-frequency density perturbations, that is, as the shock wave
moves, the perturbations spread ahead. In this example, the free boundary conditions are imposed at both ends of the computational domain $[-10,5]$.

\begin{figure}[ht!]
\centerline{\includegraphics[trim=0.8cm 0.3cm 0.8cm 0.8cm, clip, width=4.5cm]{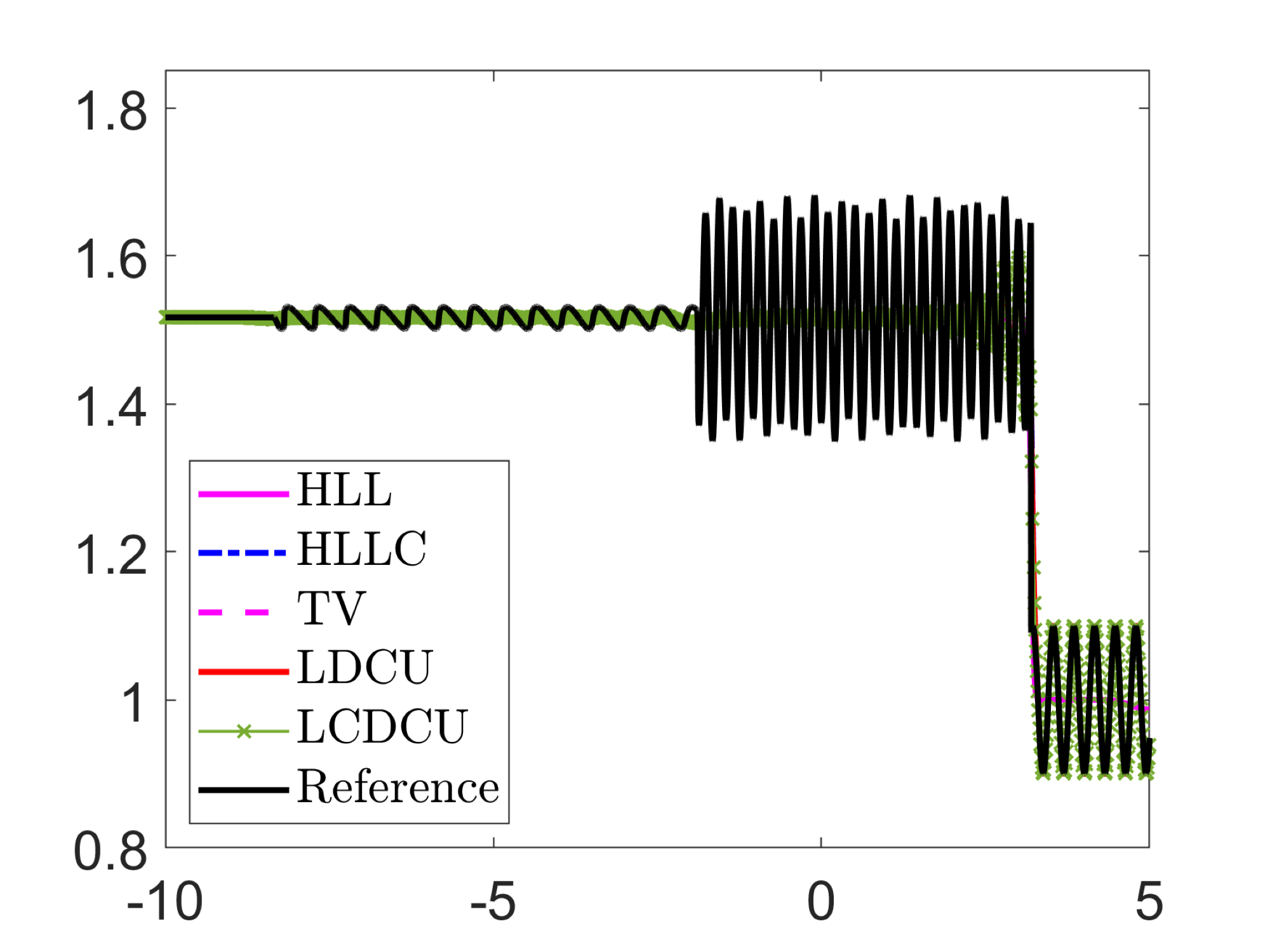}\hspace{0.5cm}
            \includegraphics[trim=0.8cm 0.3cm 0.8cm 0.8cm, clip, width=4.5cm]{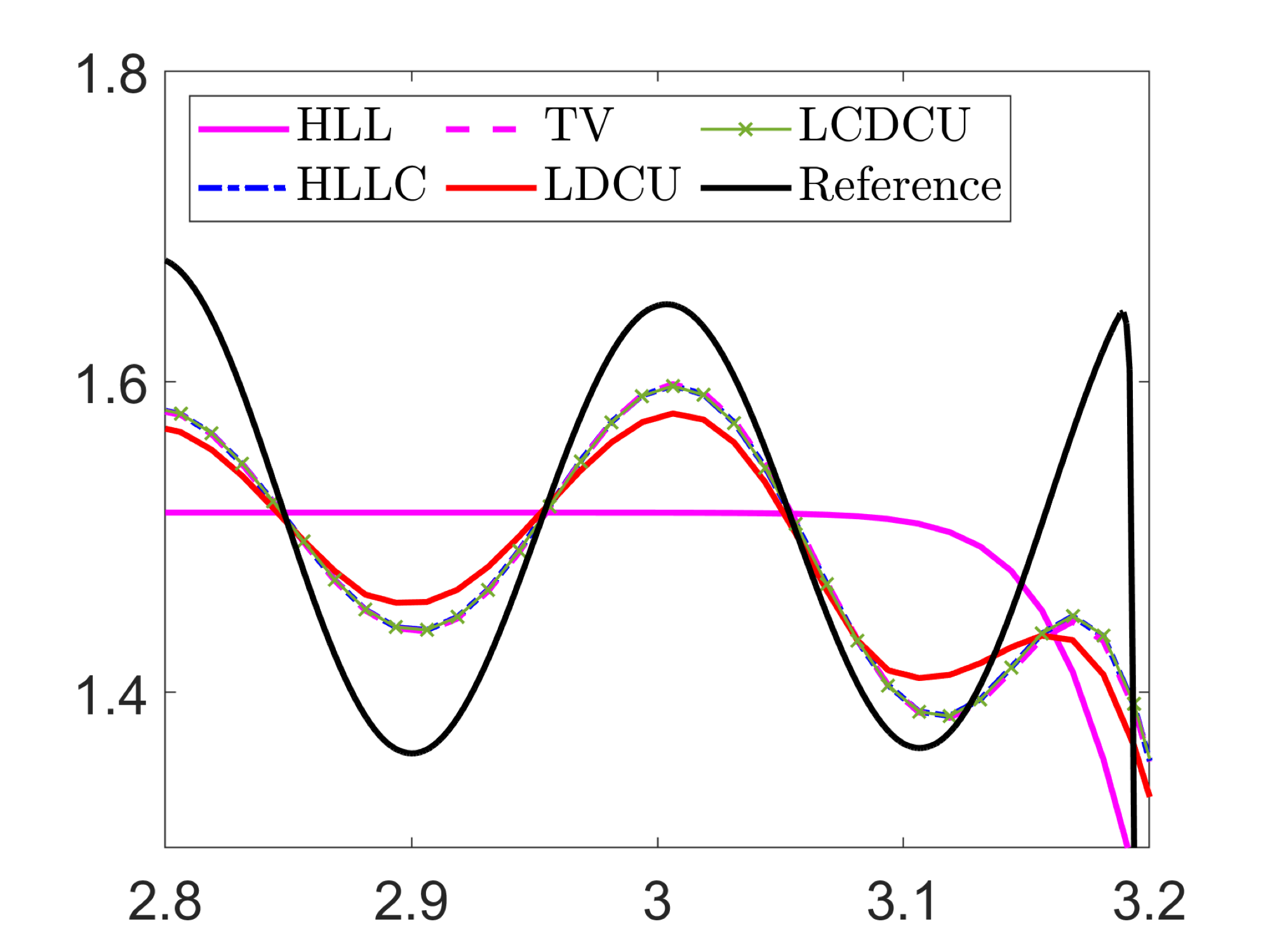}\hspace{0.5cm}
            \includegraphics[trim=0.8cm 0.3cm 0.8cm 0.8cm, clip, width=4.5cm]{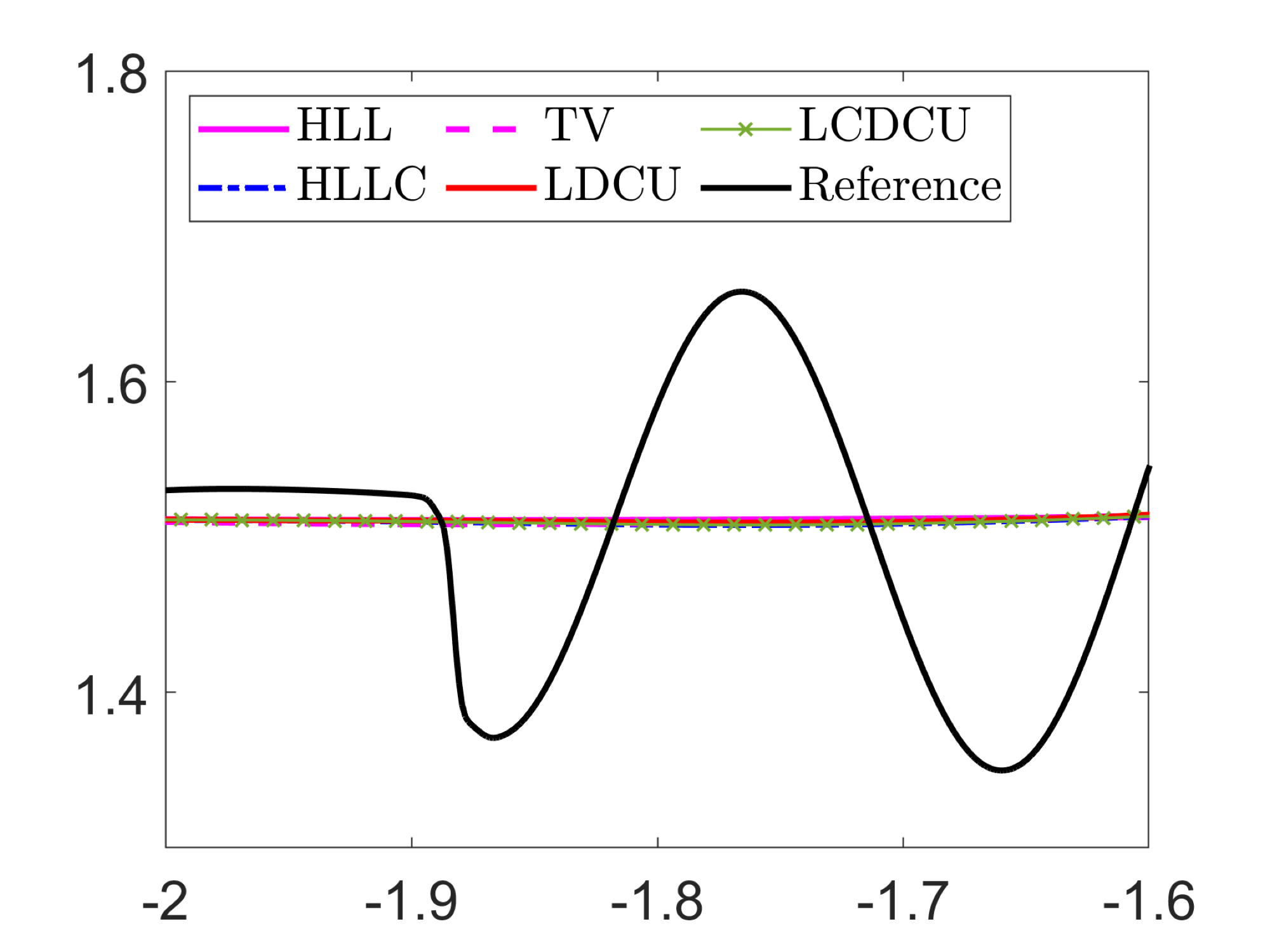}}
\vskip 12pt
\centerline{\includegraphics[trim=0.8cm 0.3cm 0.8cm 0.8cm, clip, width=4.5cm]{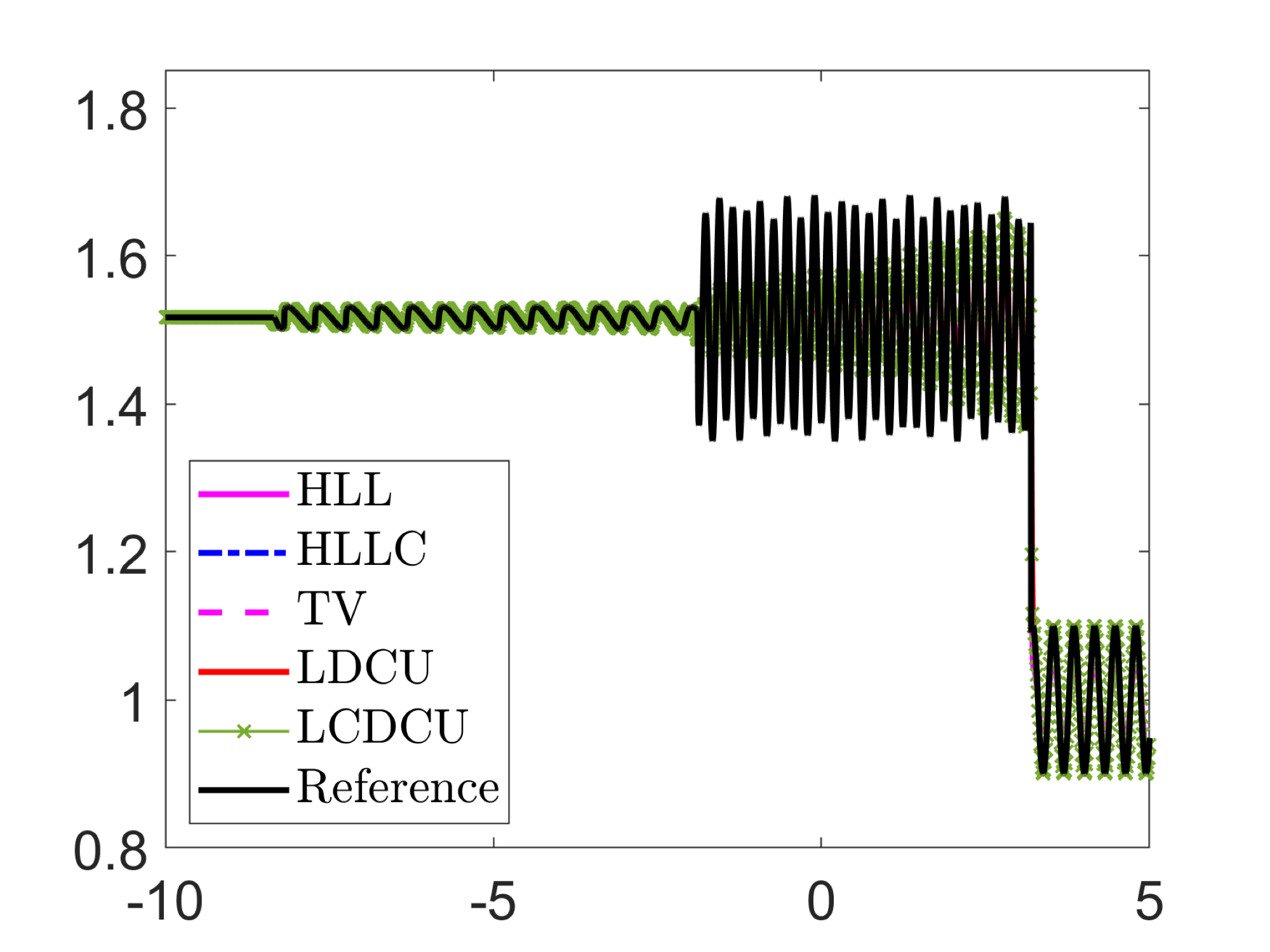}\hspace{0.5cm}
            \includegraphics[trim=0.8cm 0.3cm 0.8cm 0.8cm, clip, width=4.5cm]{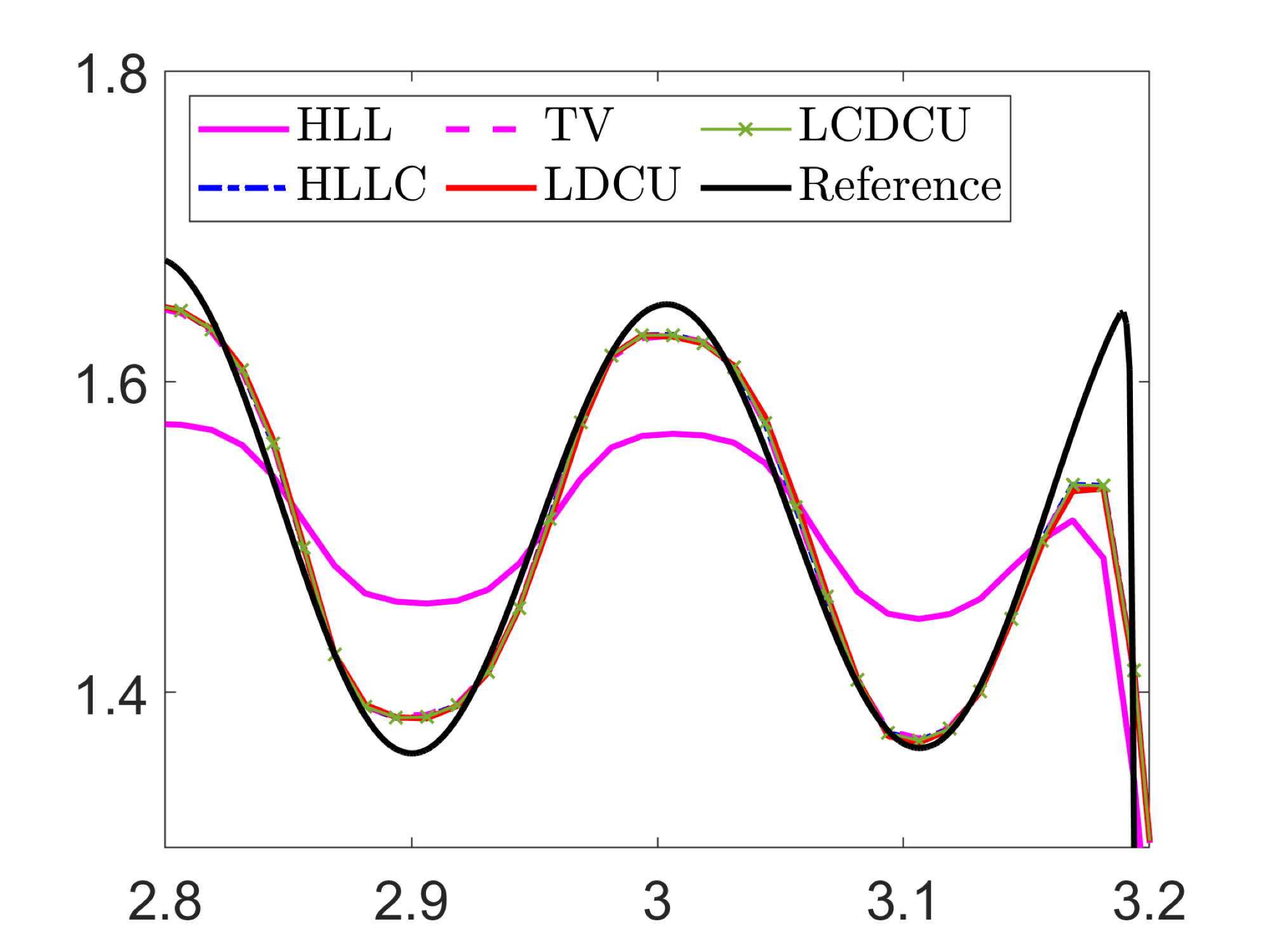}\hspace{0.5cm}
            \includegraphics[trim=0.8cm 0.3cm 0.8cm 0.8cm, clip, width=4.5cm]{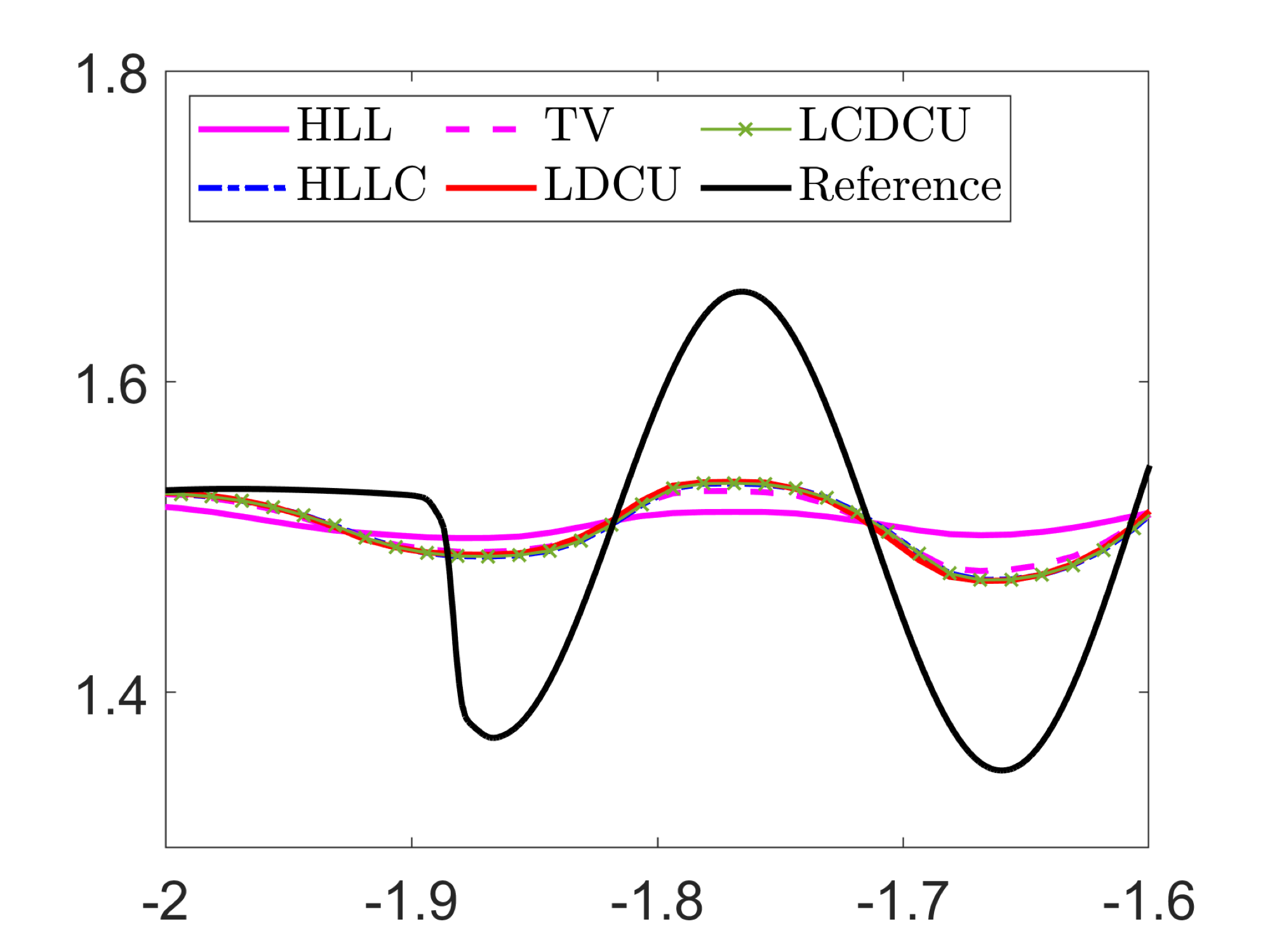}}
\caption{\sf Example 6: Density $\rho$ computed by the 1-Order and 2-Order schemes (left) and zoom at $[2.8,3.2]$ (middle) and $[-2,-1.6]$ (right).
\label{fig6a}}
\end{figure}

\begin{figure}[ht!]
\centerline{\includegraphics[trim=0.8cm 0.3cm 0.8cm 0.8cm, clip, width=4.5cm]{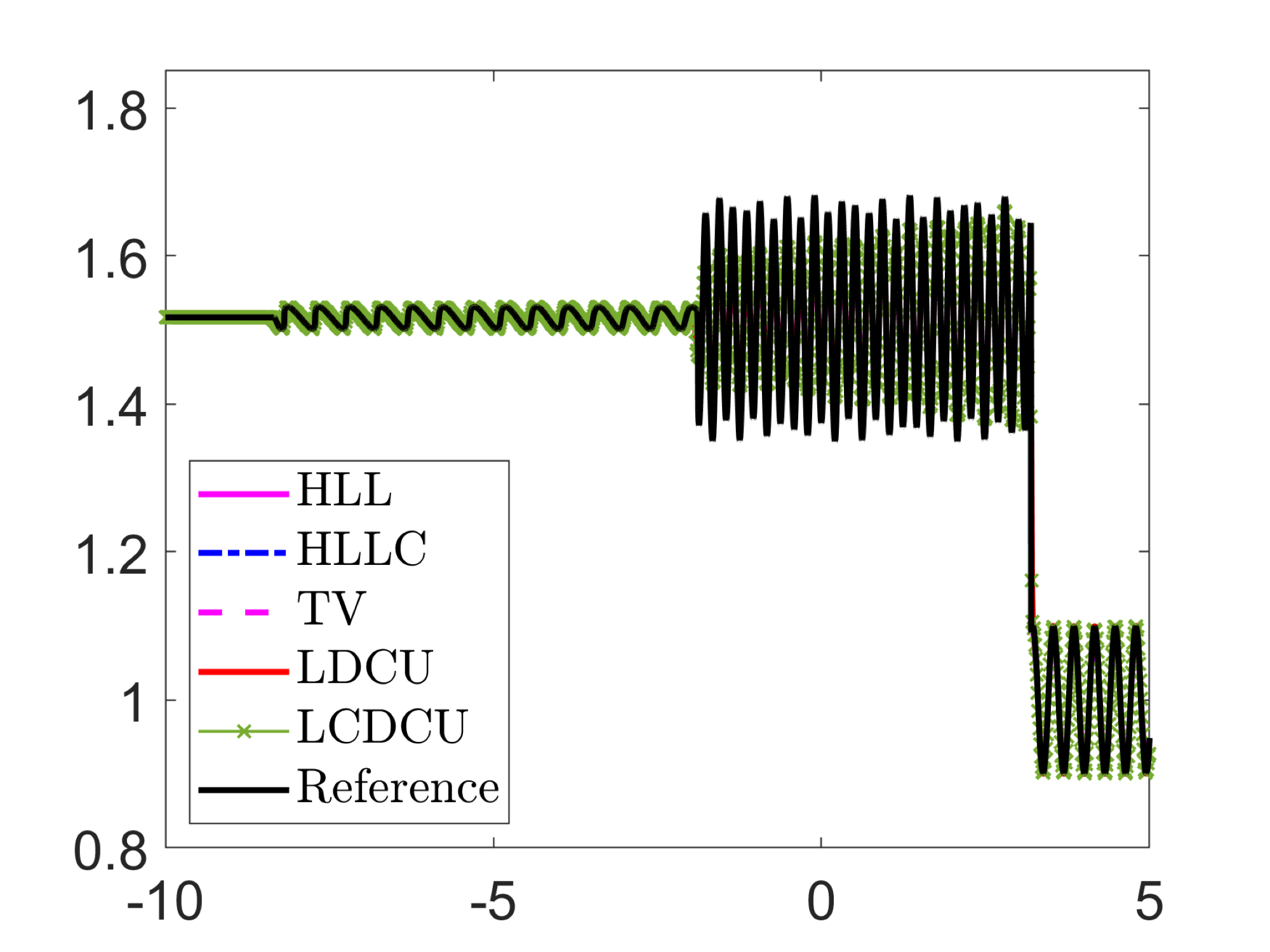}\hspace{0.5cm}
            \includegraphics[trim=0.8cm 0.3cm 0.8cm 0.8cm, clip, width=4.5cm]{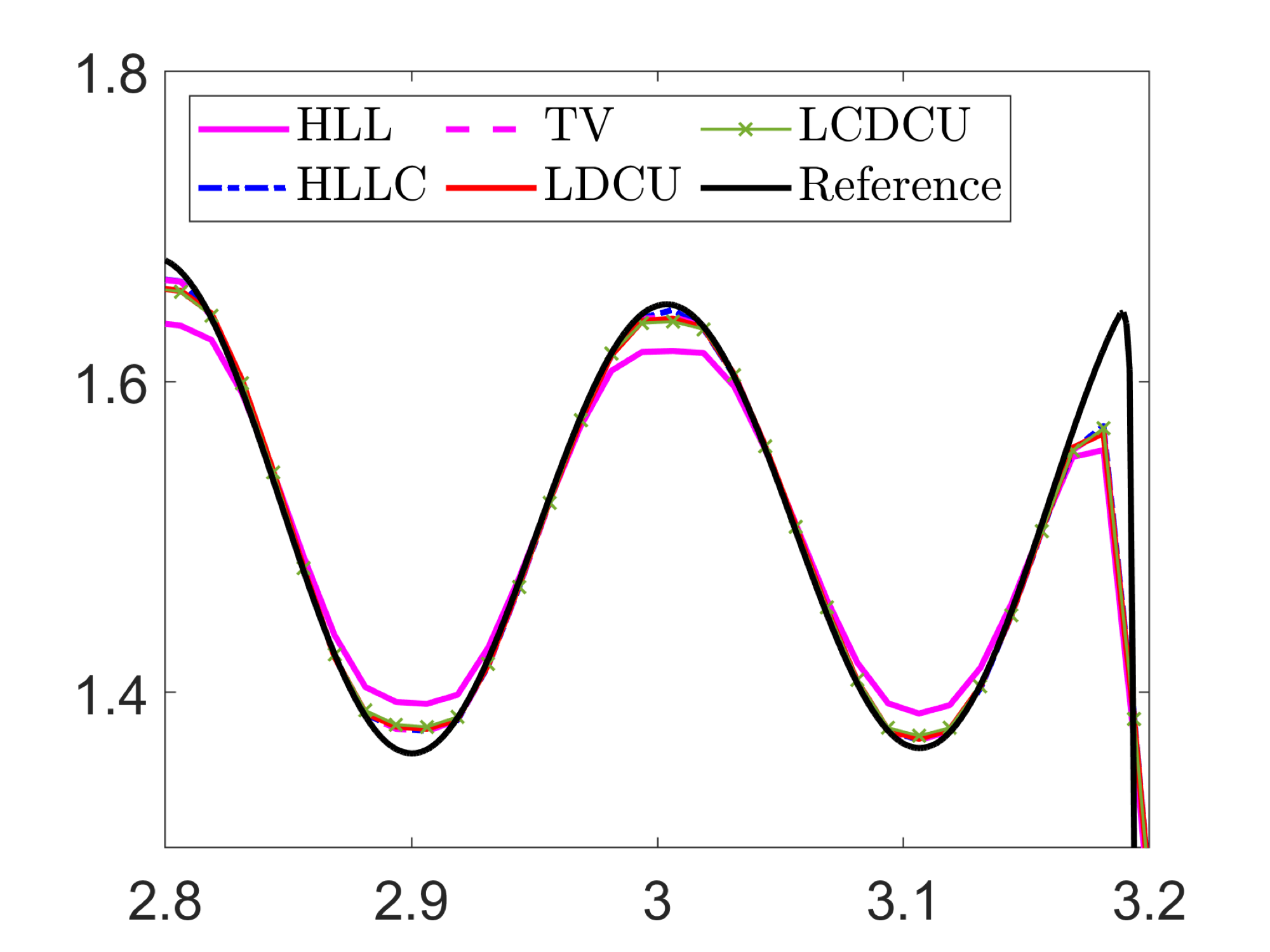}\hspace{0.5cm}
            \includegraphics[trim=0.8cm 0.3cm 0.8cm 0.8cm, clip, width=4.5cm]{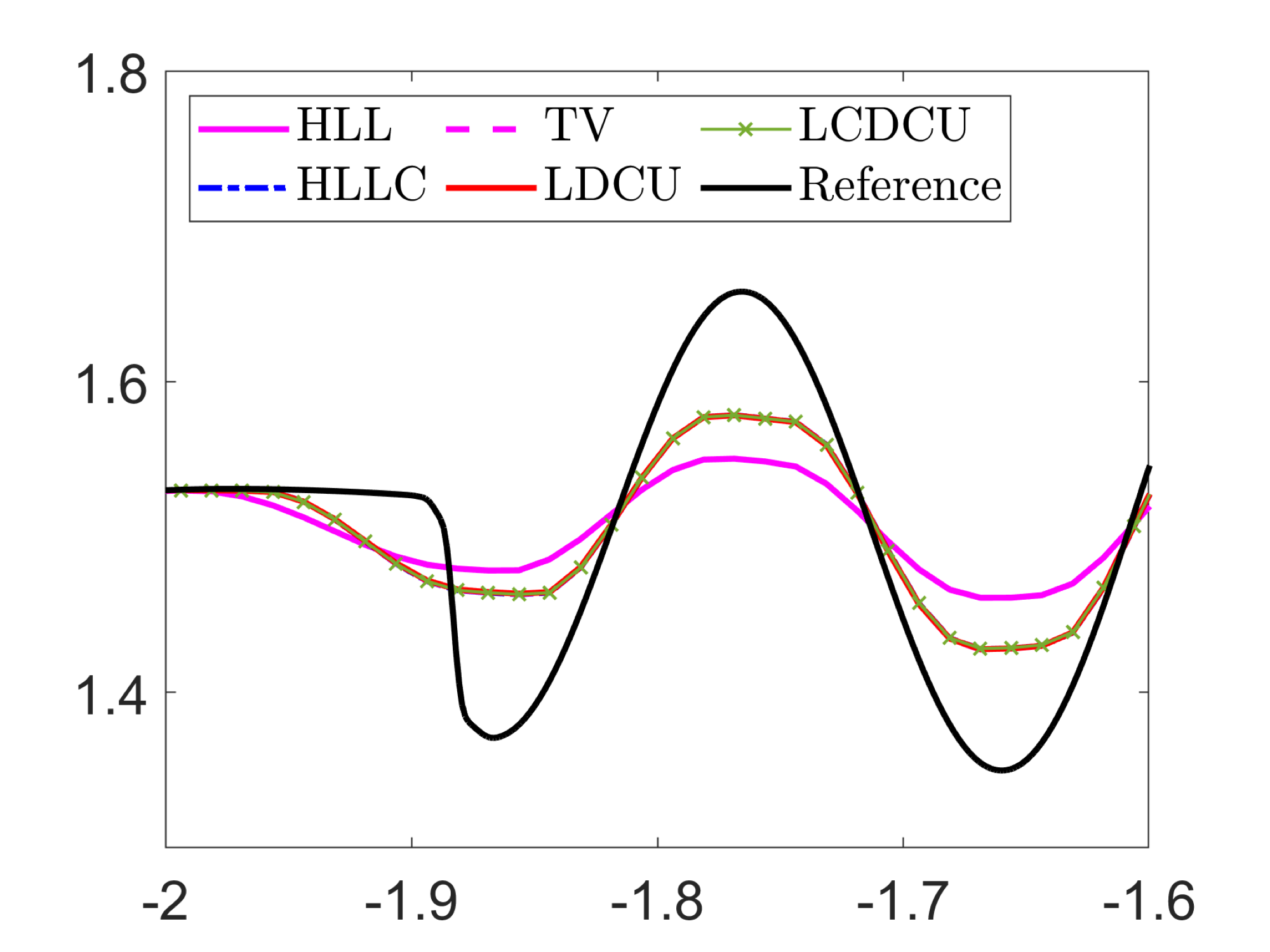}}
            \vskip 12pt
\centerline{\includegraphics[trim=0.8cm 0.3cm 0.8cm 0.8cm, clip, width=4.5cm]{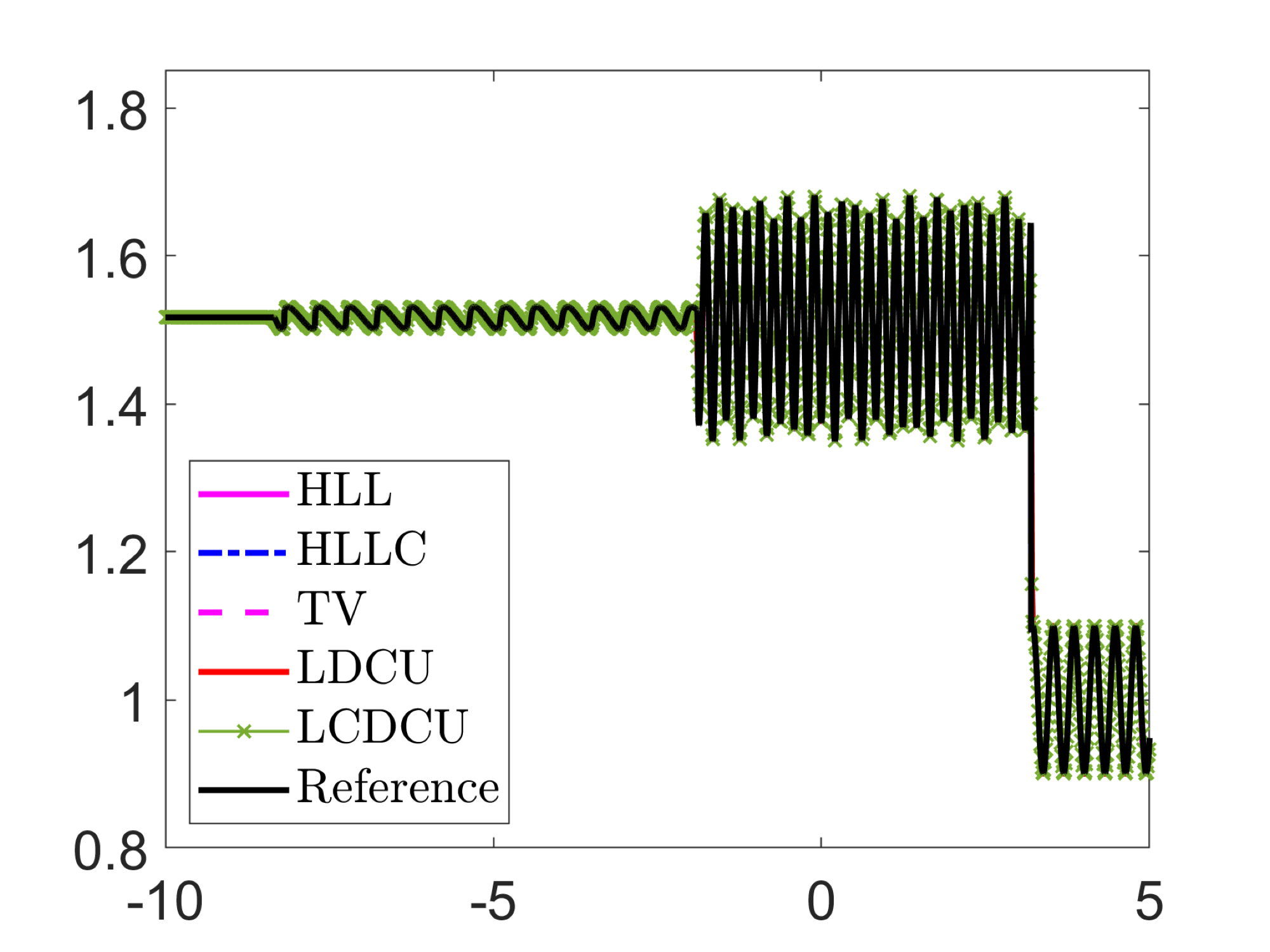}\hspace{0.5cm}
            \includegraphics[trim=0.8cm 0.3cm 0.8cm 0.8cm, clip, width=4.5cm]{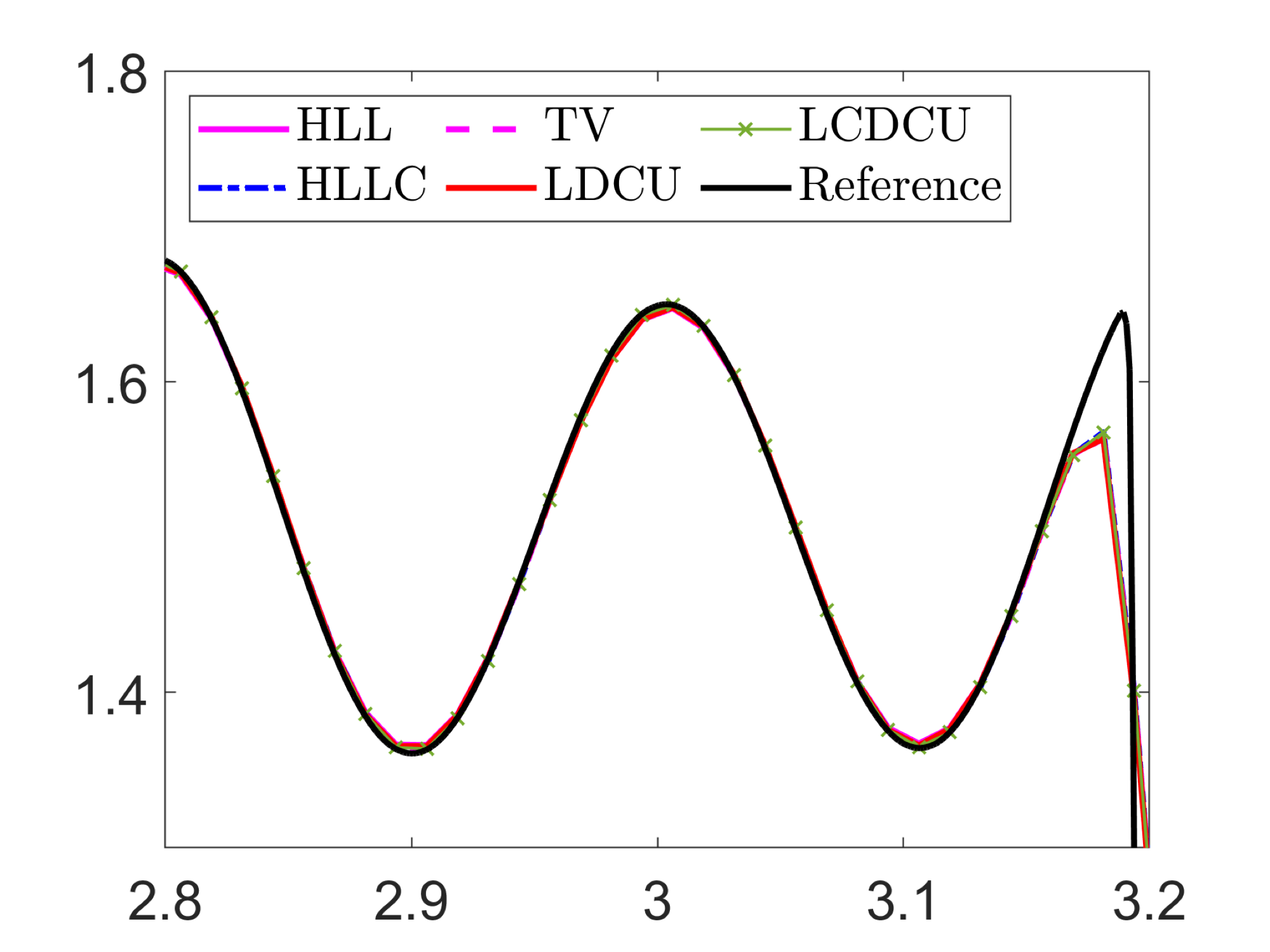}\hspace{0.5cm}
            \includegraphics[trim=0.8cm 0.3cm 0.8cm 0.8cm, clip, width=4.5cm]{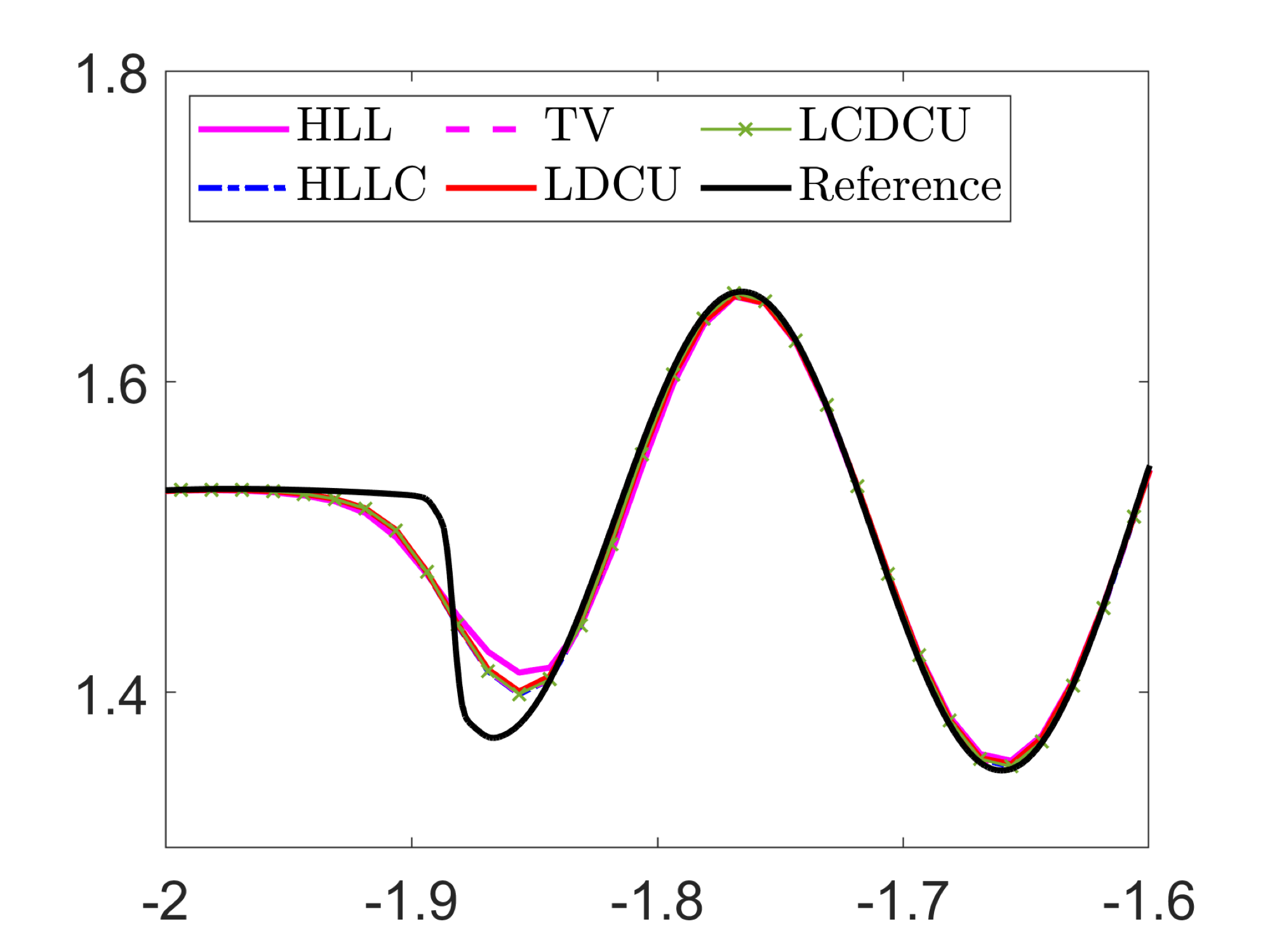}}
\caption{\sf Example 6: Density $\rho$ computed by the 3-Order and 5-Order schemes (left) and zoom at $[2.8,3.2]$ (middle) and $[-2,-1.6]$ (right).
\label{fig6}}
\end{figure}

We compute the solutions until the final time $t=5$ by the 1-Order, 2-Order, 3-Order, and 5-Order schemes on a uniform mesh of 1200 cells. The numerical results are
shown in Figures \ref{fig6a}--\ref{fig6} along with the reference solution computed by the HLL scheme on a much finer mesh of 12000 cells. The obtained results clearly demonstrate a substantial difference in the resolution computed by HLL and the four low-dissipation schemes. One can also see that in the 1-Order results, the LDCU scheme has slightly more dissipation than the TV, HLLC, and LCDCU schemes, but still significantly better than the HLL scheme. At the same time, the four low-dissipation schemes coincide when extended to higher orders; see Figures \ref{fig6a} (bottom) and \ref{fig6}.

\subsubsection*{Example 7---High-Mach Moving Shock} 
In this example, adapted from \cite{Hennemann21}, we consider the one-dimensional version of a right-moving Mach-100 shock wave. The initial conditions are given by 
\begin{equation*}
(\rho,u,p)(x,0)=\begin{cases}
(5.9970015,\,98.59147,\,11666.5),&x<0.2,\\
(1,0,1),& x>0.2,
\end{cases}
\end{equation*}
prescribed in the computational domain $[0,1]$. The initial data correspond to a right-moving Mach-100 shock wave. At the left boundary, the post-shock state is imposed, while at the right boundary, the pre-shock state is imposed. More precisely, the boundary conditions are given by
$$
(\rho,u,p)(0,t)=(5.9970015,\,98.59147,\,11666.5),
\qquad
(\rho,u,p)(1,t)=(1,0,1).
$$

\begin{figure}[ht!]
\centerline{\includegraphics[trim=0.cm 0.cm 0.cm 0.cm, clip, width=4.cm]{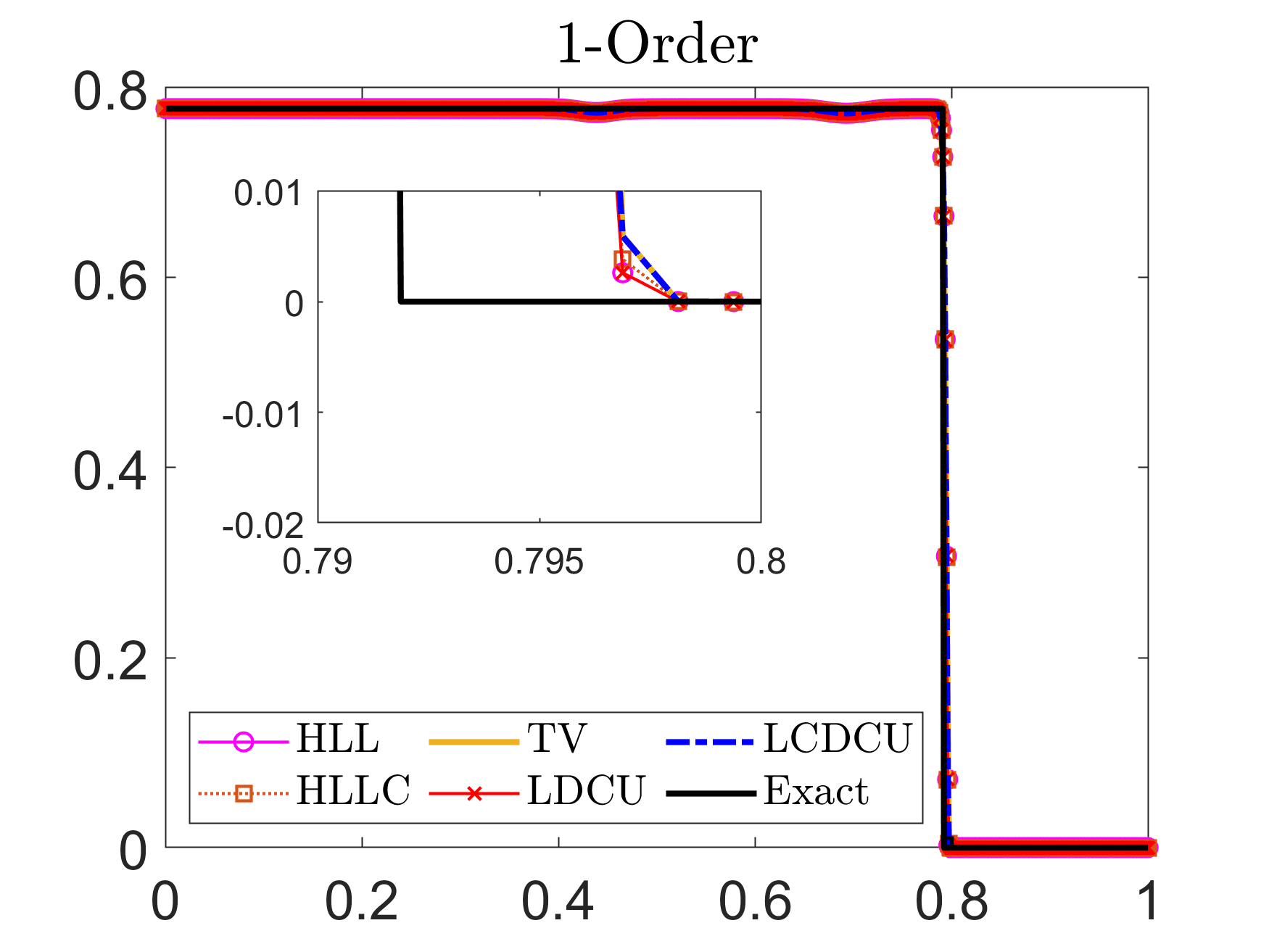}
            \includegraphics[trim=0.cm 0.cm 0.cm 0.cm, clip, width=4.cm]{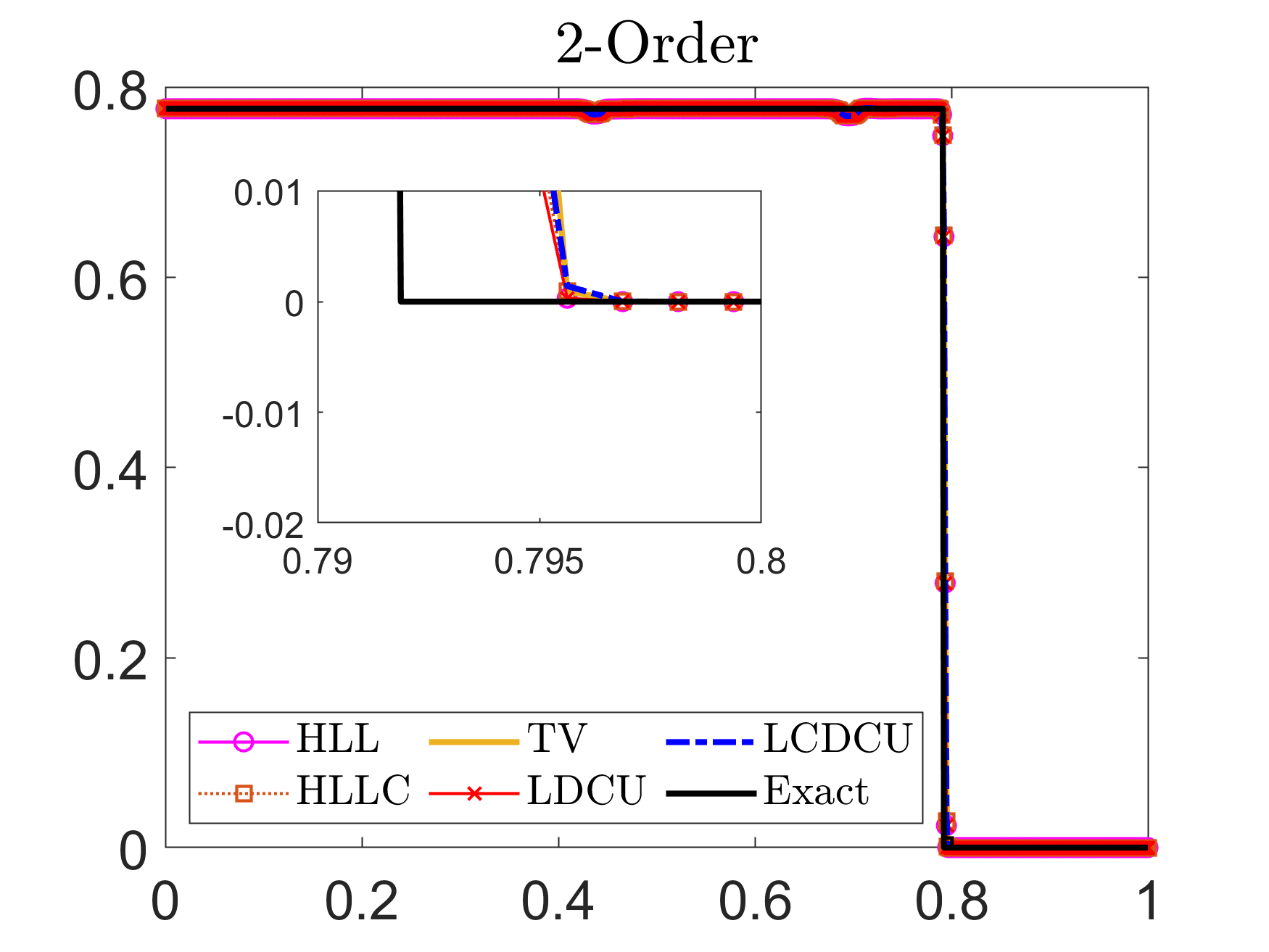}
            \includegraphics[trim=0.cm 0.cm 0.cm 0.cm, clip, width=4.cm]{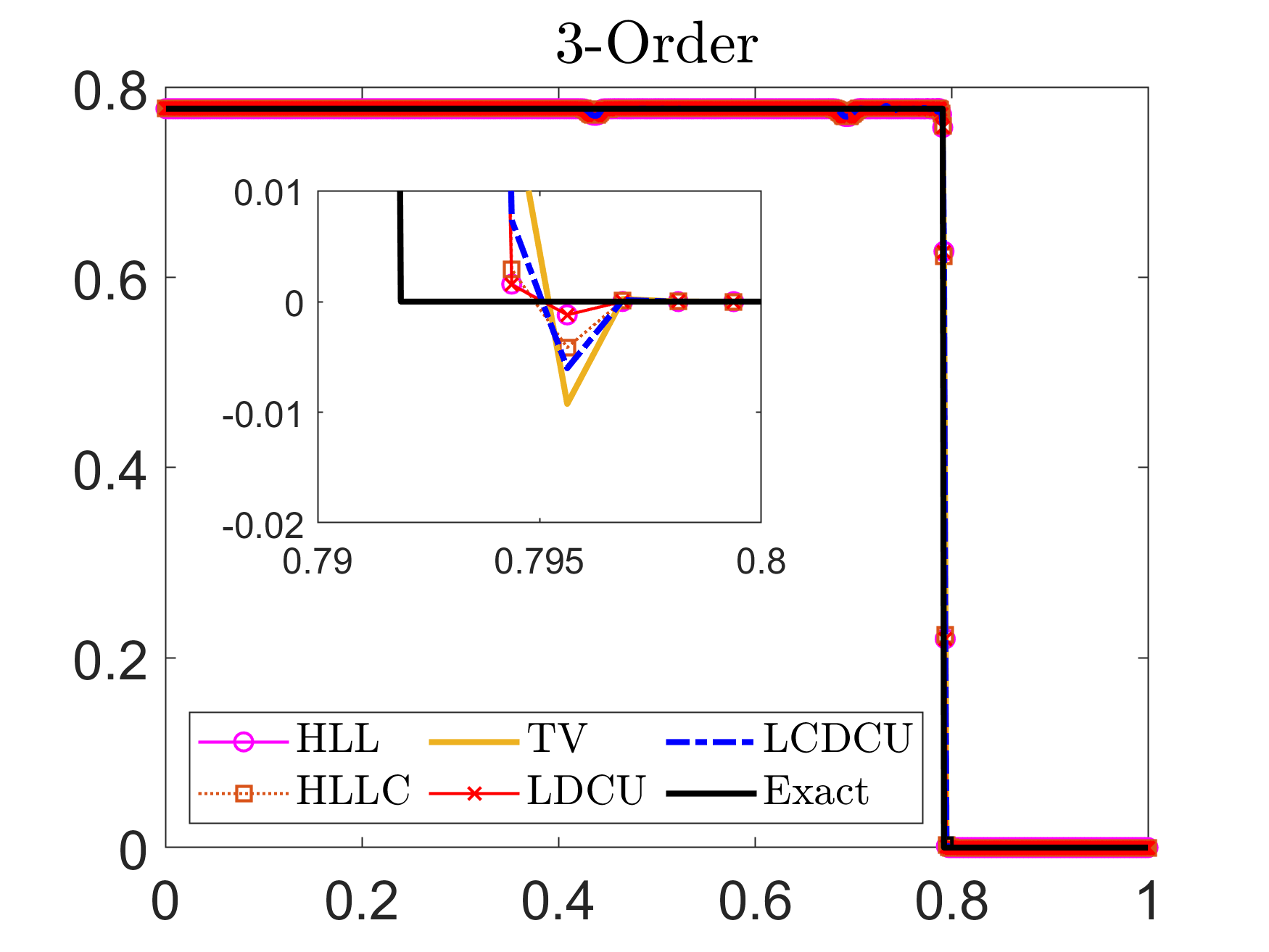}
            \includegraphics[trim=0.cm 0.cm 0.cm 0.cm, clip, width=4.cm]{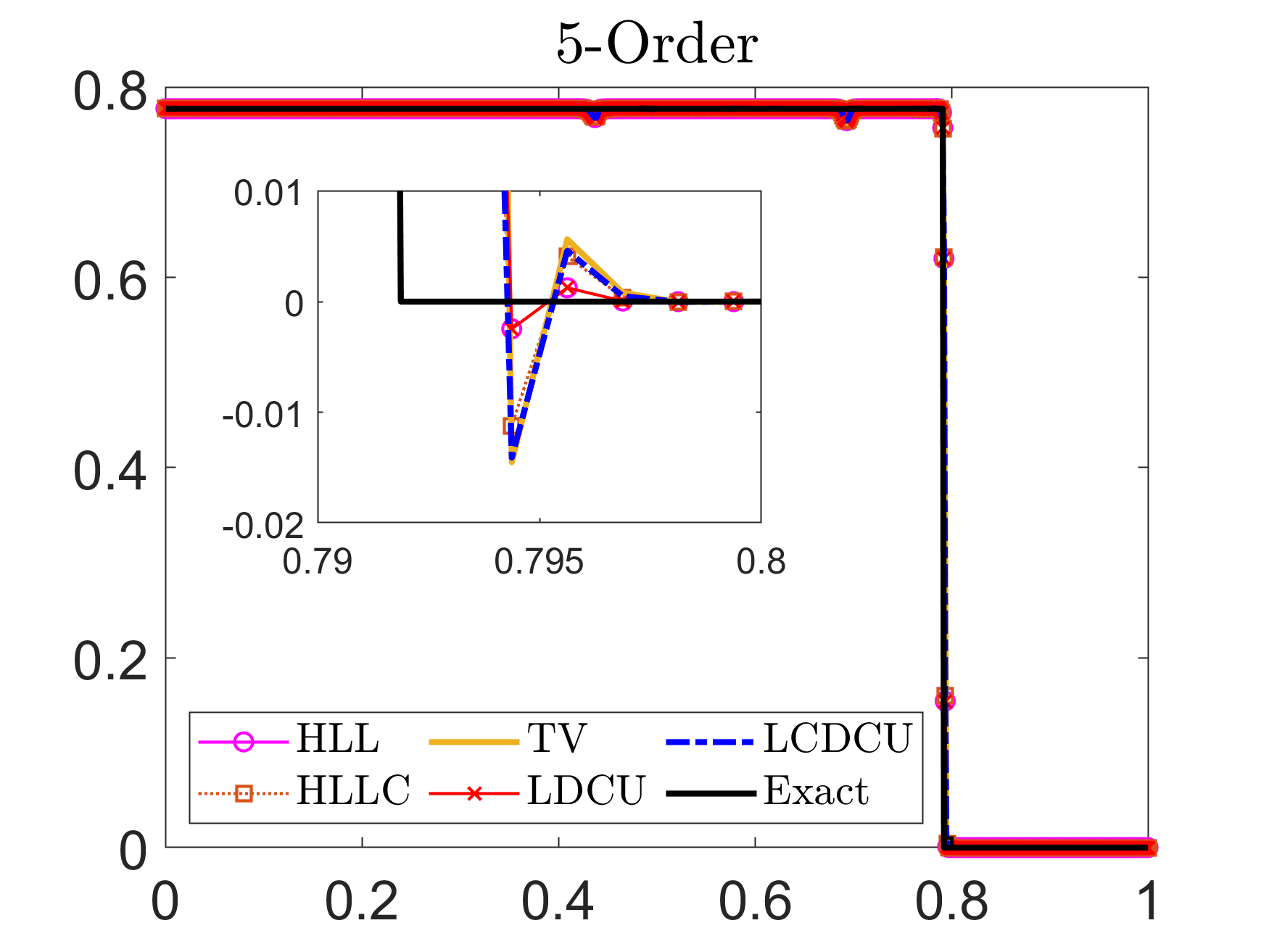}}
\caption{\sf Example 7: Decimal logarithm of the density $\log_{10}(\rho)$ computed by the HLL, HLLC, TV, LDCU, and LCDCU schemes together with the exact solution.
\label{fig6aa}}
\end{figure}

We compute the numerical solutions until the final time $t=0.005$ using the studied schemes on a uniform mesh of $800$ cells. The exact solution remains a single right-moving shock separating the two constant states, and its location is given by
$$
x_s(0.005)\approx 0.791608.
$$
The numerical results for $\log_{10}(\rho)$ are reported in Figure~\ref{fig6aa} together with the exact solution. One can see that the HLLC, TV, and LCDCU schemes exhibit spurious oscillations near the shock wave, which become more pronounced as the order increases. In contrast, the HLL and LDCU schemes do not exhibit visible numerical shock instabilities.

\subsubsection*{Example 8---Colliding Strong Shocks}
In the last 1-D example, motivated by the colliding-shocks problem in \cite{Lious2006}, we consider the following initial conditions:
$$
(\rho,u,p)(x,0) = \begin{cases}
(1,20,1),  & x<0,\\
(1,-20,1), & x>0,
\end{cases}
$$
which are prescribed in the computational domain $[-1,1]$. The two initial states move toward each other and generate two strong shock waves propagating in opposite directions.
At the left and right boundaries, the corresponding initial states are imposed. More precisely, the boundary conditions are given by
$$
(\rho,u,p)(-1,t)=(1,20,1), \qquad (\rho,u,p)(1,t)=(1,-20,1).
$$

We compute the numerical solutions until the final time $t=0.05$ using the studied schemes on a uniform mesh of $800$ cells. The exact solution consists of two shocks separating the left and right states from a constant intermediate state:
$$
(\rho,u,p)(x,t) =
\begin{cases}
(1,20,1), & x<S_Lt,\\[1mm]
(\rho_*,0,p_*), & S_Lt<x<S_Rt,\\[1mm]
(1,-20,1), & x>S_Rt,
\end{cases}
$$
where
$$
\rho_*\approx5.9283027607, \qquad p_*\approx482.1638447197,
$$
and the speeds of the left- and right-moving shocks are
$$
S_L\approx-4.0581922360, \qquad S_R\approx4.0581922360.
$$
Therefore, at the final time $t=0.05$, the two shocks are located at
$$
x_L=S_Lt\approx-0.202909612, \qquad x_R=S_Rt\approx0.202909612.
$$

\begin{figure}[ht!]
\centerline{\includegraphics[trim=0.cm 0.cm 0.cm 0.cm, clip, width=4.cm]{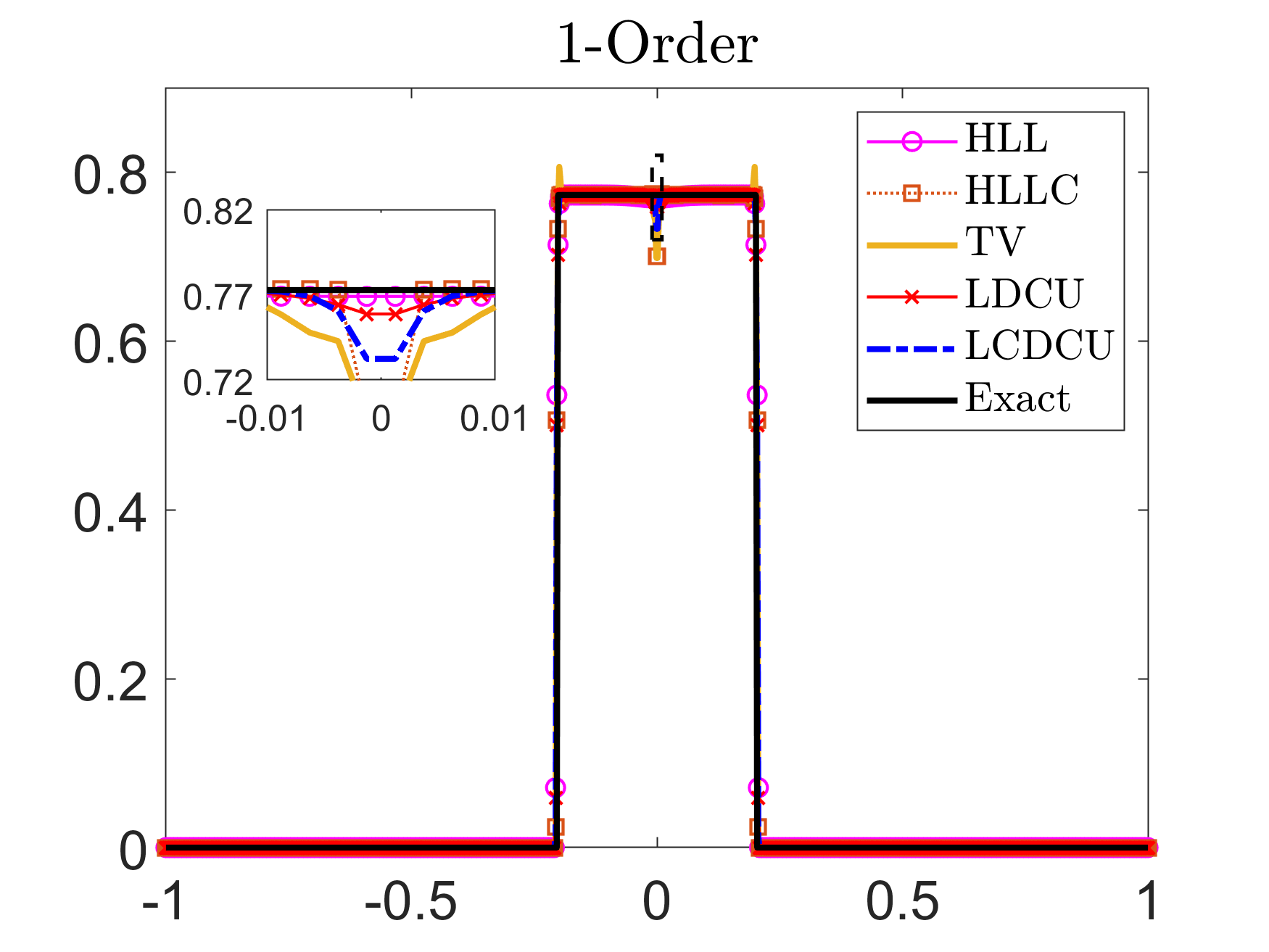}
            \includegraphics[trim=0.cm 0.cm 0.cm 0.cm, clip, width=4.cm]{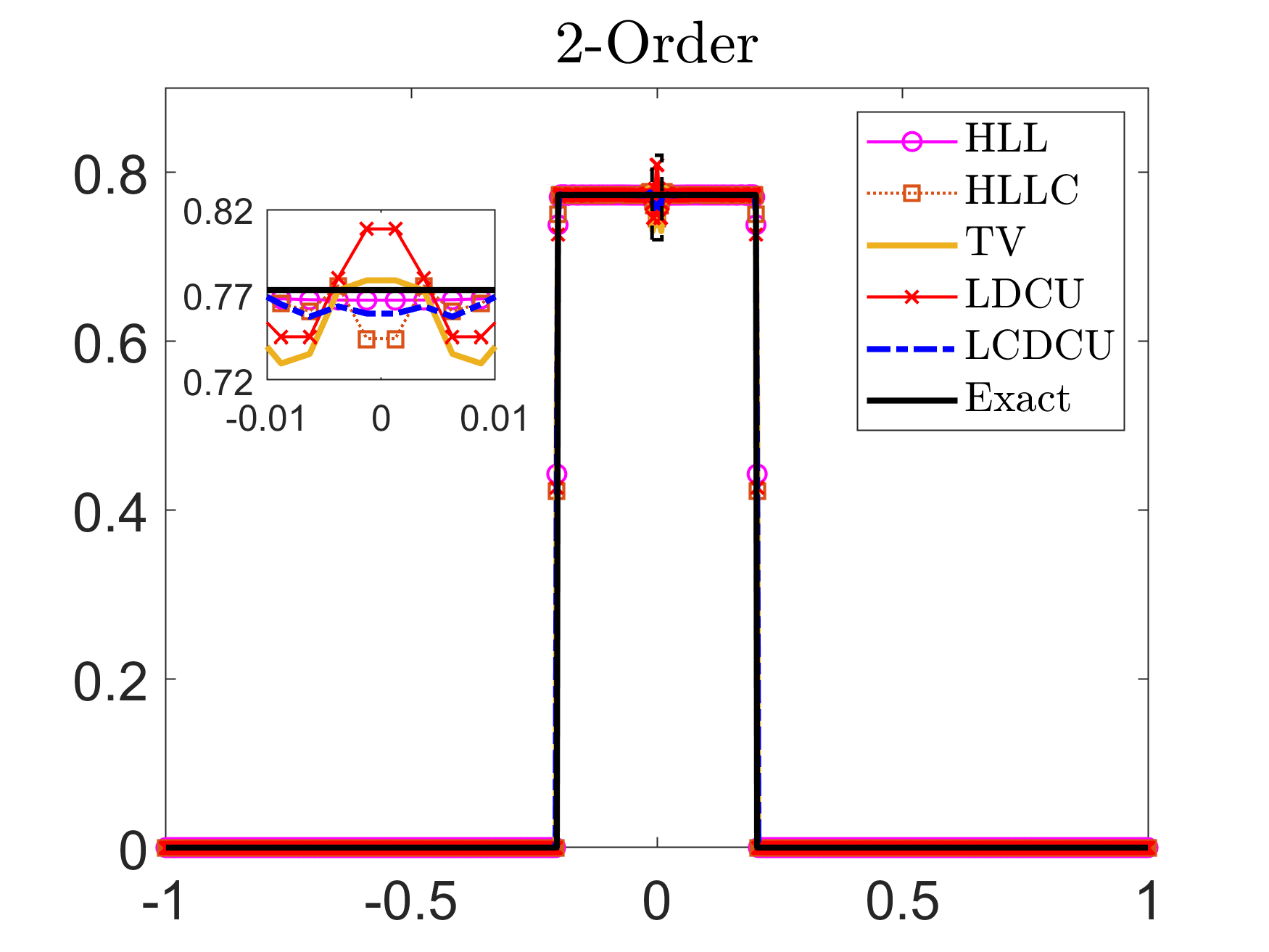}
            \includegraphics[trim=0.cm 0.cm 0.cm 0.cm, clip, width=4.cm]{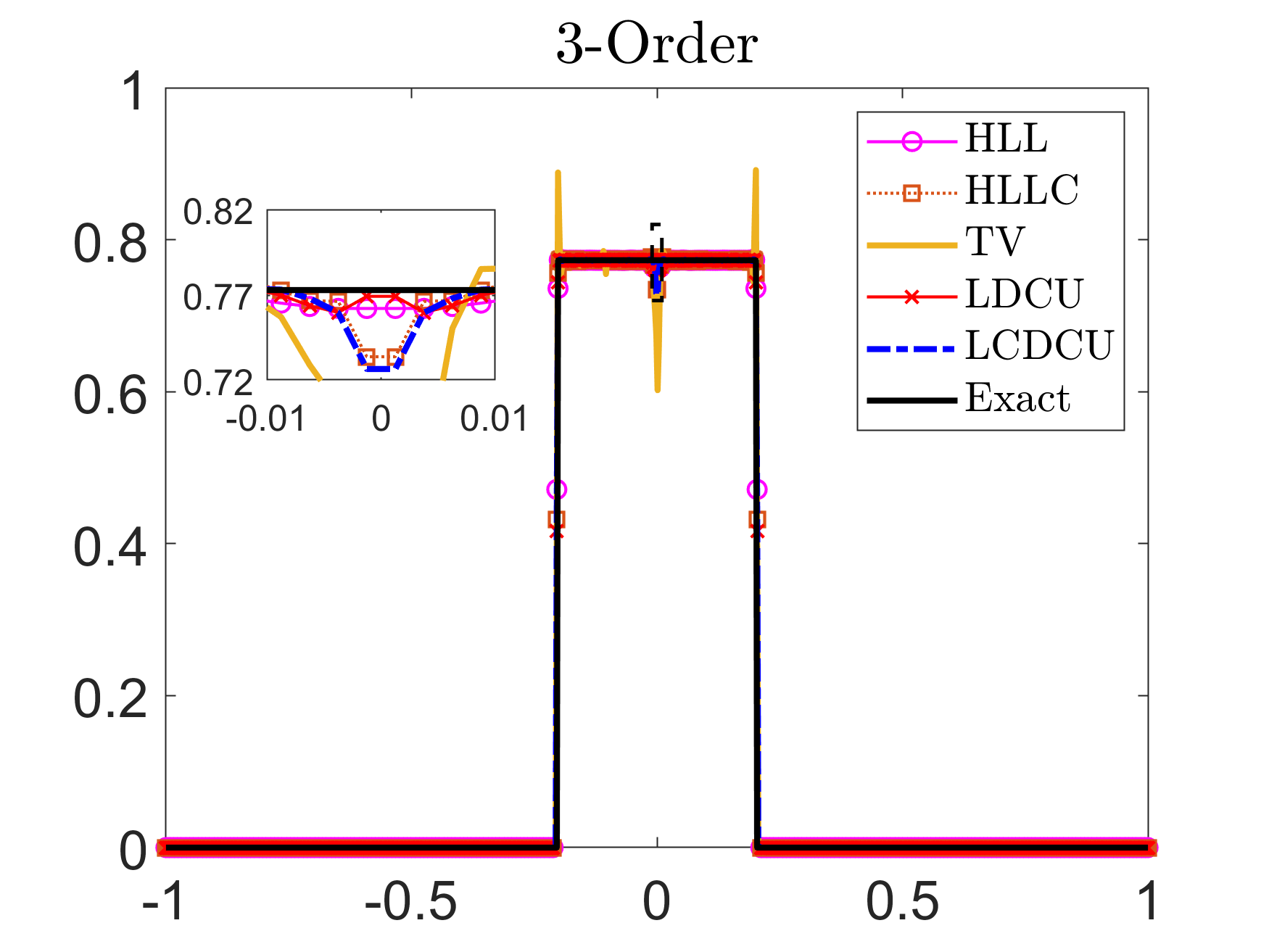}
            \includegraphics[trim=0.cm 0.cm 0.cm 0.cm, clip, width=4.cm]{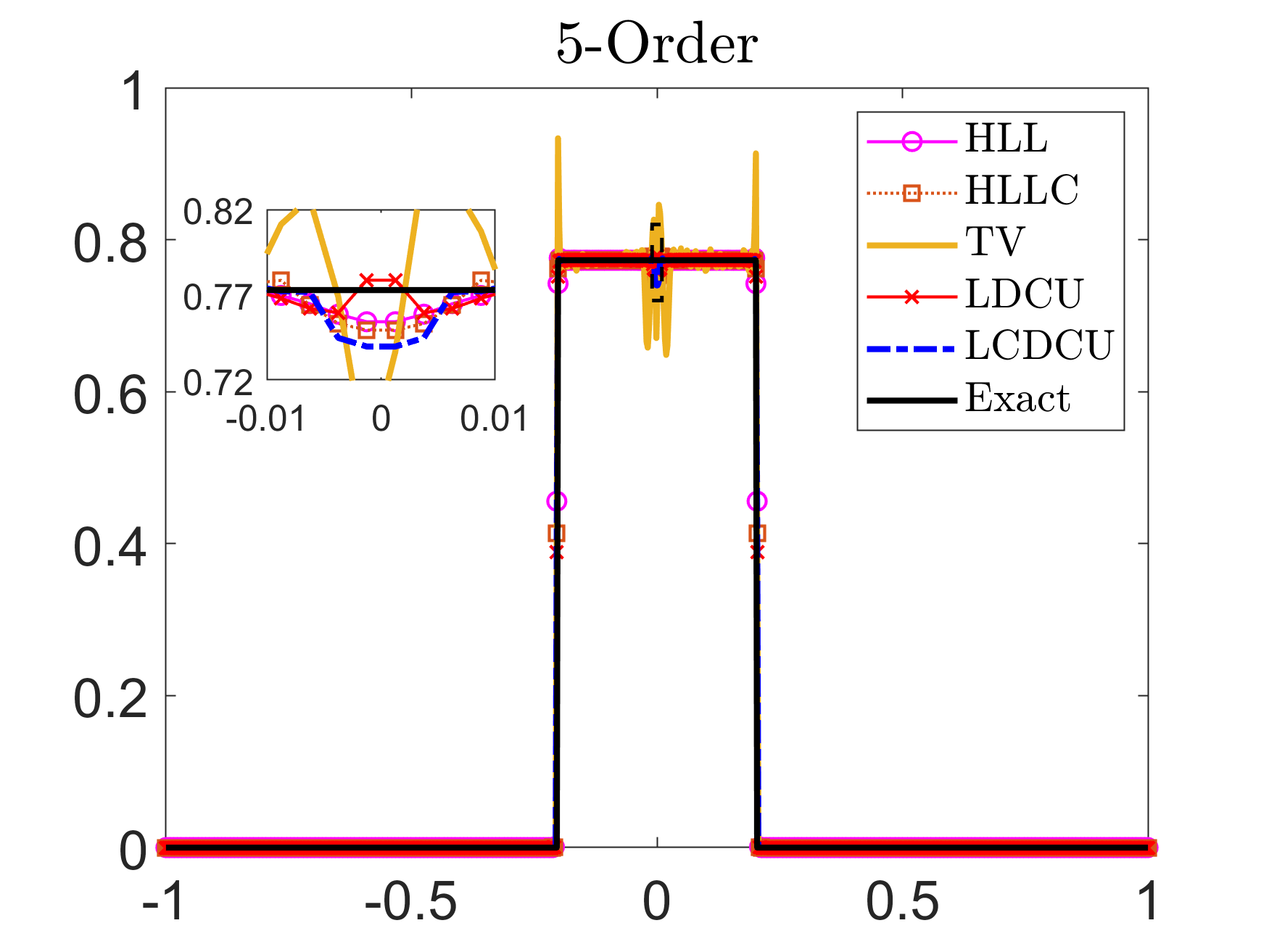}}
\caption{\sf Example 8: Decimal logarithm of the density $\log_{10}(\rho)$ computed by the HLL, HLLC, TV, LDCU, and LCDCU schemes together with the exact solution.
\label{fig6bb}}
\end{figure}

The numerical results for $\log_{10}(\rho)$ are reported in Figure~\ref{fig6bb} together with the exact solution. The four low-dissipation schemes exhibit more pronounced spurious oscillations than the corresponding HLL schemes. Among them, the TV schemes produce the strongest oscillations, especially when extended to the third and fifth orders.

\subsection{Two-Dimensional Examples}
In this section, we consider the 2-D Euler equations of gas dynamics \eref{1.2}, \eref{3.1}--\eref{3.2}.  In Examples 9--18, we take the specific heat ratio $\gamma=1.4$, while in Example 19, we take $\gamma=5/3$.

\subsubsection*{Example 9---2-D Accuracy Test}
In the first 2-D example taken from \cite{Kurganov17,CCHKL_22}, we consider the 2-D Euler equations of gas dynamics subject to the periodic initial
conditions,
\begin{equation*}
\rho(x,y,0)=1+\frac{1}{5}\sin(\pi(x+y)),\quad u(x,y,0)\equiv1,\quad v(x,y,0)\equiv-0.7,\quad p(x,y,0)\equiv1,
\end{equation*}
prescribed on $[-1,1]\times[-1,1]$. The exact solution of this initial value problem can be easily obtained and is given by
$$
\rho(x,y,t)=1+\frac{1}{5}\sin\left[\pi(x+y-0.3t)\right],\quad u(x,y,t)\equiv1,\quad v(x,y,t)\equiv-0.7,\quad p(x,y,t)\equiv1.
$$

We first compute the numerical solution until the final time $t=0.1$ using the 1-Order, 2-Order, 3-Order, and 5-Order schemes on a sequence of uniform meshes: $100\times 100$, $200\times 200$, and $400\times 400$, and then measure the $L^1$-errors and the corresponding experimental convergence rates for the density. The obtained results are presented in Table \ref{tab52}, where one can see that the studied 1-Order, 2-Order, 3-Order, and 5-Order schemes achieve the expected order of accuracy. At the same time, the four low-dissipation schemes are slightly more accurate than the HLL schemes, while the four low-dissipation schemes achieve exactly the same numerical results. As in Example 1, we had to use smaller time steps with $\dt \sim \min \big\{ (\dx)^{\frac{5}{3}}, (\dy)^{\frac{5}{3}}\big\}$ to achieve the fifth order of accuracy.

\begin{table}[ht!]
\centering
\resizebox{\linewidth}{!}{$
\begin{tabular}{|c|c|cc|cc|cc|cc|cc|cc|}
\hline
\multirow{2}{3em}{Method}&\multirow{2}{3em}{Mesh}&\multicolumn{2}{c|}{1-Order}&\multicolumn{2}{c|}{2-Order}&\multicolumn{2}{c|}{3-Order}&\multicolumn{2}{c|}{5-Order}\\
\cline{3-10}& &Error&Rate&Error&Rate&Error&Rate&Error&Rate\\
\hline
\multirow{3}{3em}{HLL}&$100\times 100$&1.18e-02 &---   &3.59e-04 &---  &6.02e-06&---     &4.36e-09 &---\\
                           &$200\times 200$&5.91e-03 &0.994 &8.38e-05 &2.10 &7.37e-07 &3.03&1.36e-10 &5.00\\
                           &$400\times 400$&2.96e-03 &0.997 &1.91e-05 &2.13 &9.22e-08 &3.00&4.42e-12 &4.95\\
\hline
\multirow{3}{3em}{HLLC}&$100\times 100$&8.39e-03 &---   &2.64e-04 &---  &4.32e-06&--- &3.12e-09 &---\\
                       &$200\times 200$&4.21e-03 &0.994 &6.16e-05 &2.10 &5.27e-07&3.03&9.76e-11 &5.00\\
                       &$400\times 400$&2.11e-03 &0.997 &1.47e-05 &2.06 &6.59e-08&3.00&3.19e-12 &4.93\\
\hline 
\multirow{3}{3em}{TV}&$100\times 100$&8.39e-03 &---   &2.64e-04 &---  &4.32e-06&--- &3.12e-09 &---\\
                       &$200\times 200$&4.21e-03 &0.994 &6.16e-05 &2.10 &5.27e-07&3.03&9.76e-11 &5.00\\
                       &$400\times 400$&2.11e-03 &0.997 &1.47e-05 &2.06 &6.59e-08&3.00&3.19e-12 &4.93\\
\hline
\multirow{3}{3em}{LDCU}&$100\times 100$&8.39e-03 &---   &2.64e-04 &---  &4.32e-06&--- &3.12e-09 &---\\
                       &$200\times 200$&4.21e-03 &0.994 &6.16e-05 &2.10 &5.27e-07&3.03&9.76e-11 &5.00\\
                       &$400\times 400$&2.11e-03 &0.997 &1.47e-05 &2.06 &6.59e-08&3.00&3.19e-12 &4.93\\
\hline 
\multirow{3}{3em}{LCDCU}&$100\times 100$&8.39e-03 &---   &2.64e-04 &---  &4.32e-06&--- &3.12e-09 &---\\
                       &$200\times 200$&4.21e-03 &0.994 &6.16e-05 &2.10 &5.27e-07&3.03&9.76e-11 &5.00\\
                       &$400\times 400$&2.11e-03 &0.997 &1.47e-05 &2.06 &6.59e-08&3.00&3.19e-12 &4.93\\
\hline 
\end{tabular}$}
\caption{\sf Example 9: The $L^1$-errors and experimental convergence rates for the density $\rho$ computed by the 1-Order, 2-Order, 3-Order, and 5-Order schemes.\label{tab52}}
\end{table}

\subsubsection*{Example 10---2-D Vortex Evolution Problem}
In this example taken from \cite{Titarev2005}; see also \cite{HS1999,CHT25}, we consider the 2-D vortex evolution problem with the following initial conditions
\begin{equation*}
(\rho(x,y,0),u(x,y,0),v(x,y,0),p(x,y,0))=\left(T^{\frac{1}{\gamma-1}}, 1-\frac{\varepsilon}{2 \pi} e^{\frac{1}{2}(1-r^2)}y, 1+\frac{\varepsilon}{2 \pi} e^{\frac{1}{2}(1-r^2)}x, \rho^\gamma \right),
\end{equation*}
where $T=1-\frac{(\gamma-1)\varepsilon^2}{8\gamma \pi^2}e^{(1-r^2)}$,  $r^2=x^2+y^2$, and $\varepsilon=5$ is the vortex strength. The initial data, prescribed in the computational domain $[-5, 5]\times [-5,5]$ subject to the periodic boundary conditions,  corresponds to a smooth vortex placed at the origin and is defined as the isentropic perturbation to the uniform flow of unit values of primitive variables and the exact solution is a vortex moving with a constant velocity at $45^{\circ}$ to the Cartesian mesh lines.

\begin{figure}[ht!]
	\centerline{\includegraphics[trim=1.cm 0.6cm 0.9cm 0.1cm, clip, width=9cm]{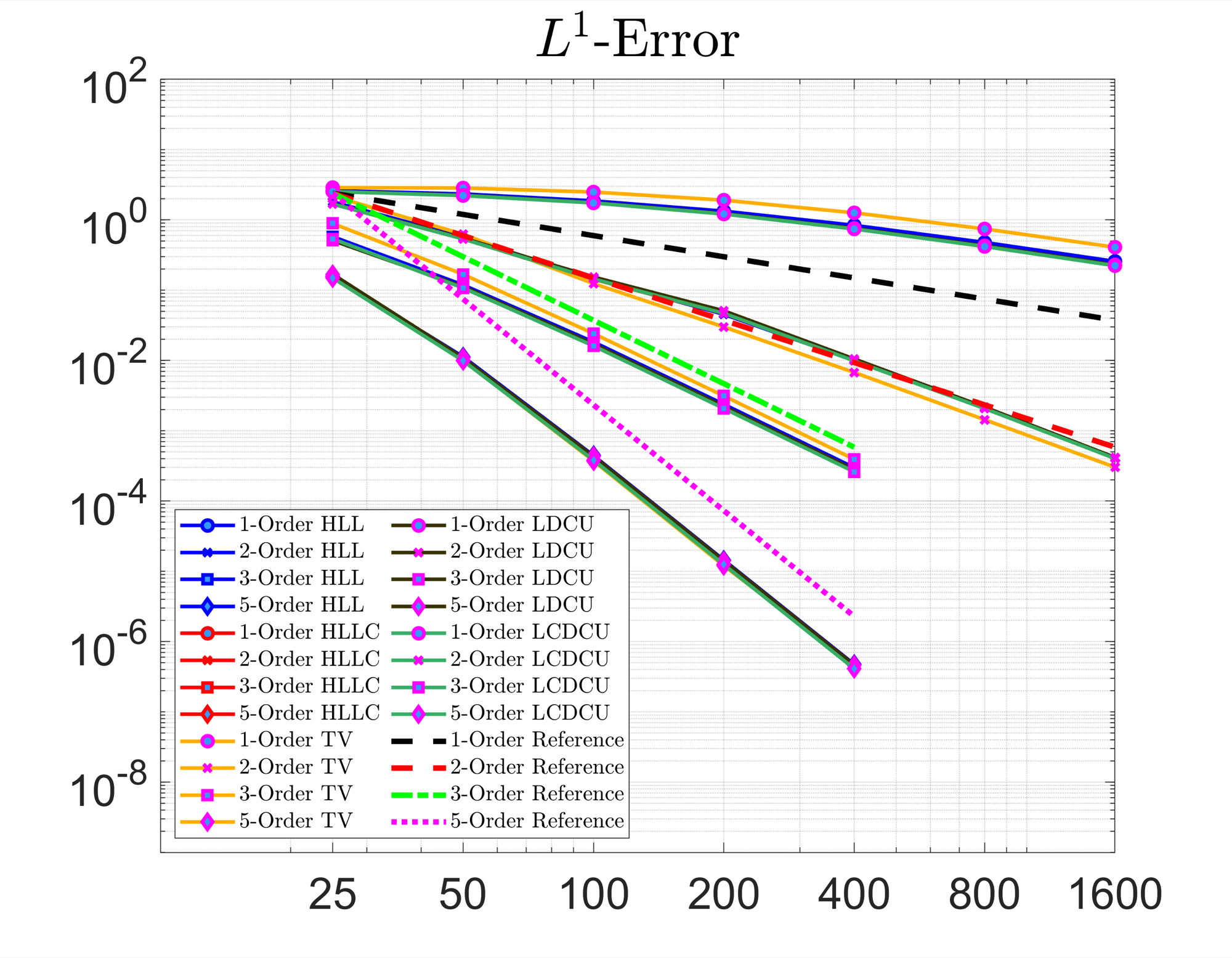}}
	\caption{\sf Example 10: $L^1$-errors for the density $\rho$ computed by the 1-Order, 2-Order, 3-Order, and 5-Order HLL, HLLC, TV, LDCU, and LCDCU schemes.}\label{fig55a}
\end{figure}

We compute the numerical solution until the final time $t=10$ using the 1-Order, 2-Order, 3-Order, and 5-Order schemes on a sequence of uniform meshes: $100\times 100$, $200\times 200$, and $400\times 400$, and then measure the $L^1$-errors between the computed solutions and the exact solutions and the corresponding experimental convergence rates for the density. The obtained results are presented in Table \ref{tab52a}, where one can see that both 1-Order and 2-Order convergence rates are observed only after significant mesh refinement, while the 3-Order and 5-Order schemes achieve the expected order of accuracy. In order to have a better view, we also show the $L^1$-errors in Figure \ref{fig55a}, where we also show the results computed by the 1-Order and 2-Order schemes on finer meshes $800\times 800$ and $1600\times 1600$ to show that the 1-Order and 2-Order schemes achieve expected convergence rates after mesh refinement. It is noticed that, in this example, the TV schemes are slightly more dissipative than the HLL schemes, but HLLC, LDCU, and LCDCU schemes are slightly more accurate than the HLL schemes.

\begin{table}[ht!]
\centering
\resizebox{\linewidth}{!}{$
\begin{tabular}{|c|c|cc|cc|cc|cc|cc|cc|}
\hline
\multirow{2}{3em}{Method}&\multirow{2}{3em}{Mesh}&\multicolumn{2}{c|}{1-Order}&\multicolumn{2}{c|}{2-Order}&\multicolumn{2}{c|}{3-Order}&\multicolumn{2}{c|}{5-Order}\\
\cline{3-10}& &Error&Rate&Error&Rate&Error&Rate&Error&Rate\\
\hline
\multirow{4}{3em}{HLL}& $100\times 100$&1.86   &---      &1.48e-01 &---  &2.45e-02&--- &7.88e-04 &---\\
                      & $200\times 200$&1.33   &4.84e-01 &4.54e-02 &1.71 &2.94e-03&3.06&2.42e-05 &5.02\\
                      & $400\times 400$&0.84   &6.71e-01 &1.00e-02 &2.18 &3.08e-04&3.25&6.16e-07 &5.30\\
\hline
 \multirow{4}{3em}{HLLC}& $100\times 100$&1.74    &---       &1.44e-01 &---  &2.26e-02&--- &7.03e-04 &---\\
                        & $200\times 200$&1.21    &5.26e-01  &4.67e-02 &1.62 &2.68e-03&3.08&2.27e-05 &4.95\\
                        & $400\times 400$&0.74    &7.04e-01  &9.88e-03 &2.24 &2.70e-04&3.31&5.60e-07 &5.34\\
 \hline
 \multirow{4}{3em}{LDCU}& $100\times 100$&1.76    &---       &1.54e-01 &---  &2.23e-02&--- &7.07e-04 &---\\
                        & $200\times 200$&1.23    &5.21e-01  &5.07e-02 &1.60 &2.71e-03&3.04&2.32e-05 &4.93\\
                        & $400\times 400$&0.76    &6.97e-01  &1.05e-02 &2.28 &2.85e-04&3.25&5.94e-07 &5.29\\
 \hline
 \multirow{4}{3em}{LCDCU}& $100\times 100$&1.74    &---       &1.44e-01 &---  &2.28e-02&--- &7.03e-04 &---\\
                         & $200\times 200$&1.21    &5.26e-01  &4.69e-02 &1.62 &2.70e-03&3.07&2.27e-05 &4.95\\
                         & $400\times 400$&0.74    &7.04e-01  &9.88e-03 &2.25 &2.70e-04&3.32&5.60e-07 &5.34\\
\hline 
\multirow{4}{3em}{TV} & $100\times 100$&2.49 &---      &1.24e-01 &---  &2.98e-02&--- &9.23e-04 &---\\
                      & $200\times 200$&1.91 &3.84e-01 &3.00e-02 &2.05 &3.58e-03&3.06&2.49e-05 &5.21\\
                      & $400\times 400$&1.26 &5.97e-01 &6.75e-03 &2.15 &3.99e-04&3.16&6.62e-07 &5.23\\
\hline
\end{tabular}
$}
\caption{\sf Example 10: $L^1$-errors and experimental convergence rates for the density $\rho$ computed by the 1-Order, 2-Order, 3-Order, and 5-Order HLL, HLLC, TV, LDCU, and LCDCU schemes.\label{tab52a}}
\end{table}

\subsubsection*{Example 11---2-D Moving Contact Waves} In this example taken from \cite{Kurganov21a} (see also \cite{CKX22,KX_22}), we consider an isolated moving contact wave with the
following initial data:
\begin{equation*}
(\rho,u,v,p)(x,y,0)=\begin{cases}
(1.4,0,0.2,1),&(x,y)\in D \\
(1,0,0.2,1),&\mbox{otherwise},
\end{cases}
\end{equation*}
where the domain
$D$ consists of the points $(x,y)$ satisfying the following conditions:
$$
\begin{aligned}
&\{-0.1<x<0.1,\,0<y<0.02\}\cup\{-0.02<x<0.02,\,0.02<y<0.1\}\cup\\
&\{(x+0.02)^2+(y-0.02)^2<0.08^2\}\cup\{(x-0.02)^2+(y-0.02)^2<0.08^2\}.
\end{aligned}
$$
The initial data are prescribed in the computational domain $[-0.2,0.2]\times[0,0.8]$ subject to the free boundary conditions.

\begin{figure}[ht!]
\centerline{\includegraphics[trim=1.8cm 2.8cm 0.9cm 2.5cm, clip, width=14.cm]{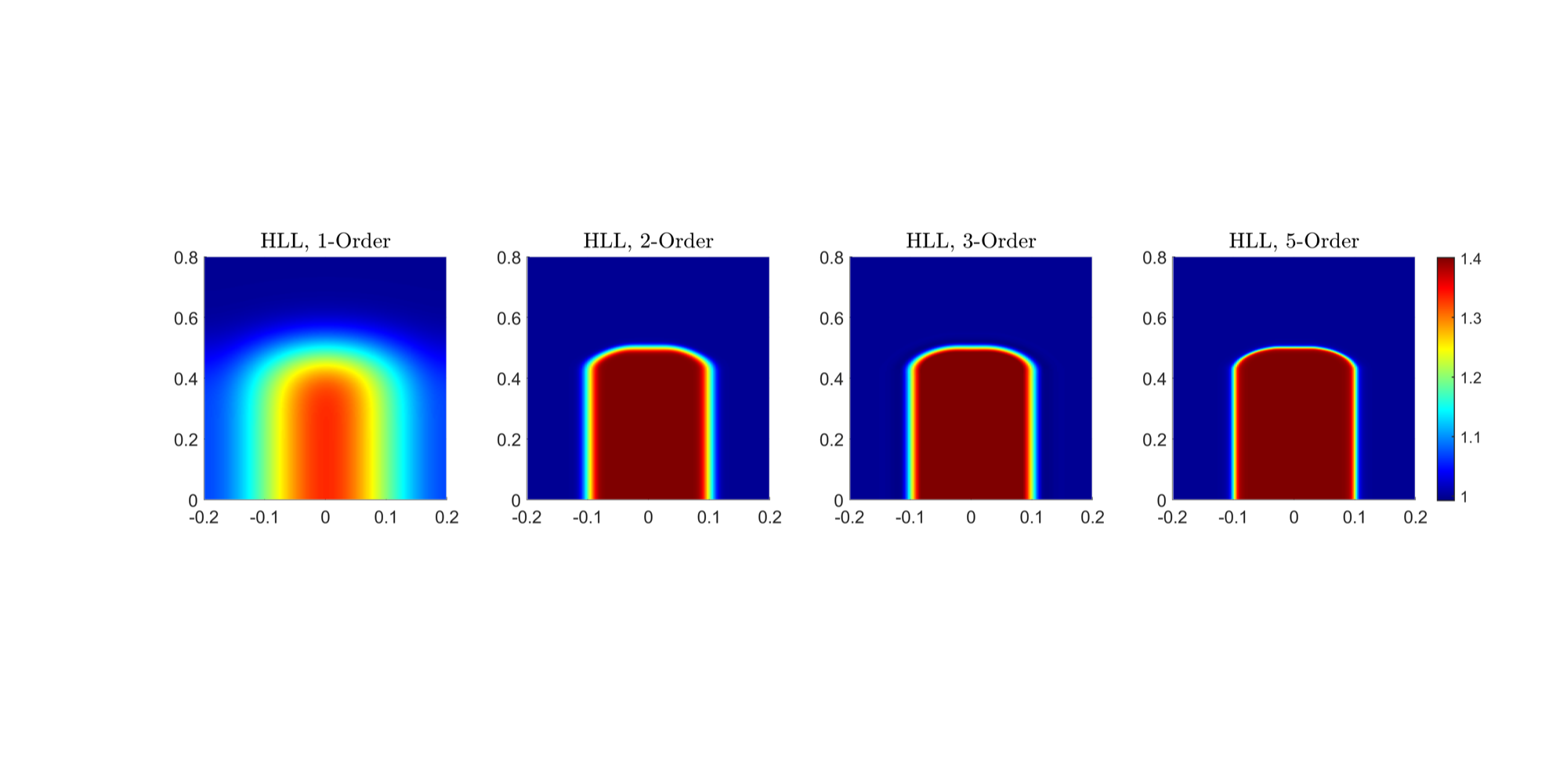}}
\vskip 10pt 
\centerline{\includegraphics[trim=1.8cm 2.8cm 0.9cm 2.5cm, clip, width=14.cm]{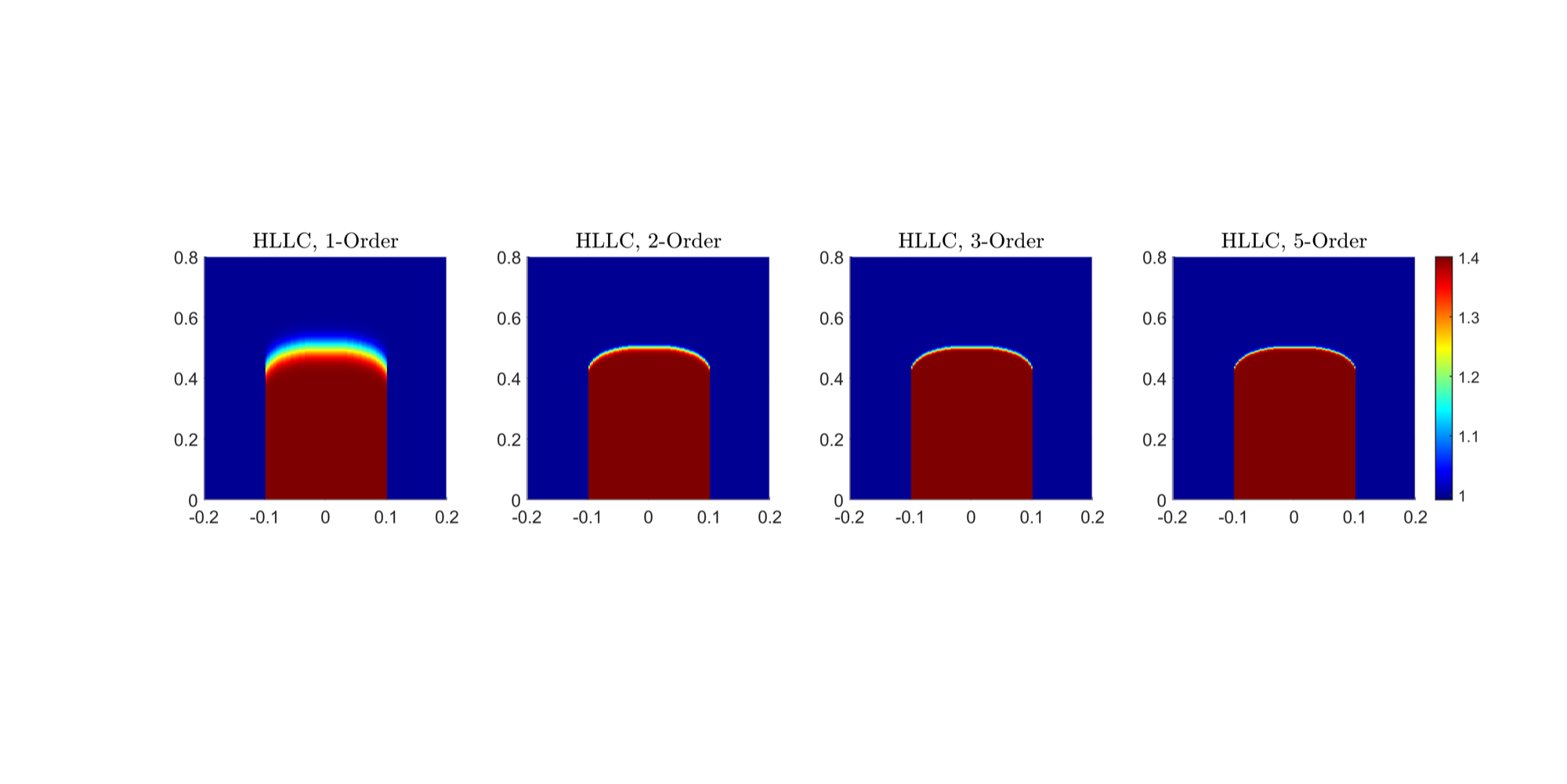}}
\vskip 10pt 
\centerline{\includegraphics[trim=1.8cm 2.8cm 0.9cm 2.5cm, clip, width=14.cm]{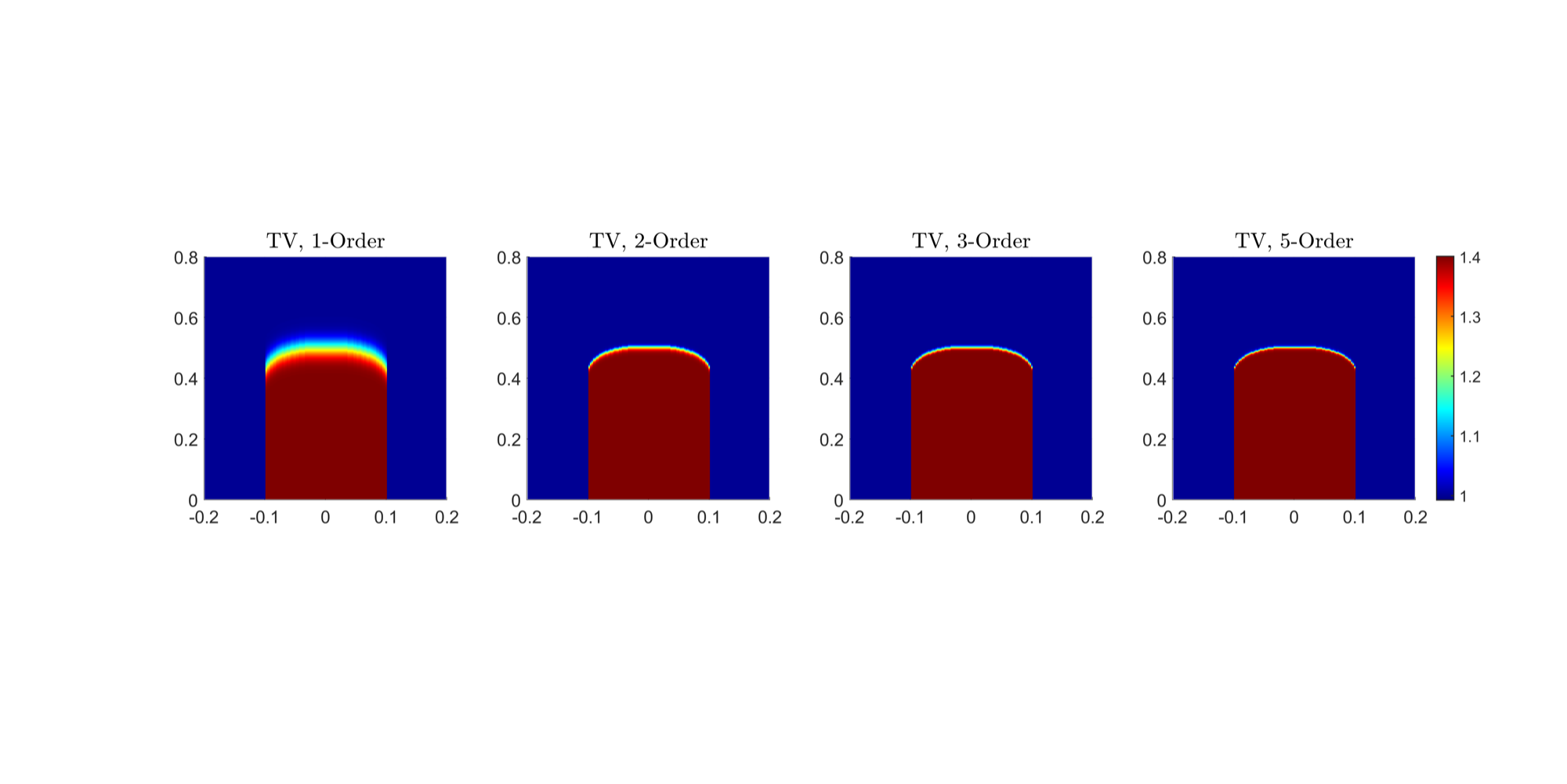}}
\vskip 10pt 
\centerline{\includegraphics[trim=1.8cm 2.8cm 0.9cm 2.5cm, clip, width=14.cm]{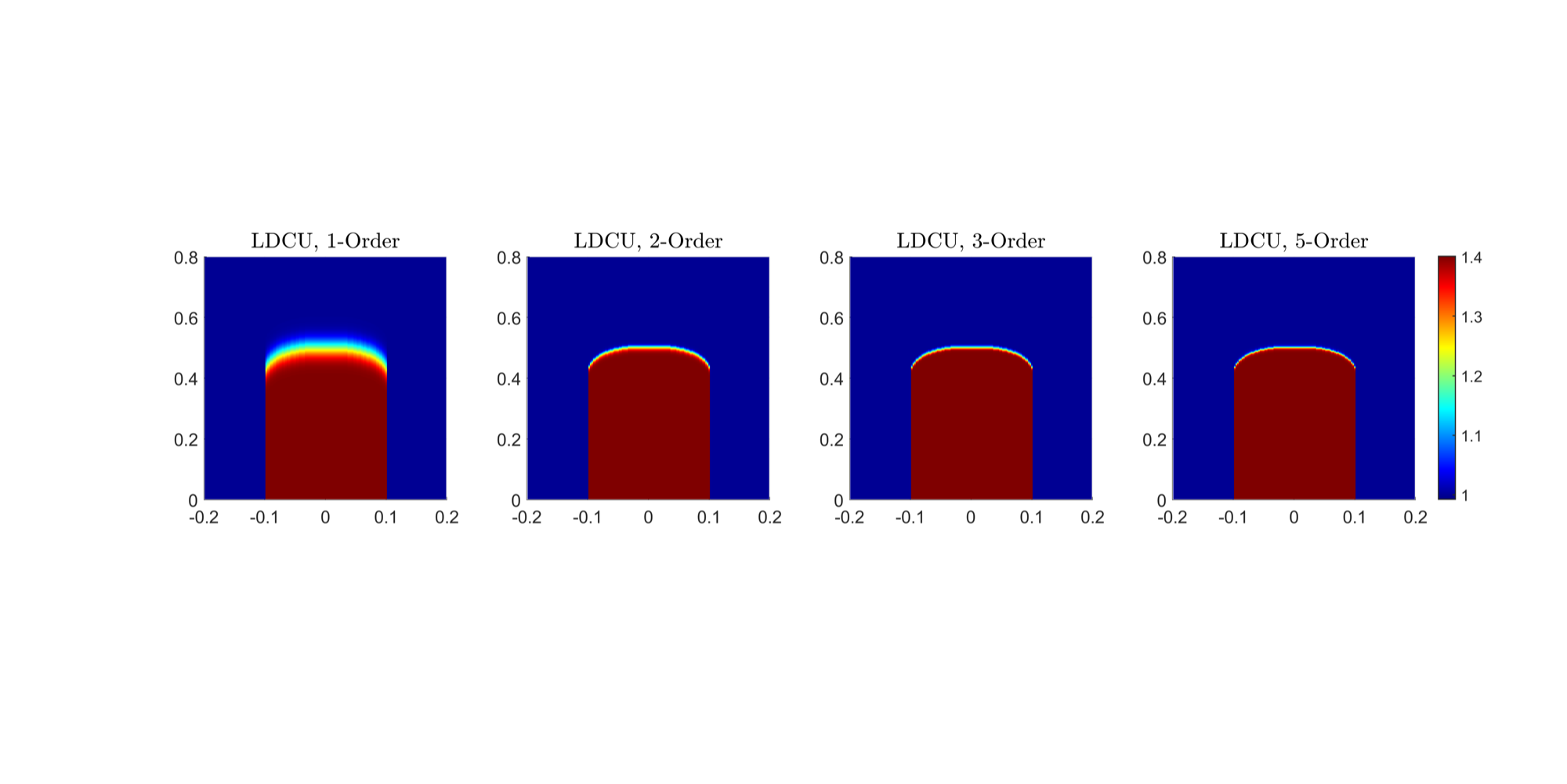}}
\vskip 10pt 
\centerline{\includegraphics[trim=1.8cm 2.8cm 0.9cm 2.5cm, clip, width=14.cm]{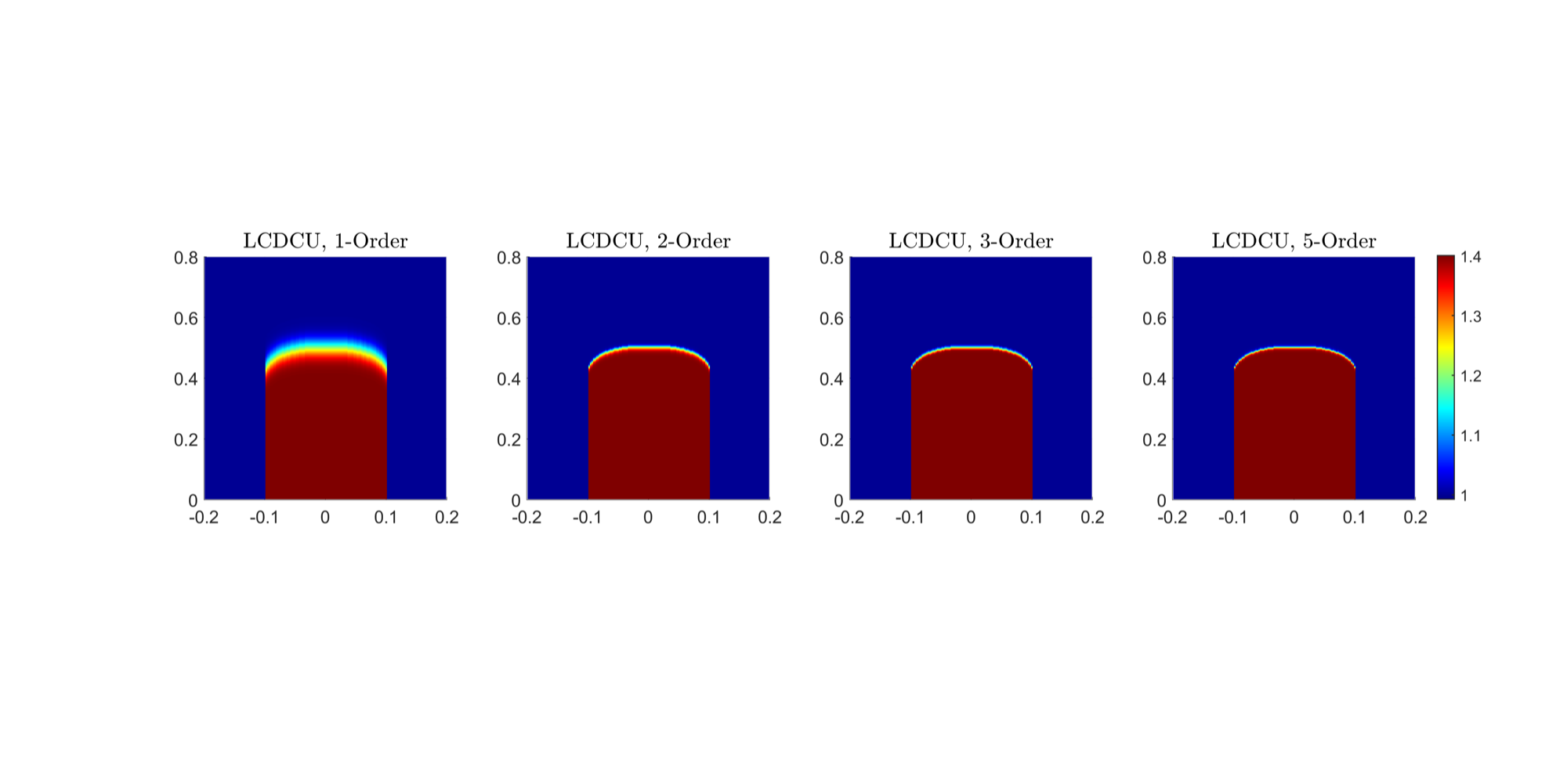}}
\caption{\sf{Example 11: Density $\rho$ computed by the 1-Order, 2-Order, 3-Order, and 5-Order HLL (top row), HLLC (second row), TV (third row), LDCU (fourth row), and LCDCU (bottom row) schemes.}\label{fig66}}
\end{figure}

We compute the numerical solutions until the final time $t=2$ by the studied 1-Order, 2-Order, 3-Order, and 5-Order schemes on a uniform mesh with $\dx=\dy=1/400$ and
plot the obtained results in Figure \ref{fig66}. As one can see, the four low-dissipation schemes produce nearly identical numerical results, which are significantly sharper than those computed by the HLL schemes, especially for the non-moving jumps across the lines $x=\pm0.1$.

\subsubsection*{Example 12---Explosion Problem}
In this example, we consider the explosion problem studied in \cite{Toro2009,Liska03} (see also \cite{KX_22,CCK23_Adaptive}). We take the following initial conditions,
\begin{equation*}
(\rho(x,y,0),u(x,y,0),v(x,y,0),p(x,y,0))=\begin{cases}
(1,0,0,1),&x^2+y^2<0.16,\\
(0.125,0,0,0.1),&\mbox{otherwise},
\end{cases}
\end{equation*}
which are prescribed in the computational domain $[-1,1]\times[-1,1]$, subject to free boundary conditions at all the four sides. The solution of this initial–boundary value problem develops circular shock and contact waves. While the shock wave is stable and requires sufficient dissipation for stable capturing, the contact wave is unstable and can only be stabilized by numerical diffusion. This makes the problem a useful benchmark for evaluating the dissipation of numerical schemes, where the goal is to minimize dissipation while retaining shock stability.

\begin{figure}[ht!]
\centerline{\includegraphics[trim=1.8cm 2.7cm 0.7cm 2.5cm, clip, width=14.cm]{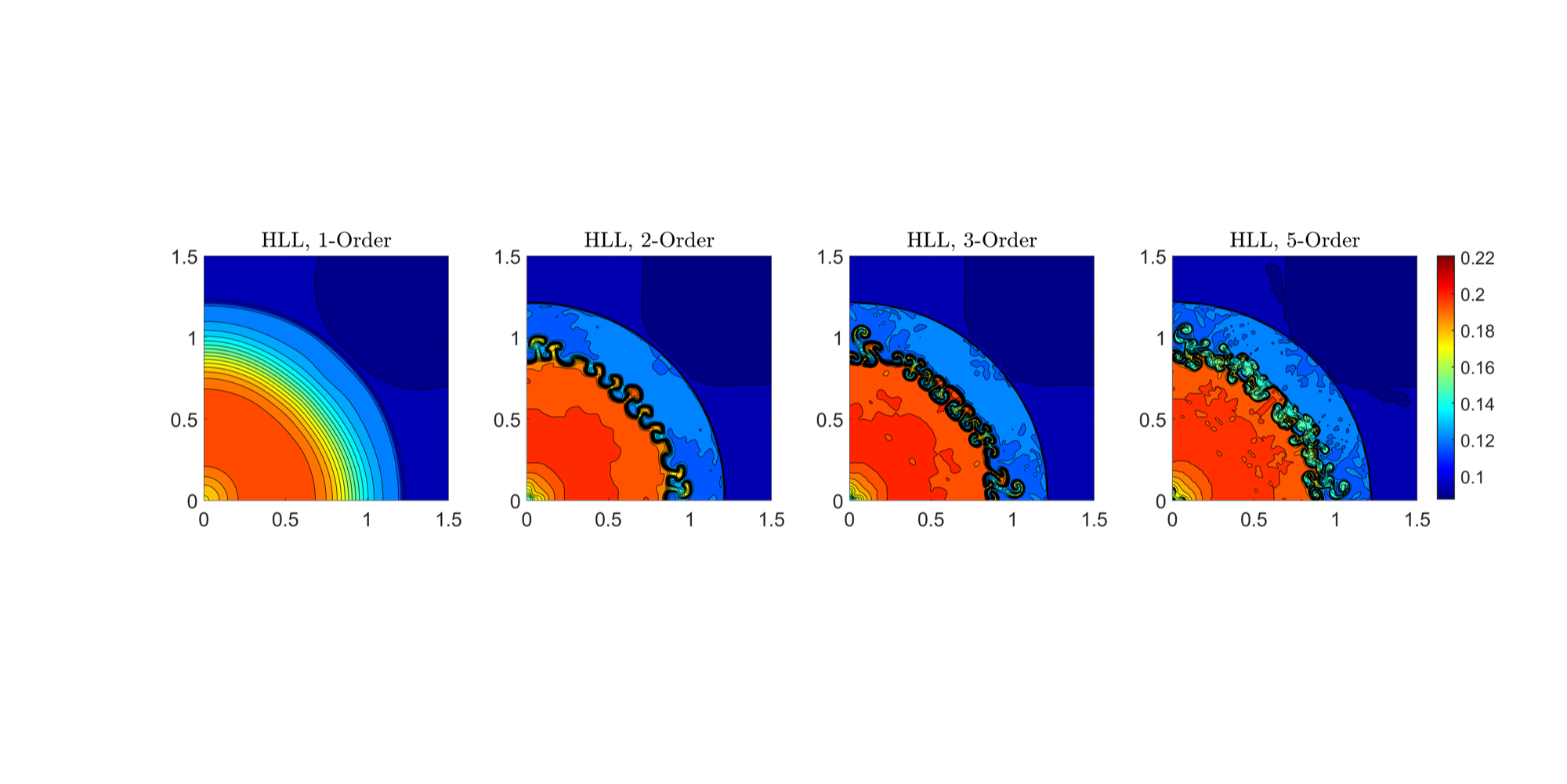}}
\vskip 10pt 
\centerline{\includegraphics[trim=1.8cm 2.7cm 0.7cm 2.5cm, clip, width=14.cm]{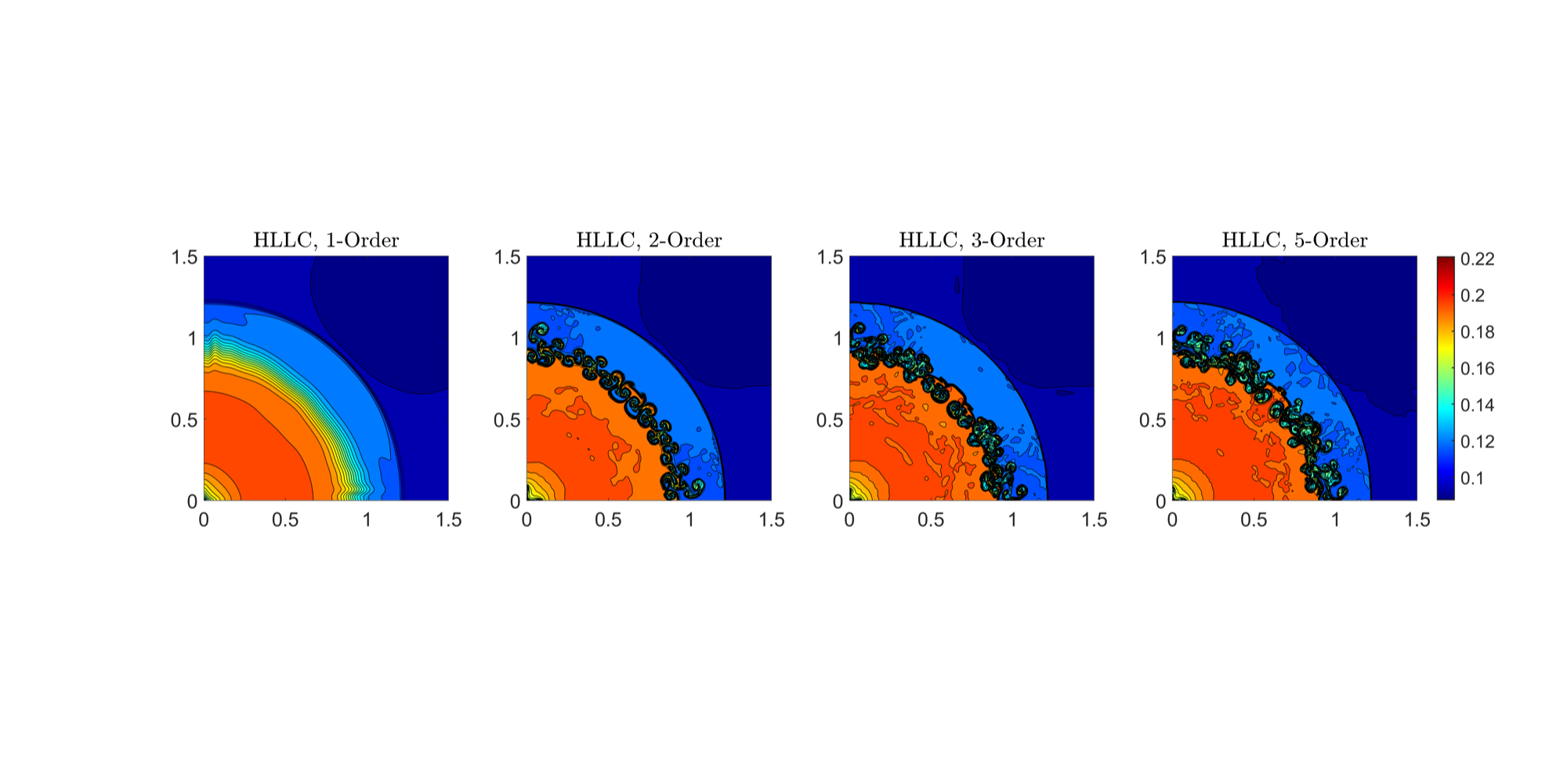}}
\vskip 10pt 
\centerline{\includegraphics[trim=1.8cm 2.7cm 0.7cm 2.5cm, clip, width=14.cm]{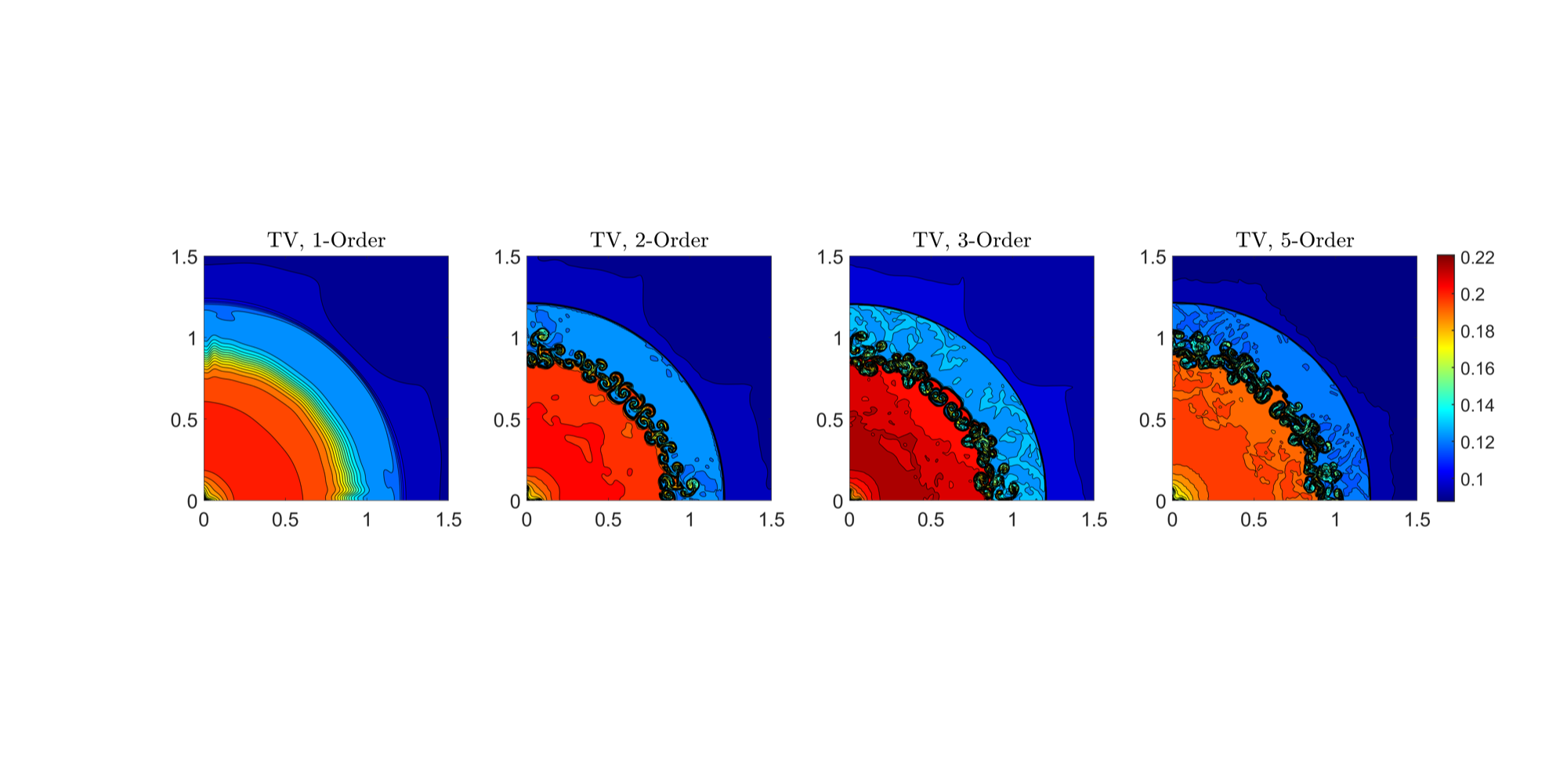}}
\vskip 10pt 
\centerline{\includegraphics[trim=1.8cm 2.7cm 0.7cm 2.5cm, clip, width=14.cm]{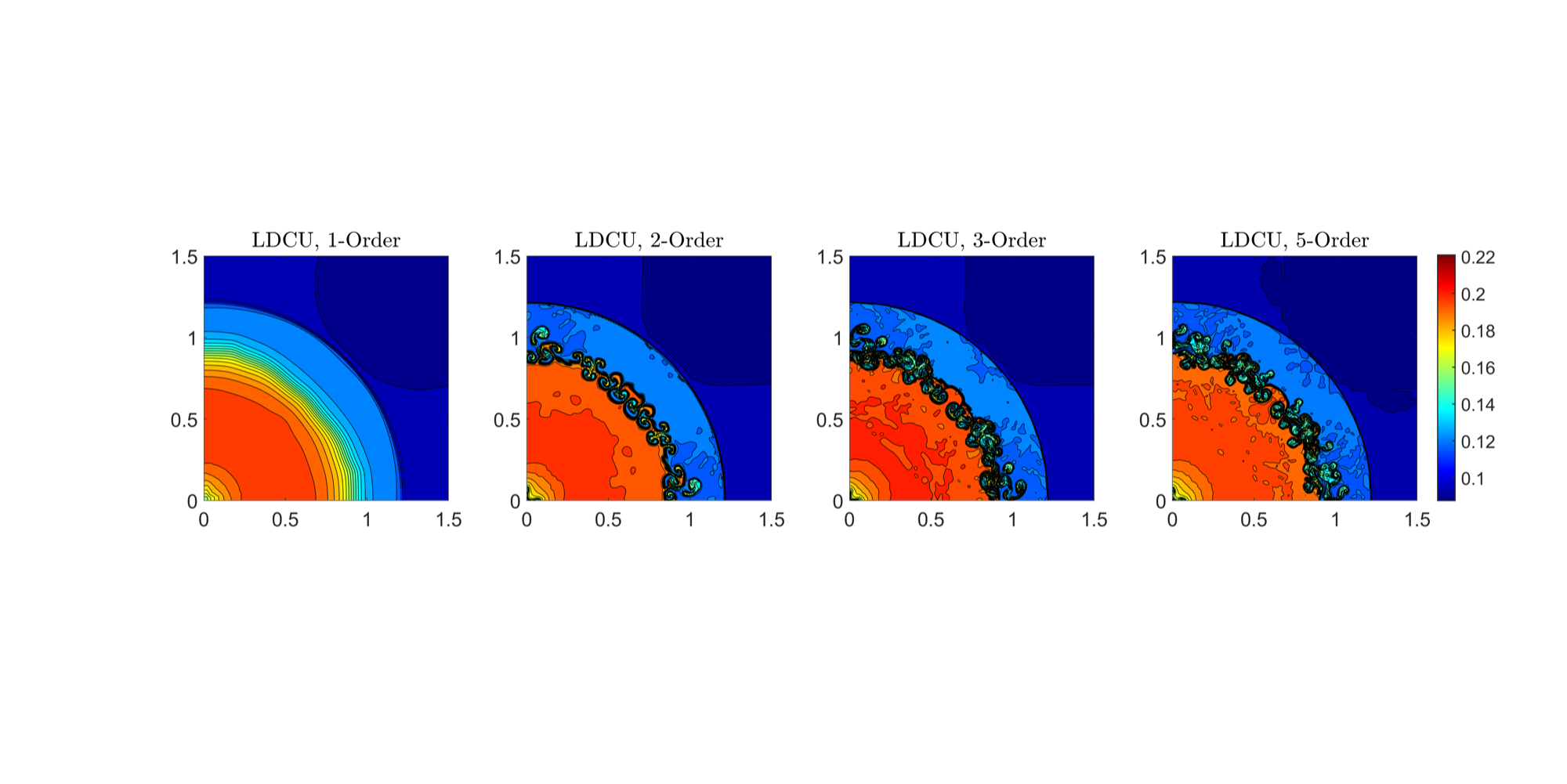}}
\vskip 10pt 
\centerline{\includegraphics[trim=1.8cm 2.7cm 0.7cm 2.5cm, clip, width=14.cm]{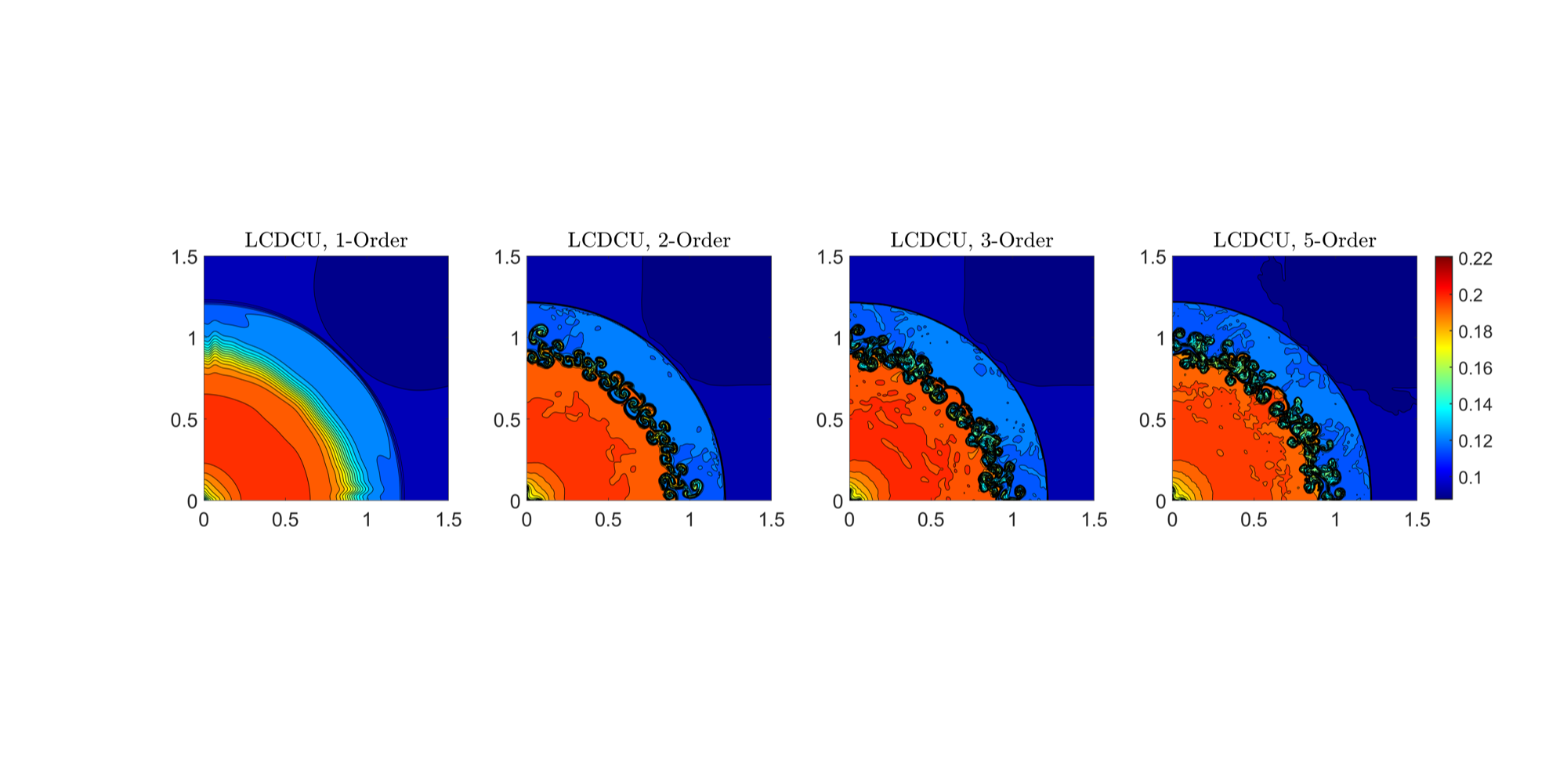}}
\caption{\sf{Example 12: Density $\rho$ computed by the 1-Order, 2-Order, 3-Order, and 5-Order HLL (top row), HLLC (second row), TV (third row), LDCU (fourth row), and LCDCU (bottom row) schemes.}\label{fig10}}
\end{figure}

We apply the studied 1-Order, 2-Order, 3-Order, and 5-Order schemes and compute the numerical solutions until the final time $t=3.2$ on a uniform mesh with $\dx=\dy=3/800$. The obtained results are presented in Figure \ref{fig10}, where one can see that the differences between the numerical results computed by the four low-dissipation schemes are limited, but much more substantially “curlier” and the mixing layer is slightly wider (indicating a more severe instability) than those computed by the HLL schemes, while the shock is still stable.

\subsubsection*{Example 13---Implosion Problem}
In this example, we consider the implosion problem taken from \cite{Liska03} (see also \cite{CCHKL_22,Kurganov07}). The initial conditions,
\begin{equation*}
(\rho(x,y,0),u(x,y,0),v(x,y,0),p(x,y,0))=\begin{cases}
(0.125,0,0,0.14),&|x|+|y|<0.15,\\
(1,0,0,1),&\mbox{otherwise},
\end{cases}
\end{equation*}
are prescribed in the computational domain $[0,0.3]\times[0,0.3]$ with solid boundary conditions imposed at all four sides. This example is designed to assess the numerical diffusion of different schemes: a jet forms near the origin and propagates along the diagonal $y=x$, and excessive diffusion may either smear the jet entirely or alter its propagation velocity.

\begin{figure}[ht!]
\centerline{\includegraphics[trim=1.8cm 2.7cm 0.7cm 2.5cm, clip, width=14.cm]{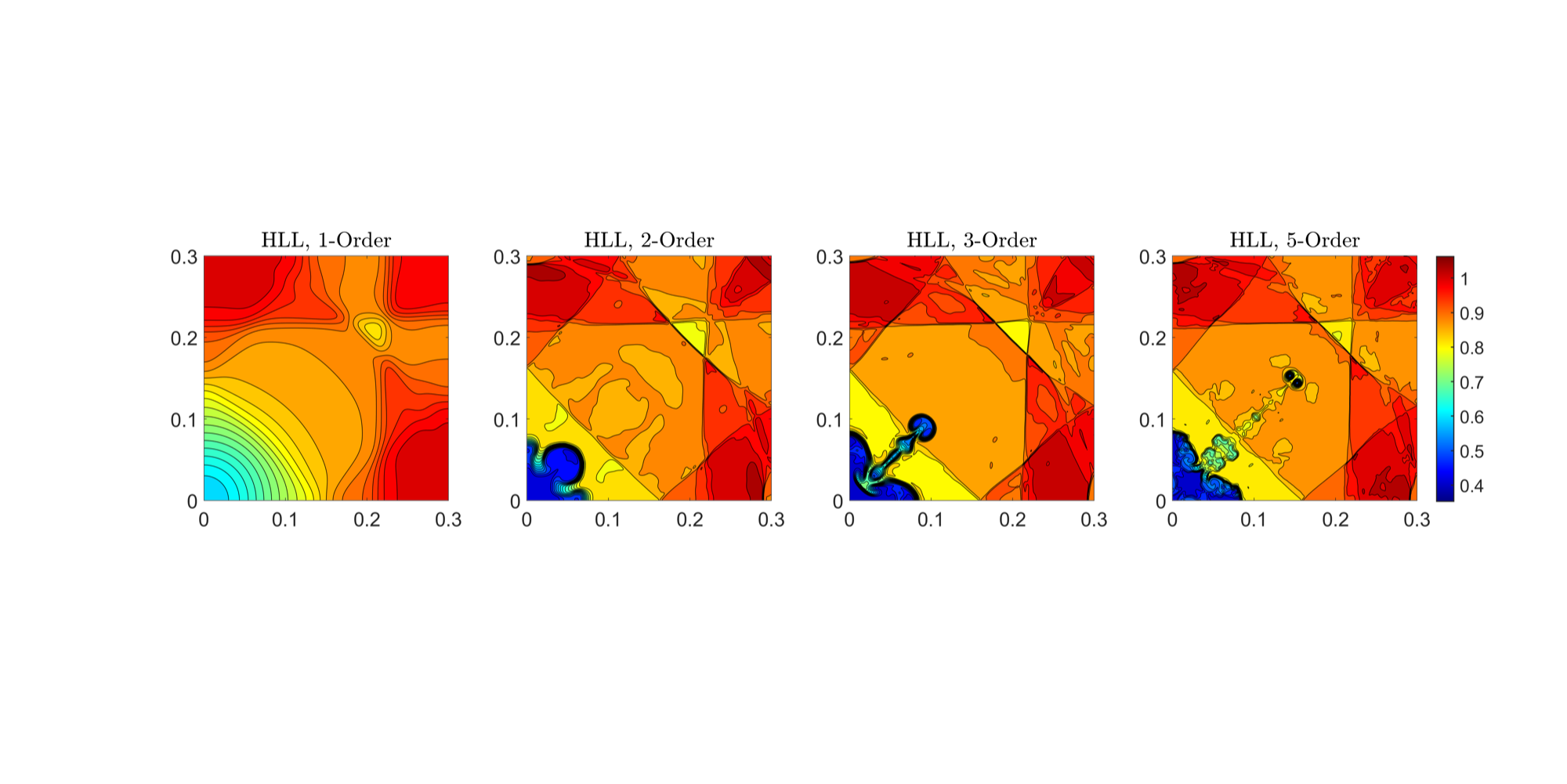}}
\vskip 10pt 
\centerline{\includegraphics[trim=1.8cm 2.7cm 0.7cm 2.5cm, clip, width=14.cm]{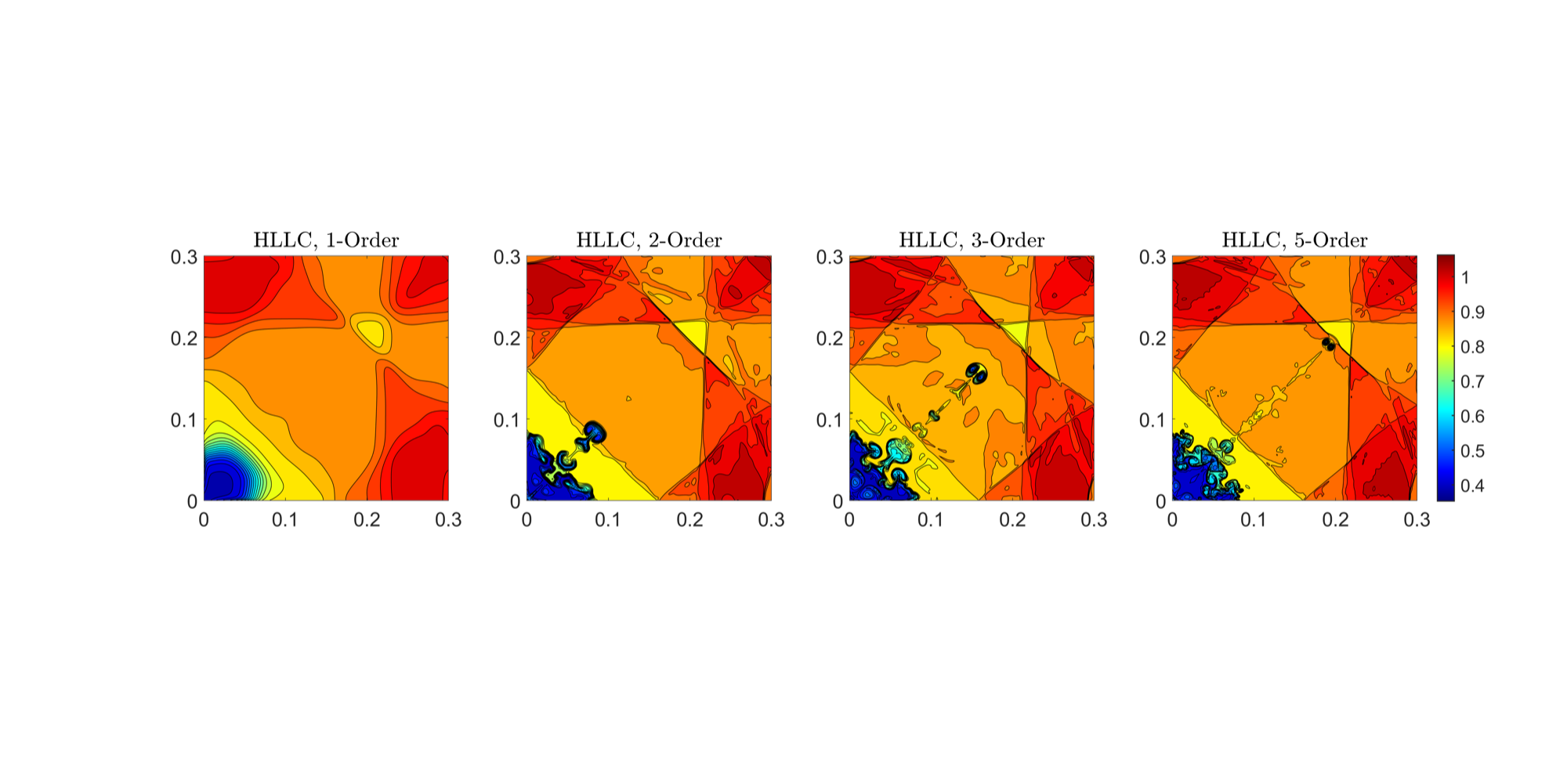}}
\vskip 10pt 
\centerline{\includegraphics[trim=1.8cm 2.7cm 0.7cm 2.5cm, clip, width=14.cm]{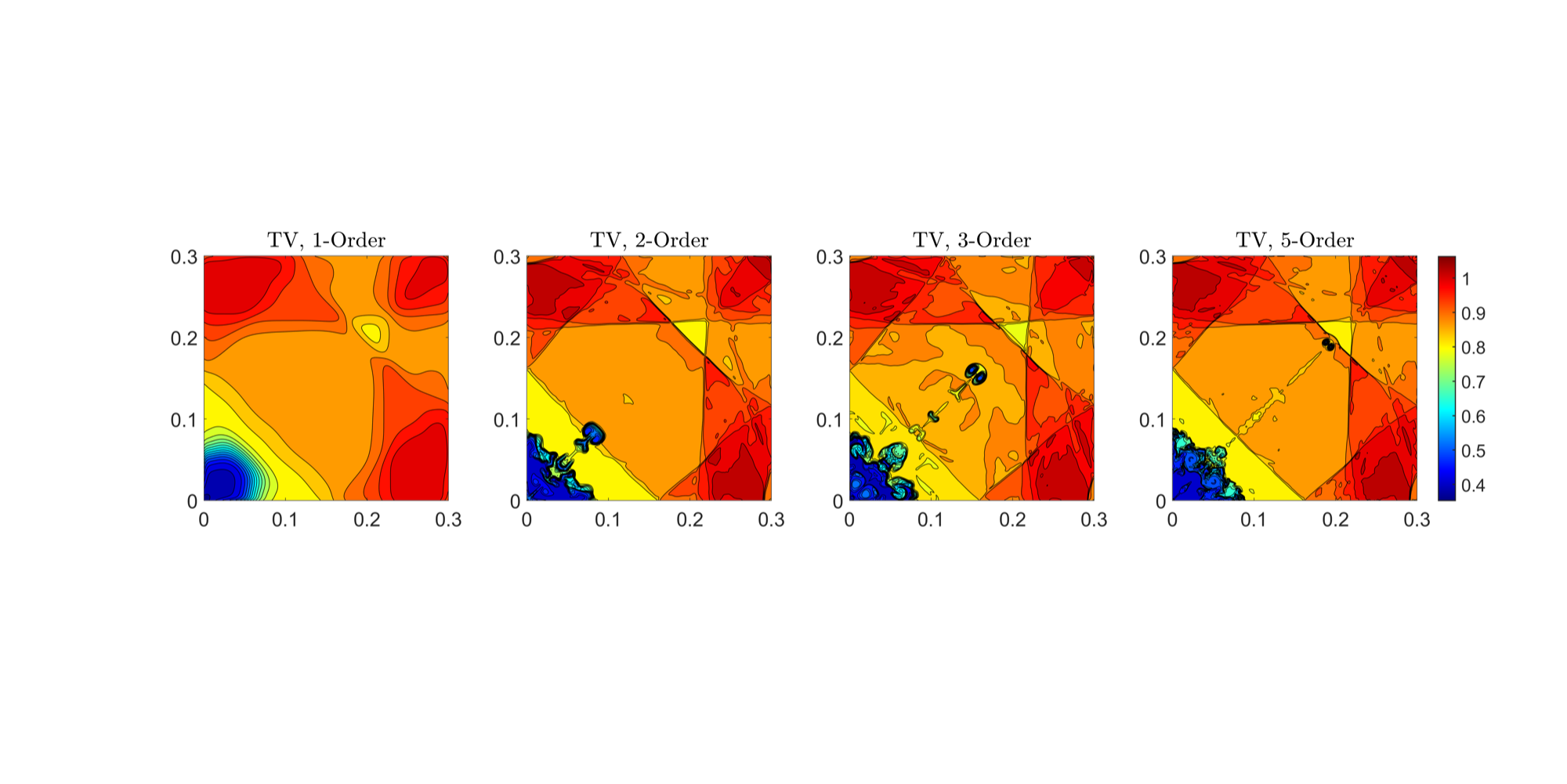}}
\vskip 10pt 
\centerline{\includegraphics[trim=1.8cm 2.7cm 0.7cm 2.5cm, clip, width=14.cm]{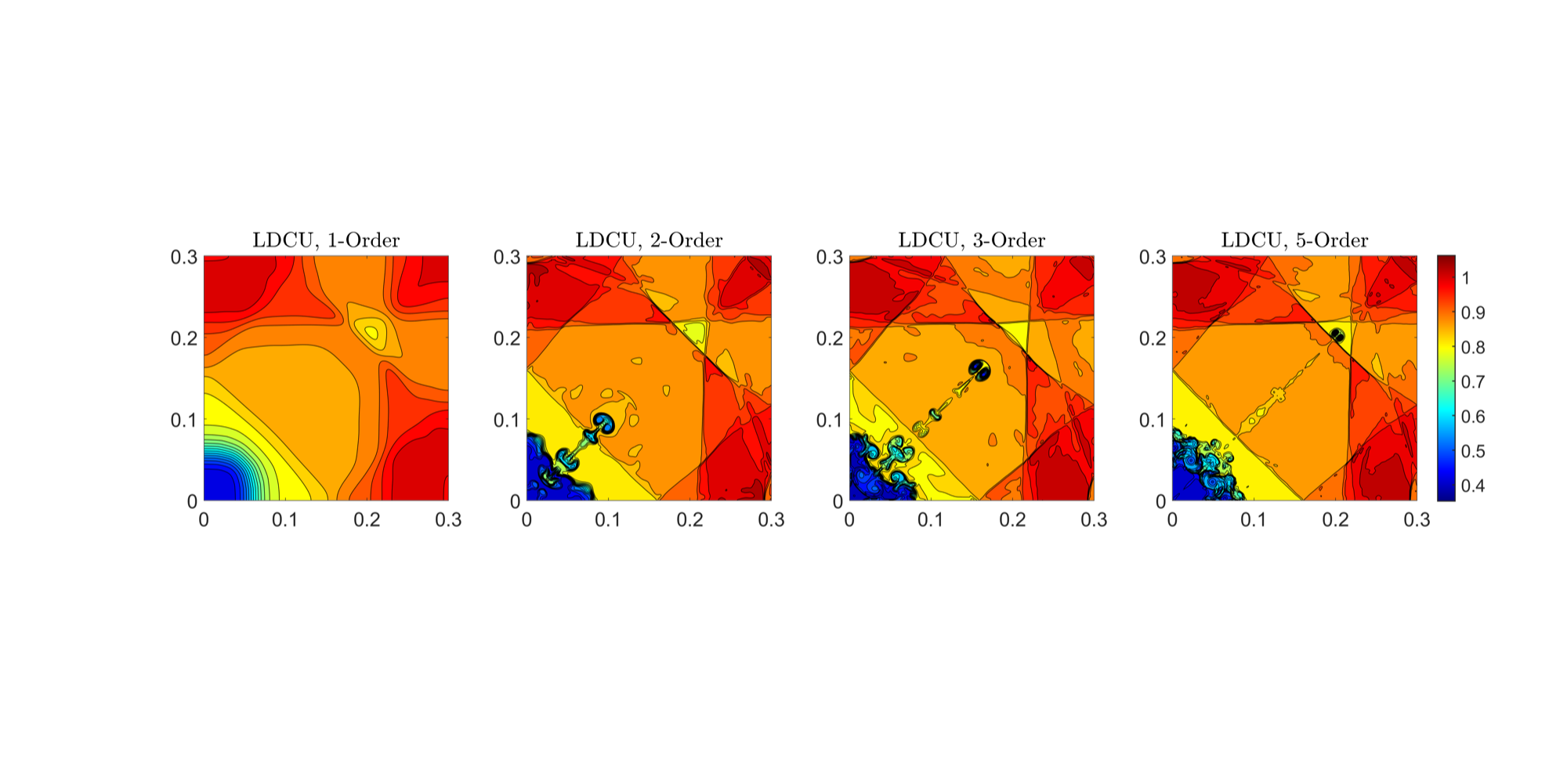}}
\vskip 10pt 
\centerline{\includegraphics[trim=1.8cm 2.7cm 0.7cm 2.5cm, clip, width=14.cm]{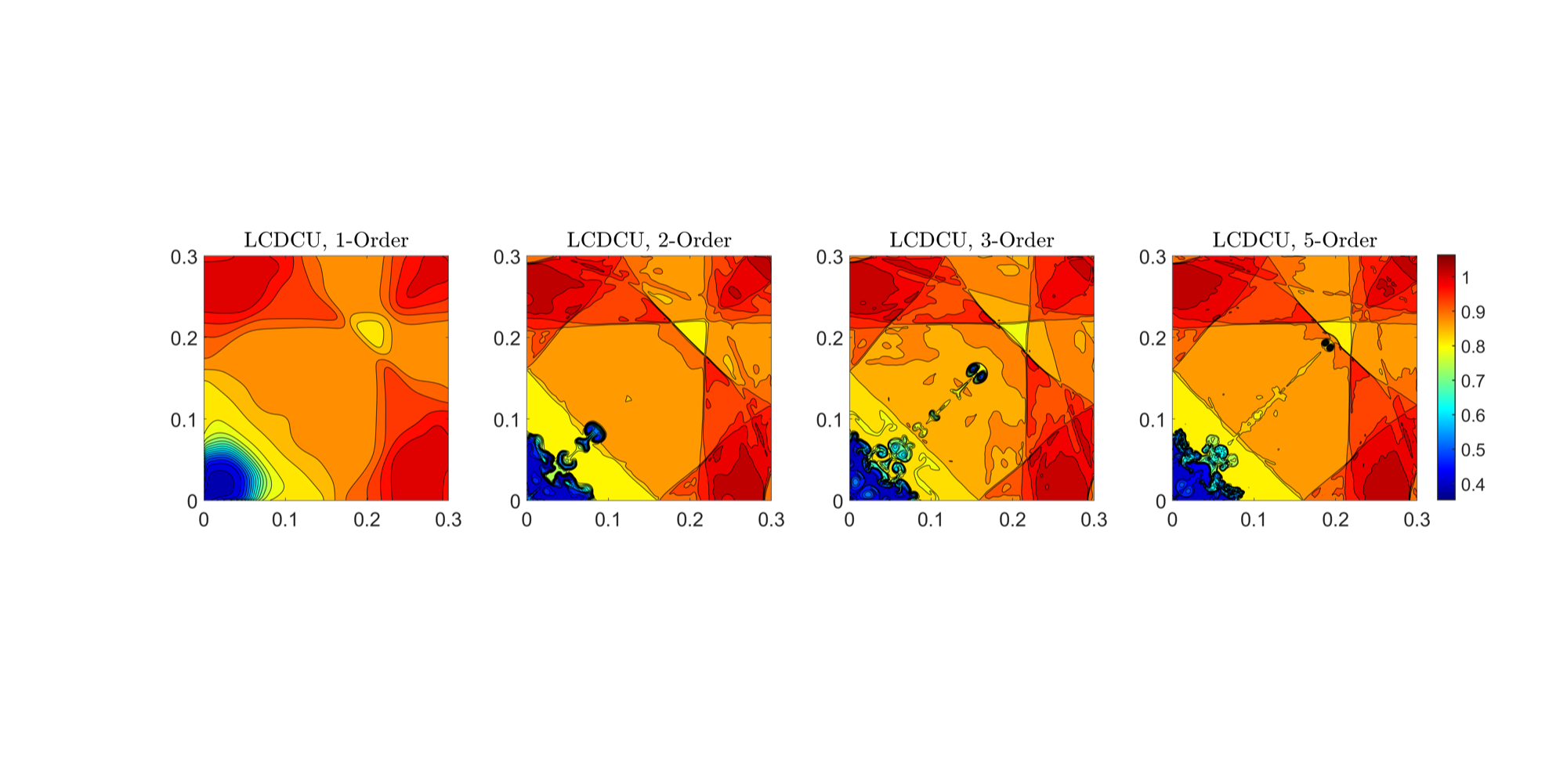}}
\caption{\sf{Example 13: Density $\rho$ computed by the 1-Order, 2-Order, 3-Order, and 5-Order HLL (top row), HLLC (second row), TV (third row), LDCU (fourth row), and LCDCU (bottom row) schemes.}\label{fig11}}
\end{figure}

We compute the numerical solutions until the final time $t=2.5$ by the 1-Order, 2-Order, 3-Order, and 5-Order schemes on a uniform mesh with $\dx=\dy=3/4000$ and plot the obtained results in Figure \ref{fig11}, where one can clearly observe that the jet propagates much further in the diagonal direction when using the four low-dissipation schemes. At the same time, the LDCU scheme is slightly less dissipative than the other three low-dissipation counterparts.

\begin{figure}[ht!]
\centerline{\includegraphics[trim=1.8cm 2.7cm 0.7cm 2.5cm, clip, width=14.cm]{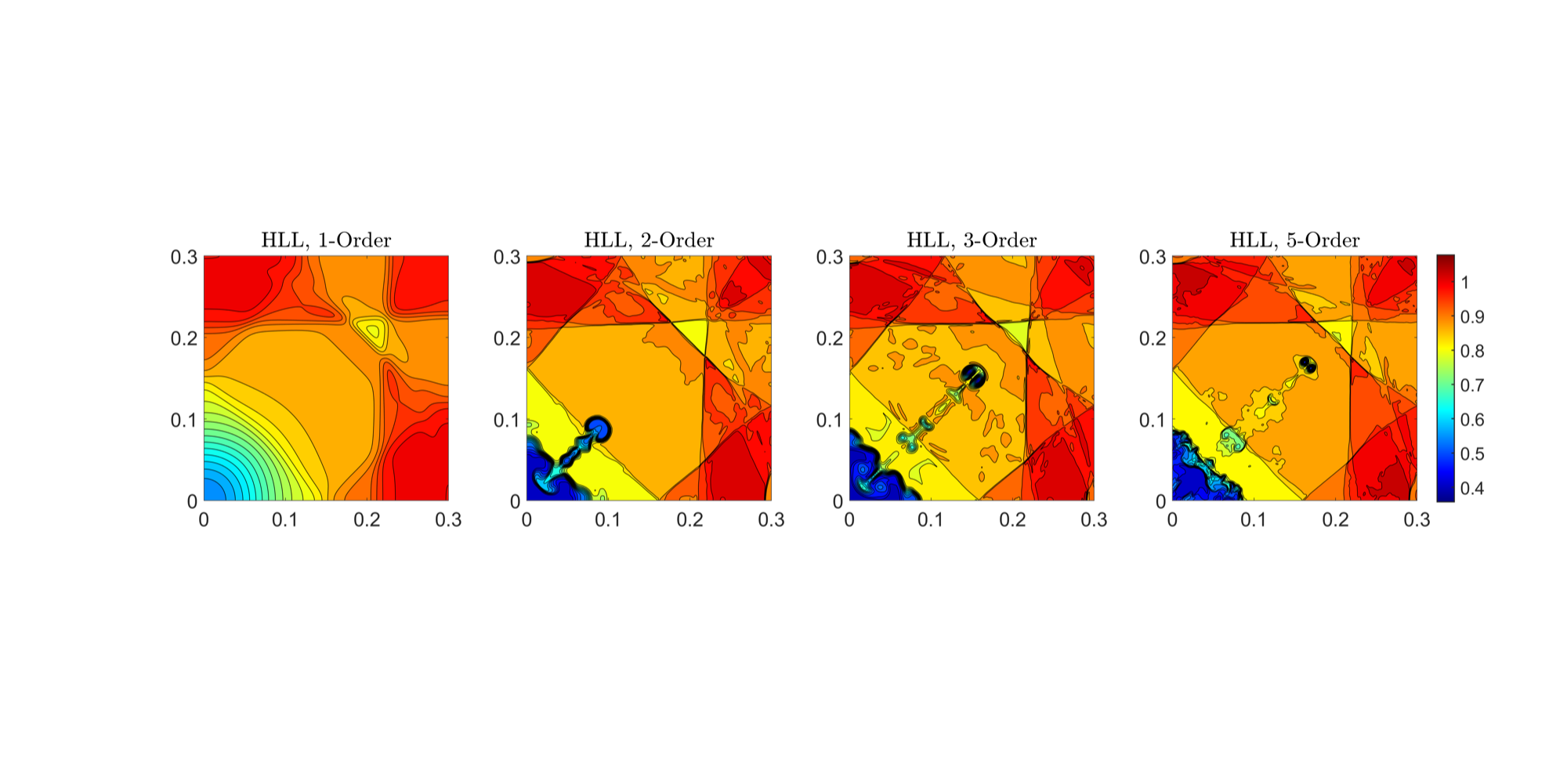}}
\vskip 10pt 
\centerline{\includegraphics[trim=1.8cm 2.7cm 0.7cm 2.5cm, clip, width=14.cm]{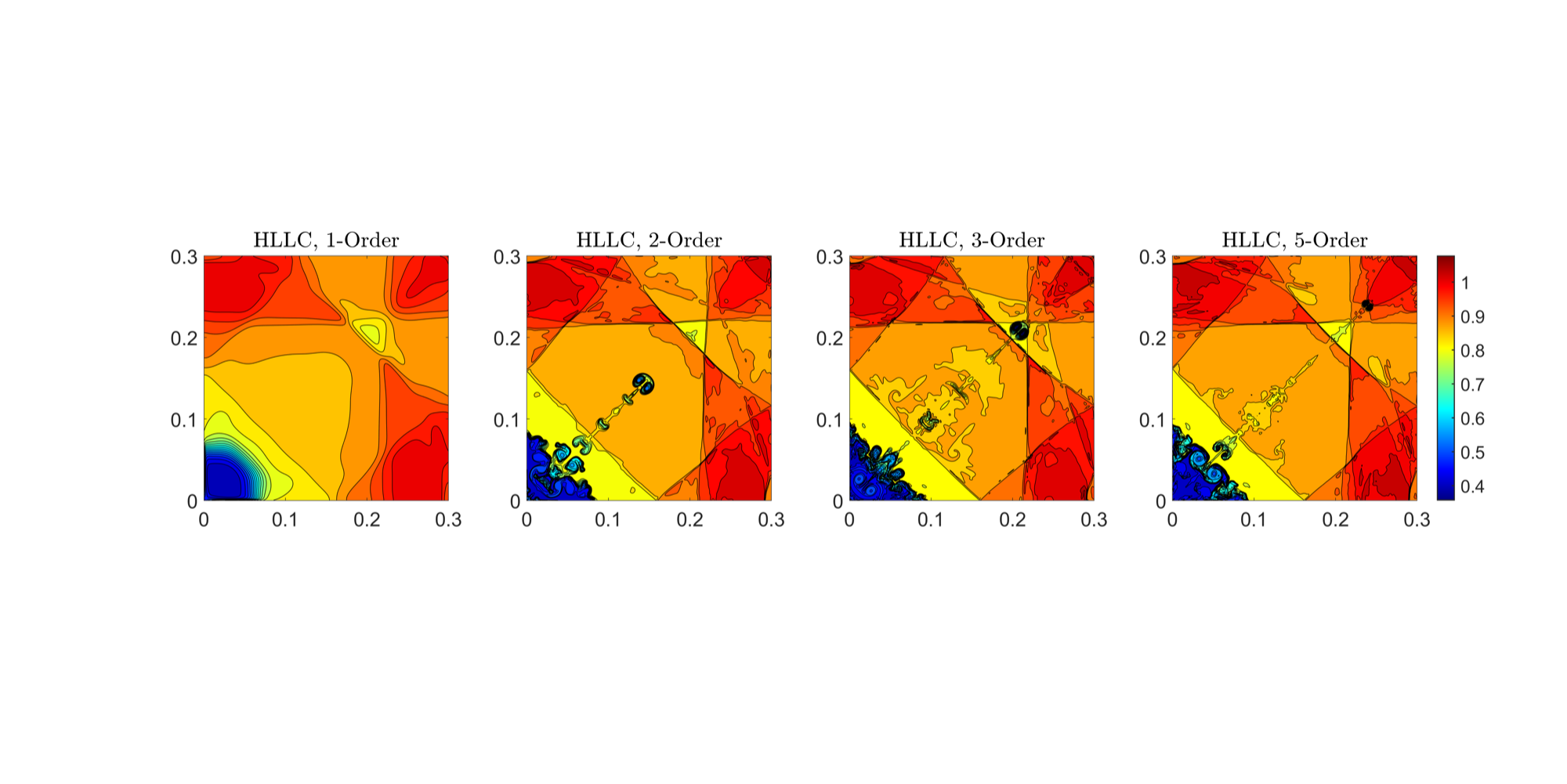}}
\vskip 10pt 
\centerline{\includegraphics[trim=1.8cm 2.7cm 0.7cm 2.5cm, clip, width=14.cm]{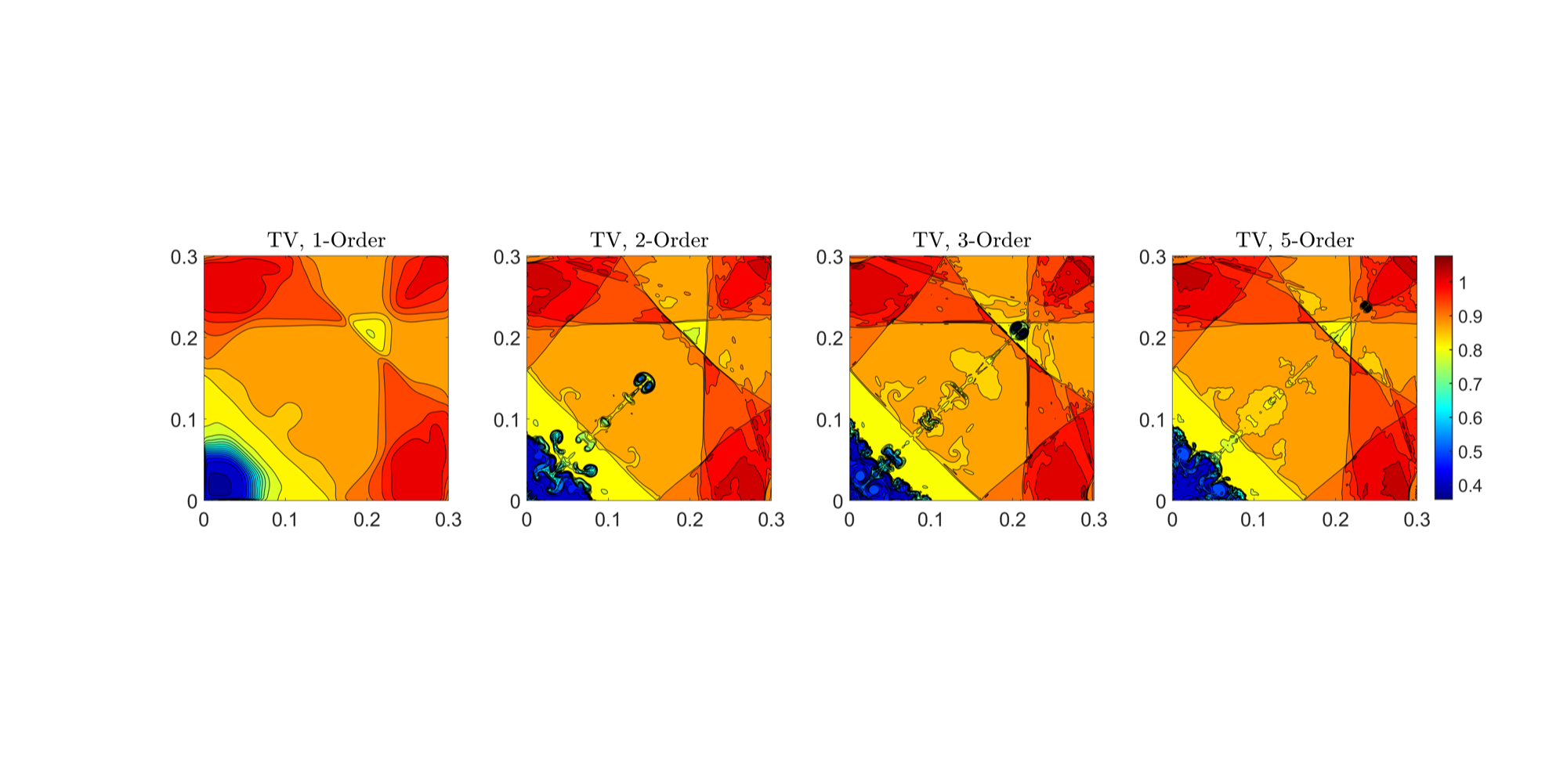}}
\vskip 10pt 
\centerline{\includegraphics[trim=1.8cm 2.7cm 0.7cm 2.5cm, clip, width=14.cm]{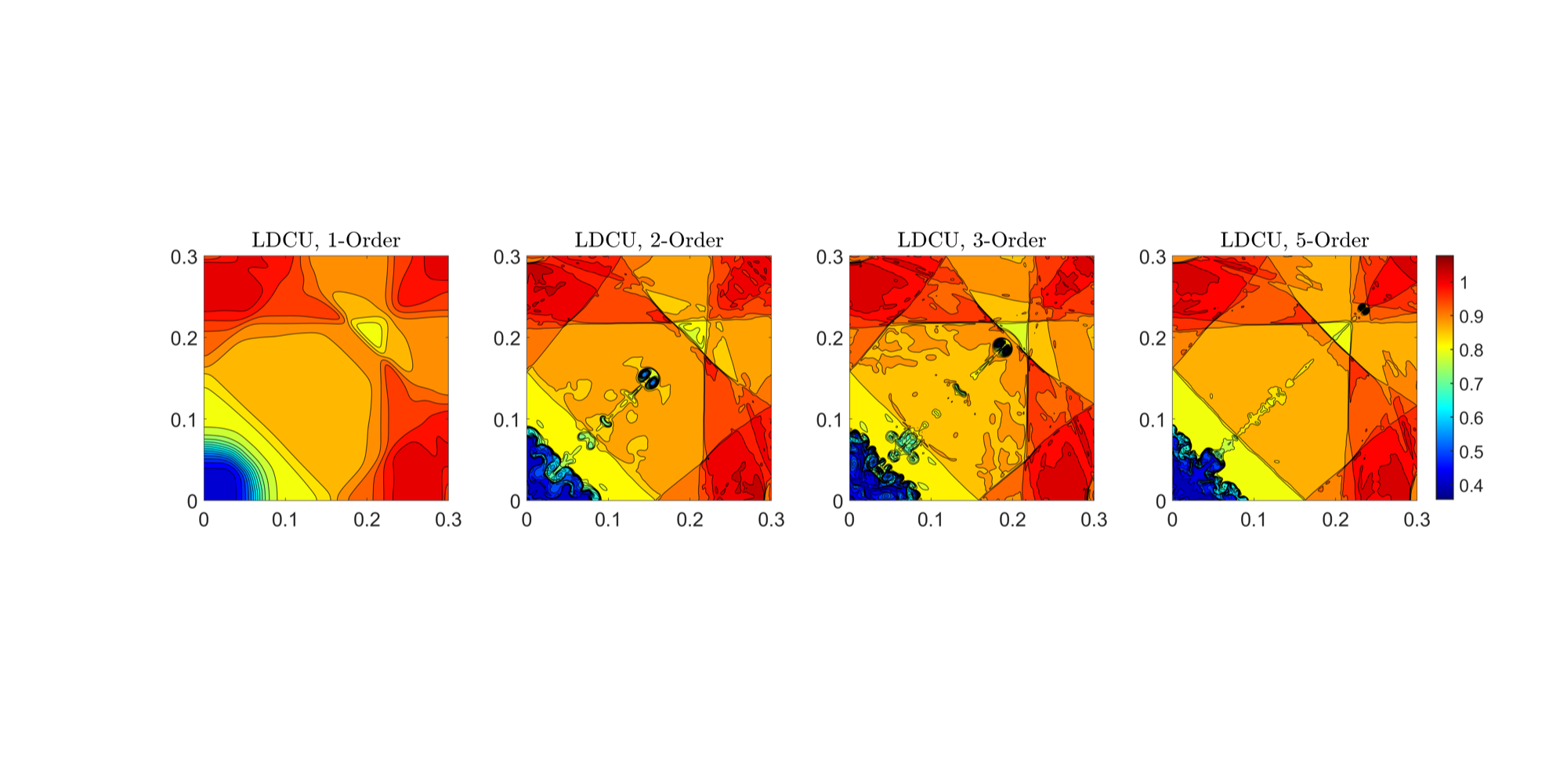}}
\vskip 10pt 
\centerline{\includegraphics[trim=1.8cm 2.7cm 0.7cm 2.5cm, clip, width=14.cm]{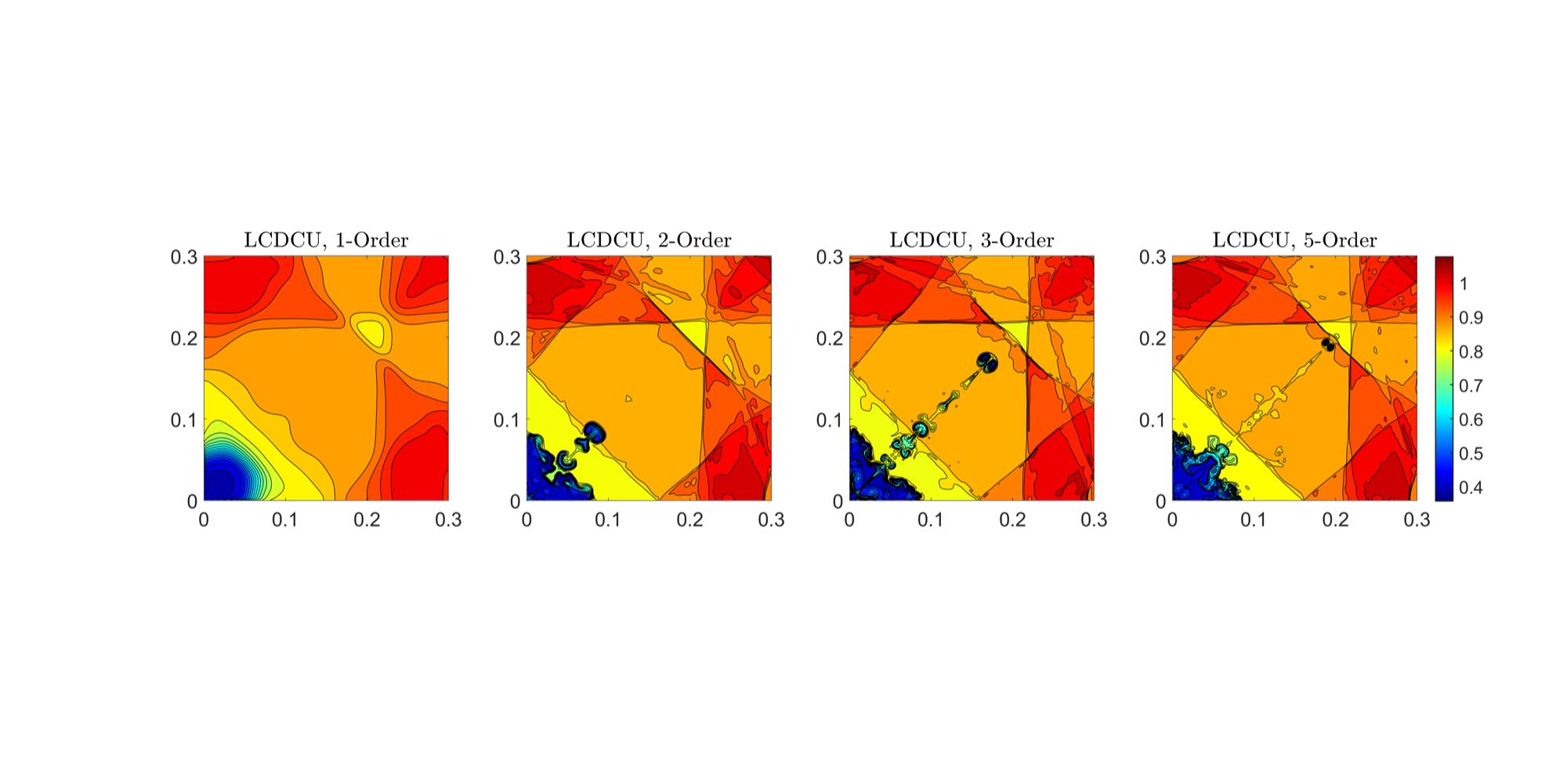}}
\caption{\sf{Example 13: Density $\rho$ computed by the 1-Order, 2-Order, 3-Order, and 5-Order HLL (top row), HLLC (second row), TV (third row), LDCU (fourth row), and LCDCU (bottom row) schemes.}\label{fig12}}
\end{figure}

It is instructive to check the efficiency of the studied schemes. To this end, we measure the CPU time consumed by the 1-Order, 2-Order, 3-Order and 5-Order LCDCU schemes and refine the mesh used by the other four schemes to the level that exactly the same CPU time is consumed to compute all four numerical solutions. The corresponding meshes are $400 \times 400$ for the 1-Order, 2-Order, 3-Order, and 5-Order LCDCU schemes, $570\times 570$, $568\times 568$, $482\times 482$, and $462\times 462$ for the HLL schemes, $568\times 568$, $565\times 565$, $479\times 479$, and $458\times 458$ for the HLLC schemes, $591\times 591$, $590\times 590$, $487\times 487$, and $466\times 466$ for the TV schemes, and $565\times 565$, $563\times 563$, $475\times 475$, and $455\times 455$ for the LDCU schemes. The obtained numerical results, presented in Figure \ref{fig12}, indicate that the four low-dissipation schemes are more efficient than the HLL counterparts. At the same time, it is noticed that the positions of the jets produced by the HLLC and TV schemes are consistent, which are slightly further than the ones produced by the LDCU schemes. The LCDCU schemes are more computationally expensive than the other three low-dissipation schemes, and the reason is that the LCDCU fluxes defined in \S \ref{sec3.5} are more complicated than the others.  

\subsubsection*{Example 14---2-D Riemann Problem (Configuration 3)} In this example, we consider Configuration 3 of the 2-D Riemann problems from
\cite{Kurganov02} (see also \cite{Schulz93,Schulz93a,Zheng01}) with the following initial conditions:
\begin{equation*}
(\rho,u,v,p)(x,y,0)=\begin{cases}
(1.5,0,0,1.5),&x>1,~y>1,\\
(0.5323,1.206,0,0.3),&x<1,~y>1,\\
(0.138,1.206,1.206,0.029),&x<1,~y<1,\\
(0.5323,0,1.206,0.3),&x>1,~y<1.
\end{cases}
\end{equation*}
which are prescribed in the computational domain $[0,1.2]\times[0,1.2]$ subject to the free boundary conditions.

\begin{figure}[ht!]
\centerline{\includegraphics[trim=1.8cm 2.7cm 0.7cm 2.5cm, clip, width=14.cm]{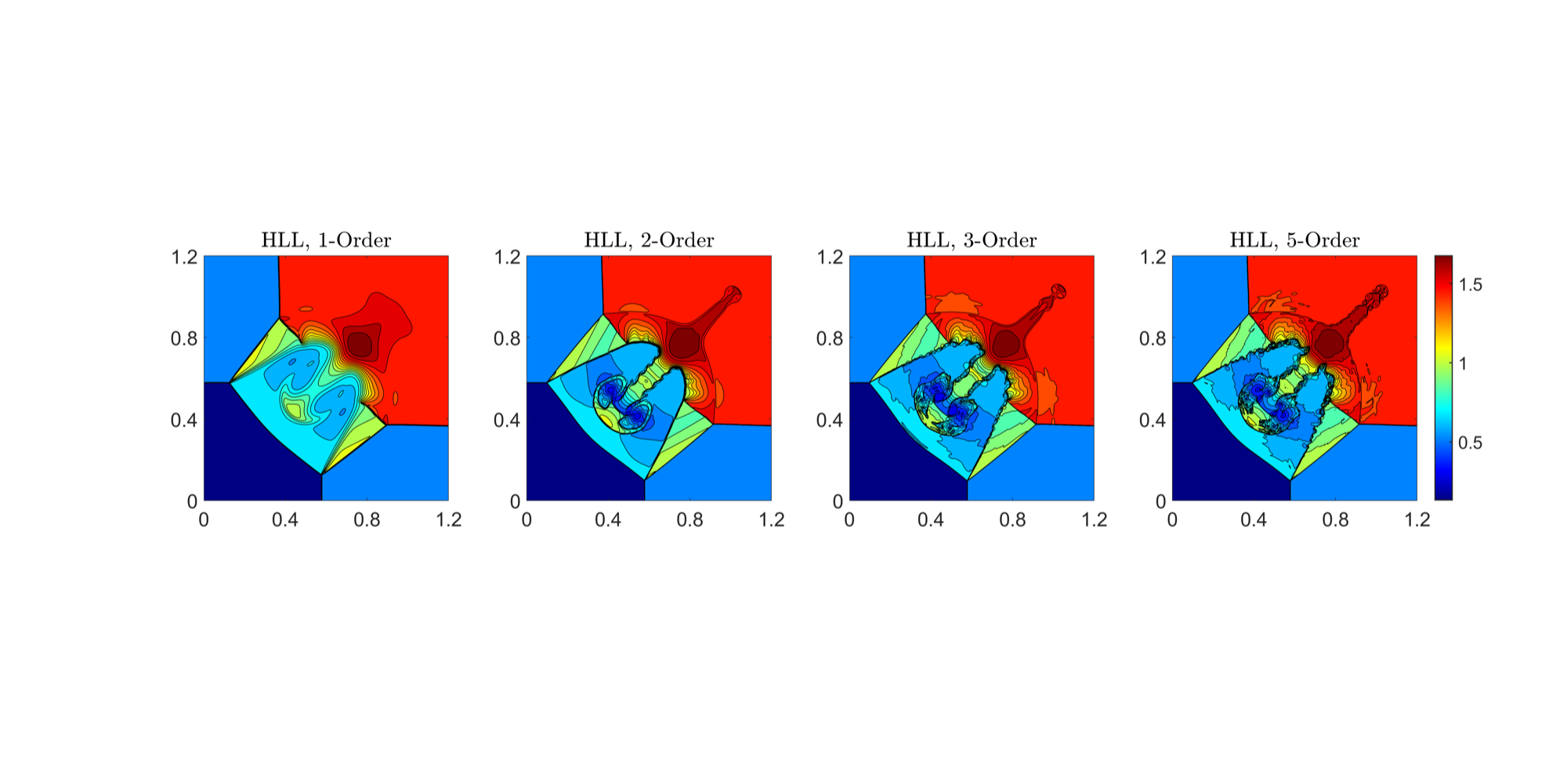}}
\vskip 10pt 
\centerline{\includegraphics[trim=1.8cm 2.7cm 0.7cm 2.5cm, clip, width=14.cm]{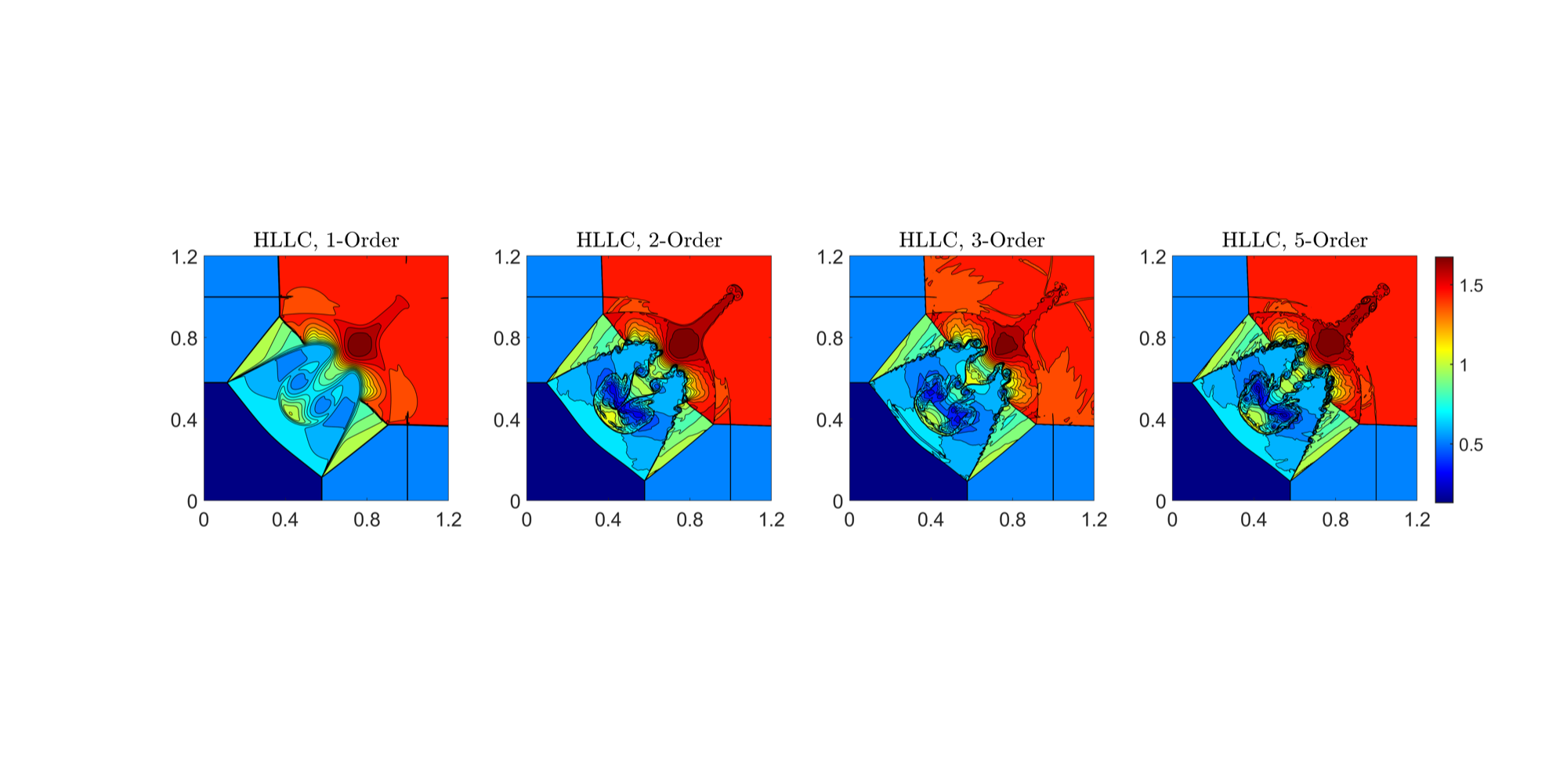}}
\vskip 10pt 
\centerline{\includegraphics[trim=1.8cm 2.7cm 0.7cm 2.5cm, clip, width=14.cm]{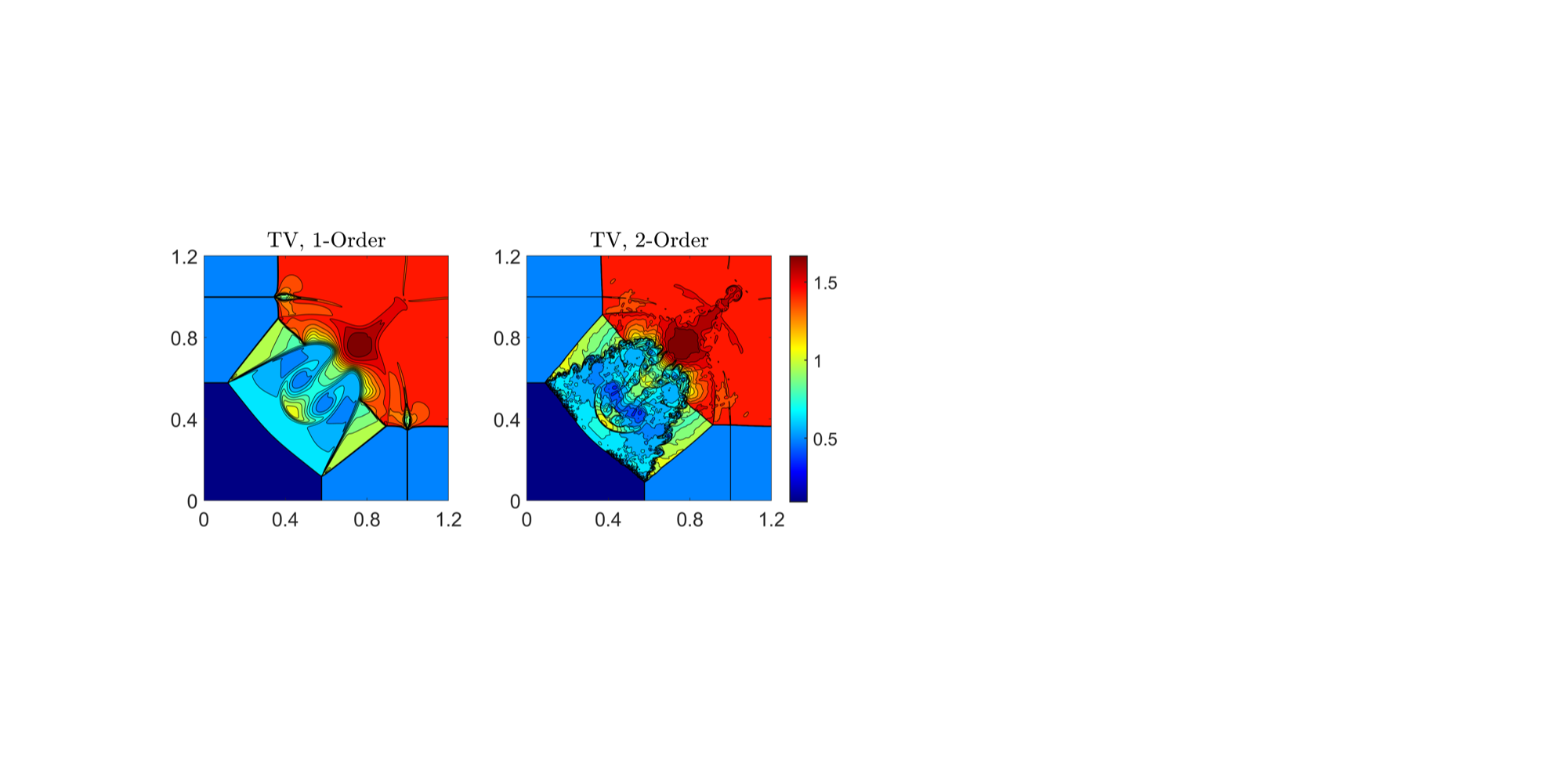}}
\vskip 10pt 
\centerline{\includegraphics[trim=1.8cm 2.7cm 0.7cm 2.5cm, clip, width=14.cm]{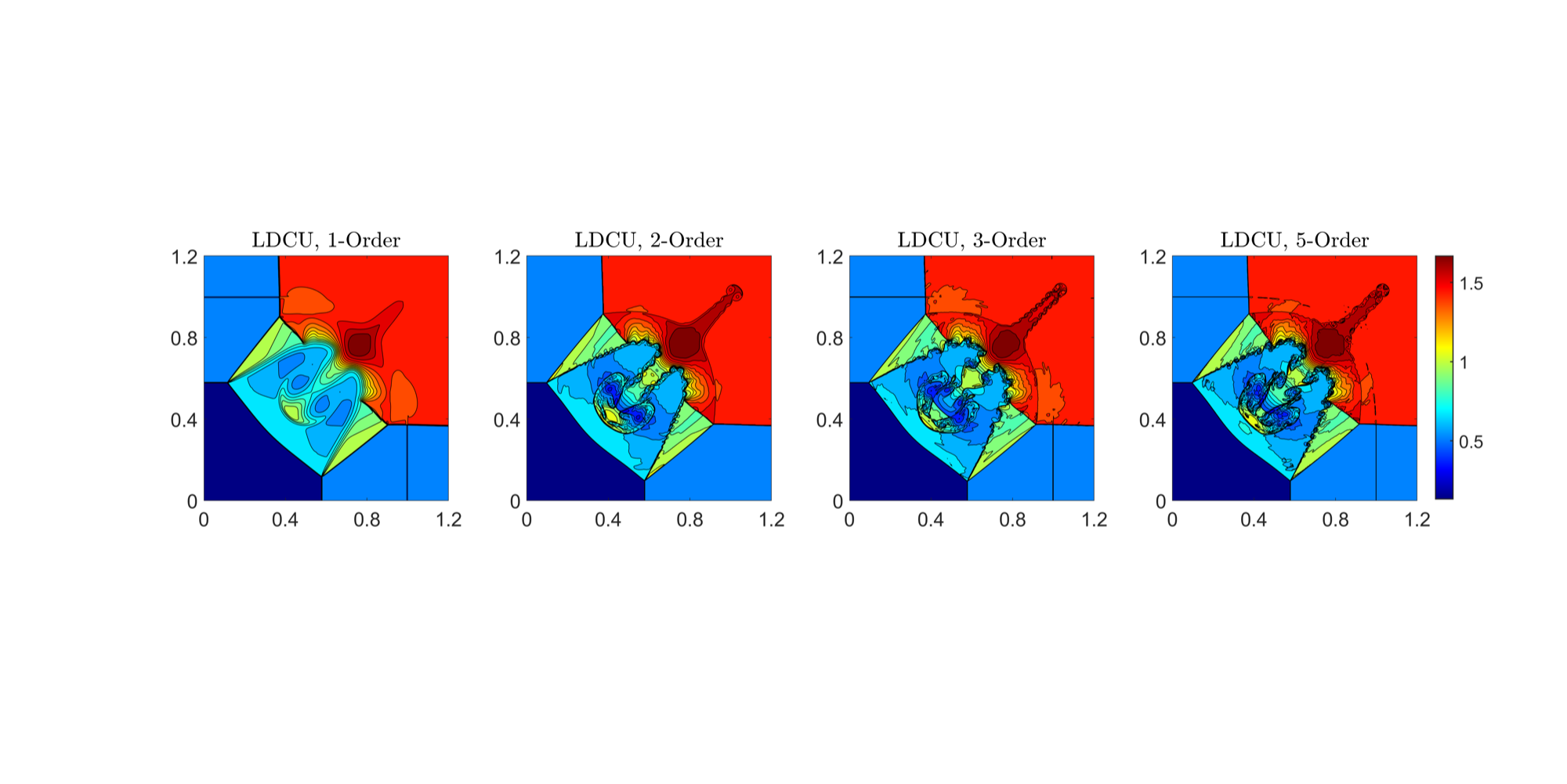}}
\vskip 10pt 
\centerline{\includegraphics[trim=1.8cm 2.7cm 0.7cm 2.5cm, clip, width=14.cm]{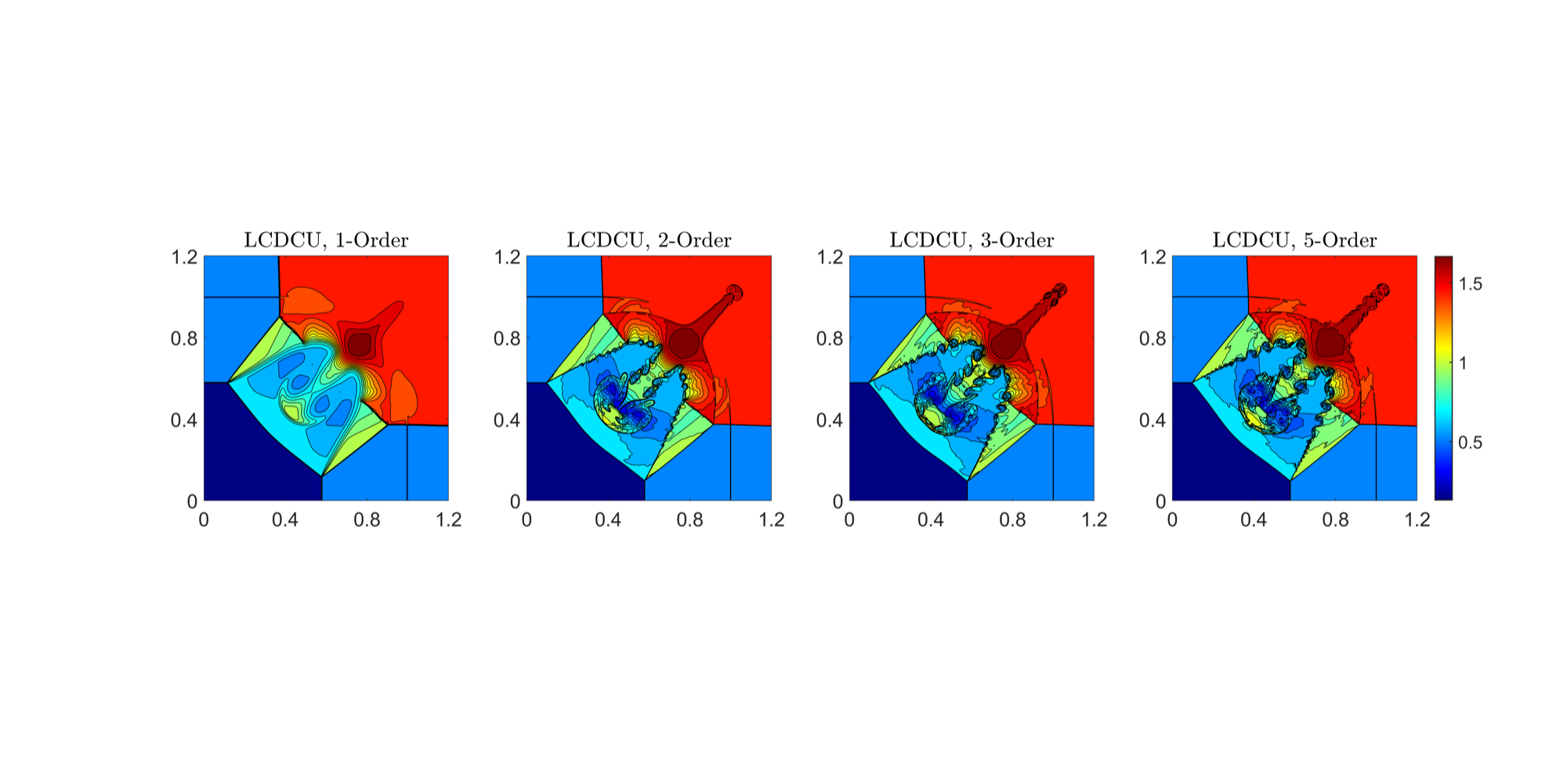}}
\caption{\sf{Example 14: Density $\rho$ computed by the 1-Order, 2-Order, 3-Order, and 5-Order HLL (top row), HLLC (second row), TV (third row; the third- and fifth-order results are omitted because of loss of positivity), LDCU (fourth row), and LCDCU (bottom row) schemes.}\label{fig11c}}
\end{figure}

We compute the numerical solution until the final time $t=1$ by the 1-Order, 2-Order, 3-Order, and 5-Order schemes on a uniform mesh with $\dx=\dy=6/5000$ and plot the obtained results in Figure \ref{fig11c}, where one can see that the HLLC, LDCU, and LCDCU schemes resolve more small-scale structures than the HLL scheme in capturing a sideband instability of the jet in the zones of strong along-jet velocity shear and the instability along the jet’s neck. However, one can see that the TV schemes can capture more details of the instability along the jet’s neck, but there are some oscillations in the numerical results.  At the same time, the 3-Order and 5-Order TV schemes fail in this simulation because numerical oscillations lead to negative pressure.

\subsubsection*{Example 15---2-D Riemann Problem (Configuration 6)} In this example, we consider Configuration 6 of the 2-D Riemann problems from
\cite{Kurganov02} (see also \cite{Schulz93,Schulz93a,Zheng01}) with the following initial conditions:
\begin{equation*}
(\rho,u,v,p)\Big|_{(x,y,0)}=\begin{cases}(1,0.75,-0.5,1),&x>0.5,~y>0.5,\\(2,0.75,0.5,1),&x<0.5,~y>0.5,\\(1,-0.75,0.5,1),&x<0.5,~y<0.5,\\
(3,-0.75,-0.5,1),&x>0.5,~y<0.5,\end{cases}
\end{equation*}
which are prescribed in the computational domain $[0,1]\times[0,1]$ subject to the free boundary conditions.

\begin{figure}[ht!]
\centerline{\includegraphics[trim=1.8cm 2.7cm 0.7cm 2.5cm, clip, width=14.cm]{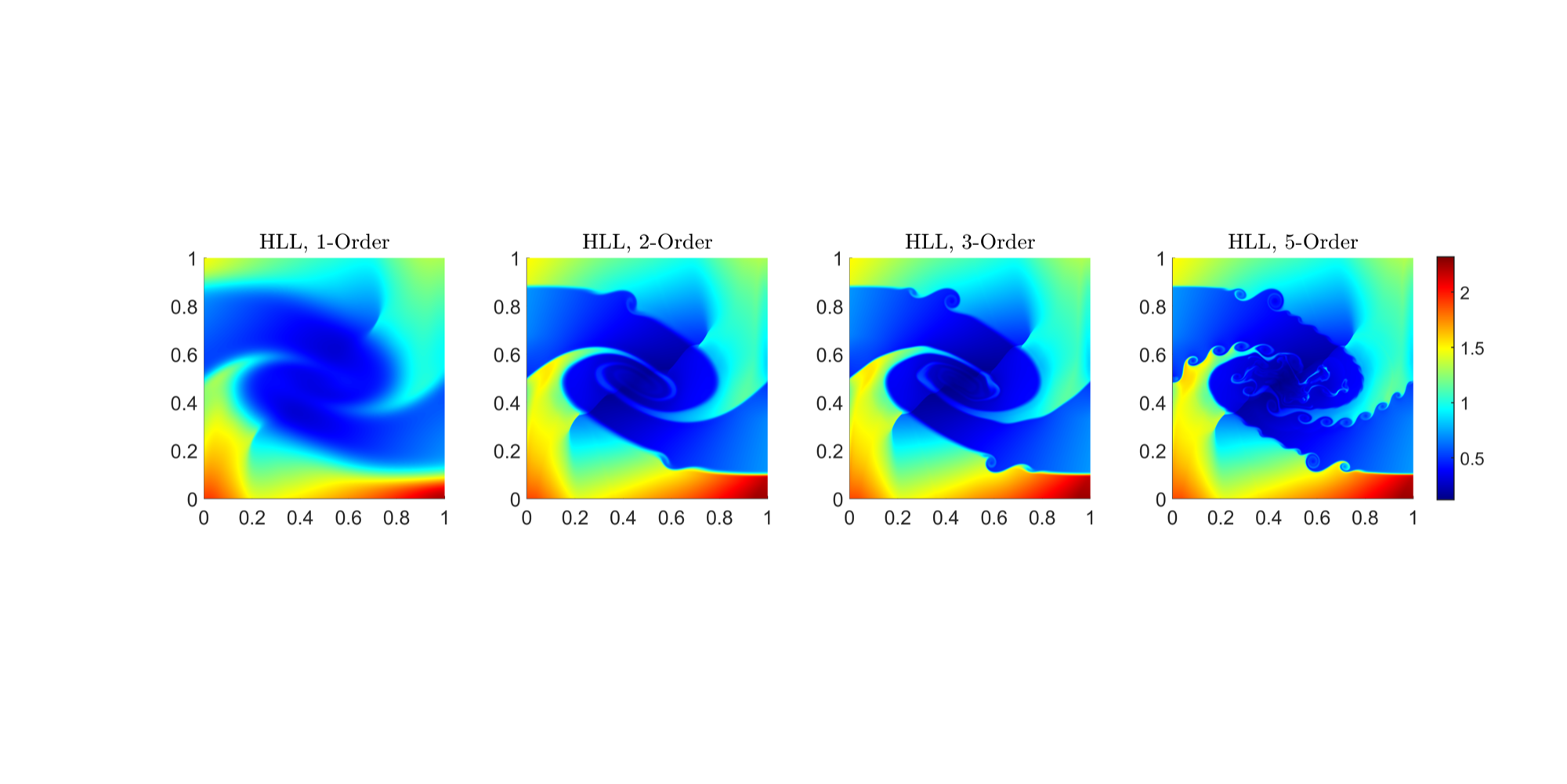}}
\vskip 10pt 
\centerline{\includegraphics[trim=1.8cm 2.7cm 0.7cm 2.5cm, clip, width=14.cm]{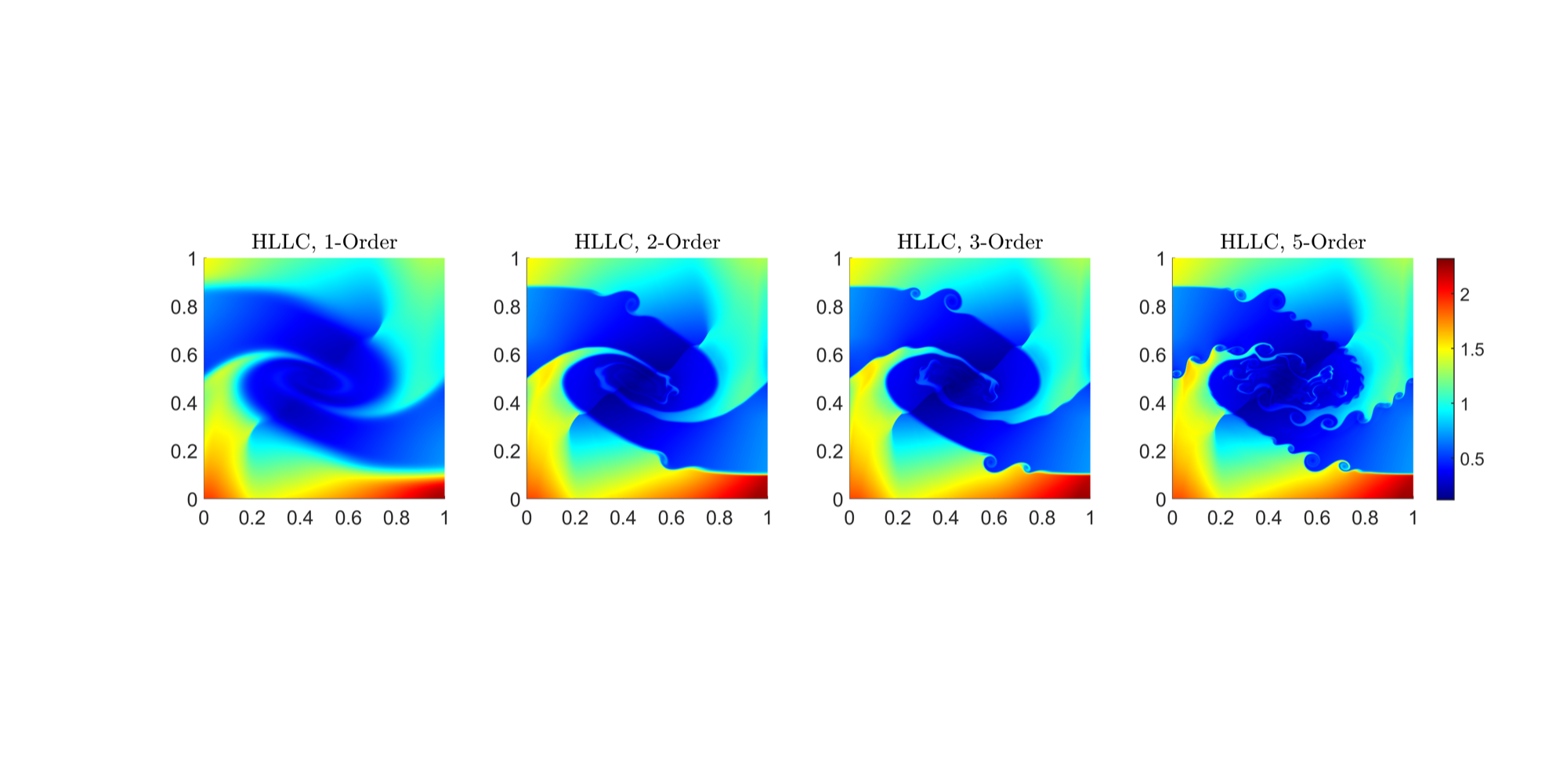}}
\vskip 10pt 
\centerline{\includegraphics[trim=1.8cm 2.7cm 0.7cm 2.5cm, clip, width=14.cm]{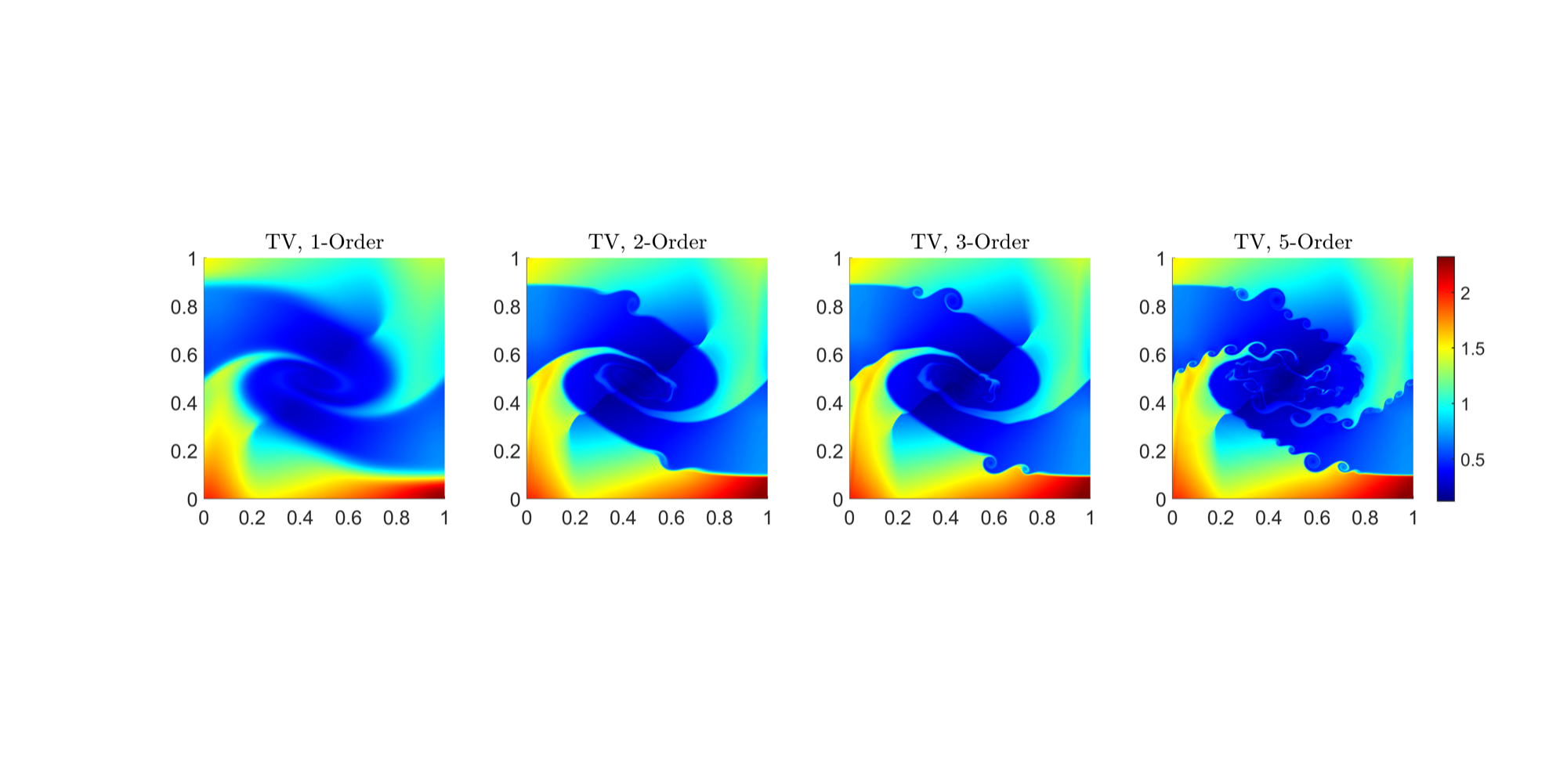}}
\vskip 10pt 
\centerline{\includegraphics[trim=1.8cm 2.7cm 0.7cm 2.5cm, clip, width=14.cm]{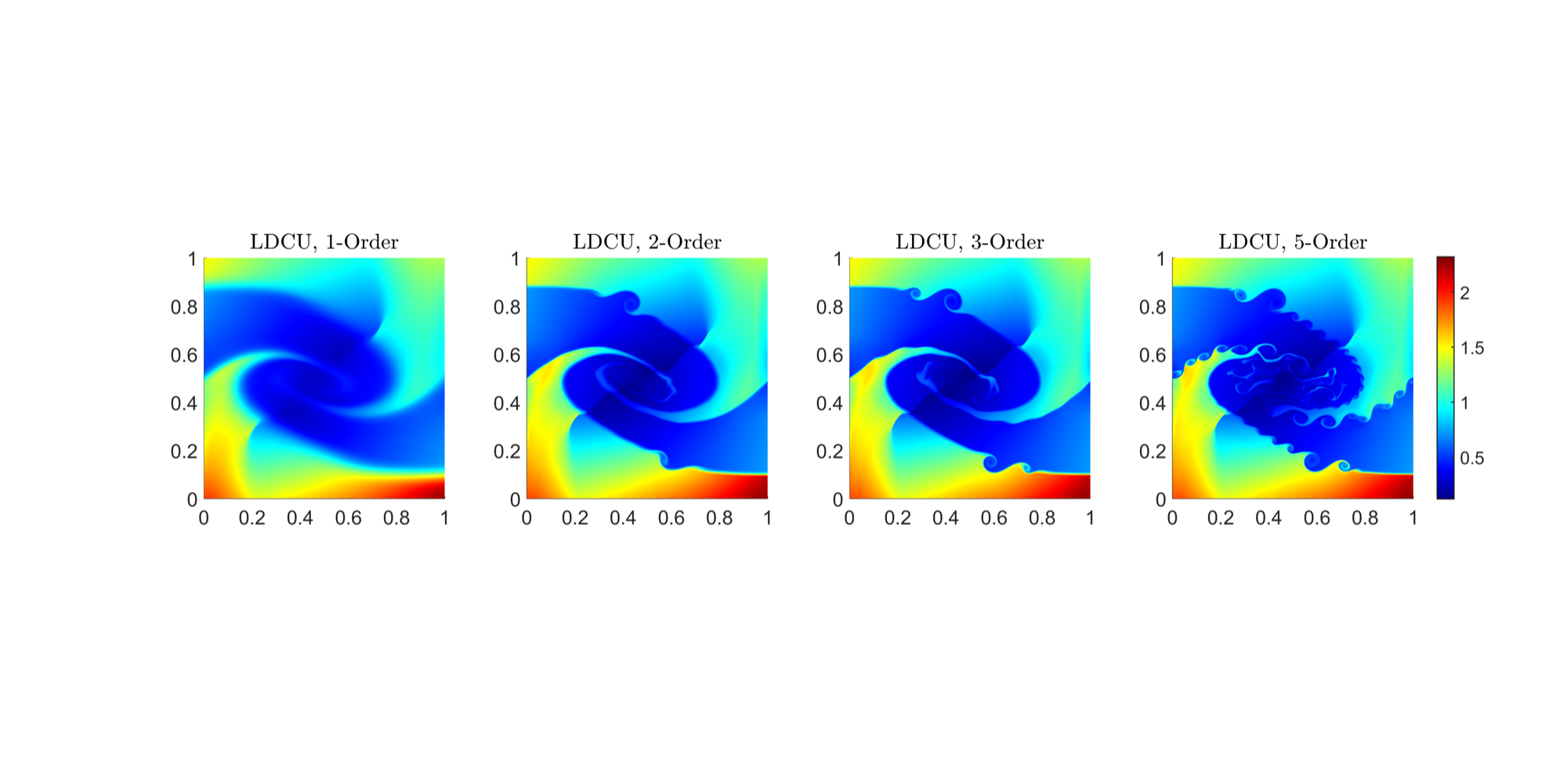}}
\vskip 10pt 
\centerline{\includegraphics[trim=1.8cm 2.7cm 0.7cm 2.5cm, clip, width=14.cm]{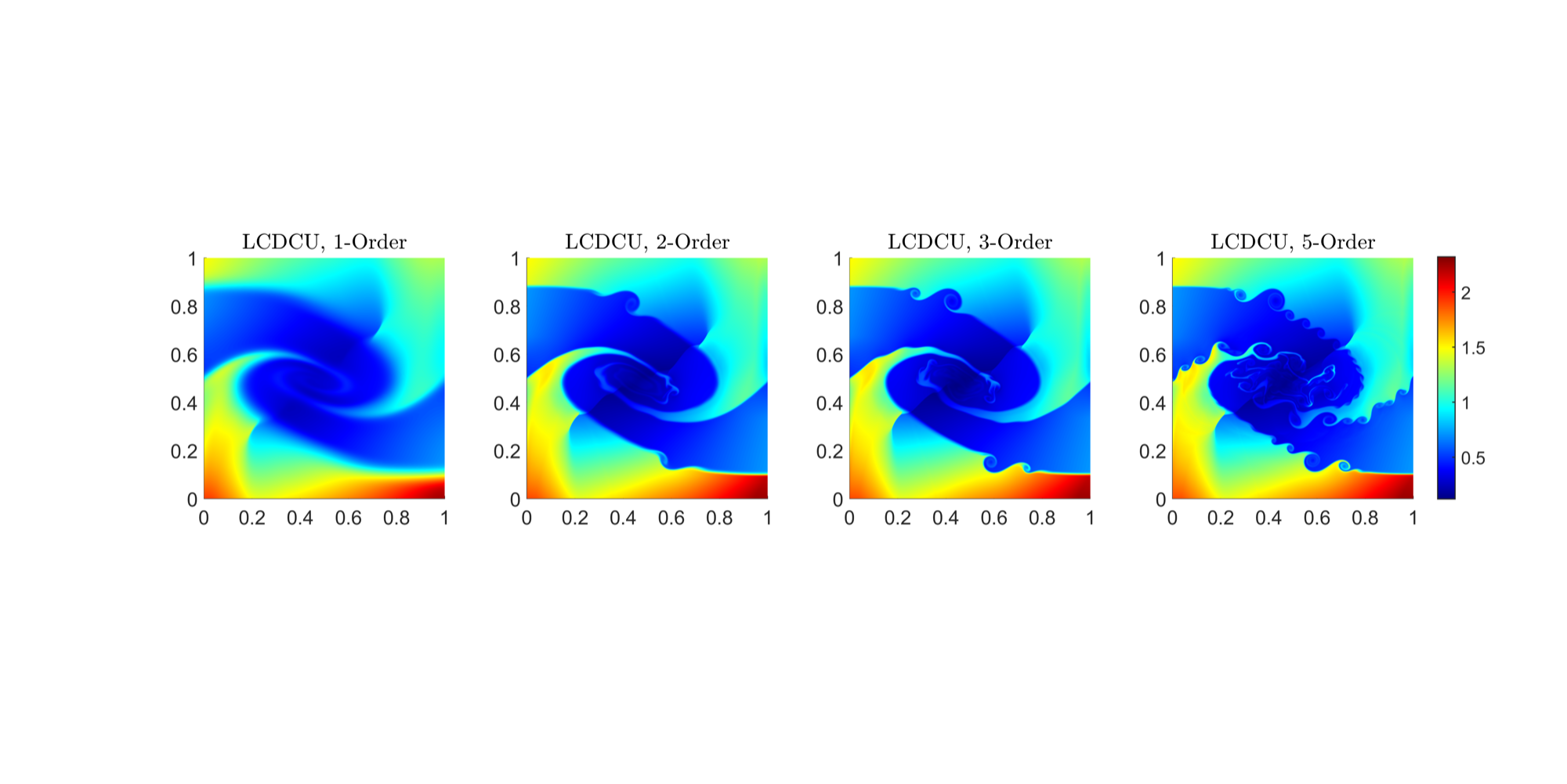}}
\caption{\sf{Example 15: Density $\rho$ computed by the 1-Order, 2-Order, 3-Order, and 5-Order HLL (top row), HLLC (second row), TV (third row), LDCU (fourth row), and LCDCU (bottom row) schemes.}\label{fig11b}}
\end{figure}

We compute the numerical solution until the final time $t=1$ by the 1-Order, 2-Order, 3-Order, and 5-Order schemes on a uniform mesh with $\dx=\dy=1/600$ and plot the obtained results in Figure \ref{fig11b}, where one can see that the four low-dissipation schemes produce more intricate vortex structures than the HLL schemes, consistently with their lower numerical dissipation. 

\subsubsection*{Example 16---Kelvin-Helmholtz (KH) Instability}
In this example, we study the KH instability taken from \cite{Fjordholm16,Panuelos20} (see also \cite{Feireisl21,CCHKL_22}). We take the following initial data:
\begin{equation*}
\begin{aligned}
&(\rho(x,y,0),u(x,y,0))=\begin{cases}
(1,-0.5+0.5e^{(y+0.25)/L}),&y<-0.25,\\
(2,0.5-0.5e^{(-y-0.25)/L}),&-0.25<y<0,\\
(2,0.5-0.5e^{(y-0.25)/L}),&0<y<0.25,\\
(1,-0.5+0.5e^{(0.25-y)/L}),&y>0.25,
\end{cases}\\
&v(x,y,0)=0.01\sin(4\pi x),\qquad p(x,y,0)\equiv1.5,
\end{aligned}
\end{equation*}
where $L$ is a smoothing parameter (we take $L=0.00625$), which corresponds to a thin shear interface with a perturbed vertical velocity field $v$ in the conducted simulations. The periodic boundary conditions are imposed on all four sides of the computational domain $[-0.5,0.5]\times[-0.5,0.5]$.

\begin{figure}[ht!]
\centerline{\includegraphics[trim=1.8cm 2.7cm 0.7cm 2.5cm, clip, width=14.cm]{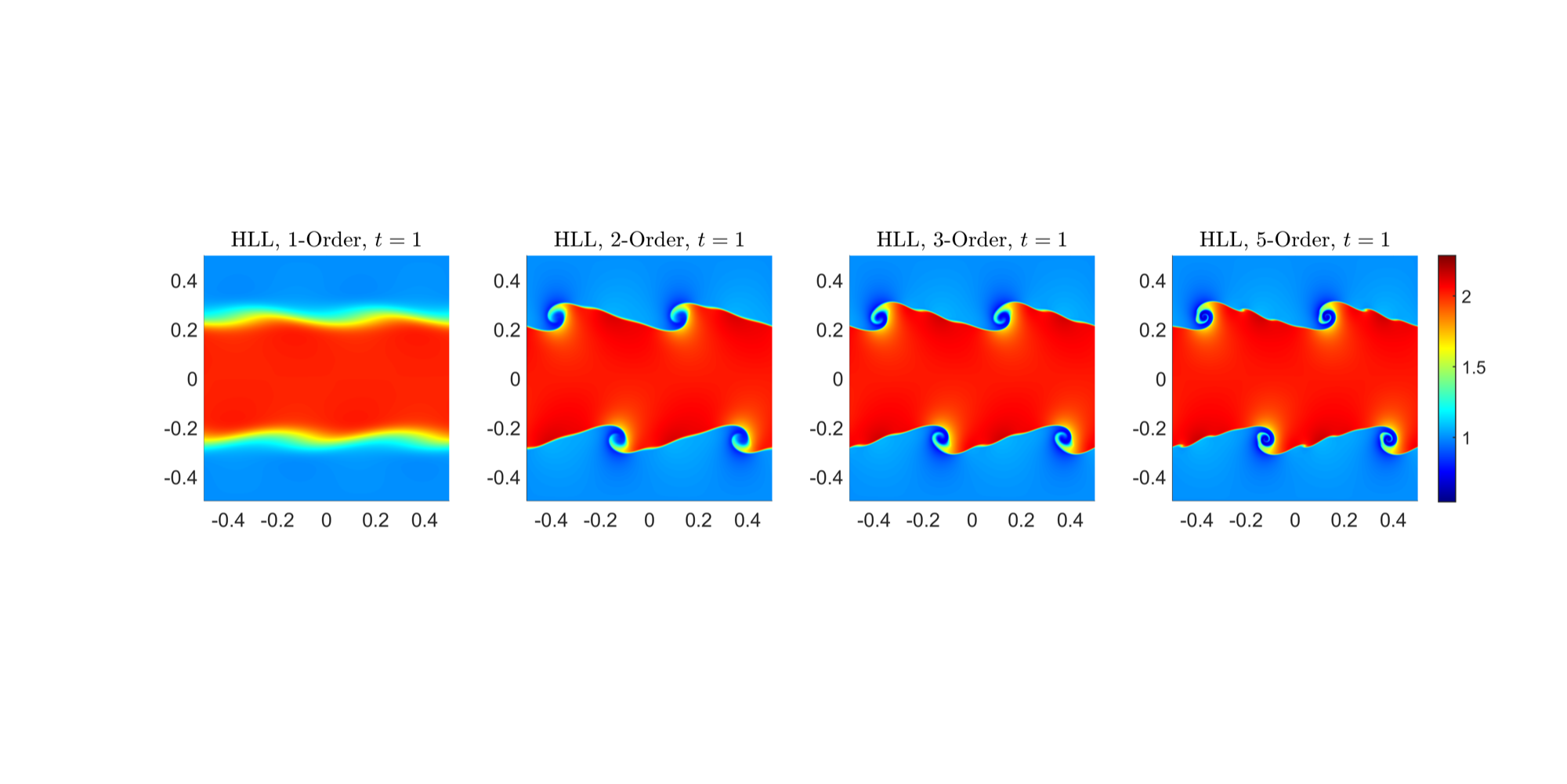}}
\vskip 15pt
\centerline{\includegraphics[trim=1.8cm 2.7cm 0.7cm 2.5cm, clip, width=14.cm]{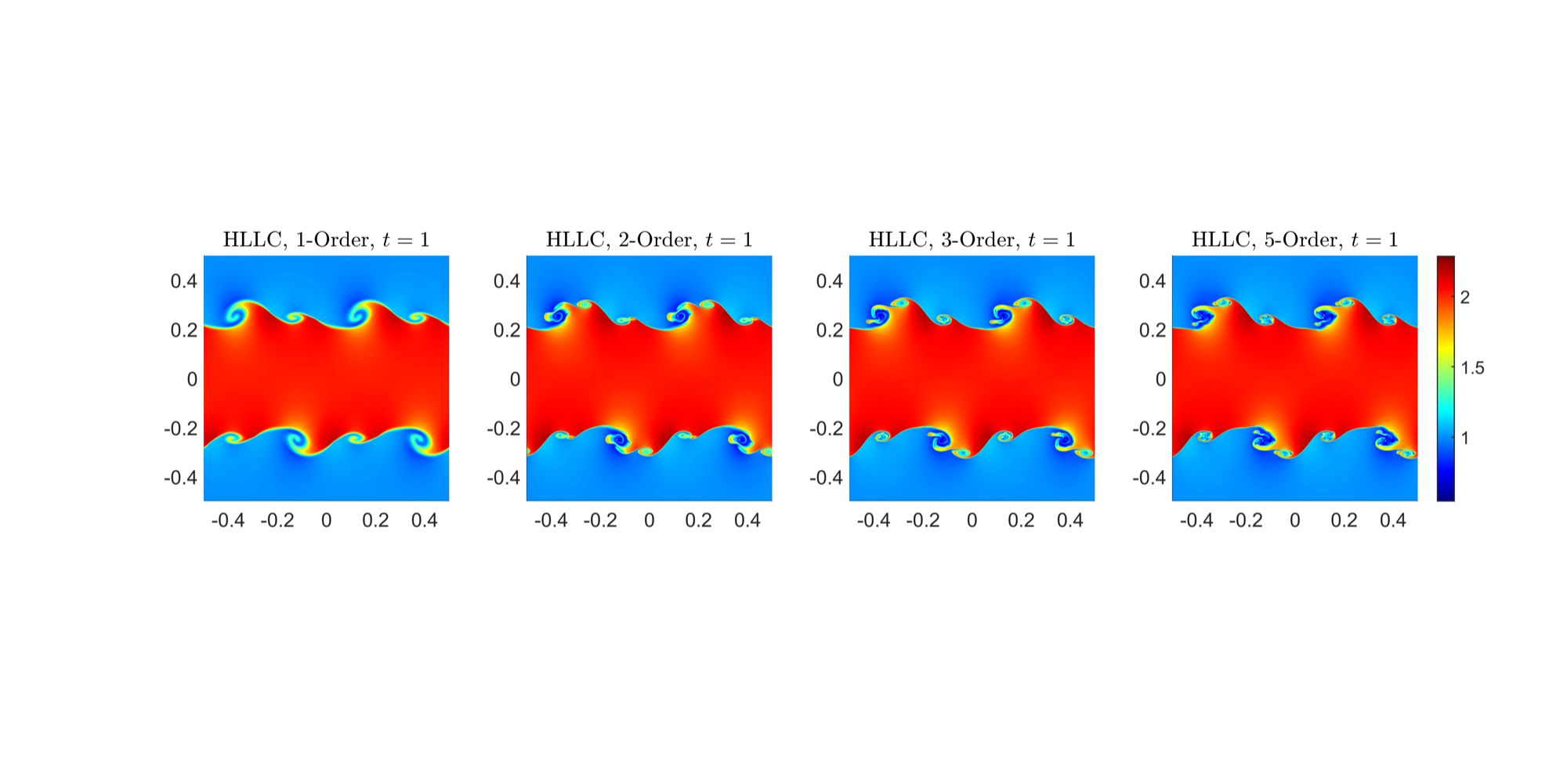}}
\vskip 15pt
\centerline{\includegraphics[trim=1.8cm 2.7cm 0.7cm 2.5cm, clip, width=14.cm]{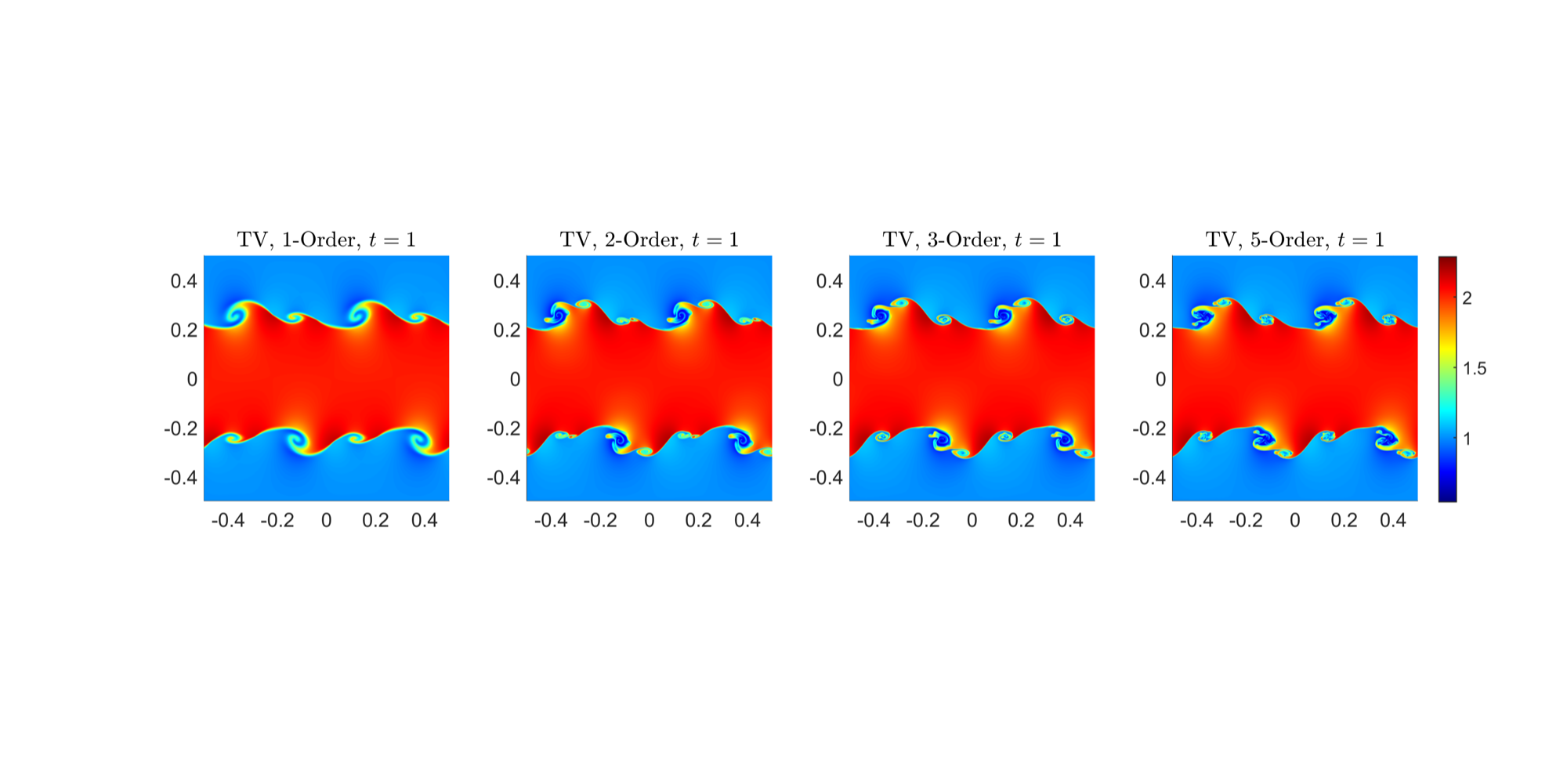}}
\vskip 15pt
\centerline{\includegraphics[trim=1.8cm 2.7cm 0.7cm 2.5cm, clip, width=14.cm]{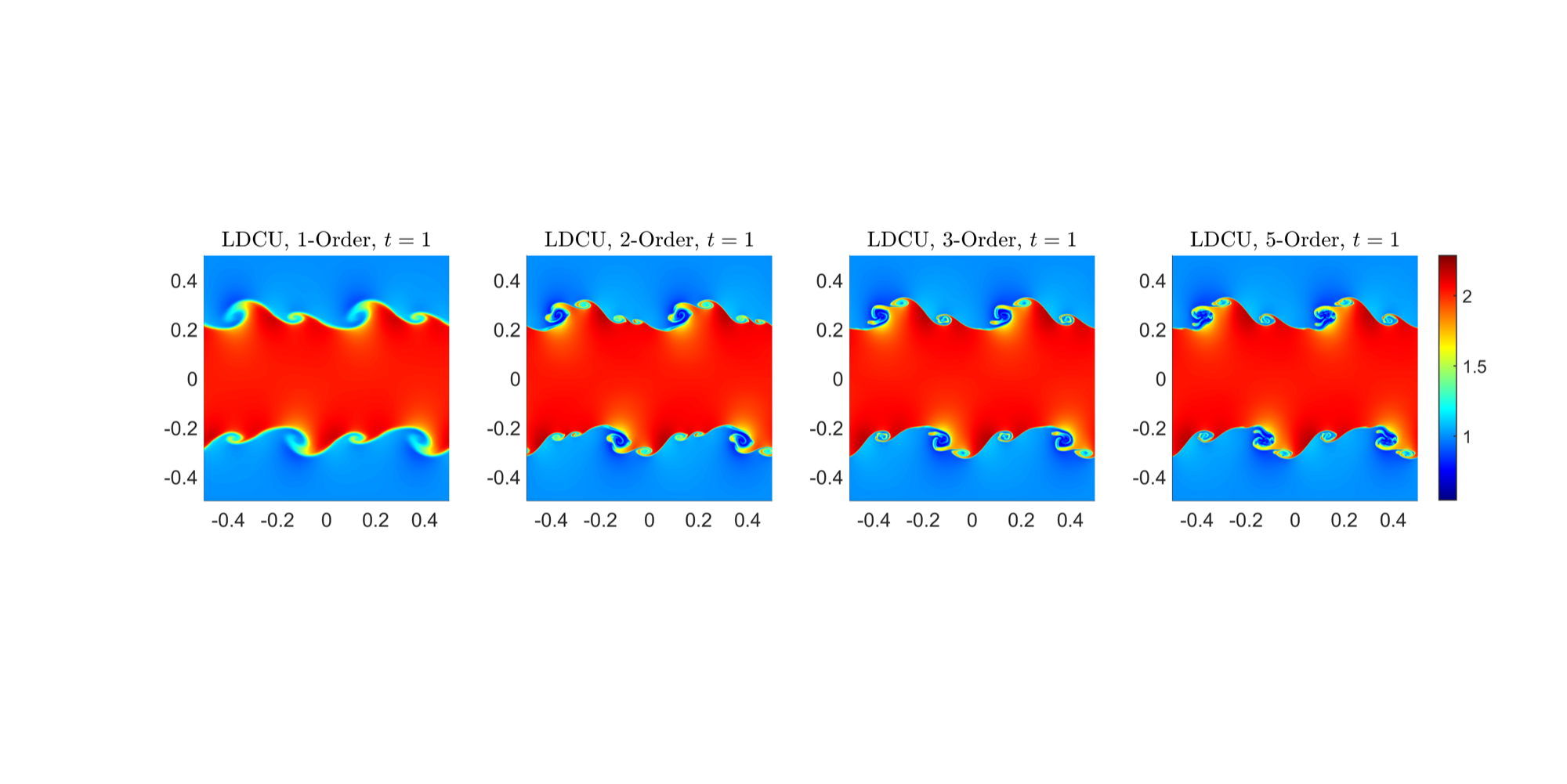}}
\vskip 15pt
\centerline{\includegraphics[trim=1.8cm 2.7cm 0.7cm 2.5cm, clip, width=14.cm]{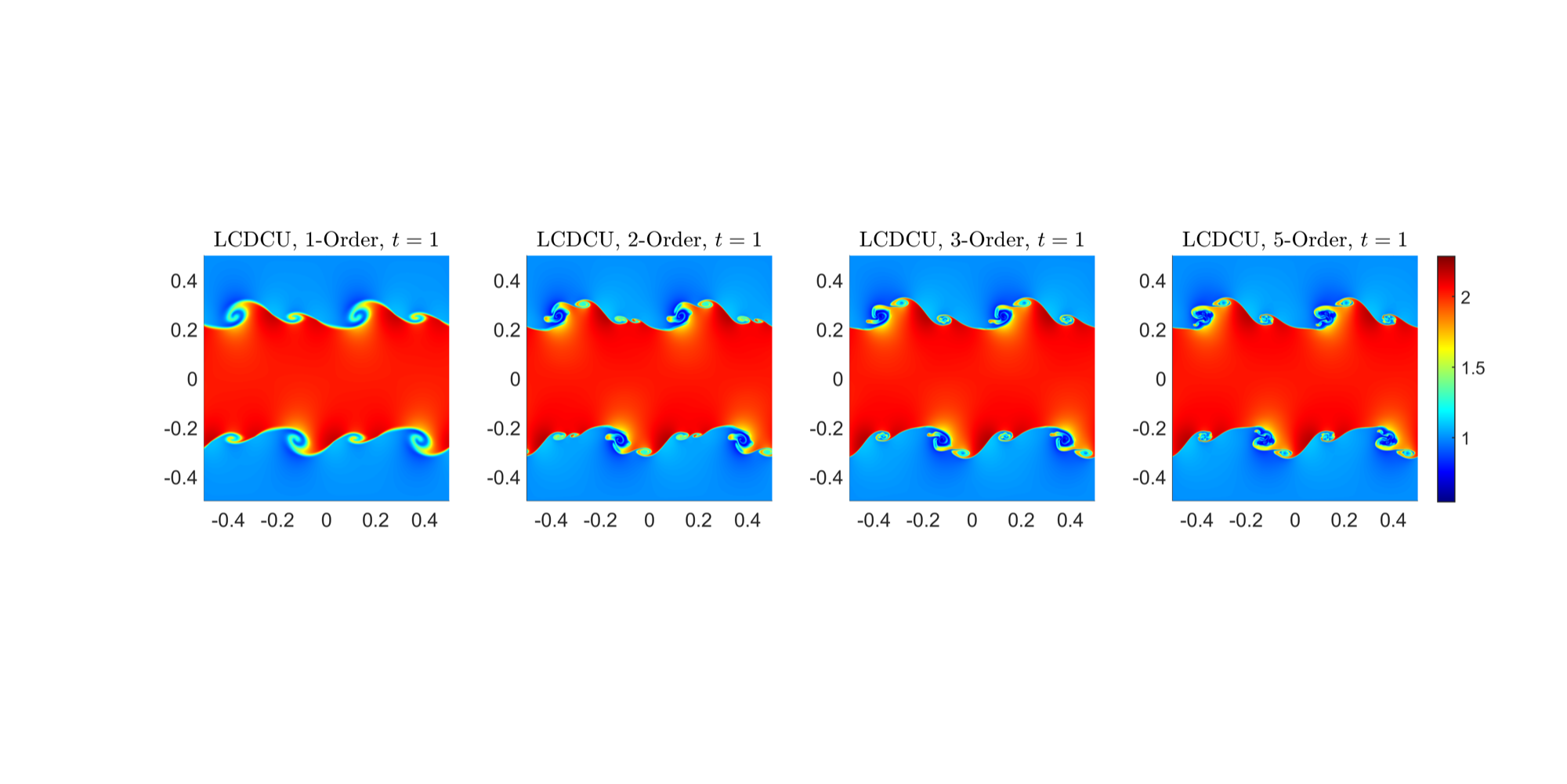}}
\caption{\sf Example 16: Time snapshots of density $\rho$ computed by the 1-Order, 2-Order, 3-Order, and 5-Order HLL (top row), HLLC (second row), TV (third row), LDCU (fourth row), and LCDCU (bottom row) schemes at $t=1$.\label{fig12a}}
\end{figure}

\begin{figure}[ht!]
\centerline{\includegraphics[trim=1.8cm 2.7cm 0.7cm 2.5cm, clip, width=14.cm]{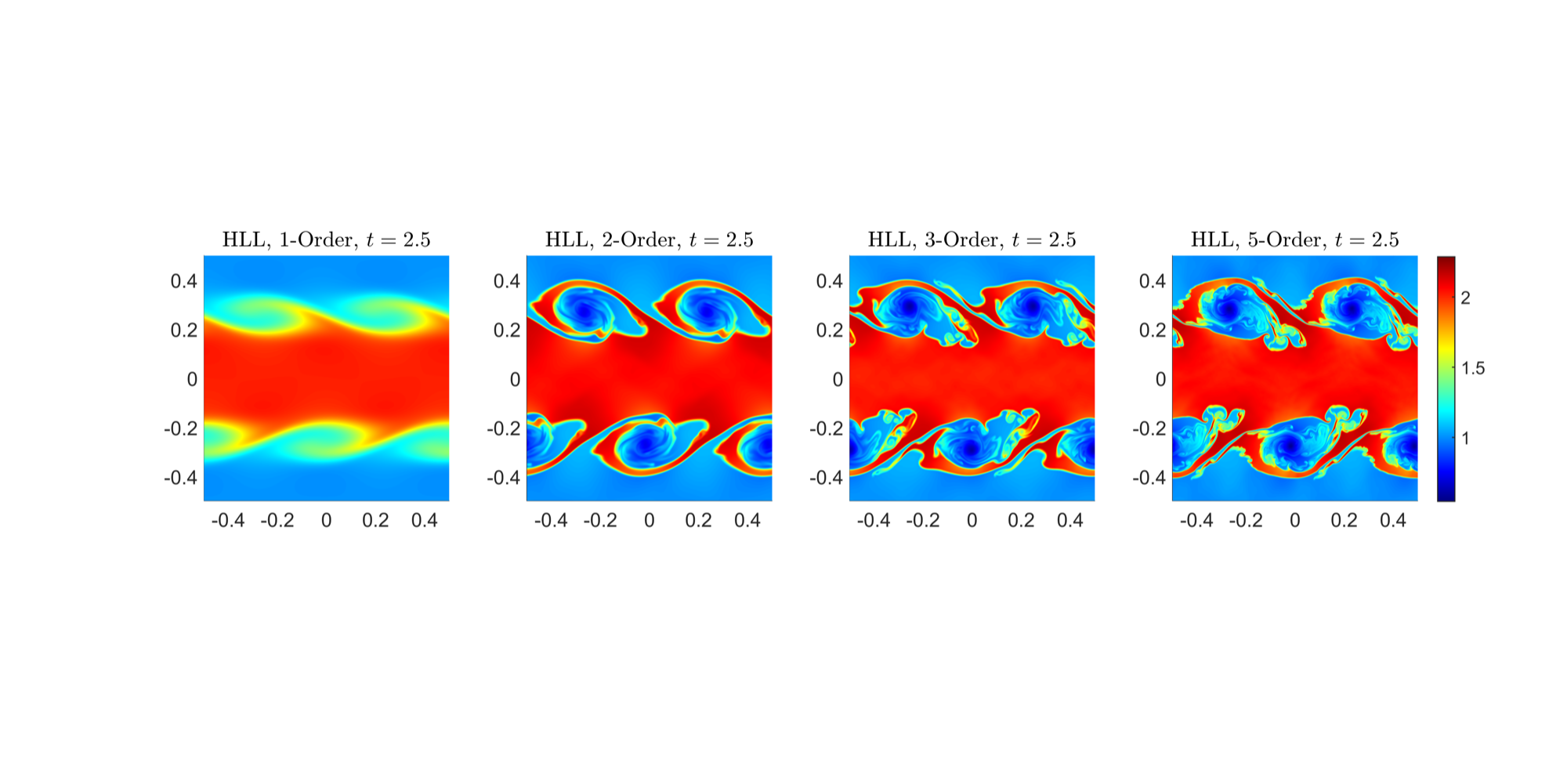}}
\vskip 15pt
\centerline{\includegraphics[trim=1.8cm 2.7cm 0.7cm 2.5cm, clip, width=14.cm]{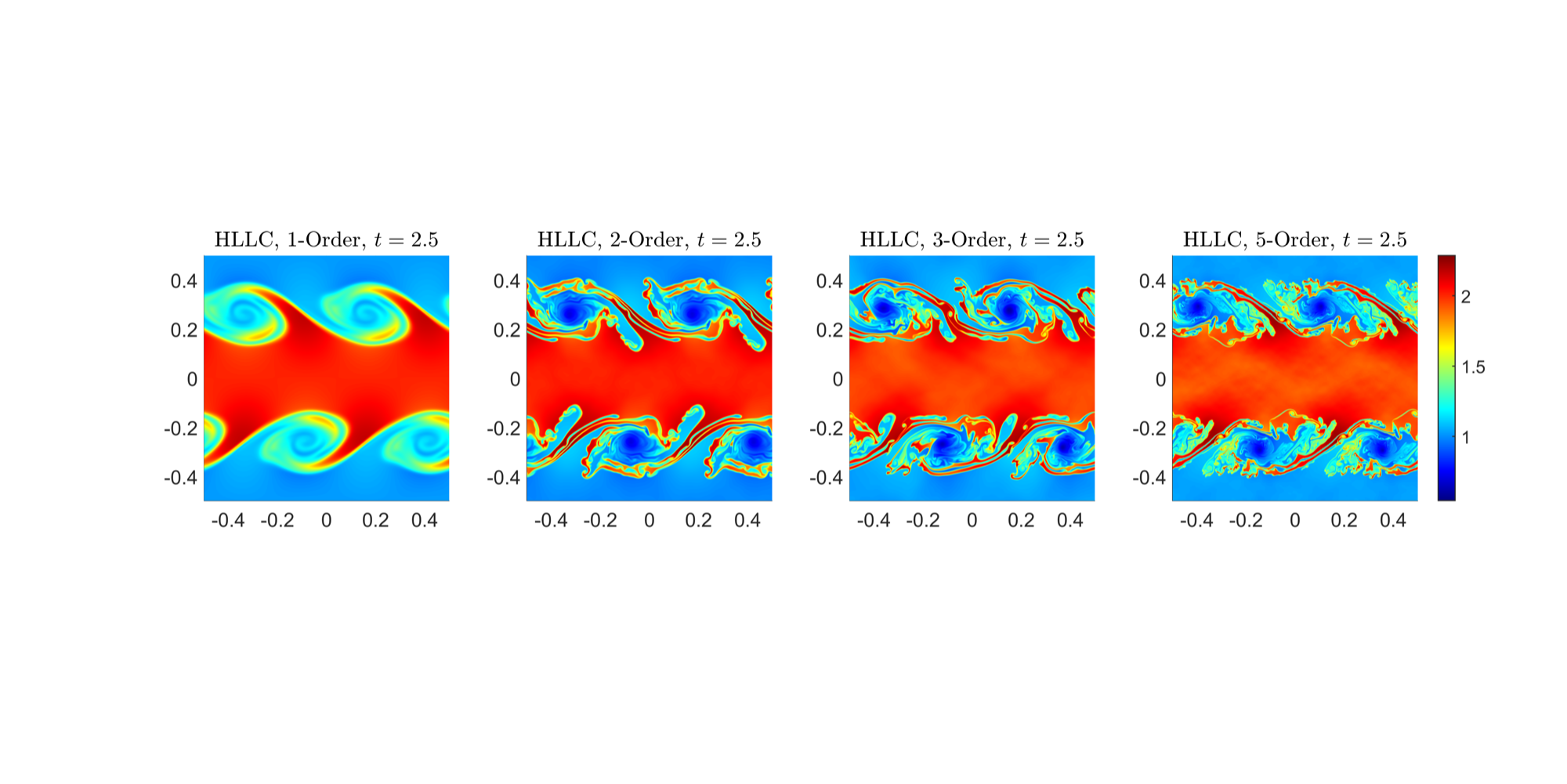}}
\vskip 15pt
\centerline{\includegraphics[trim=1.8cm 2.7cm 0.7cm 2.5cm, clip, width=14.cm]{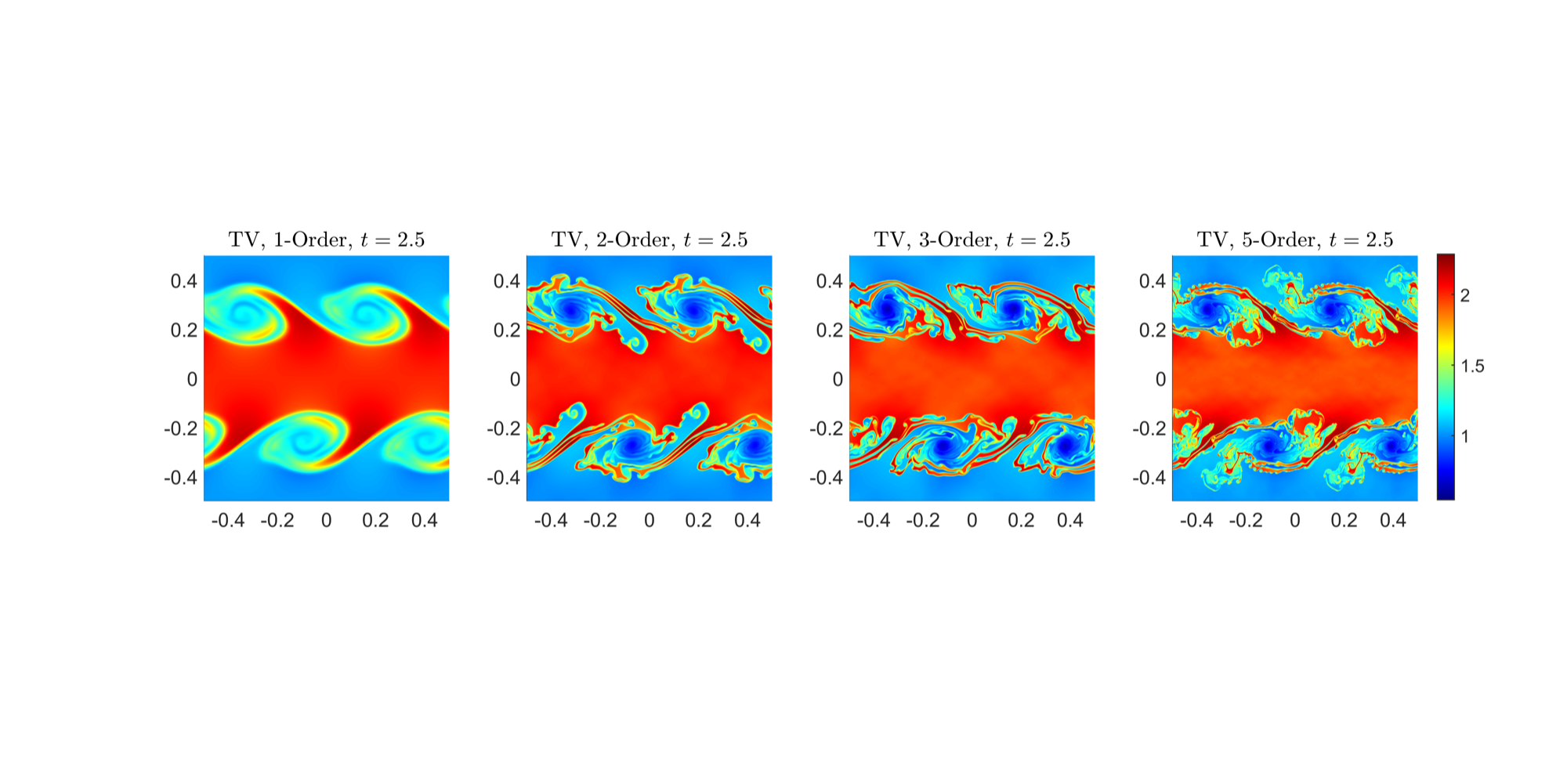}}
\vskip 15pt
\centerline{\includegraphics[trim=1.8cm 2.7cm 0.7cm 2.5cm, clip, width=14.cm]{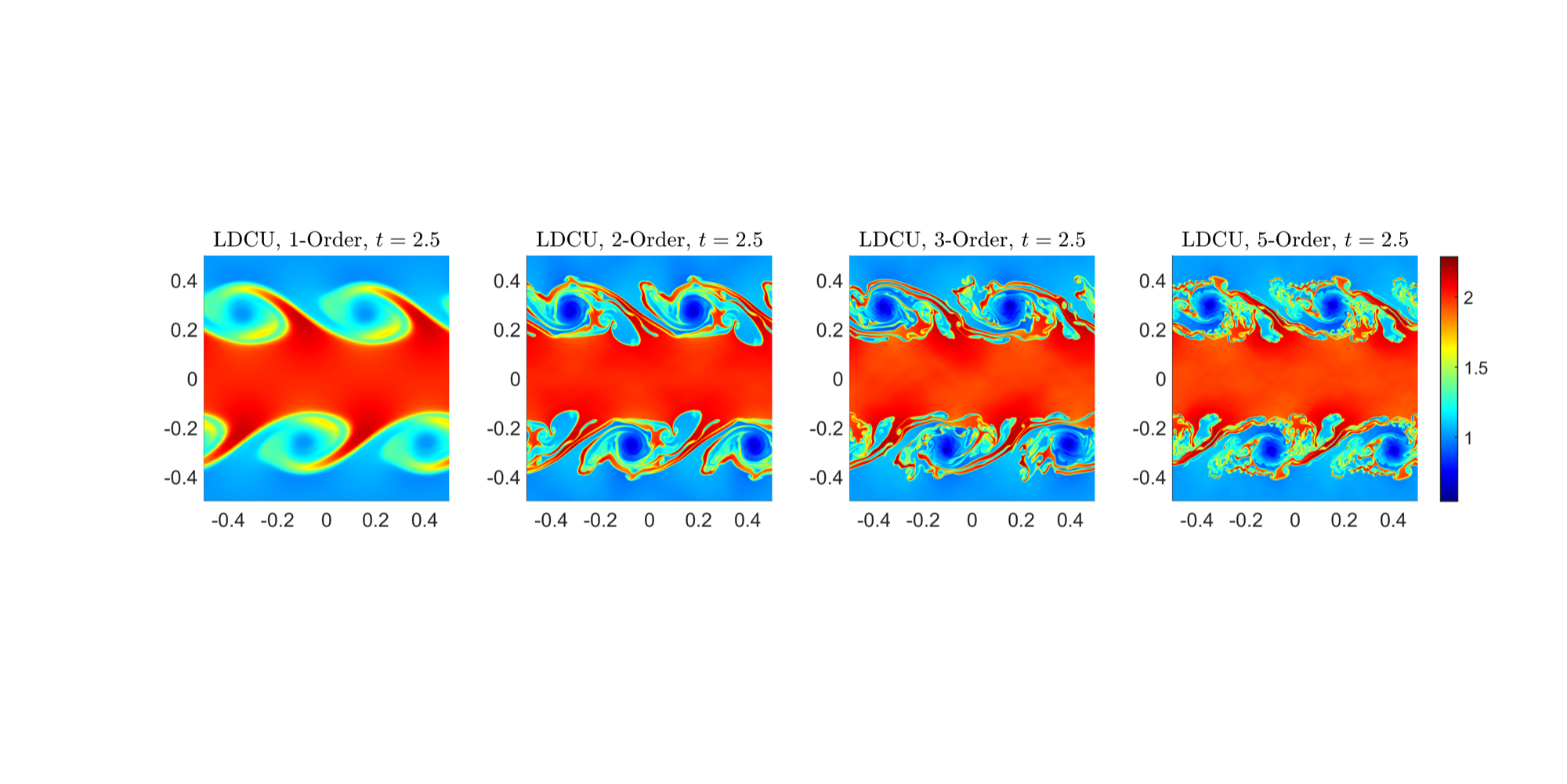}}
\vskip 15pt
\centerline{\includegraphics[trim=1.8cm 2.7cm 0.7cm 2.5cm, clip, width=14.cm]{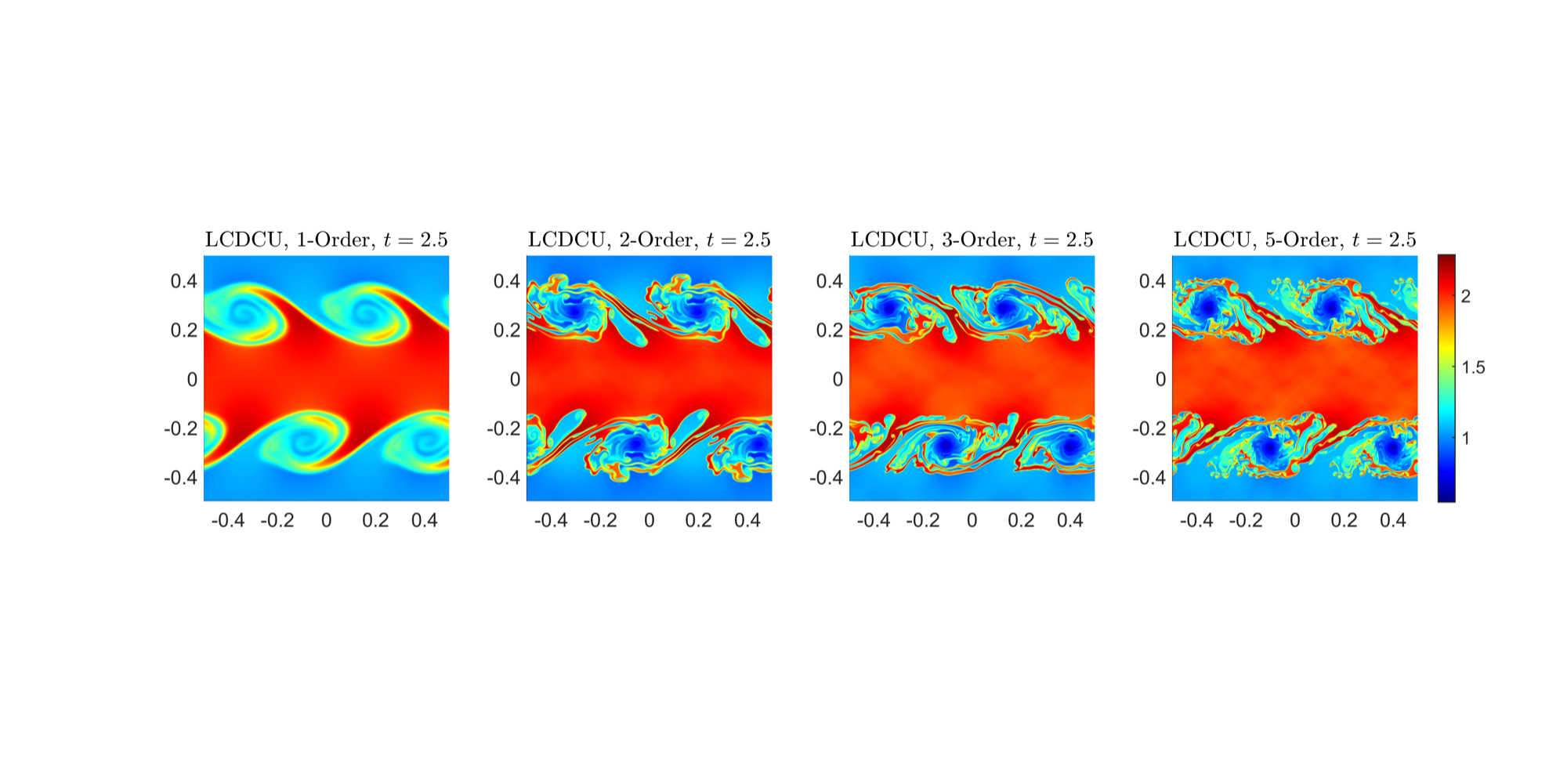}}
\caption{\sf Example 16: Time snapshots of density $\rho$ computed by the 1-Order, 2-Order, 3-Order, and 5-Order HLL (top row), HLLC (second row), TV (third row), LDCU (fourth row), and LCDCU (bottom row) schemes at $t=2.5$.\label{fig12b}}
\end{figure}

\begin{figure}[ht!]
\centerline{\includegraphics[trim=1.8cm 2.7cm 0.7cm 2.5cm, clip, width=14.cm]{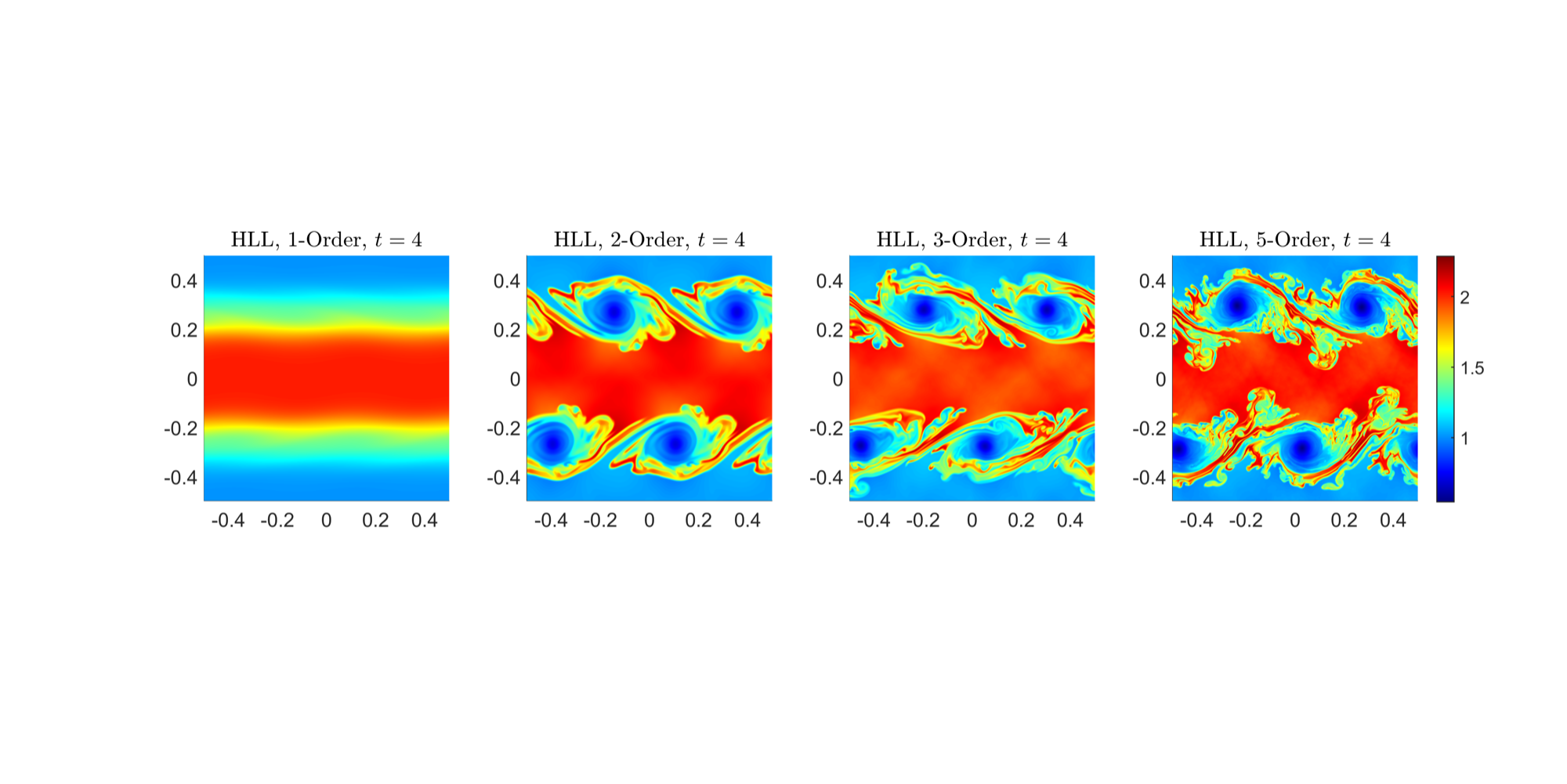}}
\vskip 15pt
\centerline{\includegraphics[trim=1.8cm 2.7cm 0.7cm 2.5cm, clip, width=14.cm]{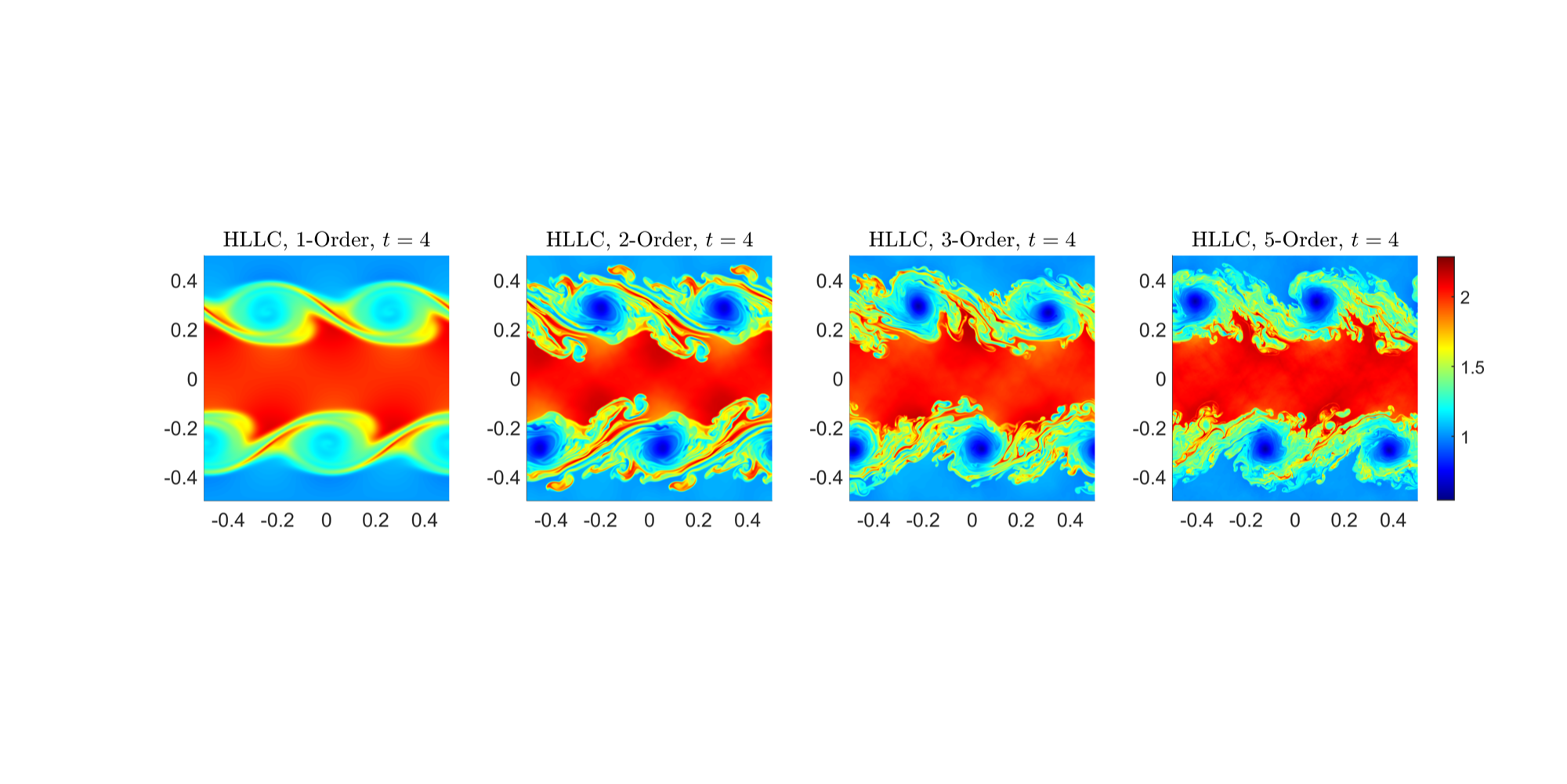}}
\vskip 15pt
\centerline{\includegraphics[trim=1.8cm 2.7cm 0.7cm 2.5cm, clip, width=14.cm]{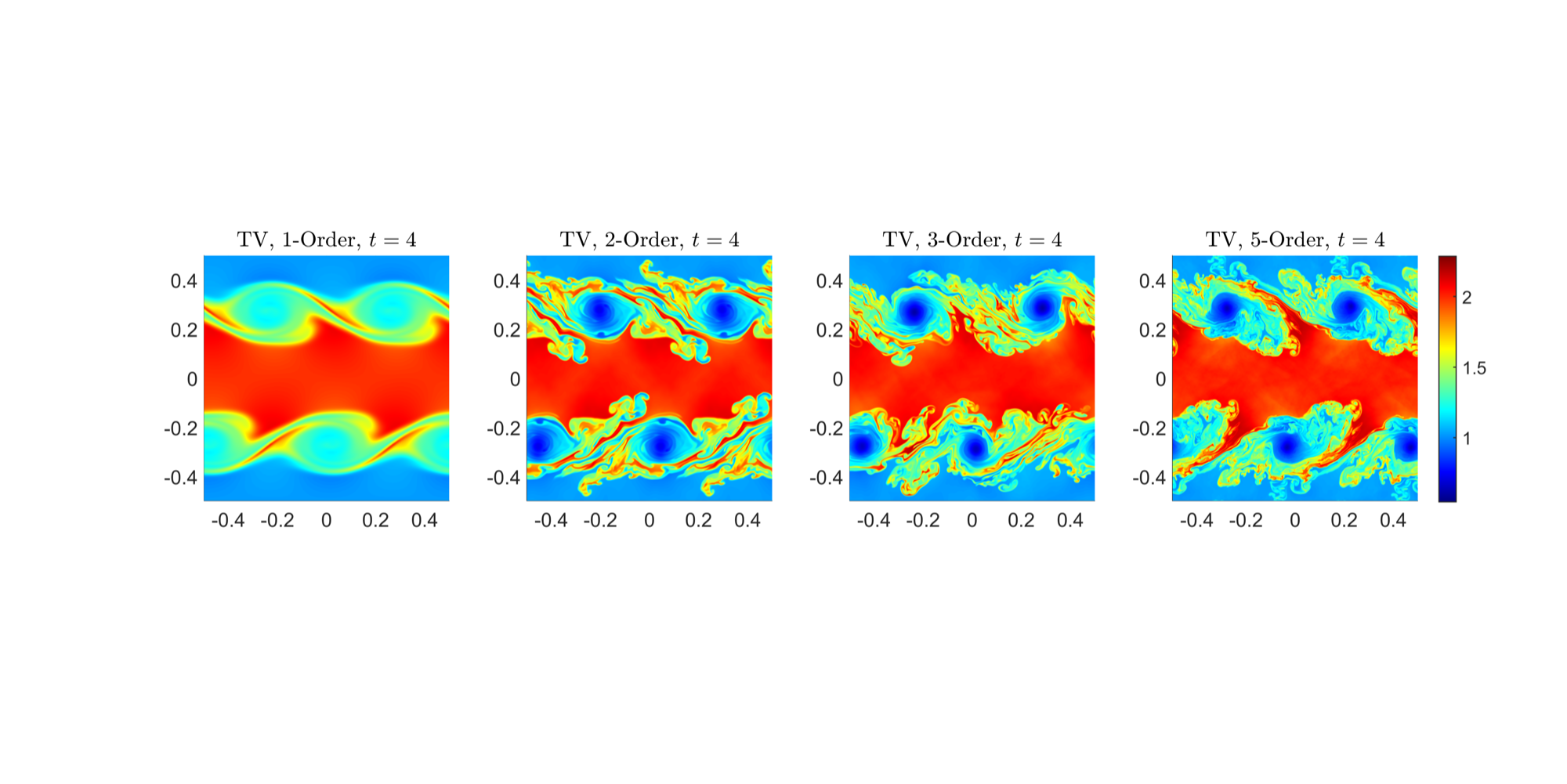}}
\vskip 15pt
\centerline{\includegraphics[trim=1.8cm 2.7cm 0.7cm 2.5cm, clip, width=14.cm]{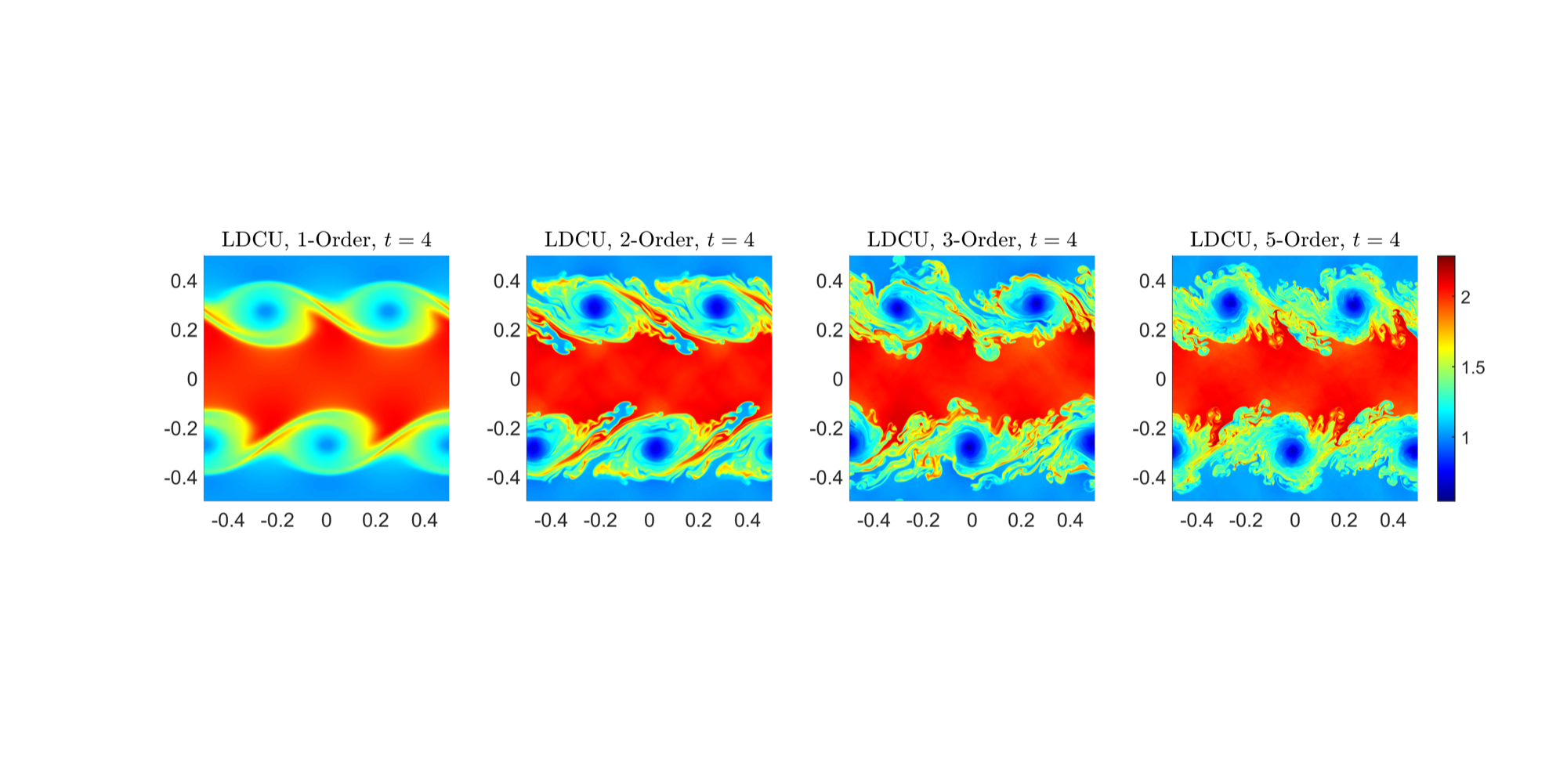}}
\vskip 15pt
\centerline{\includegraphics[trim=1.8cm 2.7cm 0.7cm 2.5cm, clip, width=14.cm]{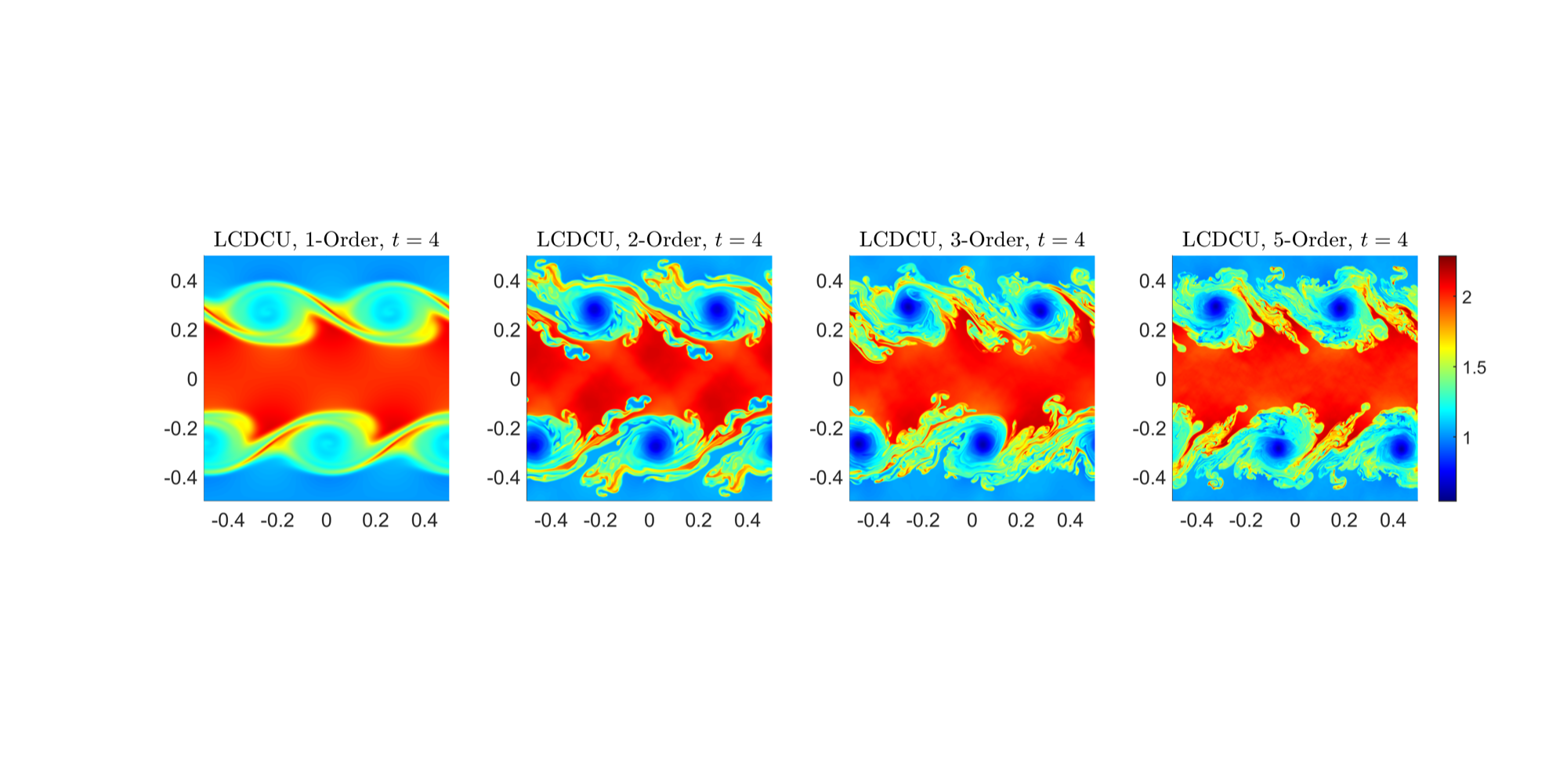}}
\caption{\sf Example 16: Time snapshots of density $\rho$ computed by the 1-Order, 2-Order, 3-Order, and 5-Order HLL (top row), HLLC (second row), TV (third row), LDCU (fourth row), and LCDCU (bottom row) schemes at $t=4$.\label{fig12c}}
\end{figure}

We compute the numerical solutions until the final time $t=4$ by the 1-Order, 2-Order, 3-Order, and 5-Order schemes on a uniform mesh of $1024\times 1024$ cells, and plot the numerical results at times $t=1$, 2.5, and 4 in Figures \ref{fig12a}--\ref{fig12c}. One can observe that at the early time $t=1$, the vortex streets generated by the high-order schemes are more pronounced; see Figure \ref{fig12a}. These structures develop over time, leading to increasingly complex vortical patterns, particularly at the later times $t=2.5$ and $4$; see Figures \ref{fig12b}--\ref{fig12c}. At the same time, it is clear that the four low-dissipation schemes produce more small-scale vortical structures than the HLL schemes.

\subsubsection*{Example 17---Quirk's Odd--Even Decoupling Problem}
In this example, we consider a variant of Quirk's odd--even decoupling problem \cite{Quirk94}, which is designed to examine the stability of numerical schemes in the presence of a strong grid-aligned shock. A right-moving Mach-6 shock is initially located near $x=0.5$. The initial conditions are given by
$$
(\rho,u,v,p)(x,y,0)
=
\begin{cases}
(7.376,\,4.861,\,0,\,41.833), & x<x_s(y),\\[1mm]
(1.4,\,0,\,0,\,1),            & x>x_s(y),
\end{cases}
$$
where the left and right states are the post-shock and pre-shock states, respectively.

To trigger the odd--even decoupling instability, we introduce a small row-dependent perturbation to the initial shock location. More precisely, for
$$
y_k=\left(k-\frac{1}{2}\right)\Delta y, \qquad \kappa_k=\min\{k,N_y+1-k\},
$$
the shock location is prescribed as
$$
x_s(y_k)=x_0+(-1)^{\kappa_k}\varepsilon, \qquad x_0=0.5+\frac{\Delta x}{2}, \qquad \varepsilon=10^{-2}\Delta x.
$$
This choice makes the perturbation symmetric with respect to the channel centerline $y=0.5$. The initial data are prescribed in the computational domain $[0,8]\times[0,1]$. At the left boundary, the post-shock state is imposed, while a free boundary condition is implemented at the right boundary. Solid-wall boundary conditions are imposed at both the bottom and top boundaries.

We compute the numerical solutions until the final time $t=1$ using the studied schemes on a uniform mesh with $\Delta x=\Delta y=1/40$ (with the CFL number $0.3$ and 0.1 for the 3-Order and 5-Order schemes, respectively, to reduce the risk of loss of positivity). The unperturbed shock propagates at the speed $s=6$ and is therefore located approximately at
$$
x_s(1)=x_0+6\approx6.5125
$$
at the final time. This test is particularly effective in revealing odd--even decoupling and other numerical shock instabilities. A shock-stable scheme should maintain a nearly planar shock front, whereas a shock-unstable scheme may produce spurious transverse oscillations.

\begin{figure}[ht!]
\centering
\includegraphics[trim=0cm 0cm 0cm 0cm,clip,width=\linewidth]{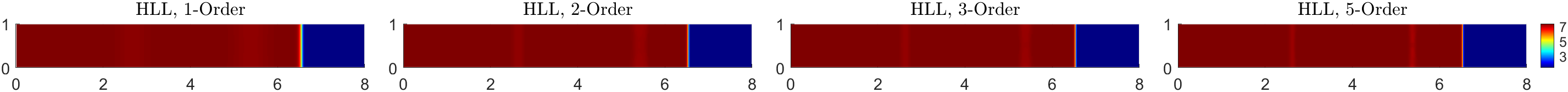}
\vspace{5pt}
\includegraphics[trim=0cm 0cm 0cm 0cm,clip,width=\linewidth]{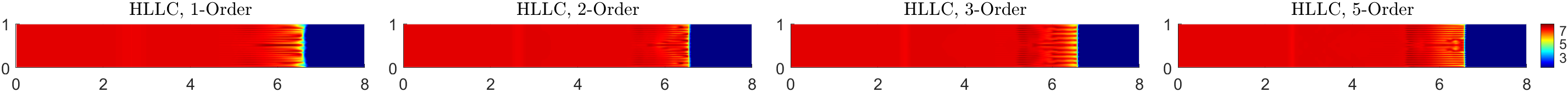}
\vspace{5pt}
\leftline{\includegraphics[trim=0cm 0cm 0cm 0cm,clip,width=0.5\linewidth]{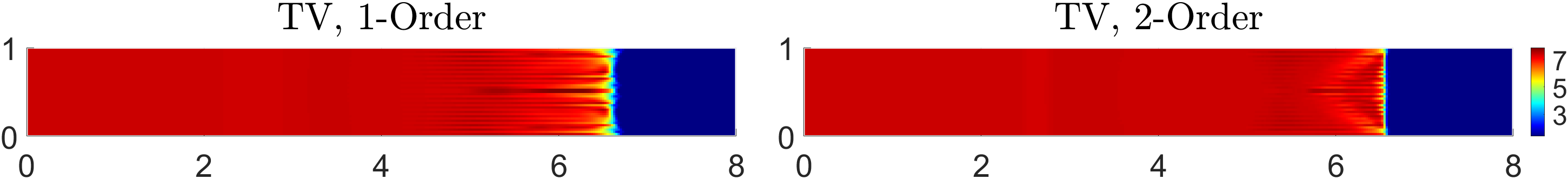}}
\vspace{5pt}
\includegraphics[trim=0cm 0cm 0cm 0cm,clip,width=\linewidth]{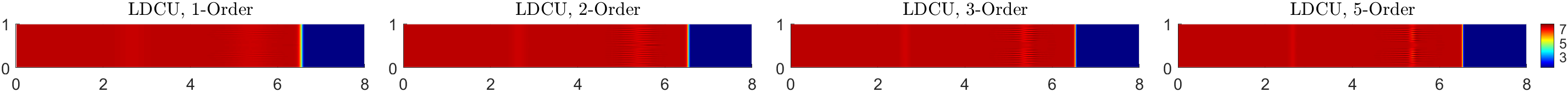}
\vspace{5pt}
\includegraphics[trim=0cm 0cm 0cm 0cm,clip,width=\linewidth]{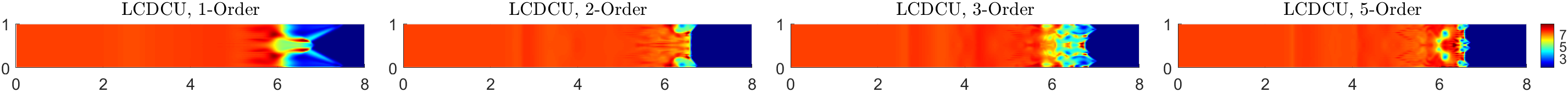}
\caption{\sf Example 17: Density $\rho$ computed by the 1-Order, 2-Order, 3-Order, and 5-Order HLL (top row), HLLC (second row), LDCU (fourth row), and LCDCU (bottom row) schemes, as well as the 1-Order and 2-Order TV (third row) schemes.}
\label{fig16aaa}
\end{figure}

The obtained numerical results are plotted in Figure~\ref{fig16aaa}. As one can see, the HLL schemes maintain essentially planar shock fronts. In contrast, the HLLC, TV, and LCDCU schemes exhibit much more pronounced numerical shock instabilities, while the LDCU schemes produce only weak instabilities compared with the HLL schemes. The 3- and 5-Order TV computations encounter a loss of positivity, and their results are therefore omitted from Figure~\ref{fig16aaa}.

\begin{figure}[ht!]
\centering
\includegraphics[trim=1.0cm 0.3cm 1.3cm 0.8cm,clip,width=4cm]{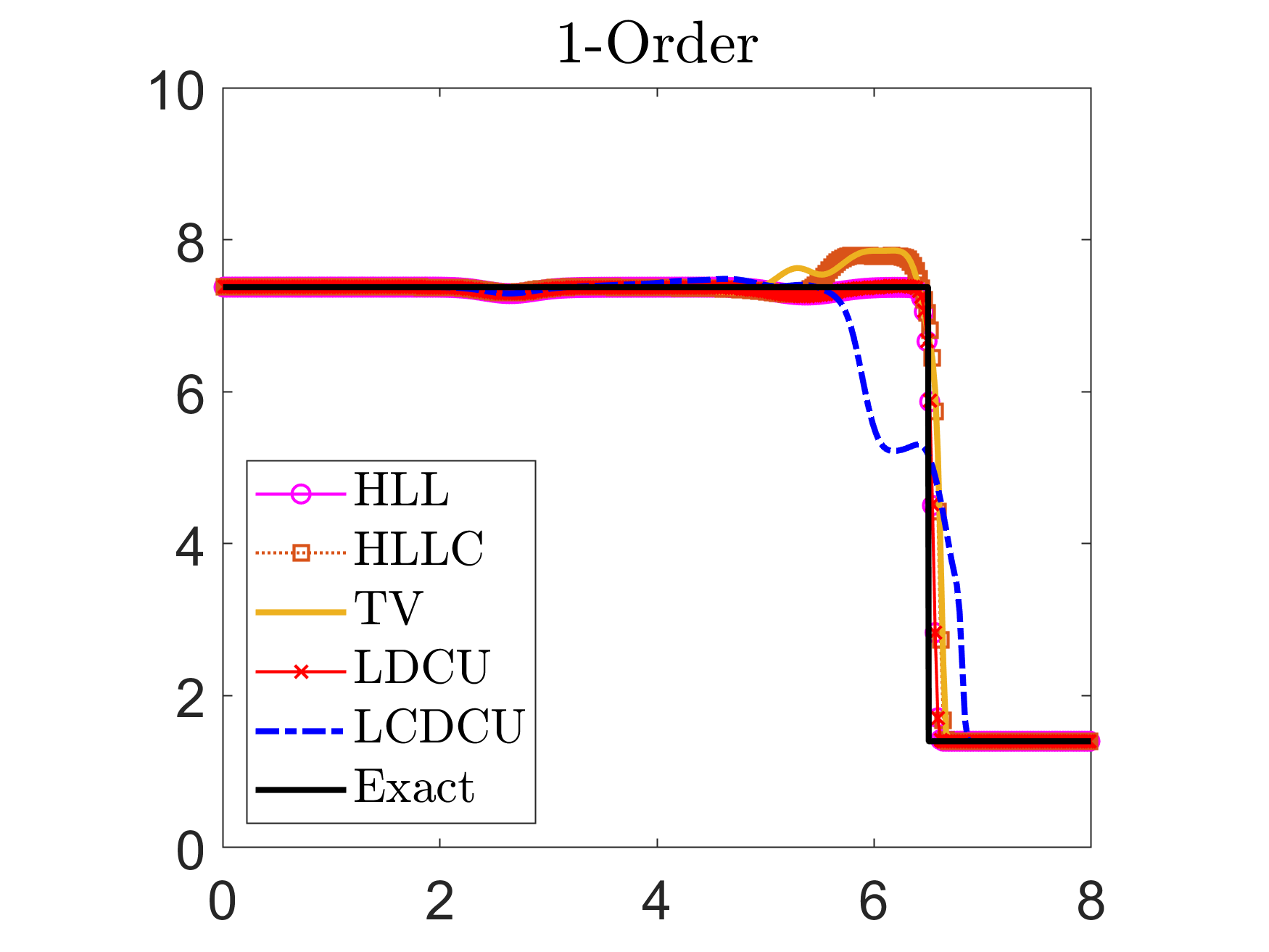}
\includegraphics[trim=1.0cm 0.3cm 0.9cm 0.8cm,clip,width=4cm]{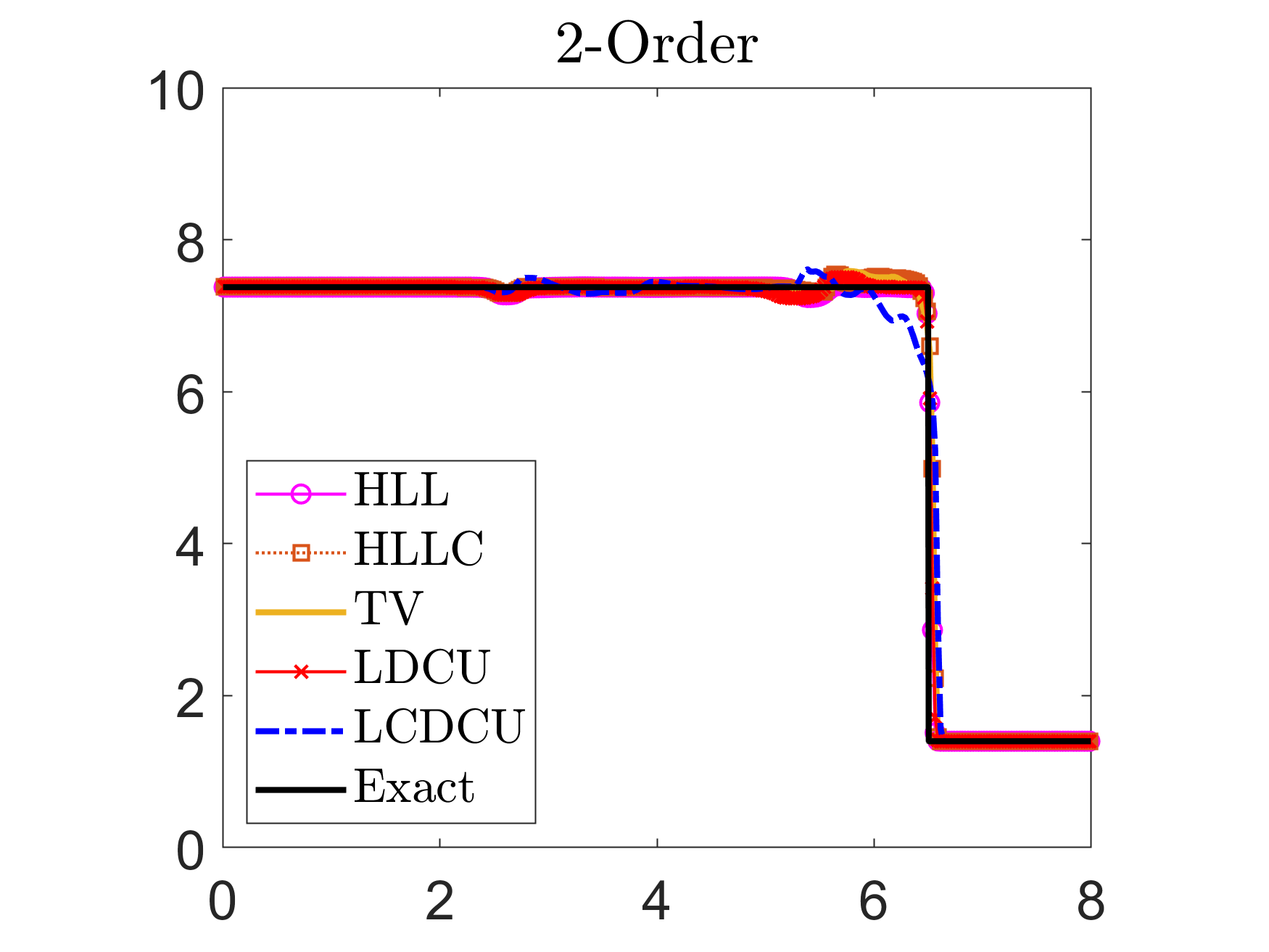}
\includegraphics[trim=1.0cm 0.3cm 0.9cm 0.8cm,clip,width=4cm]{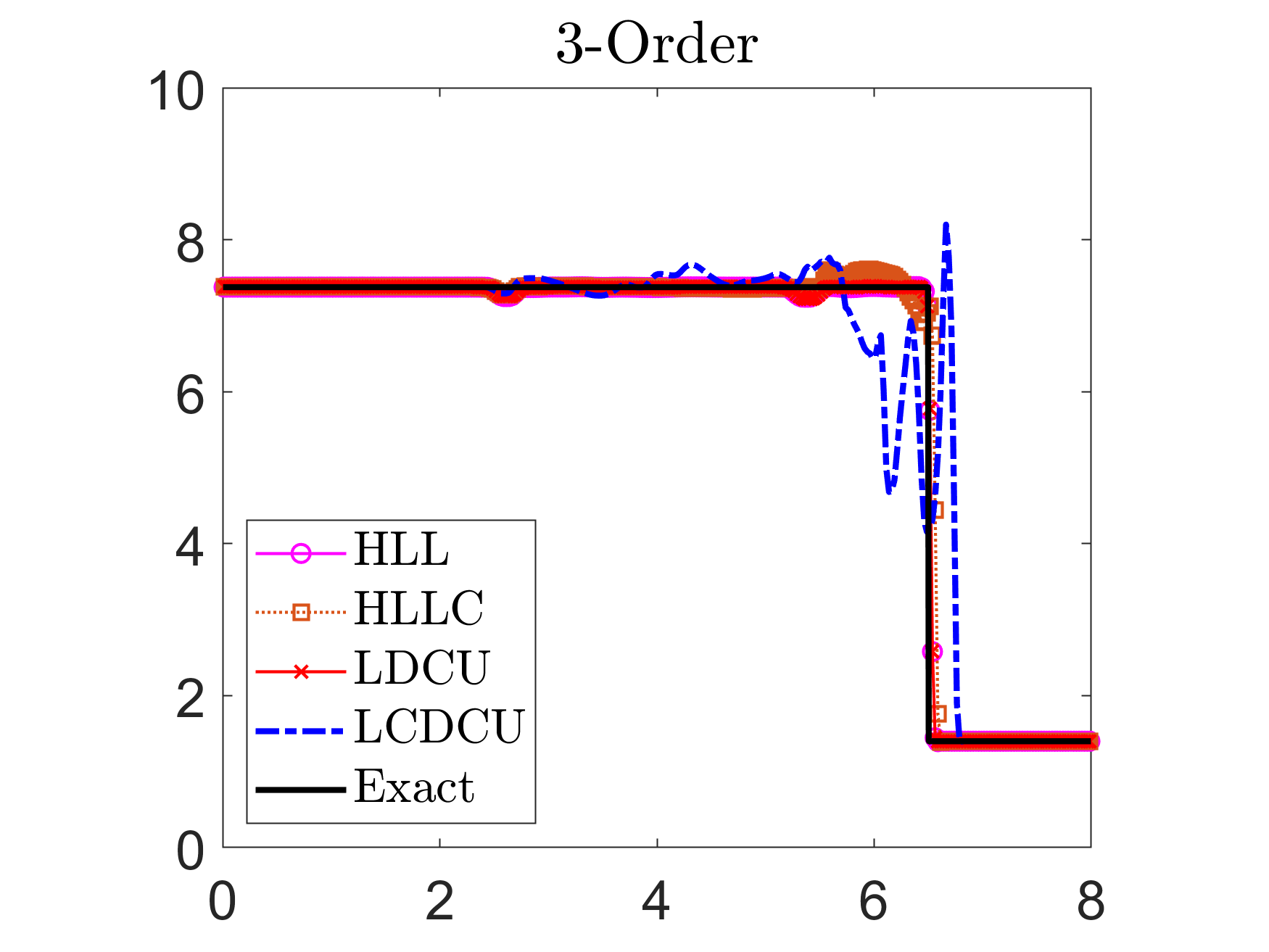}
\includegraphics[trim=1.0cm 0.3cm 0.9cm 0.8cm,clip,width=4cm]{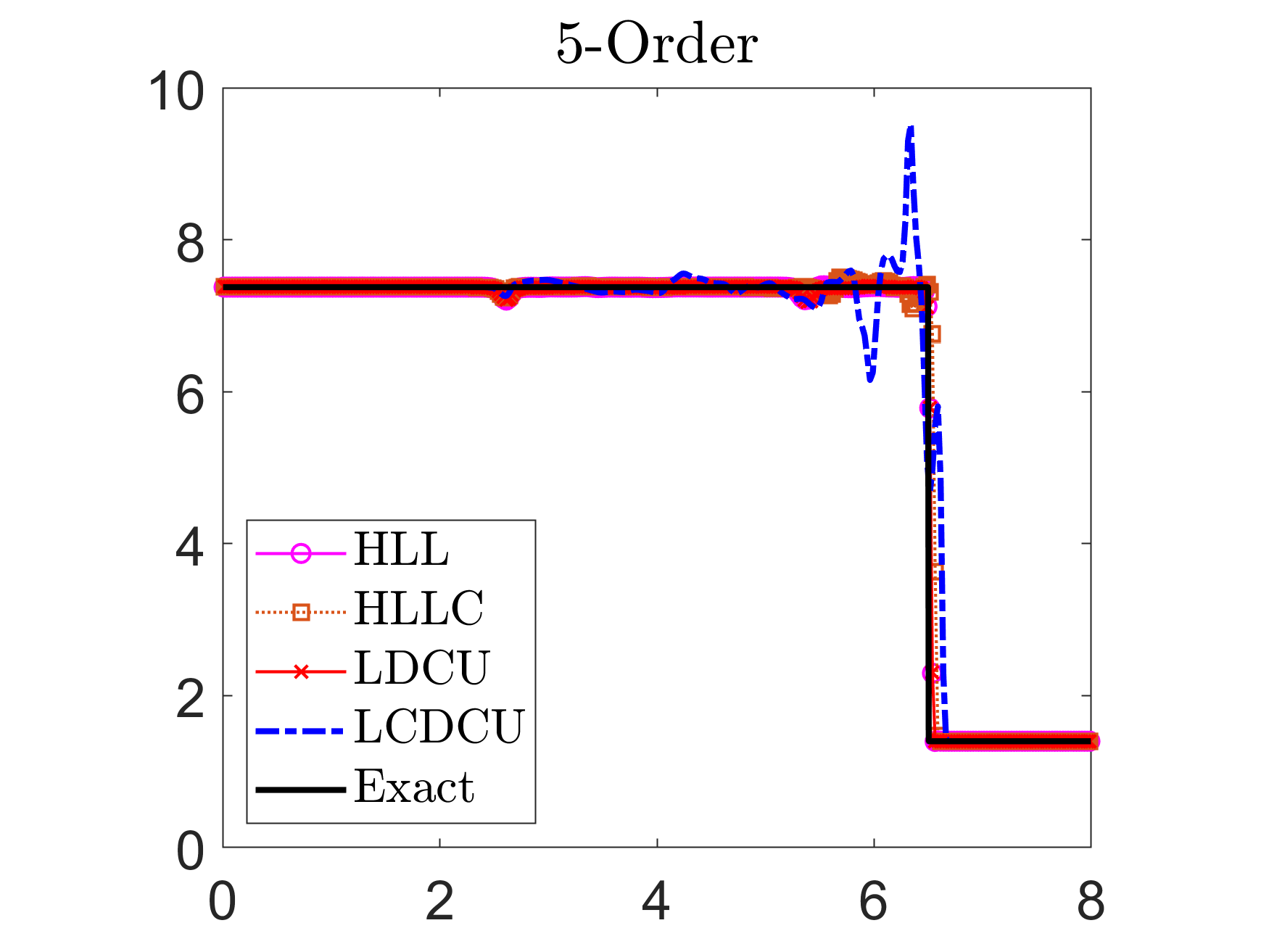}
\caption{\sf Example 17: Density $\rho$ along the centerline $y=0.5$ computed by the studied schemes.}
\label{fig16bbb}
\end{figure}

For a quantitative comparison, we also plot the density along the centerline $y=0.5$.  The obtained centerline densities are presented in Figure~\ref{fig16bbb}. The results confirm that the HLLC, TV, and LCDCU schemes develop pronounced spurious oscillations, whereas the HLL schemes remain nonoscillatory and the LDCU schemes exhibit only weak deviations from the HLL results.

\subsubsection*{Example 18---Shock Diffraction over a $90^\circ$ Corner}
In this example, we consider a Mach-2 version of the shock diffraction problem over a $90^\circ$ corner; see, e.g., \cite{Quirk94}. The flow region is the L-shaped domain
$$
\Omega=\bigl([0,1]\times[0,1]\bigr)\cup\bigl([1,2]\times[0,0.5]\bigr).
$$
The initial conditions, prescribed in the fluid region $\Omega$, are
given by
$$
(\rho,u,v,p)(x,y,0)=\begin{cases}
(1,\,0,\,0,\,1), & x<1.375,\\[1.ex]
\left(\frac{8}{3},-\frac{\sqrt{35}}{4},0,\frac{9}{2}\right),
& x>1.375,
\end{cases}
$$
which correspond to a left-moving Mach-2 shock initially located at $x=1.375$. The left and right states are the pre-shock and post-shock states, respectively. Solid-wall boundary conditions are imposed along the channel walls and the boundary of the solid block. At the open left and right boundaries, the corresponding pre-shock and post-shock states are maintained, respectively.

\begin{figure}[ht!]
\centerline{\includegraphics[trim=5.6cm 10.3cm 4.7cm 9.4cm, clip, width=16cm]{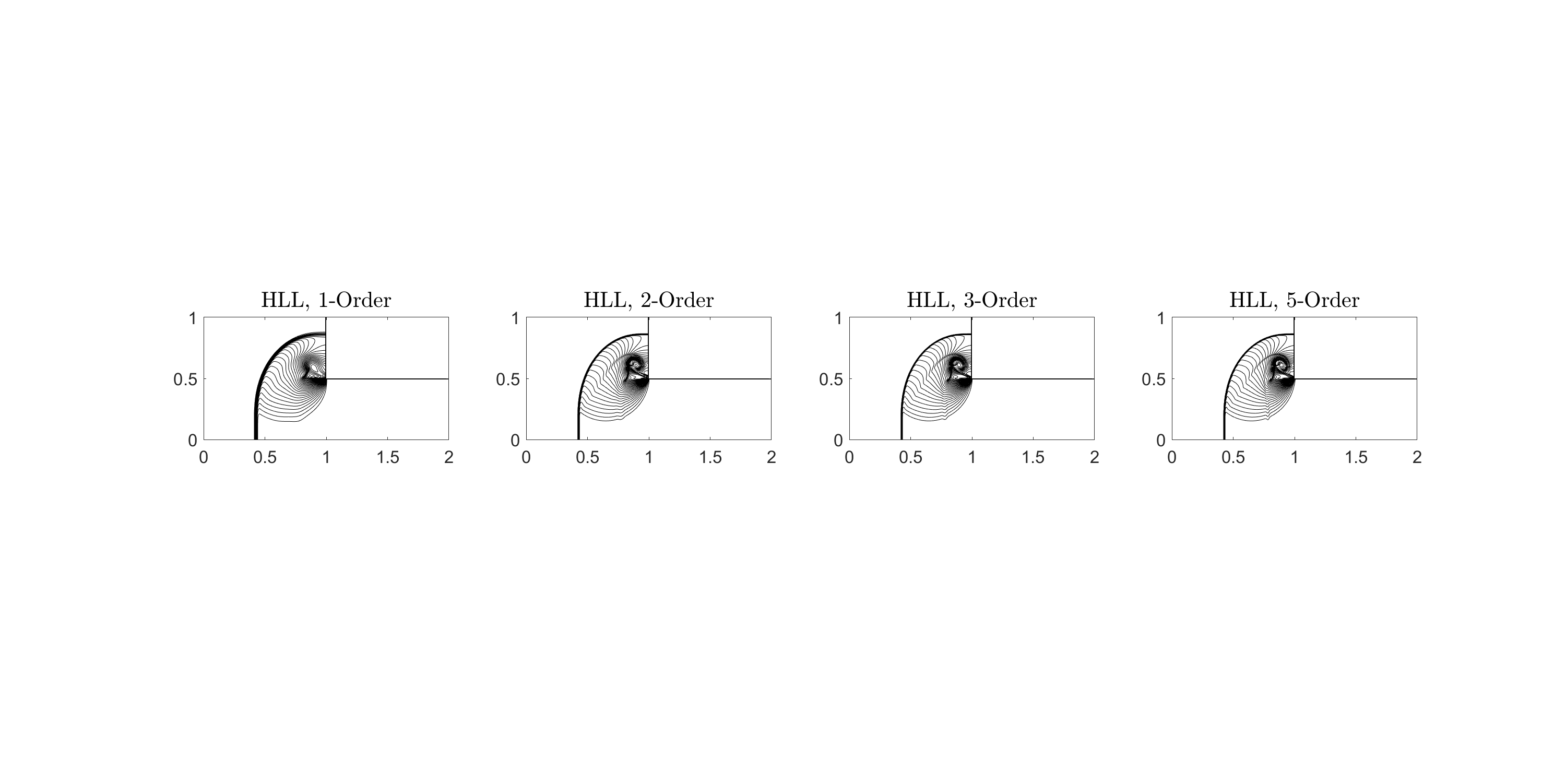}}
\vskip 15pt 
\centerline{\includegraphics[trim=5.6cm 10.3cm 4.7cm 9.4cm, clip, width=16cm]{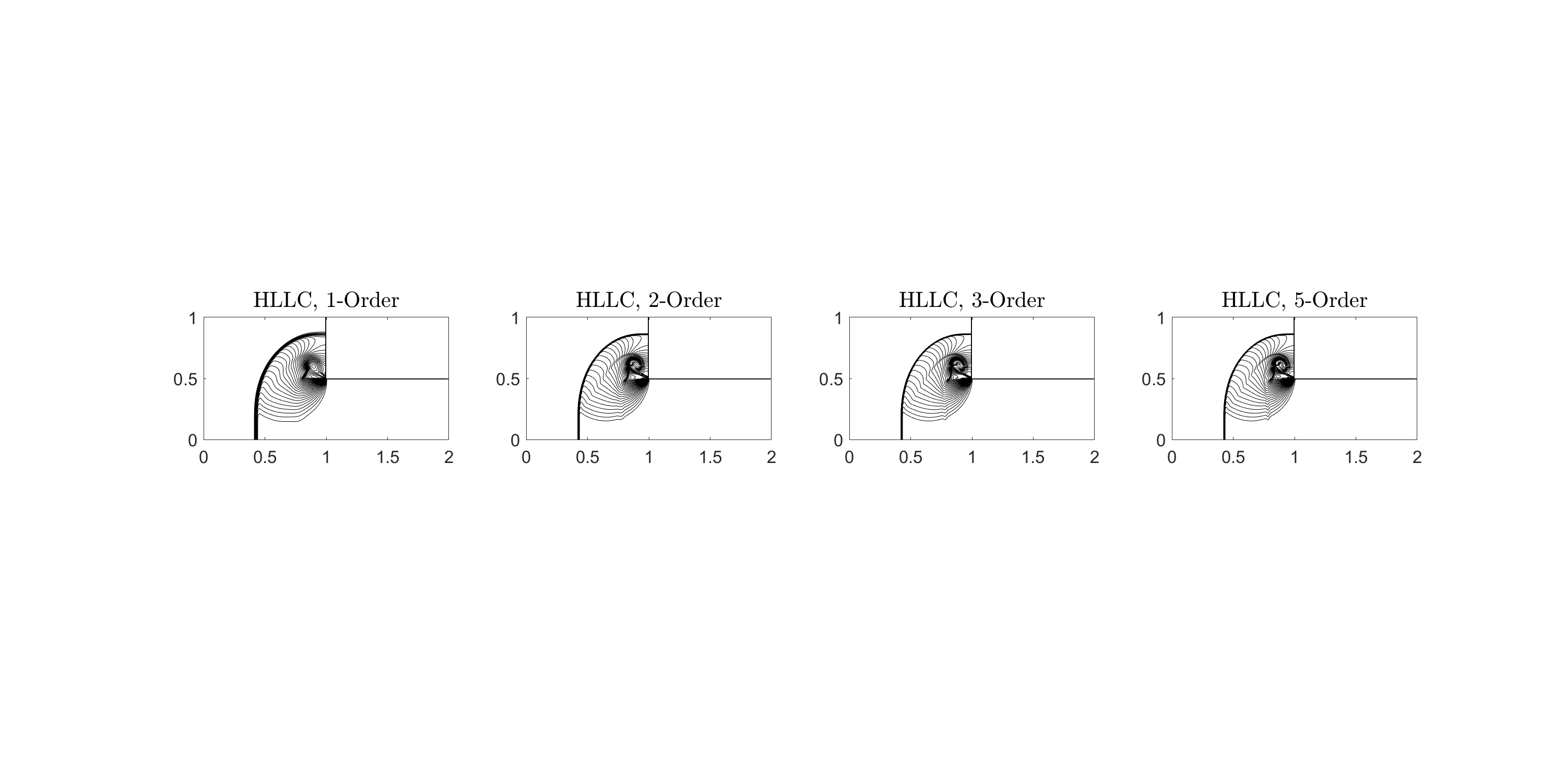}}
\vskip 15pt 
\centerline{\includegraphics[trim=5.6cm 10.3cm 4.7cm 9.4cm, clip, width=16cm]{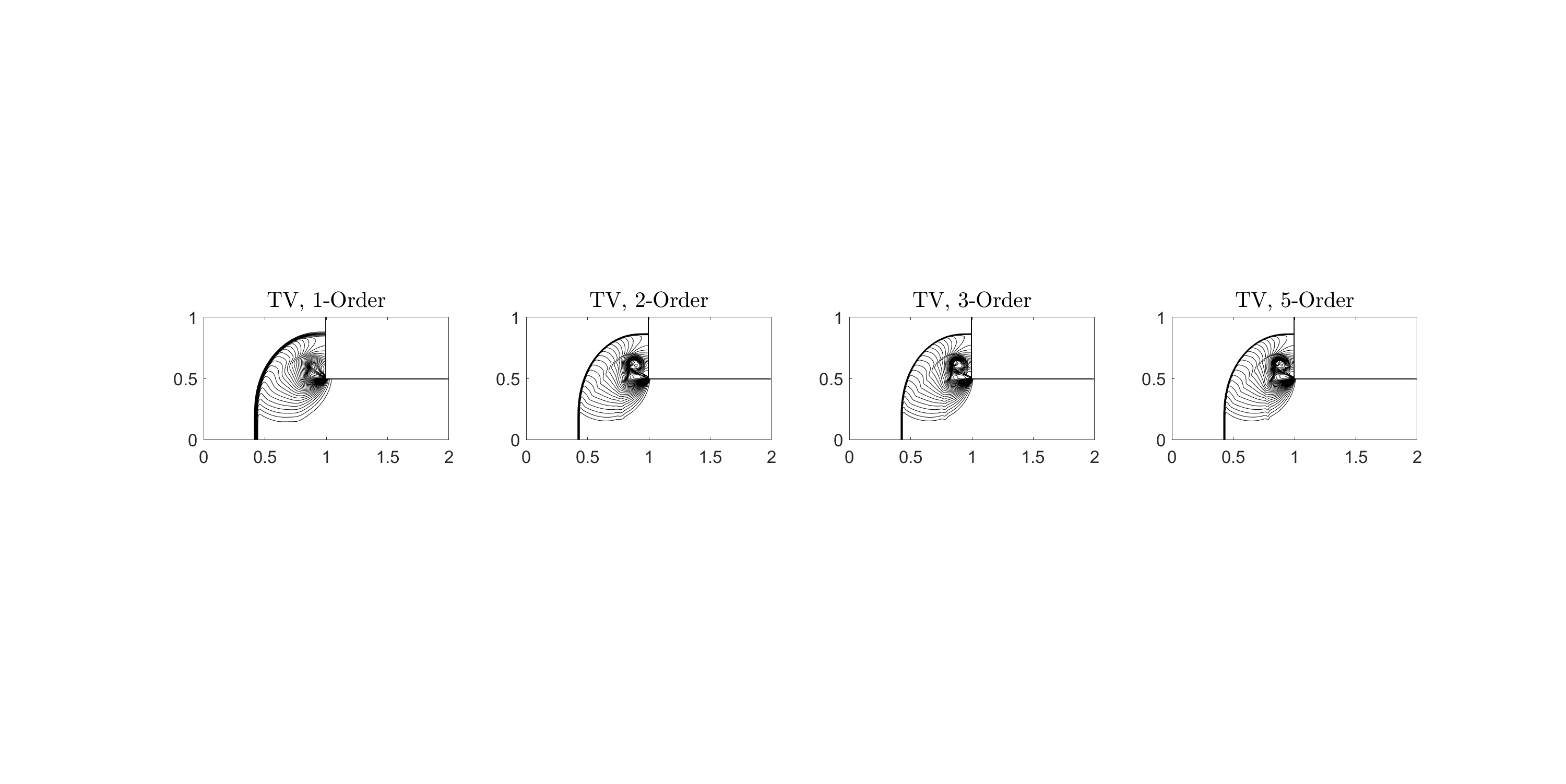}}
\vskip 15pt 
\centerline{\includegraphics[trim=5.6cm 10.3cm 4.7cm 9.4cm, clip, width=16cm]{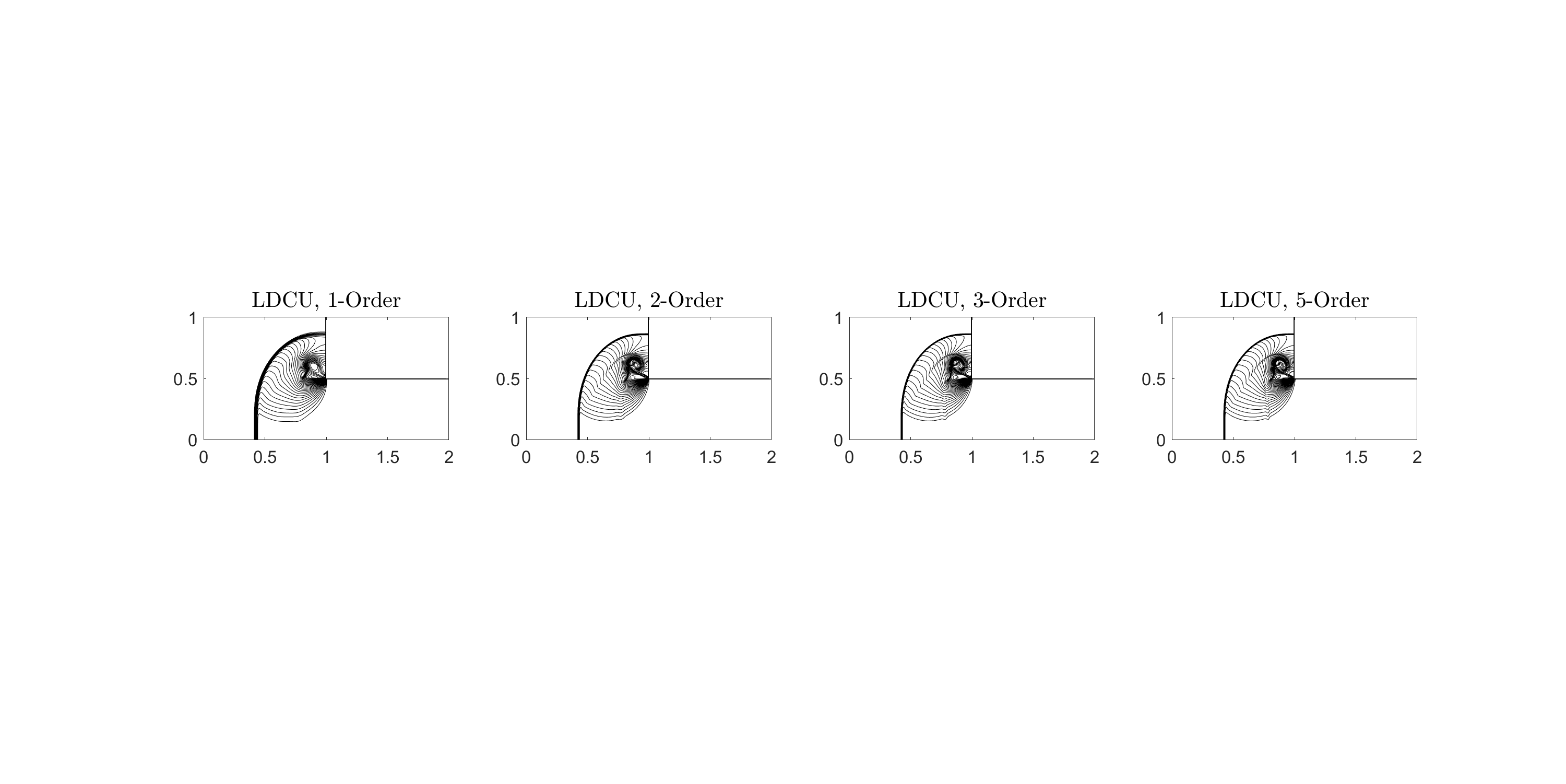}}
\vskip 15pt 
\centerline{\includegraphics[trim=5.6cm 10.3cm 4.7cm 9.4cm, clip, width=16cm]{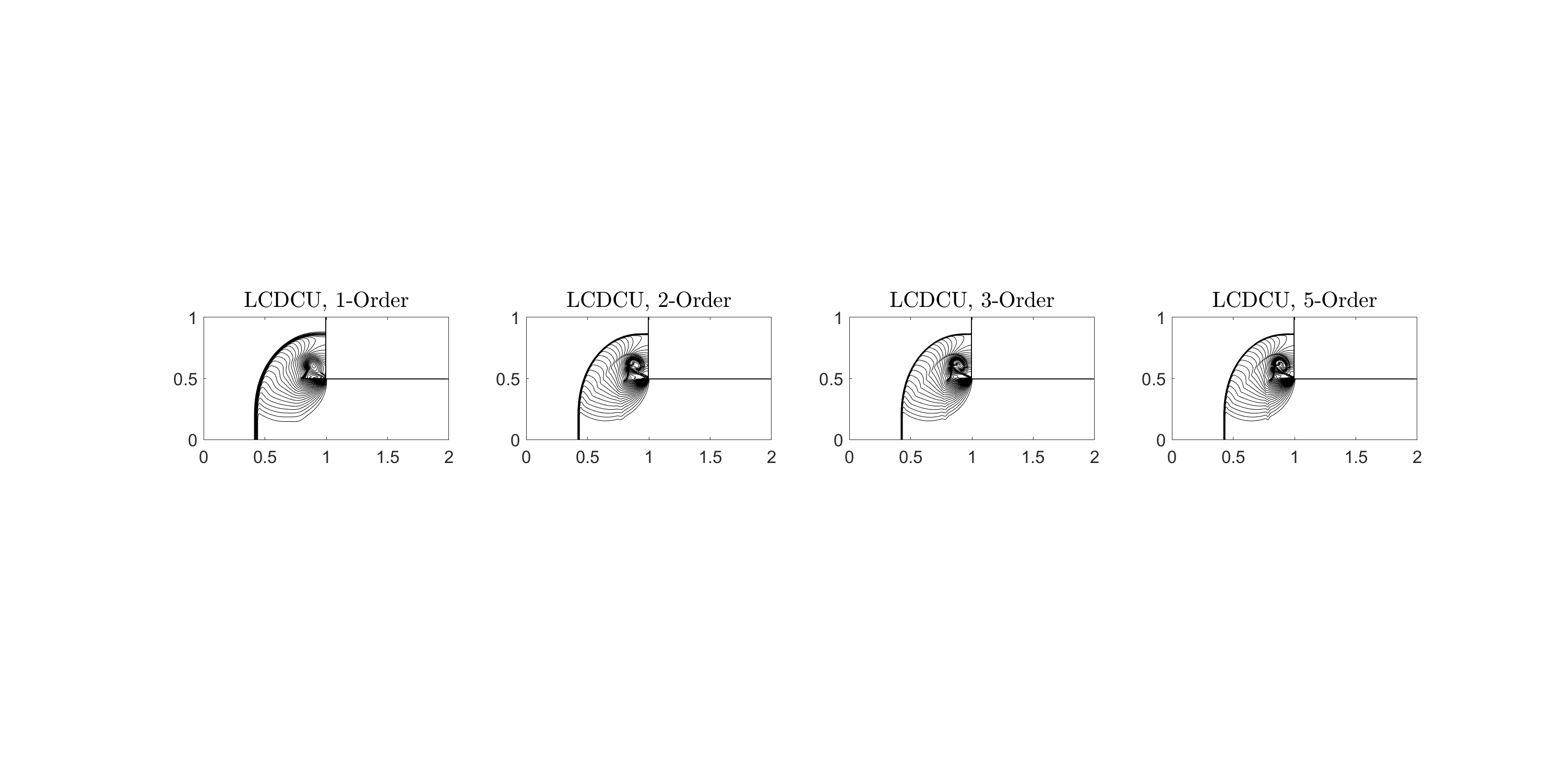}}
\caption{\sf Example 18: Density $\rho$ computed by the 1-Order, 2-Order, 3-Order, and 5-Order HLL (top row), HLLC (second row), TV (third row), LDCU (fourth row), and LCDCU (bottom row) schemes.\label{fig17aaa}}
\end{figure}

We compute the numerical solutions until the final time $t=0.4$ using the studied schemes on a uniform Cartesian mesh with $\dx=\dy=1/256$. The obtained density distributions are presented in Figure~\ref{fig17aaa}. One can see that all of the studied schemes capture the main structure of the diffracted shock wave. Compared with the HLL schemes, the four low-dissipation schemes provide sharper resolution of the complex flow structures generated around the corner. At the same time, the differences among the four low-dissipation schemes are relatively limited.

\subsubsection*{Example 19---Rayleigh-Taylor (RT) Instability}
In the last example taken from \cite{Shi03}, we investigate the RT instability, which is a physical phenomenon
occurring when a layer of heavier fluid is placed on top of a layer of lighter fluid. The model is governed by the 2-D Euler equations
\eref{3.1}--\eref{3.2} with added gravitational source terms and the modified system reads as
$$
\begin{aligned}
&\rho_t+(\rho u)_x+(\rho v)_y=0,\\
&(\rho u)_t+(\rho u^2 +p)_x+(\rho uv)_y=0,\\
&(\rho v)_t+(\rho uv)_x+(\rho v^2+p)_y=\rho,\\
&E_t+[u(E+p)]_x+[v(E+p)]_y=\rho v.
\end{aligned}
$$
We consider the following initial conditions:
\begin{equation*}
(\rho(x,y,0),u(x,y,0),v(x,y,0),p(x,y,0))=\begin{cases}
(2,0,-0.025\,c\cos(8\pi x),2y+1),&y<0.5,\\
(1,0,-0.025\,c\cos(8\pi x),y+1.5),&\mbox{otherwise},
\end{cases}
\end{equation*}
where $c:=\sqrt{\gamma p/\rho}$ is the speed of sound, prescribed in the computational domain $[0, 0.25]\times [0,1]$ with the solid wall boundary conditions imposed at $x=0$ and $x=0.25$, and the following Dirichlet boundary conditions imposed at the top and bottom boundaries:
$$
(\rho,u,v,p)|_{y=1}=(1,0,0,2.5),\qquad(\rho,u,v,p)|_{y=0}=(2,0,0,1).
$$

\begin{figure}[ht!]
\centerline{\includegraphics[trim=1.7cm 1.3cm 0.7cm 1.2cm, clip, width=14.cm]{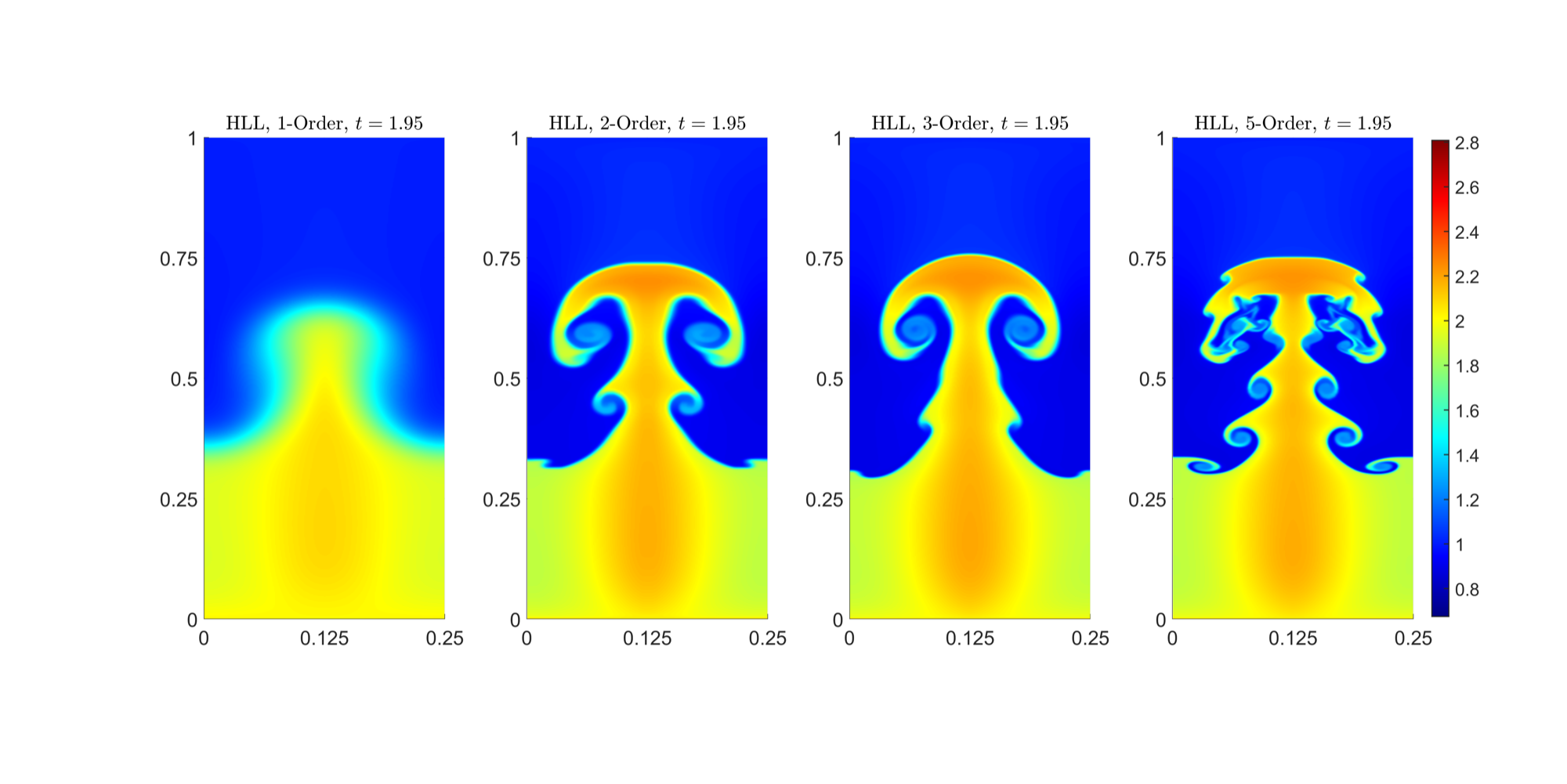}}
\vskip 15pt 
\centerline{\includegraphics[trim=1.7cm 1.3cm 0.7cm 1.2cm, clip, width=14.cm]{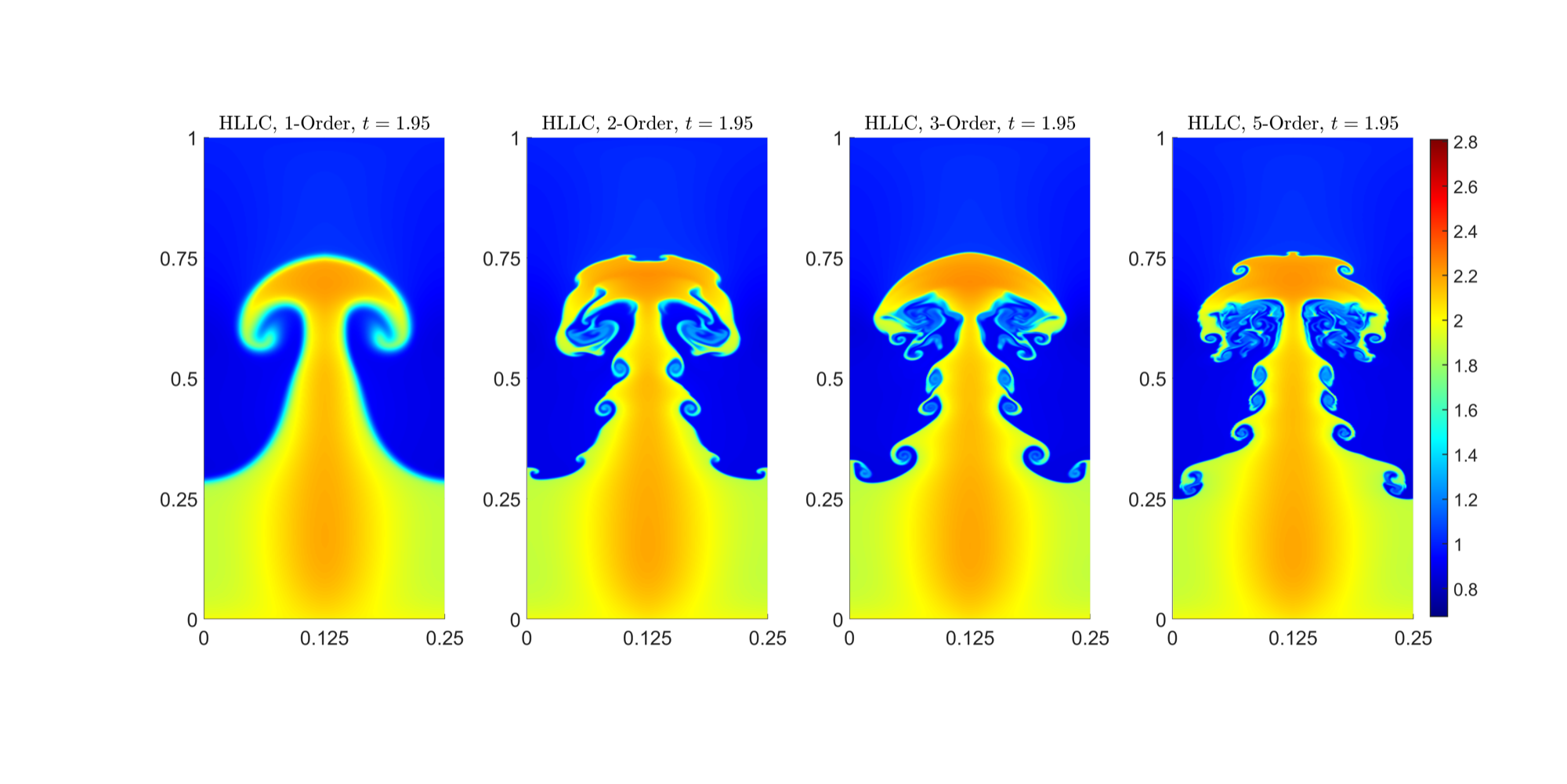}}
\vskip 15pt 
\centerline{\includegraphics[trim=1.7cm 1.3cm 0.7cm 1.2cm, clip, width=14.cm]{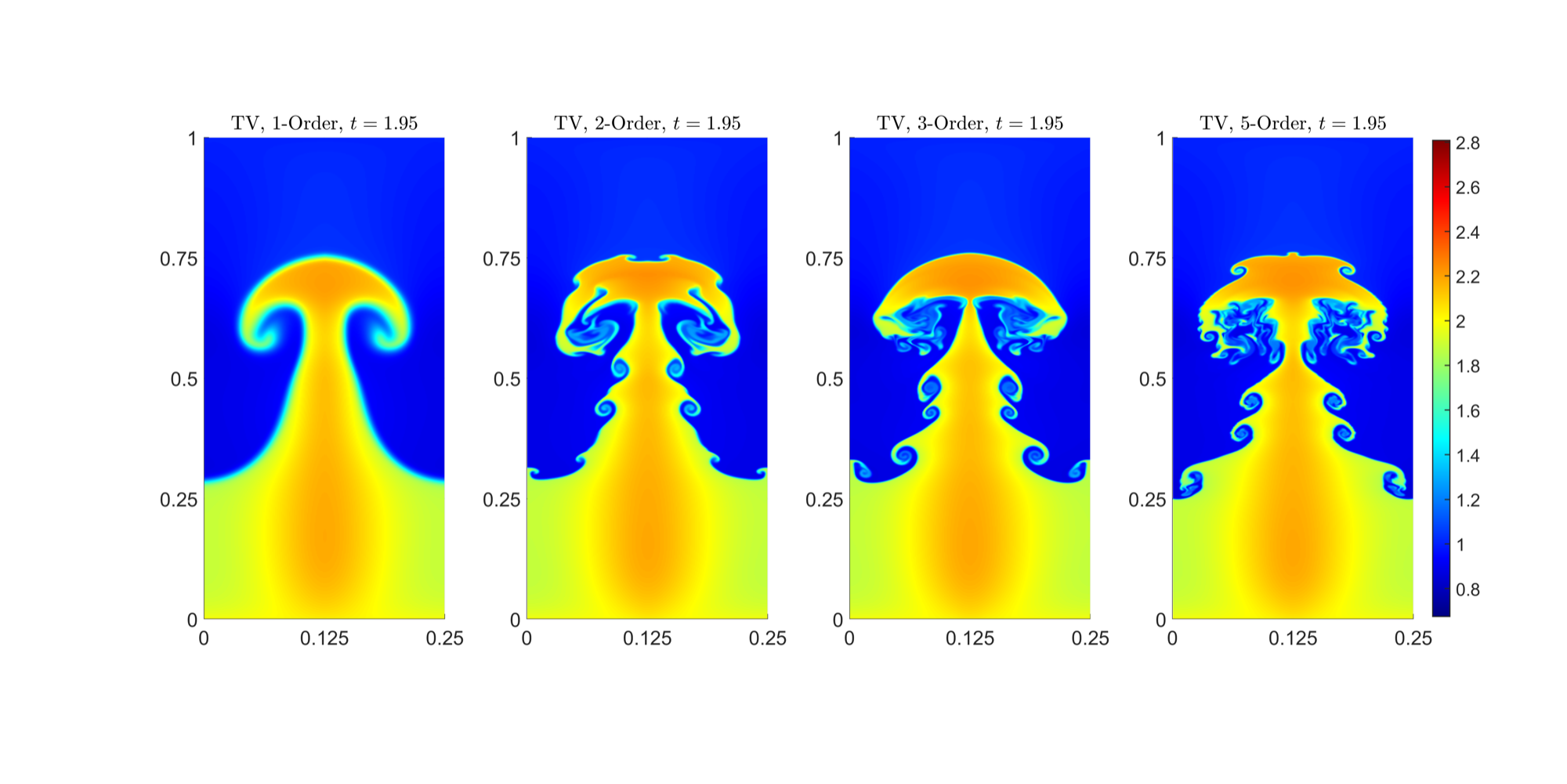}}
\caption{\sf Example 19: Density $\rho$ computed by the 1-Order (first column), 2-Order (second column), 3-Order (third column), and 5-Order (fourth column) HLL (top row), HLLC (second row), and TV (bottom row) schemes at $t=1.95$.\label{fig13a}}
\end{figure}

\begin{figure}[ht!]
\centerline{\includegraphics[trim=1.7cm 1.3cm 0.7cm 1.2cm, clip, width=14.cm]{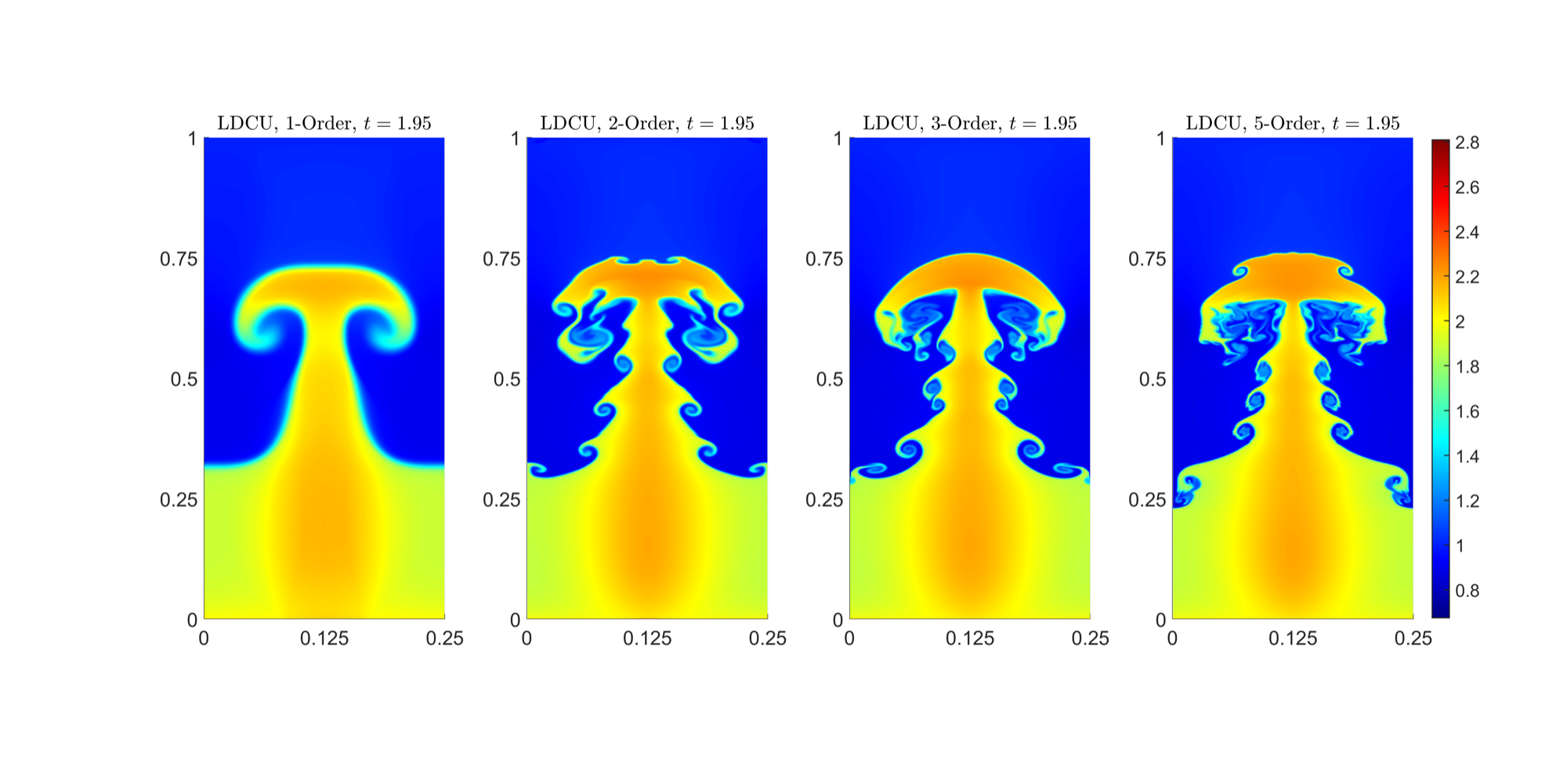}}
\vskip 15pt 
\centerline{\includegraphics[trim=1.7cm 1.3cm 0.7cm 1.2cm, clip, width=14.cm]{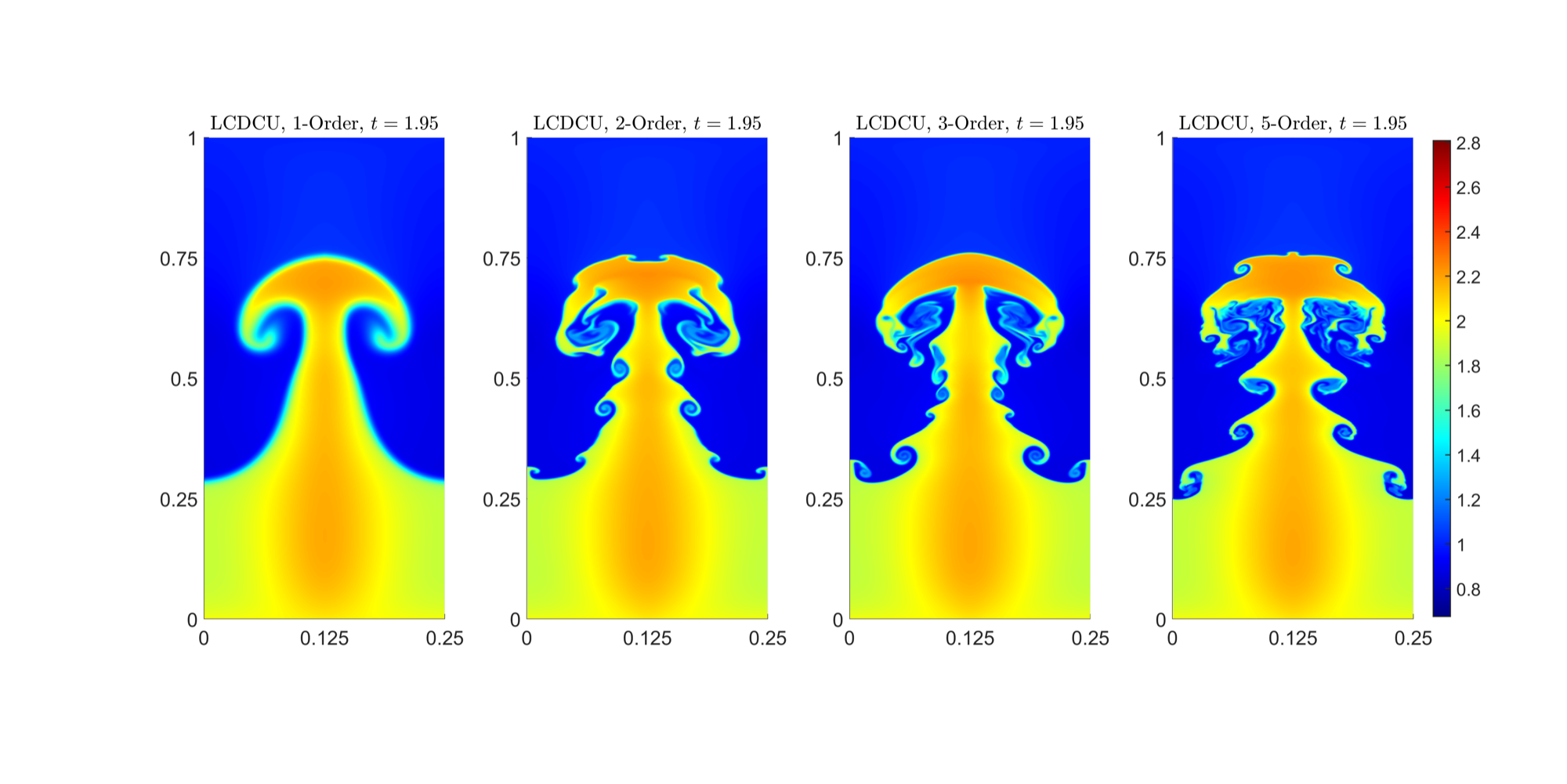}}
\caption{\sf Example 19: Density $\rho$ computed by the 1-Order (first column), 2-Order (second column), 3-Order (third column), and 5-Order (fourth column) LDCU (top row) and LCDCU (bottom row)  schemes at $t=1.95$.\label{fig13b}}
\end{figure}

\begin{figure}[ht!]
\centerline{\includegraphics[trim=1.7cm 1.3cm 0.7cm 1.2cm, clip, width=14.cm]{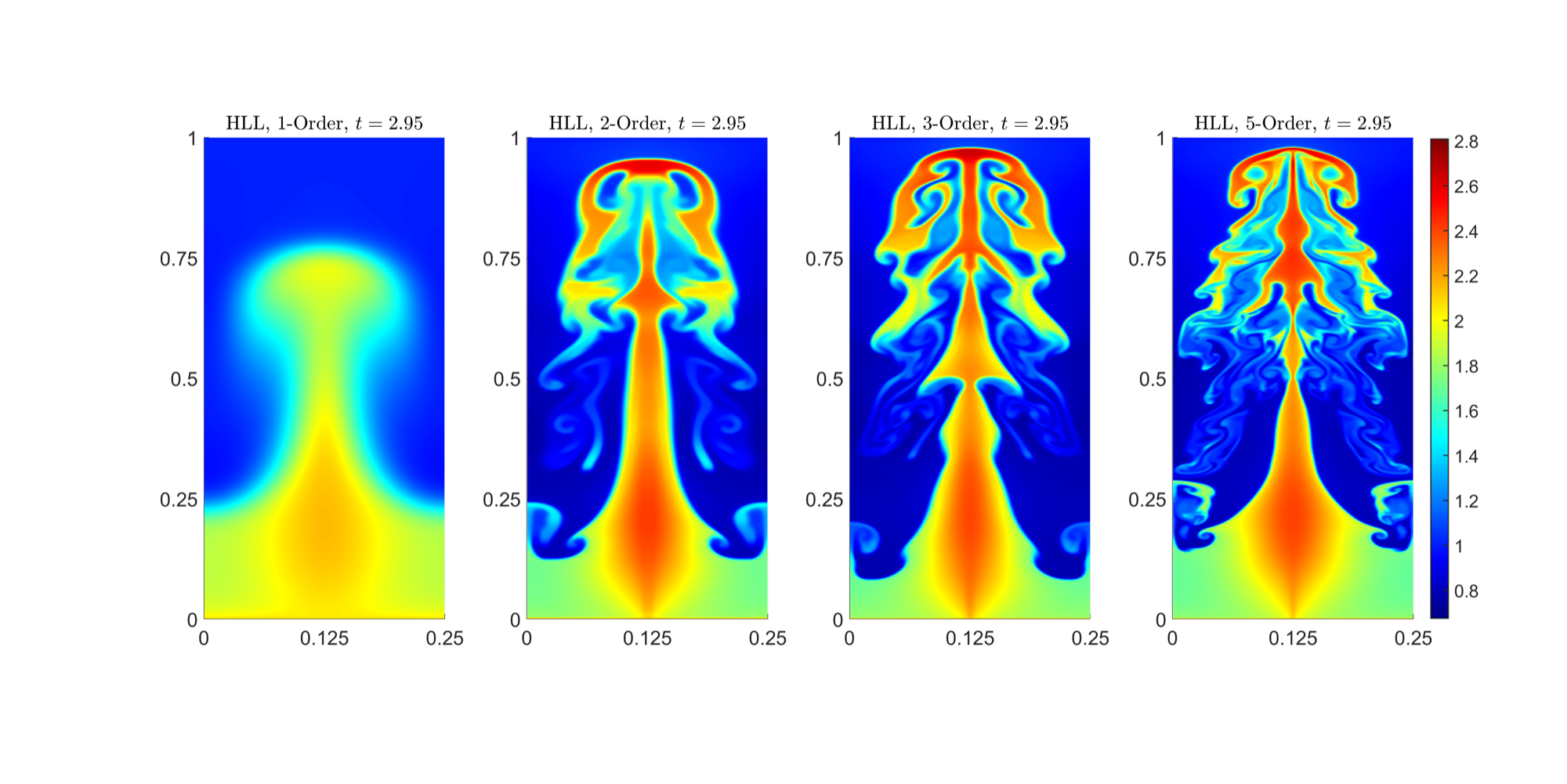}}
\vskip 15pt 
\centerline{\includegraphics[trim=1.7cm 1.3cm 0.7cm 1.2cm, clip, width=14.cm]{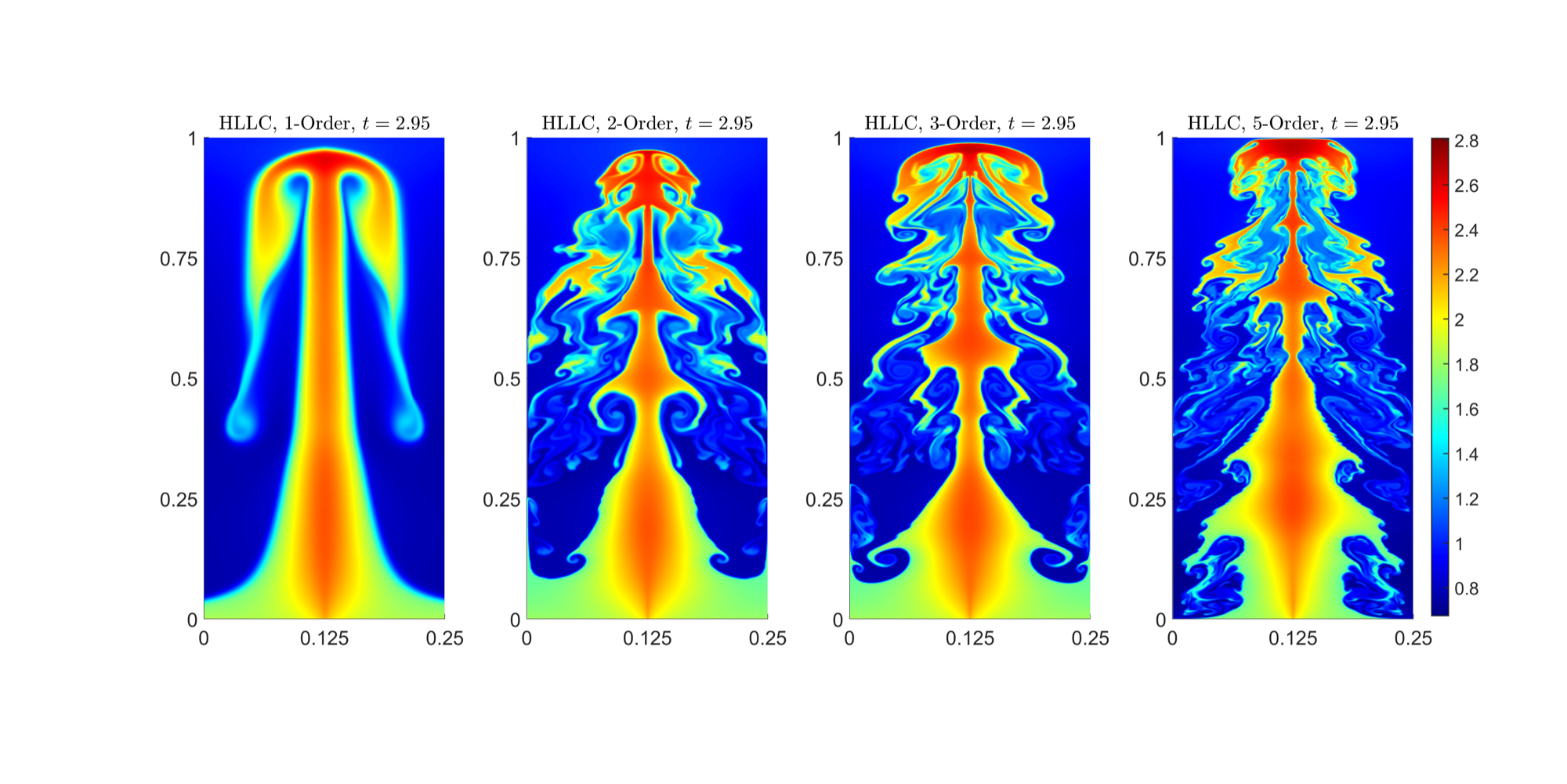}}
\vskip 15pt 
\centerline{\includegraphics[trim=1.7cm 1.3cm 0.7cm 1.2cm, clip, width=14.cm]{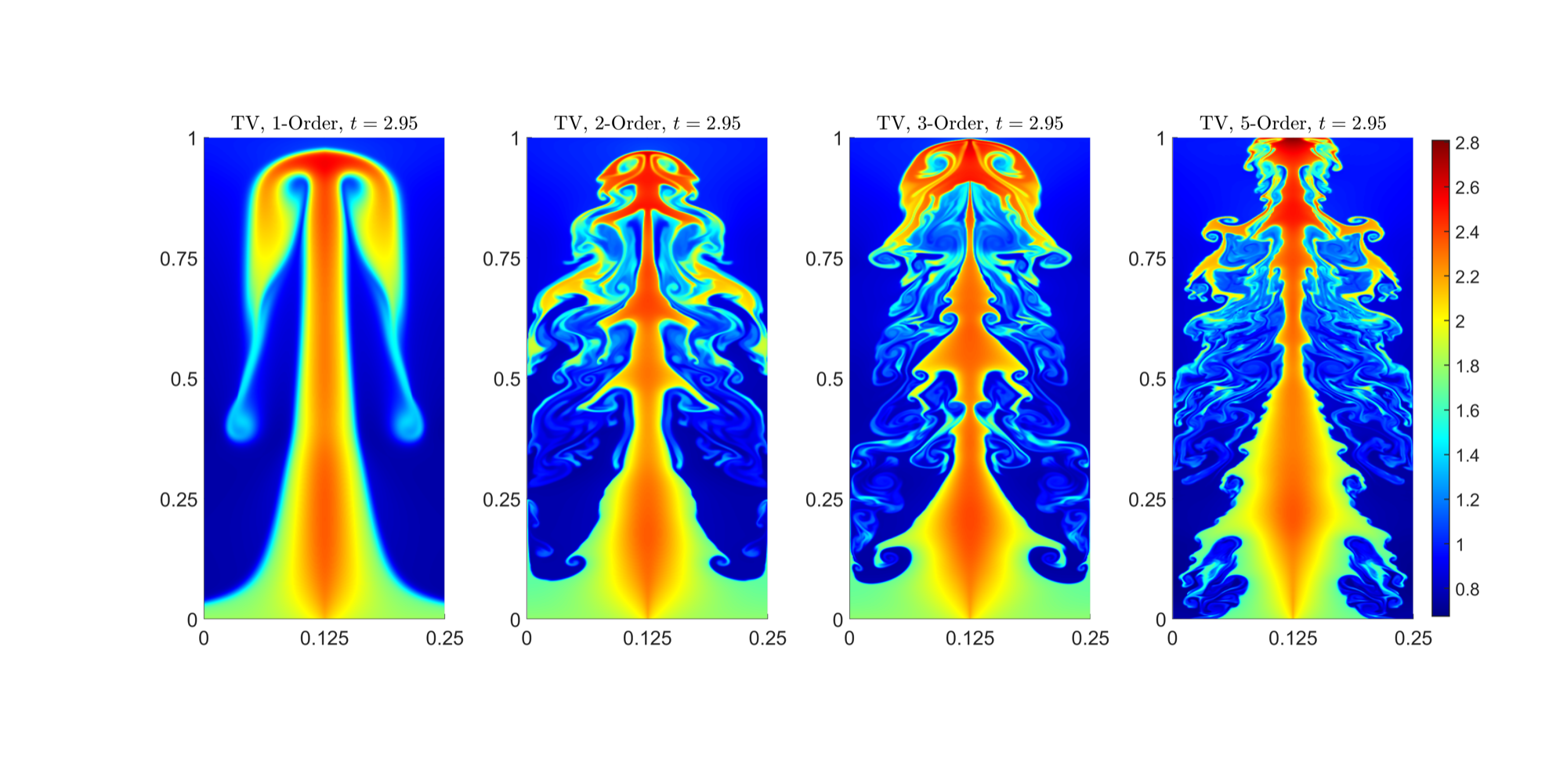}}
\caption{\sf Example 19: Density $\rho$ computed by the 1-Order (first column), 2-Order (second column), 3-Order (third column), and 5-Order (fourth column) HLL (top row), HLLC (second row), and TV (bottom row) schemes at  $t=2.95$.\label{fig13c}}
\end{figure}

\begin{figure}[ht!]
\centerline{\includegraphics[trim=1.7cm 1.3cm 0.7cm 1.2cm, clip, width=14.cm]{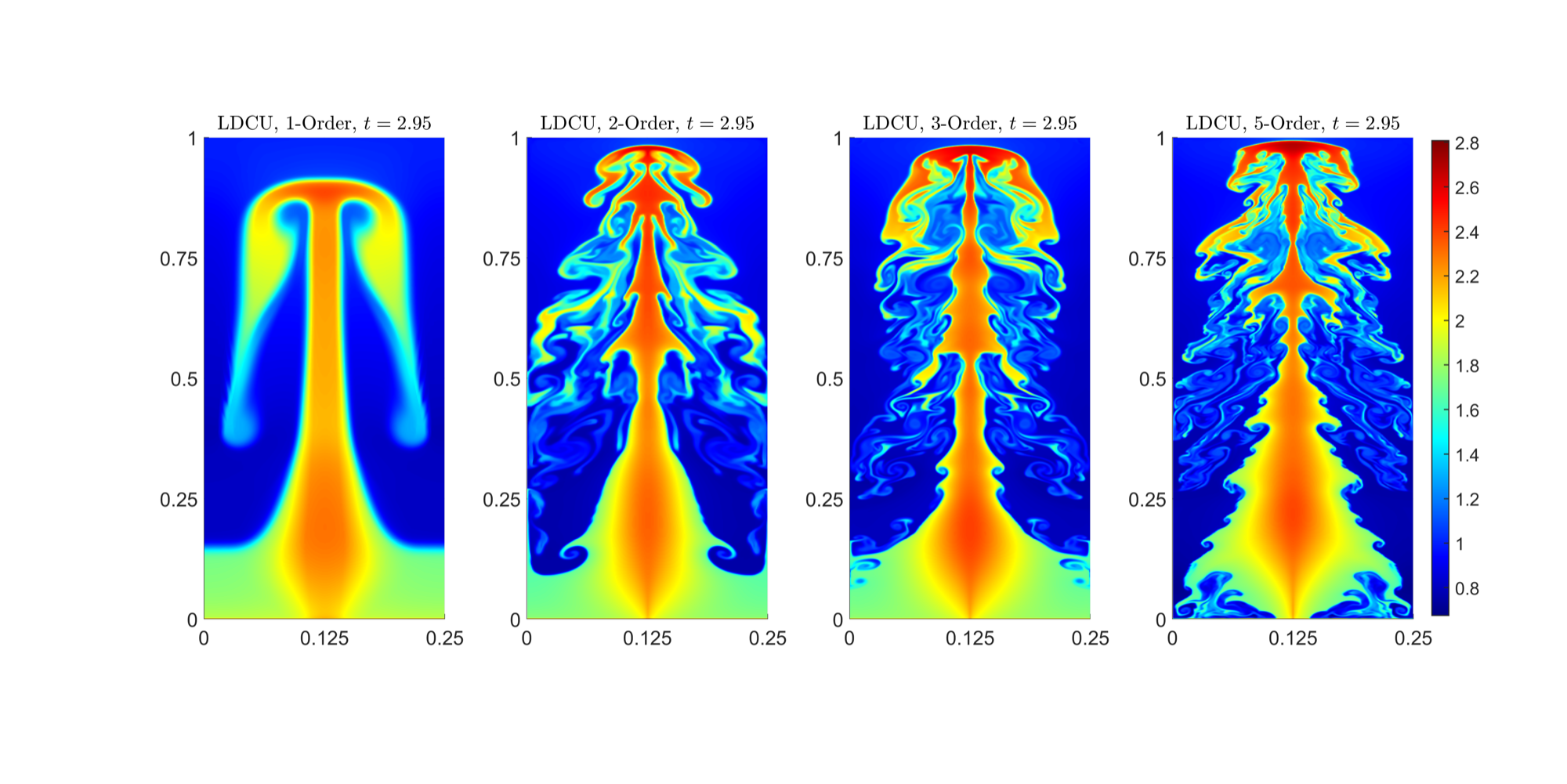}}
\vskip 15pt 
\centerline{\includegraphics[trim=1.7cm 1.3cm 0.7cm 1.2cm, clip, width=14.cm]{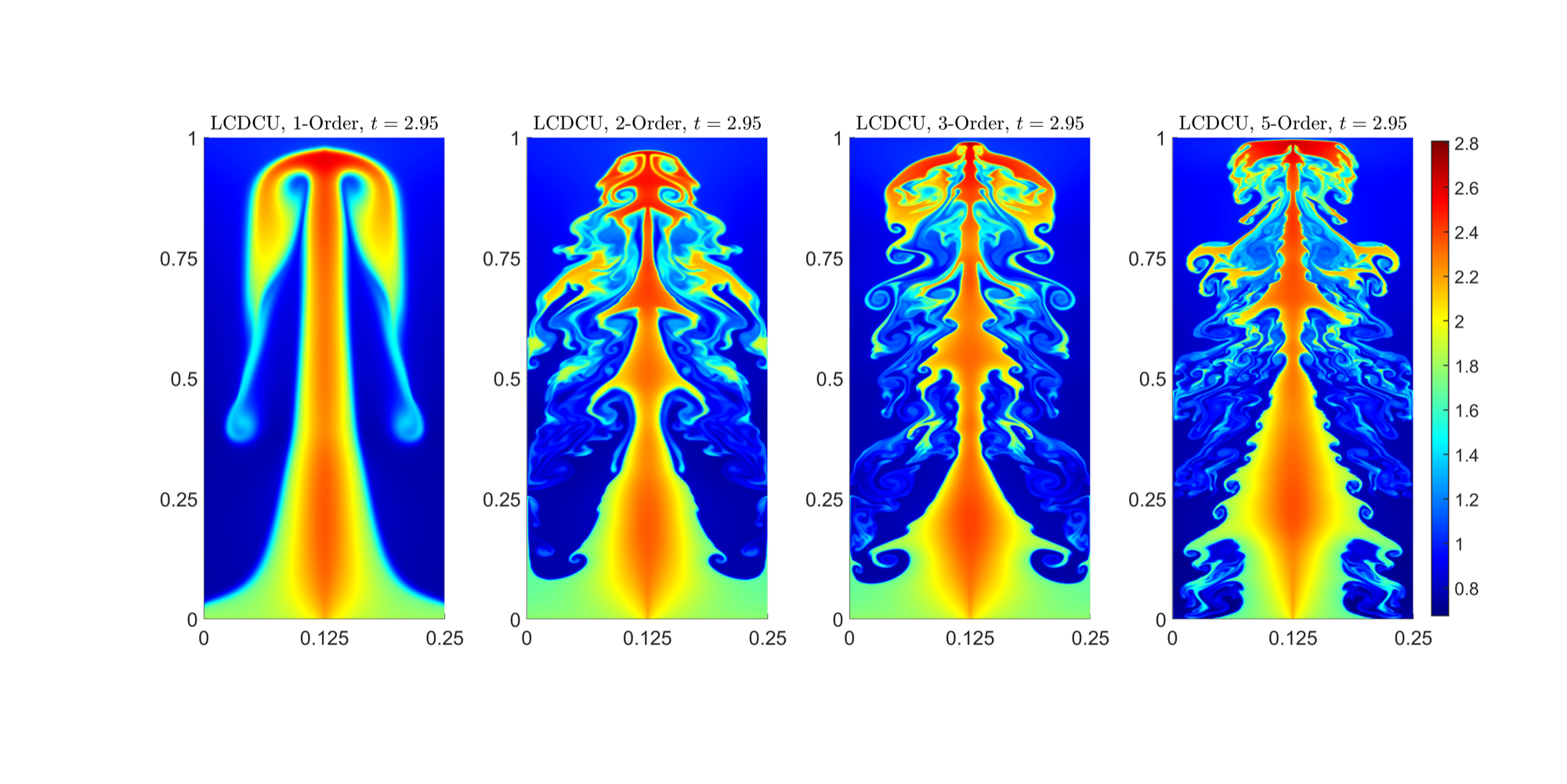}}
\caption{\sf Example 19: Density $\rho$ computed by the 1-Order (first column), 2-Order (second column), 3-Order (third column), and 5-Order (fourth column) LDCU (top row) and LCDCU (bottom row)  schemes at $t=2.95$.\label{fig13d}}
\end{figure}

We compute the numerical solutions until the final time $t=2.95$ by 1-Order, 2-Order, 3-Order, and 5-Order schemes on the uniform mesh of $256\times 1024$ cells and present the numerical results at times $t=1.95$ and 2.95 in Figures \ref{fig13a}--\ref{fig13d}. One can see that there are pronounced differences between the solutions computed by the HLL and four low-dissipation schemes since  the structures captured by the low-dissipation schemes are much more complicated, which again demonstrates that the four schemes 
contain less dissipation than the corresponding HLL counterparts. At the same time, the differences between the studied four low-dissipation schemes are limited, especially at the earlier time; see Figures \ref{fig13a}--\ref{fig13b}. At the later times, they produce different complex structures because of the unstable nature of the Rayleigh-Taylor instability problem. 

\begin{rmk}[Positivity preservation]
Preservation of positive density and pressure is an important robustness property for numerical schemes for the Euler equations. The first-order HLL and HLLC schemes are positively conservative under suitable wave-speed estimates and CFL conditions; see \cite{EMRS1991,BCLC1997}. However, these results do not automatically extend to the higher-order versions considered here. To the best of our knowledge, no general positivity-preserving result is available for the unmodified TV and LCDCU schemes used in this study. Moreover, the original LDCU construction may fail to preserve positivity; specially modified PPLDCU schemes have recently been developed to guarantee this property \cite{CGKWX2025}. Since none of our high-order implementations employs a dedicated positivity-preserving limiter, we do not claim a general positivity guarantee for them. Such a guarantee could be enforced by blending the high-order flux with a positivity-preserving first-order flux \cite{PS1996,HAS2013}. The loss of positivity observed for the third- and fifth-order TV schemes in Configuration~3 and Quirk's problem illustrates this issue.
\end{rmk}

\section{Conclusions}

In this paper, we have carried out a systematic comparative study of several low-dissipation numerical schemes for hyperbolic conservation laws, including the HLLC scheme and three recently proposed methods --- the TV flux splitting, LDCU, and LCDCU schemes. These schemes have been considered at first-, second-, third-, and fifth-order accuracy levels within finite-volume and finite-difference frameworks. Through a series of one- and two-dimensional numerical experiments for the Euler equations of gas dynamics, we have examined their accuracy, resolution, robustness, and computational efficiency. The numerical results indicate that the four low-dissipation schemes exhibit comparable levels of numerical dissipation, with noticeable differences arising only in selected test cases. The transition from first- to second-order accuracy generally produces the most significant improvement in resolution, while the third- and fifth-order versions further enhance the accuracy and resolution of complex multidimensional flow structures. Compared with their HLL counterparts, all four low-dissipation schemes provide sharper resolution of contact and shear waves. Among them, the TV splitting schemes exhibit slightly higher numerical dissipation in several examples. The newly added shock-instability tests further reveal an important trade-off between resolution and robustness. In Quirk's odd--even decoupling problem, the HLL schemes maintain essentially planar shock fronts, whereas the HLLC, TV, and LCDCU schemes exhibit more pronounced numerical shock instabilities. The LDCU schemes produce only relatively weak deviations from the HLL results. Moreover, the third- and fifth-order TV schemes lose positivity in the two-dimensional Riemann problem --- Configuration~3 and Quirk's odd--even decoupling problem. None of the high-order schemes considered in this study is equipped with a dedicated positivity-preserving limiter. Therefore, the absence of positivity failures in the other computations should be interpreted as numerical evidence of robustness rather than as a general theoretical guarantee. For problems involving near-vacuum states or extremely strong shocks, suitable positivity-preserving reconstructions or flux-limiting techniques may be required. In terms of computational cost, the LCDCU schemes are slightly more expensive than the HLLC, TV, and LDCU schemes. Nevertheless, the LCDCU framework offers greater flexibility for extension to other hyperbolic systems, since its construction primarily requires the corresponding local eigenstructure. Overall, no single scheme is uniformly superior in all aspects. The HLL schemes offer greater robustness in the presence of strong grid-aligned shocks, whereas the four low-dissipation schemes provide substantially better resolution of contact and shear waves. The comparative findings therefore provide practical guidance for selecting an appropriate numerical scheme according to the required balance among resolution, robustness, positivity, and computational cost.

\begin{DA}
\paragraph{Funding.} The work of S. Chu and  M. Herty was funded by the Deutsche Forschungsgemeinschaft (DFG, German Research Foundation) - SPP 2410 Hyperbolic Balance Laws in Fluid Mechanics: Complexity, Scales, Randomness (CoScaRa) within the Project(s) HE5386/26-1 (Numerische Verfahren für gekoppelte Mehrskalenprobleme,525842915) and (Zufällige kompressible Euler Gleichungen: Numerik und ihre Analysis, 525853336) HE5386/27-1, and  the Deutsche Forschungsgemeinschaft (DFG, German Research Foundation) - SPP 2183: Eigenschaftsgeregelte Umformprozesse with the Project(s) HE5386/19-2,19-3 Entwicklung eines flexiblen isothermen Reckschmiedeprozesses für die eigenschaftsgeregelte Herstellung von Turbinenschaufeln aus Hochtemperaturwerkstoffen (424334423). 

\paragraph{Conflicts of interest.} On behalf of all authors, the corresponding author states that there is no conflict of interest.

\paragraph{Data and software availability.} The data that support the findings of this study and Fortran codes developed by the authors and
used to obtain all of the presented numerical results are available from the corresponding author upon reasonable request.
\end{DA}

\bibliographystyle{siamnodash}
\bibliography{ref}

@article{CKX_24,
  author  = {Chu, S. and Kurganov, A. and Xin, R.},
  title   = {New low-dissipation central-upwind schemes. {P}art {II}},
  JOURNAL = {J. Sci. Comput.},
    VOLUME = {103},
      YEAR = {2025},
    NUMBER = {1},
     Note = {Paper No. 33}}

@article{Hennemann21,
  author  = {Hennemann, S. and Rueda-Ramírez, A. M. and Hindenlang, F. J. and Gassner, G. J.},
  title   = {A provably entropy stable subcell shock capturing approach for high-order split-form DG for the compressible Euler equations},
  JOURNAL = {J. Comput. Phys.},
    VOLUME = {426},
      YEAR = {2021},
    NUMBER = {1},
     Note = {Paper No. 109935}}

@article{Quirk94,
  author  = {Quirk, J. J.},
  title   = {A contribution to the great Riemann solver debate},
  JOURNAL = {Internat. J. Numer. Methods Fluids},
    VOLUME = {18},
      YEAR = {1994},
    NUMBER = {1},
    pages  = {555--574}}

@incollection {CK2023,
    AUTHOR = {Chu, S. and Kurganov, A.},
     TITLE = {Local characteristic decomposition based central-upwind scheme
              for compressible multifluids},
 BOOKTITLE = {Finite volumes for complex applications {X}. {V}ol. 2.
              {H}yperbolic and related problems},
    SERIES = {Springer Proc. Math. Stat.},
    VOLUME = {433},
     PAGES = {73--81},
 PUBLISHER = {Springer, Cham},
      YEAR = {2023},
}

@article{CKX24,
  author  = {Chu, S. and Kurganov, A. and Xin, R.},
  journal = {J. Comput. Phys.},
  title   = {Low-dissipation central-upwind schemes for compressible multifluids},
  year    = {2024},
  note    = {Paper No. 113311},
  volume  = {518}}

@article{CHK25,
  author  = {Chu, S. and Herty, M. and Kurganov, A.},
  journal = {J. Comput. Phys.},
  title   = {Novel local characteristic decomposition based path-conservative central-upwind schemes},
  year    = {2025},
  note    = {Paper No. 113692},
  volume  = {524}}

@article {GRR2017,
    AUTHOR = {Gande, N. and Rathod, Y. and Rathan, S.},
     TITLE = {Third-order {WENO} scheme with a new smoothness indicator},
   JOURNAL = {Internat. J. Numer. Methods Fluids},
    VOLUME = {85},
      YEAR = {2017},
    NUMBER = {2},
     PAGES = {90--112},
}

@article {GRR2018,
    AUTHOR = {Gande, N. and Rathod, Y. and Rathan, S.},
     TITLE = {Improved third-order weighted essentially nonoscillatory
              scheme},
   JOURNAL = {Internat. J. Numer. Methods Fluids},
    VOLUME = {87},
      YEAR = {2018},
    NUMBER = {7},
     PAGES = {329--342},
}

@article{CHT25,
  author = {Chu, S. and Herty, M. and Toro, E. F. },
  title  = {High-order flux splitting schemes for the {E}uler equations of gas dynamics},
  JOURNAL= {Comput. \& Fluids},
  VOLUME = {300},
  YEAR = {2025},
  NOTE  = {Paper No. 106738}}

@article {HS1999,
    AUTHOR = {Hu, C. and Shu, C.-W.},
     TITLE = {Weighted essentially non-oscillatory schemes on triangular
              meshes},
   JOURNAL = {J. Comput. Phys.},
    VOLUME = {150},
      YEAR = {1999},
    NUMBER = {1},
     PAGES = {97--127},
}

@article {Titarev2005,
    AUTHOR = {Titarev, V. A. and Toro, E. F.},
     TITLE = {A{DER} schemes for three-dimensional non-linear hyperbolic
              systems},
   JOURNAL = {J. Comput. Phys.},
    VOLUME = {204},
      YEAR = {2005},
    NUMBER = {2},
     PAGES = {715--736},
}

@article{CCHKL_22,
  author  = {Chertock, A. and Chu, S. and Herty, M. and Kurganov, A. and Luk\'a\v{c}ov\'a-Medvi\v{d}ov\'a, M.},
  journal = {J. Comput. Phys.},
  title   = {Local characteristic decomposition based central-upwind scheme},
  year    = {2023},
  note    = {Paper No. 111718},
  volume  = {473},
}

@article {CCK23_Adaptive,
    Author = {Chertock, A. and Chu, S. and Kurganov, A.},
    Journal = {East Asian J. Appl. Math.},
    Number = {6},
    Pages = {576-609},
    Title = {Adaptive high-order {A-WENO} schemes based on a new local smoothness indicator},
    Volume = {13},
    Year = {2023}}

@article {SO89,
	Author = {Shu, C.-W. and Osher, S.},
	Journal = {J. Comput. Phys.},
	Number = {1},
	Pages = {32--78},
	Title = {Efficient implementation of essentially nonoscillatory
              shock-capturing schemes. {II}},
	Volume = {83},
	Year = {1989}}

@article{KX_22,
    AUTHOR = {Kurganov, A. and Xin, R.},
     TITLE = {New low-dissipation central-upwind schemes},
   JOURNAL = {J. Sci. Comput.},
    VOLUME = {96},
      YEAR = {2023},
     NOTE  = {Paper No. 56}}

@article{CKX23,
  author = {Chu, S. and Kurganov, A. and Xin, R.},
  title  = {New more efficient {A-WENO} schemes},
  Journal= {J. Sci. Comput.},
  Volume = {104},
    Year = {2025},
    NOTE  = {Paper No. 53}}

@article{CH_Third,
  author = {Chu, S. and Fu, Q. and  Herty, M. and Kurganov, A.},
  title  = {Novel and efficient third-order {WENO} schemes},
  note   = {In preparation.}}

@article{Arminjon95,
	Author = {Arminjon, P. and Viallon, M. C.},
	Journal = {C. R. Acad. Sci. Paris S\'{e}r. I Math.},
	Pages = {85--88},
	Number = {3},
	Title = {G\'{e}n\'{e}ralisation du sch\'{e}ma de {N}essyahu-{T}admor pour une
              \'{e}quation hyperbolique \`a deux dimensions d'espace},
	Volume = {320},
	Year = {1995}}

@article{Lious1993,
	Author = {Liou, M. S.},
	Journal = {J. Comput. Phys.},
	Pages = {23--39},
	Number = {3},
	Title = {A new flux splitting scheme},
	Volume = {107},
	Year = {1993}}

@article{Lious1996,
	Author = {Liou, M. S.},
	Journal = {J. Comput. Phys.},
	Pages = {364--382},
	Number = {3},
	Title = {A sequel to {AUSM: AUSM+}},
	Volume = {129},
	Year = {1996}}

@incollection {CKX22,
    AUTHOR = {Chu, S. and Kurganov, A. and Xin, R.},
     TITLE = {A fifth-order A-{WENO} scheme based on the low-dissipation
              central-upwind fluxes},
 BOOKTITLE = {Hyperbolic {P}roblems: {T}heory, {N}umerics, {A}pplications.
              {V}ol. {II}},
    SERIES = {SEMA SIMAI Springer Ser.},
    VOLUME = {35},
     PAGES = {51--61},
 PUBLISHER = {Springer, Cham},
      YEAR = {2024}
}

@article{Lious2006,
	Author = {Liou, M. S.},
	Journal = {J. Comput. Phys.},
	Pages = {137--170},
	Number = {3},
	Title = {A Sequel to {AUSM}, Part {II}: {AUSM}+-up for all Speeds},
	Volume = {214},
	Year = {2006}}

@article{ZB1993,
	Author = {Zha, G.-C. and Bilgen, E.},
	Journal = { Int. J. Numer. Methods Fluids},
	Pages = {115--144},
	Number = {3},
	Title = { Numerical solution of Euler equations by a new flux vector splitting scheme},
	Volume = {17},
	Year = {1993}}

@article{DT2011,
	Author = {Dumbser, M. and Toro, E. F.},
	Journal = {Commun. Comput. Phys.},
	Pages = {635--671},
	Number = {3},
	Title = { On universal Osher-type schemes for general nonlinear hyperbolic conservation laws},
	Volume = {10},
	Year = {2011}}

@article{TSS1994,
	Author = {Toro, E. F. and Spruce, M.  and Speares, W. },
	Journal = {Shock Waves},
	Pages = {25--34},
	Number = {3},
	Title = {Restoration of the contact surface in the HLL-Riemann solver},
	Volume = {4},
	Year = {1994}}

@article {HLL1983,
    AUTHOR = {Harten, A. and Lax, P. D. and van Leer, B.},
     TITLE = {On upstream differencing and {G}odunov-type schemes for
              hyperbolic conservation laws},
   JOURNAL = {SIAM Rev.},
    VOLUME = {25},
      YEAR = {1983},
    NUMBER = {1},
     PAGES = {35--61},
}

@incollection {TT2017,
    AUTHOR = {Tokareva, S. and Toro, E. F.},
     TITLE = {A flux splitting method for the {B}aer-{N}unziato equations of
              compressible two-phase flow},
 BOOKTITLE = {Finite volumes for complex applications {VIII}---hyperbolic,
              elliptic and parabolic problems},
    SERIES = {Springer Proc. Math. Stat.},
    VOLUME = {200},
     PAGES = {127--135},
 PUBLISHER = {Springer, Cham},
      YEAR = {2017}
}

@article{LLN2012,
	Author = { Luo, H. and Luo, L. Q. and Nourgaliev, R.},
	Journal = {Commun. Comput. Phys.},
	Pages = {1495--1519},
	Number = {3},
	Title = {A reconstructed discontinuous Galerkin method for the Euler equations on arbitrary grids},
	Volume = {12},
	Year = {2012}}

@article{GVM1999,
	Author = {Gressier, J. and Villedieu, P. and Moschetta, J. M.},
	Journal = {J. Comput. Phys. },
	Pages = {199--220},
	Number = {3},
	Title = {Positivity of flux vector splitting schemes},
	Volume = {155},
	Year = {1999}}

@article{KSFW2011,
	Author = {Kitamura, K. and Shima, E. and Fujimoto, K. and Wang, Z. J.},
	Journal = {Commun. Comput. Phys.},
	Pages = {90--119},
	Number = {3},
	Title = {Performance of low-dissipation Euler fluxes and preconditioned LU-SGS at low speeds},
	Volume = {10},
	Year = {2011}}

@article{TV2012,
	Author = {Toro, E. F. and V\'azquez-Cend\'on, M. E.},
	Journal = {Comput. \& Fluids},
	Pages = {1--12},
	Number = {3},
	Title = {Flux splitting schemes for the {E}uler equations},
	Volume = {70},
	Year = {2012}}

@article {TCL2015,
    AUTHOR = {Toro, E. F. and Castro, C. E. and Lee, B. J.},
     TITLE = {A novel numerical flux for the 3{D} {E}uler equations with
              general equation of state},
   JOURNAL = {J. Comput. Phys.},
    VOLUME = {303},
      YEAR = {2015},
     PAGES = {80--94}
}

@article {BMT2016,
    AUTHOR = {Balsara, D. S. and Montecinos, G. I. and Toro, E. F.},
     TITLE = {Exploring various flux vector splittings for the
              magnetohydrodynamic system},
   JOURNAL = {J. Comput. Phys.},
    VOLUME = {311},
      YEAR = {2016},
     PAGES = {1--21}
}

@article {DBTF2019,
    AUTHOR = {Dumbser, M. and Balsara, D. S. and Tavelli, M. and Fambri, F.},
     TITLE = {A divergence-free semi-implicit finite volume scheme for
              ideal, viscous, and resistive magnetohydrodynamics},
   JOURNAL = {Internat. J. Numer. Methods Fluids},
    VOLUME = {89},
      YEAR = {2019},
    NUMBER = {1-2},
     PAGES = {16--42}
}

@article {TCVS2022,
    AUTHOR = {Toro, E. F. and Castro, C. E. and Vanzo,
              D. and Siviglia, A.},
     TITLE = {A flux-vector splitting scheme for the shallow water equations
              extended to high-order on unstructured meshes},
   JOURNAL = {Internat. J. Numer. Methods Fluids},
    VOLUME = {94},
      YEAR = {2022},
    NUMBER = {10},
     PAGES = {1679--1705}
}

@Incollection {Lious1998,
     Address = {Berlin, Verlag},
     TITLE = {Recent progress and applications of {AUSM+}},
    SERIES = {Lecture Notes in Physics},
    VOLUME = {515},
    Author = {Liou, M. S.},
 BOOKTITLE = {Sixteenth {I}nternational {C}onference on {N}umerical
              {M}ethods in {F}luid {D}ynamics},
 PUBLISHER = {Springer},
      YEAR = {1998},
}

@book{BAF,
    AUTHOR = {Ben-Artzi, M. and Falcovitz, J.},
     TITLE = {Generalized {R}iemann problems in computational fluid
dynamics},
    SERIES = {Cambridge Monographs on Applied and Computational
Mathematics},
    VOLUME = {11},
 PUBLISHER = {Cambridge University Press},
   ADDRESS = {Cambridge},
      YEAR = {2003},
     PAGES = {xvi+349},
      ISBN = {0-521-77296-6}}

@article {LXDYG2022,
    AUTHOR = {Li, X. G. and Xia, T. and Deng, Y. X. and Yang, S. Q. and
              Ge, Y. B.},
     TITLE = {A new third-order finite difference {WENO} scheme to improve convergence rate at critical points},
   JOURNAL = {Int. J. Comput. Fluid Dyn.},
    VOLUME = {36},
      YEAR = {2022},
    NUMBER = {10},
     PAGES = {857--874}
}

@incollection{VL1982a,
    Address = {Berlin, Heidelberg},
    AUTHOR = {van Leer, B.},
    TITLE = {Flux-vector splitting for the Euler equations},
    Booktitle = {Eighth International Conference on Numerical Methods in Fluid Dynamics},
    Publisher = {Springer},
    pages={507--512},
    Year = {1982}
}

@article {feireisl21,
    AUTHOR = {Feireisl, E. and Luk\'{a}\v{c}ov\'{a}-Medvi\v{d}ov\'{a}, M. and She,
              B. and Wang, Y.},
     TITLE = {Computing oscillatory solutions of the {E}uler system via
              {$\mathcal{K}$}-convergence},
   JOURNAL = {Math. Models Methods Appl. Sci.},
  FJOURNAL = {Mathematical Models and Methods in Applied Sciences},
    VOLUME = {31},
      YEAR = {2021},
    NUMBER = {3},
     PAGES = {537--576}}

@article{Fjordholm16,
	Author = {Fjordholm, U. S. and Mishra, S. and Tadmor, E.},
	Journal = {Acta Numer.},
	Pages = {567-679},
	Number = {326},
	Title = {On the computation of measure-valued solutions},
	Volume = {25},
	Year = {2016}}

@article{Fri,
        Author = {Friedrichs, K. O.},
        Fjournal = {Communications on Pure and Applied Mathematics},
        Journal = {Comm. Pure Appl. Math.},
        Pages = {345--392},
        Title = {Symmetric hyperbolic linear differential equations},
        Volume = {7},
        Year = {1954}}

@article{Garg21,
	Author = {Garg, N. K. and Kurganov, A. and Liu, Y.},
	Journal = {J. Comput. Phys.},
	Number = {326},
	Title = {Semi-discrete central-upwind {R}ankine-{H}ugoniot schemes for hyperbolic systems of conservation laws},
	Volume = {428},
	Year = {2021},
    NOTE  = {Paper No. 110078}}

@article {Godunov59,
	Author = {Godunov, S. K.},
	Journal = {Mat. Sb. (N.S.)},
	Number = {89},
	Pages = {271--306},
	Title = {A difference method for numerical calculation of discontinuous
              solutions of the equations of hydrodynamics},
	Volume = {47},
	Year = {1959}}

@book{Gottlieb11,
     AUTHOR = {Gottlieb, S. and Ketcheson, D. and Shu, C.-W.},
     TITLE = {Strong stability preserving {R}unge-{K}utta and multistep time discretizations},
     PUBLISHER = {World Scientific Publishing Co. Pte. Ltd., Hackensack, NJ},
     YEAR = {2011},
     PAGES = {xii+176}}

@book{Hesthaven18,
     AUTHOR = {Hesthaven, J. S.},
     TITLE = {Numerical methods for conservation laws: From analysis to algorithms},
     SERIES = {Comput. Sci. Eng. 18},
     PUBLISHER = {SIAM, Philadelphia},
     YEAR = {2018}}

@article {Gottlieb12,
	Author = {Gottlieb, S. and Shu, C.- W. and Tadmor, E.},
	Journal = {SIAM Rev.},
	Number = {1},
	Pages = {89--112},
	Title = {Strong stability-preserving high-order time discretization
	methods},
	Volume = {43},
	Year = {2001}}

@article{Jiang13,
    Author = {Jiang, Y. and Shu, C.- W. and Zhang, M. P. },
    Journal = {SIAM J. Sci. Comput.},
    Pages = {A1137--A1160},
    Number = {2},
    Title = {An alternative formulation of finite difference weighted {ENO} schemes with {L}ax-{W}endroff time discretization for
    	conservation laws},
    Volume = {35},
    Year = {2013}}

@article{Jiang98,
    Author = {Jiang, G. S. and Tadmor, E.},
    Journal = {SIAM J. Sci. Comput.},
    Pages = {1892--1917},
    Number = {6},
    Title = {Nonoscillatory central schemes for multidimensional hyperbolic conservation laws},
    Volume = {19},
    Year = {1998}}

@book{KLR20,
    AUTHOR = {Ketcheson, D. I. and LeVeque, R. J. and del Razo, M. J.},
     TITLE = {Riemann problems and {J}upyter solutions},
    SERIES = {Fundamentals of Algorithms},
    VOLUME = {16},
 PUBLISHER = {Society for Industrial and Applied Mathematics (SIAM), Philadelphia, PA},
      YEAR = {2020},
     PAGES = {xii+166}}

@article {Kurganov00,
    Author = {Kurganov, A. and Tadmor, E.},
    Journal = {J. Comput. Phys.},
    Number = {2},
    Pages = {720--742},
    Title = {New high-resolution semi-discrete central schemes for {H}amilton-{J}acobi equations},
    Volume = {160},
    Year = {2000}}

@article {Kurganov01,
    Author = {Kurganov, A. and   Noelle, P. and Petrova, G.},
    Journal = {SIAM J. Sci. Comput.},
    Number = {3},
    Pages = {707-740},
    Title = {Semidiscrete central-upwind schemes for hyperbolic
              conservation laws and {H}amilton-{J}acobi equations},
    Volume = {23},
    Year = {2001}}

@article {Kurganov21a,
    Author = {Kurganov, A. and  Liu, Y. and Zeitlin, Y.},
    Journal = {ESAIM Math. Model. Num. Anal.},
    Number = {24},
    Pages = {713-734},
    Title = {Numerical dissipation switch for two-dimensional central-upwind schemes},
    Volume = {55},
    Year = {2021}}

@article {Kurganov07,
	Author = {Kurganov, A. and Lin, C.-T.},
	Journal = {Commun. Comput. Phys.},
	Pages = {141--163},
	Number = {1},
	Title = {On the reduction of numerical dissipation in central-upwind
	schemes},
	Volume = {2},
	Year = {2007}}

@article{Kurganov02,
	Author = {Kurganov, A. and Tadmor, E.},
        Journal = {Numer. Methods Partial Differential Equations},
        Pages = {584--608},
        Title = {Solution of two-dimensional Riemann problems for gas dynamics without {R}iemann problem solvers},
        Volume = {18},
        Year = {2002}}

@article {Kurganov17,
    Author = {Kurganov, A. and  Prugger, M. and Wu, T.},
    Journal = {SIAM J. Sci. Comput.},
    Number = {24},
    Pages = {A947--A965},
    Title = {Second-order fully discrete central-upwind scheme for two-dimensional hyperbolic systems of conservation laws},
    Volume = {39},
    Year = {2017}}

@article{Lax,
        Author = {Lax, P. D.},
        Fjournal = {Communications on Pure and Applied Mathematics},
        Journal = {Comm. Pure Appl. Math.},
        Pages = {159--193},
        Title = {Weak solutions of nonlinear hyperbolic equations and
their numerical computation},
        Volume = {7},
        Year = {1954}}

@book{Leveque02,
     AUTHOR = {LeVeque, R. J.},
     TITLE = {Finite Volume Methods for Hyperbolic Problems},
     SERIES = {Cambridge Texts in Appl. Math.},
     PUBLISHER = {Cambridge University Press, Cambridge, UK},
     YEAR = {2002}}

@article {Levy99,
	Author = {Levy, D. and Puppo, G. and Russo, G.},
	Journal = {M2AN Math. Model. Numer. Anal.},
	Number = {3},
	Pages = {547--571},
	Title = {Central {WENO} schemes for hyperbolic systems of conservation laws},
	Volume = {33},
	Year = {1999}}

@article {Lie03,
	Author = {Lie, K.-A. and Noelle, S.},
	Journal = {SIAM J. Sci. Comput.},
	Number = {4},
	Pages = {1157--1174},
	Title = {On the artificial compression method for second-order
	nonoscillatory central difference schemes for systems of
	conservation laws},
	Volume = {24},
	Year = {2003}}

@article {Lie03a,
	Author = {Lie, K.-A. and Noelle, S.},
	Journal = {J. Sci. Comput.},
	Number = {1},
	Pages = {69--81},
	Title = {An improved quadrature rule for the flux-computation in staggered central difference schemes in multidimensions},
	Volume = {18},
	Year = {2003}}

@article{Liu17,
	Author = {Liu, H.},
	Journal = {Appl. Math. Comput.},
	Pages = {182--197},
	Title = {A numerical study of the performance of alternative weighted
	{ENO} methods based on various numerical fluxes for
	conservation law},
	Volume = {296},
	Year = {2017}}

@article {Liska03,
    Author = {Liska, R. and Wendroff, B.},
    Journal = {SIAM J. Sci. Comput.},
    Number = {24},
    Pages = {995--1017},
    Title = {Comparison of several difference schemes on 1{D} and 2{D} test problems for the Euler equations},
    Volume = {25},
    Year = {2003}}

@article {Liu98,
    Author = {Liu, X. D. and Tadmor, E.},
    Journal = {Numer. Math.},
    Number = {3},
    Pages = {397--425},
    Title = {Third order nonoscillatory central scheme for hyperbolic conservation laws},
    Volume = {79},
    Year = {1998}}

@article {Nessyahu90,
	Author = {Nessyahu, H. and Tadmor, E.},
	Journal = {J. Comput. Phys.},
	Number = {2},
	Pages = {408--463},
	Title = {Nonoscillatory central differencing for hyperbolic
	conservation laws},
	Volume = {87},
	Year = {1990}}

@article{Panuelos20,
	Author = {Panuelos, J. and Wadsley, J. and Kevlahan, N.},
	Journal = {J. Comput. Phys.},
	Title = {Low shear diffusion central schemes for particle methods},
	Volume = {414},
	Year = {2020},
    NOTE  = {Paper No. 109454}}

@article{KLX2025,
	Author = {Kurganov, A. and Liu, Z. and Xin, R.},
	Journal = {Commun. Comput. Phys.},
	Title = {Low-dissipation central-upwind schemes for elasticity in heterogeneous media},
	Volume = {38},
	Year = {2025},
    pages  = {156--180}}

@article{Schulz93,
  title={Classification of the {R}iemann problem for two-dimensional gas dynamics},
  author={Schulz-Rinne, C. W.},
  Journal = {SIAM J. Math. Anal.},
  volume={24},
  pages={76--88},
  year={1993}}

@article{Schulz93a,
  title={Numerical solution of the {R}iemann problem for two-dimensional gas dynamics},
  author={Schulz-Rinne, C. W. and Collins, J. P. and Glaz, H. M.},
  Journal = {SIAM J. Sci. Comput.},
  volume={14},
  pages={1394--1414},
  year={1993}}

@article {Shi03,
    AUTHOR = {Shi, J. and Zhang, Y.-T. and Shu, C.-W.},
     TITLE = {Resolution of high order {WENO} schemes for complicated flow structures},
   JOURNAL = {J. Comput. Phys.},
    VOLUME = {186},
      YEAR = {2003},
    NUMBER = {2},
     PAGES = {690--696}}

@article{Shu20,
	Author = {Shu, C.-W.},
	Journal = {Acta Numer.},
	Pages = {701--762},
	Number = {29},
	Title = {Essentially non-oscillatory and weighted essentially non-oscillatory schemes},
	Volume = {29},
	Year = {2020}}

@article{Shu88,
	Author = {Shu, C.-W. and Osher, S.},
	Journal = {J. Comput. Phys.},
	Pages = {439--471},
	Number = {29},
	Title = {Efficient implementation of essentially non-oscillatory shock-capturing schemes},
	Volume = {77},
	Year = {1988}}

@article {Sweby84,
	Author = {Sweby, P. K.},
	Journal = {SIAM J. Numer. Anal.},
	Number = {5},
	Pages = {995--1011},
	Title = {High resolution schemes using flux limiters for hyperbolic
	conservation laws},
	Volume = {21},
	Year = {1984}}

@book{Toro2009,
	Address = {Berlin, Heidelberg},
	Author = {Toro, E. F.},
	Edition = {Third},
	Pages = {xx+724},
	Publisher = {Springer-Verlag},
	Title = {Riemann solvers and numerical methods for fluid dynamics: {A} practical introduction},
	Year = {2009}}

@article{wang18,
    AUTHOR = {Wang, B.-S. and Li, P. and Gao, Z. and Don, W. S.},
    JOURNAL = {J. Comput. Phys.},
    PAGES = {469-477},
    NUMBER = {1},
    TITLE = {An improved fifth order alternative {WENO-Z} finite difference scheme for hyperbolic conservation laws},
    VOLUME = {374},
    YEAR = {2018}}

@article{ATvL1986,
    AUTHOR = {Anderson, W. K. and Thomas, J. L. and van Leer, B.},
    JOURNAL = {AIAA J.},
    PAGES = {1453--1460},
    NUMBER = {1},
    TITLE = {Comparison of finite volume flux vector splittings for the Euler equations},
    VOLUME = {24},
    YEAR = {1986}}

@article{ATR1989,
    AUTHOR = {Anderson, W. K. and Thomas, J. L. and Rumsey, L.},
    JOURNAL = {AIAA J.},
    PAGES = {673--674},
    NUMBER = {1},
    TITLE = {Extension and application of flux–vector splitting to calculations on dynamic meshes},
    VOLUME = {27},
    YEAR = {1989}}

@article{SW1981,
    AUTHOR = {Steger, J. L. and Warming, R. F.},
    JOURNAL = {J. Comput. Phys.},
    PAGES = {263--293},
    NUMBER = {1},
    TITLE = {Flux vector splitting of the inviscid gas dynamic equations with applications to finite difference methods},
    VOLUME = {40},
    YEAR = {1981}}

@article{VL1982,
    AUTHOR = {van Leer, B.},
    JOURNAL = {Technical Report ICASE 82-30},
    PAGES = { },
    NUMBER = {1},
    TITLE = {Flux–vector splitting for the Euler equations},
    VOLUME = { },
    YEAR = {1982}}

@article{Woodward88,
    AUTHOR = {Woodward, P. and Colella, P.},
     TITLE = {The numerical solution of two-dimensional fluid flow with strong shocks},
    JOURNAL = {J. Comput. Phys.},
    VOLUME = {54},
      YEAR = {1984},
    NUMBER = {1},
     PAGES = {115--173}}

@book{Zheng01,
     AUTHOR = {Zheng, Y.},
     TITLE = {Systems of conservation laws. Two-dimensional Riemann problems},
     SERIES = {Progress in Nonlinear Differential Equations and their Applications},
     PUBLISHER = {Birkh\"{a}user Boston, Inc., Boston, MA},
     YEAR = {2001}}

@article{Jiang96,
	Author = {Jiang, G.-S. and Shu, C.-W.},
	Journal = {J. Comput. Phys.},
	Pages = {202--228},
	Number = {5},
	Title = {Efficient implementation of weighted {ENO} schemes},
	Volume = {126},
	Year = {1996}}

@article{EMRS1991,
    AUTHOR = {Einfeldt, B. and Munz, C.-D. and Roe, P. L. and Sj{\"o}green, B.},
     TITLE = {On {G}odunov-type methods near low densities},
   JOURNAL = {J. Comput. Phys.},
    VOLUME = {92},
      YEAR = {1991},
    NUMBER = {2},
     PAGES = {273--295}}

@article{BCLC1997,
    AUTHOR = {Batten, P. and Clarke, N. and Lambert, C. and Causon, D. M.},
     TITLE = {On the choice of wavespeeds for the {HLLC} {R}iemann solver},
   JOURNAL = {SIAM J. Sci. Comput.},
    VOLUME = {18},
      YEAR = {1997},
    NUMBER = {6},
     PAGES = {1553--1570}}

@article{CGKWX2025,
    AUTHOR = {Cui, S. and Gu, Y. and Kurganov, A. and Wu, K. and Xin, R.},
     TITLE = {Positivity-preserving new low-dissipation central-upwind
              schemes for compressible {E}uler equations},
   JOURNAL = {J. Comput. Phys.},
    VOLUME = {538},
      YEAR = {2025},
      NOTE = {Paper No. 114189}}

@article{PS1996,
    AUTHOR = {Perthame, B. and Shu, C.-W.},
     TITLE = {On positivity preserving finite volume schemes for
              {E}uler equations},
   JOURNAL = {Numer. Math.},
    VOLUME = {73},
      YEAR = {1996},
    NUMBER = {1},
     PAGES = {119--130}}

@article{HAS2013,
    AUTHOR = {Hu, X. Y. and Adams, N. A. and Shu, C.-W.},
     TITLE = {Positivity-preserving method for high-order conservative
              schemes solving compressible {E}uler equations},
   JOURNAL = {J. Comput. Phys.},
    VOLUME = {242},
      YEAR = {2013},
     PAGES = {169--180}}
\end{document}